\documentclass[10pt]{amsart}

\usepackage{graphicx}        
\usepackage{multicol}        
\usepackage[bottom]{footmisc}
\usepackage{amscd,amsmath,amssymb,amsfonts}
\usepackage[cmtip, all]{xy}

\usepackage{bbm}
\usepackage{tikz-cd}
\usepackage{enumerate}
\usepackage{mathrsfs}
\usepackage{amsthm}

\usepackage[plainpages, urlcolor=blue]{hyperref}

\newtheorem{theorem}{Theorem}[section]
\newtheorem{proposition}[theorem]{Proposition}
\newtheorem{lemma}[theorem]{Lemma}
\newtheorem{corollary}[theorem]{Corollary}

\newtheorem{IntroThm}{Theorem}

\theoremstyle{definition}
\newtheorem{definition}[theorem]{Definition}
\newtheorem{IntroDef}[IntroThm]{Definition}
\theoremstyle{remark}
\newtheorem{remark}[theorem]{Remark}

\newtheorem{IntroRemark}[IntroThm]{Remark}
\newtheorem{ex}[theorem]{Example}
\newtheorem{IntroEx}[IntroThm]{Example}
\newtheorem{exs}[theorem]{Examples}
\newtheorem{Not}[theorem]{Notations}
\newtheorem{construction}[theorem]{Construction}
\numberwithin{equation}{section}

\newcommand{\CH}{{\rm CH}}

\newcommand{\red}{{\rm red}}
\newcommand{\codim}{{\rm codim}}

\newcommand{\Pic}{{\rm Pic}}

\newcommand{\Hom}{{\rm Hom}}
\newcommand{\im}{{\rm im}}
\newcommand{\Spec}{{\rm Spec\,}}

\newcommand{\Char}{{\rm char}}
\newcommand{\Tr}{{\text{Tr}}}

\newcommand{\0}{\emptyset}

\newcommand{\sC}{{\mathcal C}}

\newcommand{\sE}{{\mathcal E}}
\newcommand{\sF}{{\mathcal F}}
\newcommand{\sG}{{\mathcal G}}
\newcommand{\sH}{{\mathcal H}}
\newcommand{\sI}{{\mathcal I}}

\newcommand{\sK}{{\mathcal K}}
\newcommand{\sL}{{\mathcal L}}
\newcommand{\sM}{{\mathcal M}}
\newcommand{\sN}{{\mathcal N}}
\newcommand{\sO}{{\mathcal O}}

\newcommand{\sT}{{\mathcal T}}
\newcommand{\sU}{{\mathcal U}}
\newcommand{\sV}{{\mathcal V}}
\newcommand{\sW}{{\mathcal W}}
\newcommand{\sX}{{\mathcal X}}
\newcommand{\sY}{{\mathcal Y}}

\newcommand{\SF}{{\mathcal S\mathcal F}}
\newcommand{\A}{{\mathbb A}}

\newcommand{\G}{{\mathbb G}}

\newcommand{\N}{{\mathbb N}}
\renewcommand{\P}{{\mathbb P}}

\newcommand{\R}{{\mathbb R}}
\newcommand{\T}{{\mathbb T}}

\newcommand{\V}{{\mathbb V}}
\newcommand{\W}{{\mathbb W}}

\newcommand{\Z}{{\mathbb Z}}

\newcommand{\BM}{{\operatorname{B.M.}}}

\renewcommand{\det}{\operatorname{det}}

\newcommand{\id}{{\operatorname{\rm Id}}}
\newcommand{\Zar}{{\text{\rm Zar}}} 
\newcommand{\Nis}{{\text{\rm Nis}}} 

\newcommand{\Sch}{{\operatorname{\mathbf{Sch}}}} 
\newcommand{\colim}{{\mathop{\rm colim}}}

\newcommand{\holim}{\mathop{{\rm holim}}}
\newcommand{\op}{{\text{\rm op}}}
\newcommand{\<}{\langle}
\renewcommand{\>}{\rangle}

\newcommand{\mov}{{\mathfrak{m}}}

\renewcommand{\dim}{{\operatorname{\rm dim}}}

\newcommand{\sep}{{\operatorname{sep}}}

\newcommand{\fr}{{\operatorname{fr}}}

\newcommand{\del}{\partial}

\newcommand{\Spc}{{\mathbf{Spc}}}

\newcommand{\Sm}{{\mathbf{Sm}}}

\newcommand{\Ab}{{\mathbf{Ab}}}

\newcommand{\Ind}{{\operatorname{ind}}}
\newcommand{\Sym}{{\operatorname{Sym}}}

\newcommand{\Gr}{{\operatorname{\rm Gr}}} 
\newcommand{\rnk}{{\operatorname{\text{rnk}}}} 
\newcommand{\Bl}{\text{Bl}}

\newcommand{\GW}{{\operatorname{GW}}} 
\newcommand{\sGW}{{\mathcal{GW}}} 
\newcommand{\SH}{{\operatorname{SH}}} 
\newcommand{\Th}{{\operatorname{Th}}} 
 
\renewcommand{\th}{{\operatorname{th}}} 
\newcommand{\sHom}{\mathcal{H}om}
 
\newcommand{\Aut}{{\operatorname{Aut}}}

\newcommand{\GL}{\operatorname{GL}}
\newcommand{\SL}{\operatorname{SL}}

\newcommand{\BSL}{\operatorname{BSL}}

\newcommand{\sq}{{\operatorname{sq}}}

\newcommand{\perf}{\text{\it perf}}

\newcommand{\et}{\text{\'et}}

\newcommand{\Map}{{\operatorname{Maps}}}

\newcommand{\ind}[1]{}

\newcommand{\inp}[1]{}

\newcommand{\res}{{\operatorname{res}}}

\newcommand{\gen}{{\operatorname{gen}}}

\newcommand{\gr}{{\operatorname{gr}}}
\newcommand{\Tot}{{\operatorname{Tot}}}

\newcommand{\EM}{{\operatorname{EM}}}

\newcommand{\coker}{{\operatorname{coker}}}

\makeatletter
 \@namedef{subjclassname@2020}{\textup{2020} Mathematics Subject Classification}
\makeatother

\begin{document}

\title{Atiyah-Bott localization in equivariant Witt cohomology}

\date{ \today}

\author[M.~Levine]{Marc~Levine}
\address{Universit\"at Duisburg-Essen,
Fakult\"at Mathematik, Campus Essen, 45117 Essen, Germany}
\email{marc.levine@uni-due.de}


\setcounter{tocdepth}{1}

\begin{abstract} Let $N$ be the normalizer of the diagonal torus $T_1\cong \G_m$ in $\SL_2$. We prove localization theorems for $\SL_2^n$ and $N^n$ for equivariant cohomology with coefficients in the (twisted) Witt sheaf, along the lines of the classical   localization theorems for equivariant cohomology for a torus action. We also have an analog of the Bott residue formula for $\SL_2^n$ and $N$. In the case of an $\SL_2^n$-action, there is a rather serious restriction on the orbit type. For an $N$-action, there is no restriction for the localization result, but for the Bott residue theorem,  one requires a certain type of decomposition of the fixed points for the $T_1$-action, which is always available if the subscheme of $T_1$ fixed points has dimension zero. \end{abstract}

\thanks{The author was partially supported by the DFG through the grant  LE 2259/7-2 and by the ERC through the project QUADAG.  This paper is part of a project that has received funding from the European Research Council (ERC) under the European Union's Horizon 2020 research and innovation programme (grant agreement No. 832833).\\ 
\includegraphics[scale=0.08]{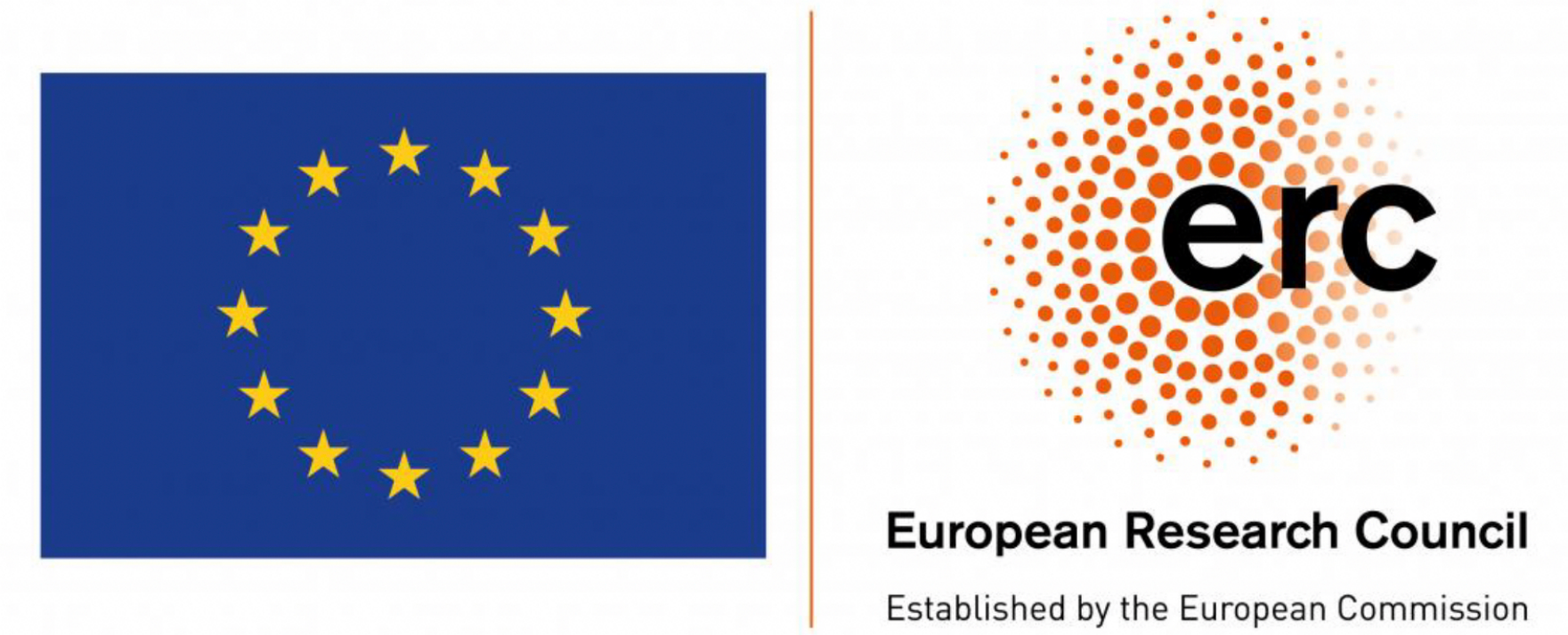}}

\maketitle

\tableofcontents

\section*{Introduction}  Localization to the fixed point set in equivariant cohomology and the Bott residue theorem give  very effective tools for computing cohomology and characteristic classes for manifolds with an action of a finite group or a torus. See for example \cite{AtiyahBott}, \cite{BerlineVergne}, \cite{Borel} \cite{Bott}, \cite{Quillen}. In the algebraic setting, this yields similar computations for Chow rings and $\CH^*$-valued Chern classes; this has been developed by  Edidin-Graham \cite{EGLoc} and Thomason \cite{Thomason2} (for algebraic $K$-theory). As part of a larger program for refining such invariants to motivic cohomology theories on schemes over more general bases, yielding  for example invariants in the Grothendieck-Witt ring of the base scheme, we would like to add such localization methods to the toolbox.

The simplest theory that yields such quadratic information is giving by the cohomology of the sheaf of Witt rings. When working over a field $k$, this will yield invariants in the Witt ring $W(k)$; combined with the classical numerical invariants in $\Z=\CH^0(k)$, this lifts canonically to invariants in the Grothendieck-Witt ring $\GW(k)$. One can often achieve this type of ``glued'' invariant directly by working in the Chow-Witt ring (see e.g. \cite{Fasel}), but it turns out that the very fact that the Chow-Witt ring mixes quadratic and classical invariants often makes for rather  complicated formulas in the concrete applications. For this reason, we will concentrate our efforts in this paper on the cohomology of the Witt sheaf. 

However, one cannot use a torus action to achieve a Witt ring version of the classical localization theorems. The reason for this is very simple:   localization  for the action of a torus $T\cong \G_m^n$ is based on the isomorphism
\[
H^*_{BT}(pt,\Z)\cong \Z[t_1,\ldots, t_n],
\]
where $t_i=c_1(L_i)$ and $L_i$ is the line bundle on $BT$ corresponding to the character $\chi_i{\colon}T\to \G_m$ given by the $i$th projection. In the algebraic setting, one uses instead the $T$-equivariant Chow groups of Totaro \cite{Totaro} and Edidin-Graham \cite{EGEquivIntThy}, with similar result
\[
\CH^*(BT)=\Z[t_1,\ldots, t_n].
\]
Localization to the fixed point subscheme of a $T$-action on a $k$-scheme $X$ requires inverting a non-zero element $P_X\in \CH^d(BT)$ for some $d\ge1$. The Witt-sheaf analog would require inverting some $P_X\in H^d(BT, \sW)$ for some $d\ge1$, but the Witt-sheaf cohomology of $BT$ is given by
\[
H^d(BT, \sW)=\begin{cases} W(k)&\text{ for }d=0,\\0&\text{ else,}\end{cases}
\]
and thus the localization method fails. 

One can use instead the group-scheme $\SL_2^n$. Ananyevskiy \cite[Introduction]{Anan} has computed $H^*(\BSL_2^n,\sW)$ as
\[
H^*(\BSL_2^n,\sW)=W(k)[e_1,\ldots, e_n],
\]
where $e_i\in H^2(\BSL_2^n,\sW)$ is the Euler class of the rank 2 bundle corresponding to the representation given by the $i$th projection $\SL_2^n\to \SL_2$ followed by the usual inclusion $\SL_2\subset \GL_2$. This is  perfectly analogous to the Chow ring of $B\G_m^n$ and one can hope for a corresponding localization theorem. Unfortunately, the situation is not so simple and one finds some rather strong restrictions for an $\SL_2^n$-action on a $k$-scheme $X$ that are necessary to permit an analog of the classical localization results. We now describe in more detail the situation for an $\SL_2^n$-localization theorem.

For a group-scheme $G$ over $k$, let $\Sch^G/k$ denote the category of quasi-projective $k$-schemes with a $G$-action. For $X\in \Sch^G/k$, the $G$-fixed subscheme $X^G$ is the maximal closed subscheme of $X$ on which $G$ acts by the identity.

 For $X\in \Sch^G/k$ with a $G$-linearized invertible sheaf $\sL$,  there are well-defined equivariant cohomology and Borel-Moore homology for the  $\sL$-twisted sheaf of Witt groups,
\[
H^*_G(X, \sW(\sL)),
H^\BM_{G,*}(X, \sW(\sL));
\]
see  Definition~\ref{def:EquivCoh}  and Definition~\ref{def:BMWIttHomology}.

\begin{IntroDef}  For $X$ in $\Sch^{\SL_2^n}/k$, we call the $\SL_2^n$-action {\em localizing} if the following conditions hold:\\[2pt]
Let $O\subset X$ be a $\SL_2^n$-orbit, with quotient $\SL_2^n\backslash O=\Spec k_O$, that is, $k_O$ is the field of $\SL_2^n$-invariant functions on $O$. We consider $O$ as an object in $\Sch^{\SL_2^n}/k_O$. Then there is a morphism $O\to (\SL_2^n\times_kk_O)/H$ in $\Sch^{\SL_2^n}/k_O$, where $H\subset \SL_2^n\times_kk_O$ is a closed subgroup-scheme of one of the two following types:\\[2pt]
i) $H$ is a maximal parabolic subgroup-scheme\\
ii) $H$ is $k_O$-conjugate to a ``diagonal'' subgroup-scheme $G_{ij}\subset \SL_2^n$,
\[
G_{ij}=\{(g_1,\ldots, g_n)\in \SL_2^n\mid g_i=g_j\}
\]
for some $i\neq j$.
\end{IntroDef}

\begin{IntroThm}[Localization for $\SL_2^n$: Theorem~\ref{thm:MainSLLoc}]\label{IntroThm1} Let $k$ be a field of characteristic $\neq 2$ and let $X$ be in $\Sch^{\SL_2^n}/k$. Let $\sL$ be a $\SL_2^n$-linearized invertible sheaf on $X$ and let $i{\colon}X^{\SL_2^n}\to X$ be the inclusion of the $\SL_2^n$-fixed points. Suppose that the $\SL_2^n$-action is localizing. Let
\[
e_*=\prod_{i=1}^ne_i \cdot\prod_{1\le i<j\le n}e_i-e_j\in H^{*\ge2}(\BSL_2^n, \sW)
\]
Then  the push-forward map
\[
i_*{\colon}H^\BM_{\SL_2^n, *}(X^{\SL_2^n}, \sW(i^*\sL))[1/e_*]\to H^\BM_{\SL_2^n, *}(X, \sW(\sL))[1/e_*]
\]
is an isomorphism. Moreover
\[
H^\BM_{\SL_2^n, *}(X^{\SL_2^n}, \sW(i^*\sL))\cong 
H^\BM_{*}(X^{\SL_2^n}, \sW(i^*\sL))\otimes_{W(k)}W(k)[e_1,\ldots, e_n]
\]
\end{IntroThm}

\begin{IntroEx} 1. Let $\SL_2^n$ act on $\P^{2n-1}$ via the block-diagonal embedding 
$\SL_2^n\hookrightarrow \SL_{2n}$ and the action of $\SL_{2n}$ on $\P^{2n-1}$ through its standard action on $\A^{2n}$. Then this action is localizing, as is the induced action on the Grassmann varieties $\Gr(r, 2n)$.  
\\[5pt]
2. Letting $F$ be the rank two representation of $\SL_2$  corresponding to the inclusion $\SL_2\subset \GL_2$, the $\SL_2$-action on $\P(\Sym^2F)$ is {\em not} localizing. Moreover, the $\SL_2$-fixed point locus is empty, and $H^*_{\SL_2}(\P(\Sym^2F), \sW)$ is a free $H^*(\BSL_2,\sW)$-module (see \cite[Lemma 5.2]{LevineBG}), so localization fails in this case.
\end{IntroEx}

The situation is somewhat nicer for an action by $N^n$, where $N\subset \SL_2$ is the normalizer of the diagonal torus $T_1\subset \SL_2$. $N$ is generated by the diagonal torus $T_1\subset\SL_2$ and an additional element $\sigma$
\[
\sigma=\begin{pmatrix} 0&1\\-1&0\end{pmatrix}.
\]
We have computed $H^*(BN, \sW)$ in \cite{LevineBG}. This is slightly larger than $H^*(\BSL_2, \sW)=W(k)[e]$, but agrees with $H^*(\BSL_2, \sW)$ in positive degrees and  the two are isomorphic after inverting $e$; similarly, 
\[
H^*(BN^n, \sW)[1/\prod_{i=1}^n e_i]= H^*(\BSL_2^n, \sW)[1/\prod_{i=1}^n e_i]=
W(k)[e_1^{\pm1},\ldots, e_n^{\pm 1}],
\]
so $N^n$ is another suitable candidate for a localization theorem. Here we have a close analogy with the classical case, although there is a minor complication in localizing to the $N^n$-fixed points, arising from the fact that $N^n$ is not connected. Recognizing this, our first localization theorem localizes the Borel-Moore homology of a scheme $X$ with $N^n$ action to the $T_1^n$-fixed points. In the statement below, we use the canonical ring homomorphism
\[
\Z[e_1,\ldots, e_n]\to H^*(BN^n, \sW)
\]
sending $e_i$ to the image of $e_i\in H^*(\BSL_2^n, \sW)$ under the natural map  $H^*(\BSL_2^n, \sW)\to H^*(BN^n, \sW)$.

\begin{IntroThm}[Localization for $N^n$: Theorem~\ref{thm:NnLocalization}]\label{IntroThm2} Let $k$ be a field of characteristic $\neq2$ and let $X$ be in $\Sch^{N^n}/k$. Let $\sL$ be an $N^n$-linearized invertible sheaf on $X$ and let $i{\colon}X^{T_1^n}\to X$ be the inclusion of the $T_1^n$-fixed points. Then there is a non-zero homogeneous element $P_X\in \Z[e_1,\ldots, e_n]$ such that, after inverting the image of $P_X$ in $H^*(BN^n, \sW)$, the push-forward map
\[
i_*{\colon}H^\BM_{G, *}(X^{T_1^n}, \sW(i^*\sL))\to H^\BM_{G, *}(X, \sW(\sL))
\]
is an isomorphism. 
\end{IntroThm}

For the case $n=1$, we have a finer result, under an additional assumption.

\begin{IntroDef}\label{def:StrictIntro}   Let $X$ be a quasi-projective $k$-scheme with $N$-action. Let $\bar\sigma$ be the image of $\sigma$ in $N/T_1$. We let $X^{T_1}_N\subset X^{T_1}$ denote the union of the irreducible components  $Z\subset X^{T_1}$ such that  $\bar\sigma\cdot Z=Z$. Let $X^{T_1}_\Ind$ be the union of the irreducible components $Z\subset X^{T_1}$ such that $(\bar\sigma\cdot Z)\cap Z=\0$. 
 \\[5pt]
1. We call the $N$-action {\em semi-strict} if 
$X^{T_1}_\red=X^{T_1}_N\cup X^{T_1}_\Ind$.
\\[5pt]
2. If the $N$-action on $X$ is semi-strict, we say the $N$-action on $X$ is {\em strict} if 
$X^{T_1}_N\cap X^{T_1}_\Ind=\0$ and  we can decompose $X^{T_1}_N$ as a disjoint union of two $N$-stable  closed  subschemes 
\[
X^{T_1}_N=X^N\amalg X^{T_1}_\fr 
\]
where the $N/T_1$-action on $X^{T_1}_\fr$ is free.   
 \end{IntroDef}
 
 \begin{IntroEx}\label{IntroEx:SStrictEtc} Let $X\subset \A^2_k$ be the closed subscheme $xy=0$, with trivial $T_1$-action and with $\bar\sigma$ acting by $\bar\sigma(x,y)=(y,x)$. Then $X^{T_1}_N=\0=X^{T_1}_\Ind$, and  the $N$-action is not semi-strict.

Let $X=\A^n_k$, with $T_1$ acting trivially and with $\bar\sigma$ acting by $\bar\sigma(t)=-t$; in this case $X=X^{T_1}_N$, $X^{T_1}_\Ind=\0$, $X^N=\{0\}$. Moreover,   $\bar\sigma$ acts freely on $X\setminus X^N$, but this is not a closed subscheme of $X^{T_1}_N$. Thus, the $N$-action is semi-strict, but not strict.
 
 If  $k(\sqrt{a})/k$ is a non-trivial degree two extension of $k$, and $X=\Spec k(\sqrt{a})$ with $T_1$ acting trivially and $\bar\sigma$ acting by conjugation over $k$, then the action is strict, with $X=X^{T_1}_\fr$, and $X^{T_1}_\Ind=X^N=\0$. 
 
 If $X$ is a quasi-projective $k$ scheme with $N$-action and $X^{T_1}$ has dimension zero, then the $N$-action is strict. Moreover, for $x\in X^{T_1}_\fr$, either $N\cdot x=x\amalg x'$ with $\sigma(x)=x'$, or $\sigma(x)=x$ and $\sigma^*$ acts on $k(x)$ by an involution, identifying $k(x)$ with $k(x)^\sigma(\sqrt{a})$ for some $a\in (k(x)^\sigma)^\times$, not a square, and with $\sigma^*(\sqrt{a})=-\sqrt{a}$.
  \end{IntroEx}

In the case of a semi-strict $N$-action, we have
\begin{IntroThm}[Localization for an $N$-action: Theorem~\ref{thm:AtiyahBottLocalizationN}]\label{IntroThm3} Let $k$ be a field of characteristic $\neq2$ and let $X$  be in $\Sch^N/k$. Suppose the $N$-action is semi-strict. Let $\sL$ be an $N$-linearized invertible sheaf on $X$ and let $i{\colon}X^{T_1}_N \to X$ be the inclusion. Then there is an integer $M>0$ such that the push-forward map
\[
i_*{\colon}H^\BM_{N, *}(X^{T_1}_N, \sW(i^*\sL))[1/Me]\to H^\BM_{G, *}(X, \sW(\sL))[1/Me]
\]
is an isomorphism. 
\end{IntroThm}

\begin{IntroRemark}\label{IntroRem:SemiStrict} The reader will have noticed that the integer $M$ in this last result might be even, in which case $M$ might be a nilpotent in $W(k)$ and the theorem tells us nothing.  However, if the base-field $k$ admits a real embedding, this will not be the case, even if one would lose 2-primary information by inverting $M$. Similar remarks hold for the polynomial $P_X$ in Theorem~\ref{IntroThm2}.
\end{IntroRemark}

We have an analog of the Bott residue formula in the case of a localizing  $\SL_2^n$-action. In case $k$ has characteristic zero, let $Y$ be a reduced connected $k$-scheme with trivial $\SL_2^n$-action and let $V$ be a $\SL_2^n$-linearized vector bundle on $Y$. Then there is a ``generic representation type'' $[V^{gen}]$, defined as an isomorphism class of a $\SL_2^n$-representation over $k$, and an associated generic Euler class $e^{gen}_{\SL_2^n}\in H^*(\BSL_2^n, \sW)$, well-defined up to multiplication by a unit in $W(k)$. 

\begin{IntroThm}[Bott Residue Theorem for $\SL_2^n$: Theorem~\ref{thm:BottResidue}] Let $G=\SL_2^n$ and suppose $k$ has characteristic zero. Take $X\in \Sch^{\SL_2^n}/k$ with $X$ reduced,  and let $\sL$ be  a $\SL_2^n$-linearized invertible sheaf on $X$. We suppose that the action is localizing, and that the inclusion $i{\colon}X^{\SL_2^n}\to X$ is a regular embedding. Let $i_j{\colon}X^{\SL_2^n}_j\to X$, $j=1,\ldots, s$,  be the connected components of $X^{\SL_2^n}$, with normal bundle $N_{i_j}$. Suppose in addition that $e^{gen}_{\SL_2^n}(N_{i_j})\neq 0$ for each $j$.  Let $P$ be as in Theorem~\ref{IntroThm1}. Then 
\[
e_{\SL_2^n}(N_{i_j})\in H^*_{\SL_2^n}(X^{\SL_2^n}_j,\sW(\det^{-1} N_{i_j}))[(Pe^{gen}_{\SL_2^n}(N_{i_j}))^{-1}]
\]
is invertible. Letting $e^{gen}_{\SL_2^n}(N_i) =\prod_{j=1}^s e^{gen}_{\SL_2^n}(N_{i_j})$,  the inverse of the isomorphism 
\[
i_*{\colon}H^\BM_{G*}(X^{\SL_2^n},\sW(i^*\sL))[(Pe^{gen}_{\SL_2^n}(N_i))^{-1}]\xrightarrow{\sim} H^\BM_{G*}(X,\sW(\sL)[Pe^{gen}_{\SL_2^n}(N_i))^{-1}] 
\]
is the map 
\[
x\mapsto \prod_{j=1}^s i_j^!(x)\cap e_{\SL_2^n}(N_{i_j})^{-1},
\]
where we use the evident isomorphism
\[
H^\BM_{G*}(X^{\SL_2^n},\sW(i^*\sL))=\prod_{j=1}^sH^\BM_{G*}(X_j^{\SL_2^n},\sW(i_j^*\sL)).
\]
\end{IntroThm}
 
In the case of a reduced $X\in \Sch^N/k$, we need to assume that the action is semi-strict, and that $i{\colon}X^{T_1}_N\to X$ is a regular embedding. 

\begin{IntroThm}[Bott Residue Theorem for $N$: Theorem~\ref{thm:BottResidue}] Let $k$ be a field of characteristic $\neq2$. Take $X\in \Sch^N/k$ with $X$ reduced, and let $\sL$ be an $N$-linearized invertible sheaf on $X$. We suppose the $N$-action is semi-strict, and the inclusion $i{\colon}X^{T_1}_N\to X$ is a regular embedding.  For each connected component $i_j{\colon}(X^{T_1}_N)_j\to X$ of $X^{T_1}_N$, let $N_{i_j}$ be the normal bundle of $i_j$   Then there is an $M>0$ such that 
\[
e_N(N_{i_j})\in H^*_N((X^{T_1}_N)_j,\sW)[1/Me]
\]
is invertible,  and $i_*$ defines an isomorphism
\[
i_*{\colon}H^\BM_{N, *}(X^{T_1}_N,\sW(i^*\sL))[1/Me)]\xrightarrow{\sim} H^\BM_{N*}(X,\sW(\sL)[1/Me] 
\]
with inverse the map 
\[
x\mapsto \prod_j i_j^!(x)/e_N(N_{i_j}).
\]
\end{IntroThm}

Just as in the classical case, for $G=\SL_2^n, N^n$, we have the usual description of the $G$-equivariant Witt Borel-Moore homology (or cohomology for $X$ smooth) for a trivial $G$-action,
\[
H^\BM_{G, *}(X, \sW)\cong H^\BM_*(X, \sW)\otimes_{W(k)}H^{-*}(BG, \sW).
\]
See Corollary~\ref{cor:Trivial}. 

If $X^{T_1}_N \neq X^N$, one would also want to compute $H^\BM_{N, *}(X^{T_1}_N,\sW(i^*\sL))$. By using localization for Borel-Moore homology, one can reduce to the case of a strict $N$-action, with $X^{T_1}_N=X^N\amalg X^{T_1}_\fr$, and in this case,   one  would  need to compute 
$H^\BM_{N,*}(X^{T_1}_\fr, \sW(\sL))$. In general, I found no nice answer for this, but using localization again, one should be able to reduce to the case  of dimension zero.  This is just the case of  $\Spec F(\sqrt{a})$, $F$ a field, with $\bar\sigma$ acting by the usual conjugation and $T_1$ acting trivially, as in Example~\ref{IntroEx:SStrictEtc}. This case is discussed in \S\ref{sec:TwistedQuadrics}, with the end result given in Theorem~\ref{thm:Nontrivial0DimlCase}.  We also give an explicit computation (Corollary~\ref{cor:Nontrivial0DimlCasePushforward}) of the push-forward
\[
H^*_N(\Spec F(\sqrt{a}), \sW)\to H^*_N(\Spec F, \sW)
\]
which suffices for computing the push-forward
\[
H^\BM_{N, 0}(X, \sW)\to W(k)
\]
for $X$ proper over $k$.

As to the organization of this paper, we recall in  \S\ref{sec:Quot}  some basic concepts on group actions, and review some important structural results on quotients, especially Thomason's theorem on torus actions. We assemble some useful results from representation theory in \S\ref{sec:RepThy}, and in \S\ref{sec:CohBMHom} we provide the reader with the relevant background on the ``motivic'' framework of cohomology and Borel-Moore homology of algebraic varieties, following \cite{DJK, Hoyois6, LR}, that we will be using throughout the paper. We review the equivariant theory and specialize to the case of equivariant Borel-Moore homology for the sheaf of Witt groups in \S\ref{sec:EquivBMWitt}, and in \S\ref{sec:WittCohBG} we recall computations of the Witt-sheaf cohomology of the classifying spaces $BN$ and $\BSL_2$. We then turn to our main applications, the localization theorems and residue formulas. We handle localization for the case of $\SL_2^n$ in \S\ref{sec:SLLOc}, the case of $N^n$ in \S\ref{sec:NnLoc}, and specialize to the case of an $N$-action in \S\ref{sec:NLoc}. Once we have the localization results, the corresponding residue theorems are easy to prove with the help of the excess intersection formula of \cite{DJK}; this is done in \S\ref{sec:BottRes}. 

We compute the $N$-equivariant cohomology of the scheme $\Spec k(\sqrt{a})$, with $N$ acting through its $\Z/2$-quotient $N/T_1$ by conjugation over $k$, in \S\ref{sec:TwistedQuadrics}. 

In  \cite{LevineBG}, we had computed the structure of $H^*(BN, \sW)$ as a $W(k)$-$\Z$-graded algebra, and also the structure of $H^{*\ge2}(BN, \sW(\gamma))$ as a module over $H^*(\BSL_2, \sW)$,  where   $\gamma$ is the non-trivial generator of $\Pic(BN)\cong \Z/2$; there we had incorrectly stated this as a computation of the entire cohomology $H^{*}(BN, \sW(\gamma))$, ignoring the contribution of $H^0(BN, \sW(\gamma))$\footnote{This was later corrected and extended to more general cohomology theories in the thesis of A. D'Angelo \cite{DAngelo}.}. We take the opportunity here in \S\ref{sec:SpectSeqEulerClass} to correct this error and to give a complete description of the $\Z\times \Z/2$-graded $W(k)$-algebra  $H^*(BN, \sW)\oplus H^*(BN, \sW(\gamma))$ (Theorem~\ref{thm:CohBN}),  The $H^*(\BSL_2,\sW)$-module structure suffices for most applications, but we compute the $W(k)$-algebra structure for the sake of completeness and general interest. Some results from this section are needed in 
\S\ref{sec:TwistedQuadrics}, and some of the definitions and notations from \S\ref{sec:TwistedQuadrics} are used in \S\ref{sec:SpectSeqEulerClass}, but we avoid any circularity of the arguments.

Our description of $H^*_N(\Spec k(\sqrt{a}), \sW)$ and  $H^*(BN, \sW)\oplus H^*(BN, \sW(\gamma))$ is in terms of explicit generators and relations.  These final two sections  involve detailed computations and as such are rather technical, but they also contain numerous computational methods that we hope will find use elsewhere.

I would like to thank Philippe Gille for his very helpful remarks on $\SL_2$-homo\-gen\-eous spaces and for supplying the statement and proof of Lemma~\ref{lem:LocalTrivGRep}. Thanks are also due to Detlev Hoffmann and Raman Parimala for telling me about Lam's exact triangle, used in the proof of Proposition~\ref{prop:EquivWittCohBarXa}. Many thanks as well to Aravind Asok for his very helpful comments on $\pi_1^{\A^1}(\SL_2)$, which guided many of the computations in \S\ref{sec:SpectSeqEulerClass},  to Lorenzo Mantovani for explaining the theory of equivariant Borel-Moore homology as developed by him and A. Di Lorenzo \cite{DiLorenzoMantovani}, and to Jochen Heinloth for his  explanations about the representability of quotients. Discussions with Alessandro D'Angelo clarified many points of equivariant cohomology and Borel-Moore homology for $\SL$-oriented, $\eta$-local theories. Finally, I am grateful to the referee for their very helpful comments and suggestions.

\section{Quotients and fixed points}\label{sec:Quot}
In this section, we work over a noetherian separated base-scheme $B$ of finite Krull dimension and let $G$ be a linear algebraic group scheme over $B$, flat over $B$, with a fixed embedding  $i{\colon}G\to \GL_N/B$ realizing $G$ as a closed subgroup-scheme of some $\GL_N/B$.  

We let $\Sch^G/B$ denote the category of $G$-quasi-projective $B$-schemes, that is, an object is a quasi-projective $B$-scheme $X$ with a $G$-action $G\times X\to X$ that extends to a $G$-action $G\times \P(V)_B\to \P(V)_B$ with respect to some vector bundle $V\to B$ with a $G$-action over the trivial action on $B$,  together with a locally closed immersion $i{\colon}X\to \P(V)_B$; that is, the $G$-action on $X$ admits a $G$-linearized very ample invertible sheaf, or in other words, $X$ admits a $G$-linearization.  Morphisms are $G$-equivariant morphisms of $B$-schemes (not necessarily respecting the $G$-linearizations).
We write $\Sm^G/B$ for the full subcategory of $\Sch^G/B$ with objects those $X$ that are smooth over $B$, and write  $\Sch/B$, resp., $\Sm/B$  for $\Sch^{\id}/B$, resp. $\Sm^{\id}/B$ .   We write $\Sch_B$ for the category of finite-type $B$-schemes.  

Given a vector bundle $V$ or coherent sheaf $\sF$ on some $X\in \Sch^G/B$, we refer to a $G$-action on $V$ or $\sF$ over the given $G$-action on $X$ as a {\em $G$-linearization} of $V$ or $\sF$, and refer to $V$ or $\sF$ with such an action as $G$-linearized.

We recall some definitions from \cite{Fogarty}.
For $X\in \Sch^G/B$, we say that $G$ acts trivially on $X$ if the action morphism $G\times_BX\to X$ is the projection to $X$.  For general $X\in \Sch^G/B$, we have the $G$-fixed subscheme $X^G\hookrightarrow X$, this being the maximal $G$-stable closed subscheme $Y$ of $X$ on which $G$ acts trivially. 

We recall some results due to Thomason \cite{Thomason}. Recall that a {\em geometric quotient} $q{\colon}X\to G\backslash X$ for some $X\in \Sch^G/B$ is a $G$-equivariant  morphism in $\Sch_B$, where $G$ acts trivially on $G\backslash X$, such that $q$ represents the quotient as fppf sheaf (i.e., $q$ is a categorical quotient in $\Sch_B$), and such that each geometric fiber of $q$ is a single orbit (after passing to the associated reduced closed subscheme).

\begin{proposition}\label{prop:Quotient}  Let $X\in \Sch^G/B$ be reduced.\\[5pt]
1.  There is a decomposition of $X$ into locally closed  reduced $G$-stable subschemes, $X=\amalg_{i=1}^r X_i$, such that for each $i$ the quotient $q_i{\colon}X_i\to G\backslash X_i$ (as an fppf sheaf) exists, with $G\backslash X_i\in\Sch/B$. Moreover  $X_i\to G\backslash X_i$ is a geometric quotient as well. If $G$ is smooth over $B$, then we may take the decomposition as above with each $q_i$ smooth.  \\[2pt]
2. Suppose that $G$ is a split torus over $B$, $G\cong \G_m^n/B$. Then  there is a decomposition of $X$ as in (1) such that for each $i$, $X_i$ is affine, and there is a diagonalizable subgroup-scheme $T'_i$ of $G$, acting trivially on $X_i$, such that the quotient torus $T''_i:=G/T'_i\cong \G_m^{n_i}/B$ acts freely on $X_i$, making $X_i$ into a  trivial  $T''_i$-torsor over $G\backslash X_i=T_i''\backslash X_i$. If $B$ is excellent, we may take the decomposition so that in addition, $X_i$  is regular.
 \end{proposition}

\begin{proof} For (1), it follows from \cite[Proposition 4.7]{Thomason} that there is a non-empty open $G$-stable subscheme $U$ of $X$ such that the  fppf quotient sheaf $G\backslash U$ is  represented by a separated scheme of finite type over $B$, with quotient map $q{\colon}U\to G\backslash U$ smooth if $G$ is smooth over $B$. By \cite[Proposition 1(2), Th\'eor\`eme 3(1)]{Raynaud}, this implies that  $q{\colon}U\to G\backslash U$ is a geometric quotient.  Replacing $G\backslash U$ with a non-empty  affine open subscheme $V$ and replacing $U$ with $q^{-1}(V)$, we may assume that $G\backslash U$ is affine, hence in $\Sch/B$. Replacing $X$ with the complement $X\setminus U$ with its induced $G$-action, the result follows by noetherian induction. 

Part (2) follows from Thomason's generic slice theorem for a torus action, \cite[Proposition 4.10 and Remark 4.11]{Thomason}.
\end{proof}

\begin{remark}  Suppose  $B=\Spec k$, and that $G$ is linearly reductive over $k$. Take $X\in \Sm^G/k$.  By \cite[Theorem 5.4]{Fogarty},  $X^G$ is smooth over $k$. Examples include $G=\SL_2^n$ and $k$ of characteristic zero, or $G=\G_m^n$, and $k$ arbitrary. Letting $N\subset \SL_2$ be the normalizer of the diagonal torus, $N^n$ is linearly reductive over $k$ if $\Char k\neq 2$.
\end{remark}

\begin{definition} Let $G$ be a smooth,  finite type  group-scheme over  $B$, and let $X\in \Sm^G/B$ be a smooth $B$-scheme with $G$-action $\rho{\colon}G\times_B X\to X$. Consider $G\times_B X$ as a group scheme over $X$ via the projection,   acting on the $X$-scheme
\[
p_2{\colon}X\times_BX\to X
\]
via 
\[
\rho_\delta{\colon}G\times_BX\to X\times_BX,\ \rho_\delta(g,x)=(gx,x).
\]
Let $G_X\subset G\times_BX$ be the isotropy group-scheme of the diagonal section $\Delta_X\subset X\times_BX$.  We call $X$ a {\em homogeneous space} for $G$ (or a $G$-homogeneous space) if the induced morphism
\[
\bar\rho_\delta{\colon}G\times_BX/G_X\to X\times_BX
\]
is an isomorphism. 
  \end{definition}

\begin{remark} Let $F$ be a field, and $G$ an affine group-scheme over $F$, and  let 
$X\in\Sm^G/F$ be a $G$-homogeneous space. 

Let $F_X\subset \Gamma(X, \sO_X)$ be the subring of $G$-invariant functions on $X$,
\[
F_X:=\Gamma(X, \sO_X)^G.
\]
Then $F_X$ is a finite separable field extension of $F$ and $G\backslash X=\Spec F_X$. Indeed, if $X=G/H$ for some closed subgroup-scheme $H$ of $G$, then clearly $G\backslash X=G\backslash (G/H)=G\backslash G=\Spec F$ and  $F_X=F$ in this case. In general, since $X$ is smooth over $F$, $X$ has an $F'$-point $x$ over some finite separable extension $F'$ of $F$. Then $X_{F'}\cong G_{F'}/G_x$, where $G_x$ is the isotropy group of  $x$. Thus, $G_{F'}\backslash X_{F'}=\Spec F'$, so $F'\supset F_X\supset F$ and $F_X$ is thus a finite separable extension of $F$. 

Since $X$ is a $G$-homogeneous space,  every $G$-equivariant map $f{\colon}X\to Y$ in $\Sch_F$, where $G$ acts trivially on $Y$, must factor through an irreducible, zero-dimensional closed subscheme $y$ of $Y$. Again, this is clear if $X=G/H$. Letting $F'$ be as above, we have the surjective map $X_{F'}\to X$, and $f$ induces the $G_{F'}$-equivariant map $f'{\colon}X_{F'}\to Y_{F'}$, which must factor through a map $\Spec F'\to Y_{F'}$, so the image of $f$ is contained in the image of $\Spec F'\to Y_{F'}\to Y$.

In particular, $f$ factors through an affine neighborhood $U=\Spec A$ of $y$ in $Y$, and is thus induced by an $F$-algebra map $A\to F_X$.  This shows that  $X\to \Spec F_X$ is the categorical quotient of $X$ by $G$, i.e., $G\backslash X=\Spec F_X$.

To conclude, we note that $X$ is smooth over $F_X$.
\end{remark}

\begin{definition}\label{def:Orbit}  Let $G$ be a smooth,  finite type  group-scheme over  $B$ and take $X\in \Sch^G/B$ with structure morphism $\pi{\colon}X\to B$. Form a locally closed decomposition $X=\amalg_{i=1}^r X_i$  
as in Proposition~\ref{prop:Quotient} and let $q_i{\colon}X_i\to G\backslash X_i$ be the corresponding quotient map. For each $y\in G\backslash X_i$, we have the subscheme $q_i^{-1}(y)$ of $X$, with its induced $G_{k(y)}$-action, defining an object $q_i^{-1}(y)\in \Sm^{G_{k(y)}}/k(y)$; we refer to 
$q_i^{-1}(y)\in \Sch^{G_{k(y)}}/k(y)$ as an {\em orbit} of $G$ in $X$. For $O:=q_i^{-1}(y)$, we denote the subfield $q_i^*k(y)\subset \Gamma(O, \sO_O)$ by $k_O$,   giving us the orbit $O\in \Sm^{G_{k_O}}/k_O$.
\end{definition}
We note that the notion of orbit $O$ and corresponding field $k_O$ are  independent of the choice of decomposition $\amalg_{i=1}^r X_i$ of $X$. Also, an orbit $O\in \Sm^{G_{k_O}}/k_O$ is a homogeneous space for $G_{k_O}$. 

\section{Some representation theory}\label{sec:RepThy}

\begin{lemma}\label{lem:LocalTrivGRep} Let $G$ be an affine linearly reductive algebraic group scheme over a field $k$. We suppose  that $G$ is either reductive and $k$-split, or, for  each finitely generated extension field $L$ of $k$, with separable closure $L^{sep}$,  that each $L^{sep}$-representation of $G$ is $L^{sep}$-conjugate to a $k$-representation of $G$.   Let $\sO$ be a local $k$-algebra, essentially of finite type over $k$ and let $f{\colon}G_\sO\to \GL_n/\sO$ be an $\sO$-representation. Then there is a $k$-representation $f_0{\colon}G\to \GL_n/k$ such that  $f$ is $\sO$-conjugate to the base-extension of $f_0$. Moreover, $f_0$ is unique up to $k$-conjugation.
\end{lemma}

I am indebted to Philippe Gille for the argument.

\begin{proof} Let $L$ denote the residue field of $\sO$, with separable closure $L^{sep}$, and let $\bar{f}{\colon} G_L\to \GL_n/L$ be the $L$-representation given by the reduction of $f$. In the case of a $k$-split, reductive group scheme, the  theory of Tits algebras \cite{Tits}  shows that $\bar{f}$ is $L^{sep}$-conjugate to a $k$-representation $f_0{\colon}G\to \GL_n/k$; if $G$ is not of this type, then this follows from our assumption on $G$.

Let $C_{f,f_0}^\sO$ be the transporter scheme over $\sO$, 
\[
C_{f,f_0}^\sO(R)=\{g\in \GL_n(R)\mid gfg^{-1}=f_0\}
\]
and define the group-scheme $C^\sO_{f_0, f_0}$ over $\sO$ similarly. Let $\sO^{sh}$ be the strict henselization of $\sO$. Margaux's theorem  \cite[Corollary 4.9]{Margaux} implies that $f$ is $\sO^{sh}$-conjugate to $f_0$. This shows that $C^\sO_{f,f_0}$ is a $C^\sO_{f_0,f_0}$-torsor for the \'etale topology, so is classified by an element of $H^1_\et(\sO, C^\sO_{f_0,f_0})$. By Schur's lemma, $C^\sO_{f_0,f_0}$ is a product of Weil restrictions of group-schemes $\GL_{n_i}$, so $H^1_\et(\sO, C^\sO_{f_0,f_0})=\{*\}$. Thus 
$C^\sO_{f,f_0}$ admits a section over $\sO$, showing that $f$ is $\sO$-conjugate to the base-extension of $f_0$ to $\sO$. 

Applying this argument to the $k$-representation $f_0$ gives the uniqueness: if  $\bar{f}$ is $L^{sep}$-conjugate to a second $k$-representation $f_0'$, then $C^k_{f_0',f_0}$  is a trivial $C^k_{f_0,f_0}$-torsor, so $f_0$ and $f_0'$ are $k$-conjugate.
\end{proof} 

\begin{definition} Take $Y\in \Sch^G/k$ with trivial $G$-action, and let $\sV$ be a $G$-linearized locally free coherent sheaf on $Y$. We say that $\sV$ is {\em $G$-trivialized} on some open subscheme $j{\colon}U\to Y$ if there is a $k$-representation $V(f_0)$,  $f_0{\colon}G\to \GL_n$,  and an isomorphism $\phi{\colon}j^*\sV\xrightarrow{\sim}\sO_U\otimes_k V(f_0)$. We call $U$ a $G$-trivializing open subscheme and $\phi$ a $G$-trivialization. 
\end{definition}

\begin{remark}\label{rem:G-Triv}
It follows from Lemma~\ref{lem:LocalTrivGRep} that if $G$ satisfies the hypotheses of that Lemma, then each point $y\in Y$ admits a $G$-trivializing open neighborhood $ U\ni y$. 
\end{remark}

Let $N\subset \SL_2$ be the normalization of the diagonal torus $T_1\subset \SL_2$ 
\[
T_1(R)=\{\begin{pmatrix}t&0\\0&t^{-1}\end{pmatrix}\mid t\in R^\times\}
\]
The group scheme $N\subset \SL_2$ is a $\Z/2$-extension of   $T_1\cong \G_m$:
\[
1\to T_1\to N\to \Z/2\to 0,
\]
and the $1\in \Z/2$ lifts to the order four element
\[
\sigma:=\begin{pmatrix}0&1\\-1&0\end{pmatrix}\in N.
\]
We let $\bar\sigma$ denote the image of $\sigma$ in $N/T_1$.

Suppose that $\Char k\neq 2$. Then $N$ is linearly reductive and  has the family of irreducible representations $\rho_m{\colon}N\to \GL_2$, $m=1,2,\ldots$, defined by
\[
\rho_m(\begin{pmatrix}t&0\\0&t^{-1}\end{pmatrix})=
\begin{pmatrix}t^m&0\\0&t^{-m}\end{pmatrix}
\]
\[
\rho_m(\sigma)=\begin{pmatrix}0&1\\(-1)^m&0\end{pmatrix}
\]
We let  $\rho_0{\colon}N\to \G_m$ denote the trivial character and $\rho_0^-{\colon}N\to \G_m$ the sign character, $N\to N/T_1\cong \{\pm1\}$.   

The set of representations $\{\rho_m\mid m\ge0\}\cup \{\rho_0^-\}$ gives a complete set of representatives for all irreducible $k$-representations of $N$.  Each $\rho_m$ gives the bundle $\tilde{\sO}(m):=V(\rho_m)\times^NEN$ on $BN$, where $V(\rho_m)$ is the vector space of the representation $\rho_m$.  We define $\tilde{\sO}(0)^-:= V(\rho_0^-)\times^NEN$, $\tilde{\sO}(0):= V(\rho_0)\times^NEN$.

For $m>0$, the basis $e_1\wedge e_2$ of $V(\rho_m)$ defines an isomorphisms $\det\rho_m\cong \rho_0$, $\det\tilde{\sO}(m)\cong \tilde{\sO}(0)$ for $m$ odd and $\det\rho_m\cong \rho^-_0$,  $\det\tilde{\sO}(m)\cong \tilde{\sO}(0)^-$ for $m$ even.   

\begin{remark}\label{rem:NnRep}  For $k$ a field of characteristic $\neq2$, the group $N^n$ is linearly reductive and the classification of the irreducible representations of $N$ described above shows that the hypotheses of Lemma~\ref{lem:LocalTrivGRep} are satisfied. 
\end{remark}

We extend the notion of an $N$-trivialization to the case where only $T_1\subset N$ acts trivially on the scheme $Y$.

\begin{definition}\label{def:NSemiLinear} Take $Y\in \Sch^N/k$ and let $\sV$ be an $N$-linearized locally free coherent sheaf on $Y$. Suppose that $T_1\subset N$ acts trivially on $Y$, giving the action of $N/T_1\cong \Z/2$ on $Y$, with $\sigma\in N$ acting by an involution $\tau{\colon}Y\to Y$. \\[5pt]
1. We have the decomposition of $\sV$ into weight subspaces for the $T_1$-action
\[
\sV=\oplus_m\sV_m
\]
where $T_1$ acts on $\sV_m$ via the character $\chi_r(t)= t^m$. Define the {\em moving part} of $\sV$, $\sV^\mov$, by
\[
\sV^\mov:=\oplus_{m\neq0}\sV_m
\]
2. Given a $k$-representation $f{\colon}N\to \GL_n$ and $U\subset Y$ an $N$-stable open subscheme, the {\em $\sO_U$-semi-linear extension  of $f$} is the   sheaf $\sO_U\otimes_kV(f)$ on $U$, where $T_1$ acts $\sO_U$-linearly via $f$
\[
t\cdot (x\otimes v):=x\otimes f(t)(v)
\]
and $\sigma$ acts as
\[
\sigma\cdot(x\otimes v)=\tau^*(x)\otimes f(\sigma)(v).
\]
We   write this $N$-linearized sheaf as $\sO_U\otimes^{\tau,\sigma}_kV(f)$, and call it  the $\sO_U$ semi-linear extension of $f$.\\[2pt]
3.  We say that $\sV$ is $N$-trivialized on an $N$-stable open subscheme $j{\colon}U\to Y$ if there is a $k$-representation $f{\colon}N\to \GL_n$ such that $j^*\sV$ is isomorphic to the $\sO_U$-semi-linear extension  of $f$, 
\[
j^*\sV\cong \sO_U\otimes^{\tau,\sigma}_kV(f).
\]
\end{definition}

\begin{lemma}\label{lem:NLocalTriv} Suppose $k$ has characteristic $\neq2$, let $Y$, $\sV$,  and $\tau{\colon}Y\to Y$ be as in Definition~\ref{def:NSemiLinear}. Take $y\in Y$. \\[5pt]
1. There exists an $N$-stable  open neighborhood $j{\colon}U\to Y$ of $y$ and an $N$-trivialization of $j^*\sV^\mov$,  $j^*\sV^\mov\cong \sO_U\otimes^{\tau,\sigma}_kV(f^\mov)$.  Moreover,  we may take $f^\mov$ to be a sum of copies of the representations $\rho_m$, $m>0$.\\[2pt]
2. Suppose that $N/T_1$ acts freely on $Y$. Then there exists an $N$-stable  open neighborhood $j{\colon}U\to Y$ of $y$ and an $N$-trivialization of $j^*\sV$,  $j^*\sV\cong \sO_U\otimes^{\tau,\sigma}_kV(f)$.  Moreover,  we may take $f$ to be a sum of copies of the representations $\rho_m$, $m\ge0$.
\end{lemma}

\begin{proof}
Decompose $\sV$ into a finite sum of weight spaces for the $T_1$,
\[
\sV=\oplus_m\sV_m
\]
where $t$ acts on $\sV_m$ by the character $\chi_m(t)=t^m$. We can find an $N$-stable open neighborhood $U$  of $y$ with each $\sV_m\cong \sO_U^{m_r}$ as $\sO_U$-module; since $\sigma\cdot t=t^{-1}\sigma$, the action by $\sigma$ restricts to $\sO_Y$-linear isomorphisms
\[
\phi_{m,\sigma}{\colon} \tau^*\sV_m\xrightarrow{\sim}\sV_{-m}
\]
Thus we have the decomposition of $\sV$ into $N$-stable subsheaves
\[
\sV=\sV_0\oplus\oplus_{m>0}(\sV_m\oplus \sV_{-m})
\]

Suppose $m>0$. Choose a frame $e_1,\ldots, e_{r_m}$ for $j^*\sV_m$ and let $e_j':=\phi_{m,\sigma}(\tau^*e_j)$. Then $e_1',\ldots, e_{r_m}'$ gives a framing for $\sV_{-m}$, and 
$\phi_{-m,\sigma}(\tau^*e'_j)=(-1)^me_j$. This gives the isomorphism of $j^*(\sV_m\oplus \sV_{-m})$ with the $\sO_U$-semi-linear extension of $\rho_m^{n_m}$, with $T_1$ acting $\sO_U$-linearly. Putting these together gives the isomorphism of $N$-linearized sheaves on $U$
\[
j^*\sV^\mov\cong \sO_U\otimes^{\tau,\sigma}_kV(\oplus_{m>0}\rho_m^{r_m})
\]
proving (1).

We now take $m=0$; we may assume that $U=\Spec R$ is affine and that the $N/T_1$-action on $U$ is free. Let $q{\colon}U\to \bar{U}:=N\backslash U$ be the quotient map, which is an \'etale morphism of degree two, with $\bar{U}=\Spec \bar{R}$, where $\bar{R}\subset R$ is the ring of $\sigma$-invariants in $R$.

Choose a framing $e_1,\ldots, e_{r_0}$ for $j^*\sV_0$, and write
\[
\phi_{0,\sigma}(\tau^*e_j)=\sum_ia_{ij}e_i
\]
This gives us the matrix $A=(a_{ij})$ with $A\cdot A^\tau$ the $r_0\times r_0$ identity matrix. By Hilbert's Theorem 90 for $\GL_n$ with respect to the \'etale extension $\bar{R}\hookrightarrow R$, after shrinking $U$ if necessary, there is a $B\in \GL_n(R)$ with 
$A=B\cdot (B^\tau)^{-1}$. Changing our framing by $B=(b_{ij})$, by taking $f_j=\sum_ib_{ij}e_i$, we have 
\[
\phi_{0,\sigma}(\tau^*f_j)=f_j
\]
giving the isomorphism  $j^*\sV_0\cong \sO_U\otimes^{\tau,\sigma}V(\rho_0^{r_0})$. We take $f:=\oplus_{m\ge0}\rho_m^{m_r}$, proving (2).
\end{proof}

 \begin{construction}[Generic representation class]\label{const:GenEuler} Let $G$ be  group scheme over $k$, satisfying the hypotheses of Lemma~\ref{lem:LocalTrivGRep}. Let $\sV$ be a locally free coherent sheaf with $G$-linearization on some connected $Y\in \Sch^G/k$. \\[5pt]
Case 1: $G$ acts trivially on $Y$.  Let $y$ be a point of $Y$. By Remark~\ref{rem:G-Triv}, there is a $G$-stable open neighborhood $j_{U_y}{\colon}U_y\hookrightarrow Y$ of $y\in Y$ and a $G$-trivialization $\psi_y{\colon}j_{U_y}^*\sV\xrightarrow{\sim} \sO_U\otimes_kV(f)$ for some $k$-representation $f$ of $G$. Taking the canonical decomposition of $V(f)$ into isotypical components gives a corresponding decomposition of $j_{U_y}^*\sV$. By the uniqueness part of Lemma~\ref{lem:LocalTrivGRep},  the fact that $Y$ is connected implies we have a global decomposition of $\sV$ into isotypical components, indexed by the irreducible $k$-representations of $G$,
\[
\sV=\oplus_{\phi}\sV_\phi
\]
such that on each $G$-trivializing open subscheme $j_U{\colon}U\to Y$, we have
\[
j_U^*\sV_\phi\cong \sO_U\otimes_kV(\phi)^{n_\phi}
\]
where $V(\phi)$ is the irreducible $k$-representation indexed by $\phi$. 

Again using the connectivity of $Y$, we see that the $k$-representation $\oplus_\phi V(\phi)^{n_\phi}$ is uniquely determined by $\sV$, up to isomorphism of $k$-representations of $G$.  We denote the isomorphism class of the $k$-representation $\oplus_\phi V(\phi)^{n_\phi}$ of $G$ by $[\sV^{gen}]$.
\\[5pt]
\\[5pt]
Case 2. $G=N$, $T_1$ acts trivially on $Y$ and $\sV=\sV^\mov$.  In this case we first decompose $\sV$ into $T_1$ isotypical components $\sV=\oplus_{m\neq0}\sV_m$ as in the proof of Lemma~\ref{lem:NLocalTriv}(1). For each $N$-trivializing open subscheme $j_U{\colon}U\to Y$, we have the isomorphism 
\[
j_U^*\sV\cong \sO_U\otimes^{\tau,\sigma}_k\oplus_{m>0}\sO_U\otimes_k^{\tau,\sigma}V(\rho_m)^{n_m}.
\]
Since this local description follows uniquely from the decomposition into $T_1$-isotypical components, the isomorphism class of the $k$-representation $\oplus_{m>0}V(\rho_m)^{n_m}$ is uniquely determined by $\sV$. Again, we denote the isomorphism class $\oplus_{m>0}V(\rho_m)^{n_m}$ by $[\sV^{gen}]$.\\[2pt]
Case 3. $G=N$, $T_1$ acts trivially on $Y$ and the map $q{\colon}Y\to N\backslash Y=\<\bar\sigma\>\backslash Y$ is an \'etale degree two cover.  In this case we first decompose $\sV$ into $T_1$ isotypical components $\sV=\oplus_m\sV_m$ as in the proof of Lemma~\ref{lem:NLocalTriv}(2), and proceed as in Case 2. For each $N$-trivializing open subscheme $j_U{\colon}U\to Y$, we have the isomorphism 
\[
j_U^*\sV\cong \sO_U\otimes^{\tau,\sigma}_kV(\rho_0)^{n_0}\oplus\oplus_{m>0}\sO_U\otimes_k^{\tau,\sigma}V(\rho_m)^{n_m}.
\]
Since this local description follows uniquely from the decomposition into $T_1$-isotypical components, the isomorphism class of the $k$-representation $V(\rho_0)^{n_0}\oplus\oplus_{m>0}V(\rho_m)^{n_m}$ is uniquely determined by $\sV$. Again, we denote the isomorphism class $V(\rho_0)^{n_0}\oplus\oplus_{m>0}V(\rho_m)^{n_m}$ by $[\sV^{gen}]$. \end{construction}

\begin{remark}\label{rem:DetTrivialization} Let $k$ a field of characteristic zero. Then $G=\SL_2^n$ is linearly reductive over $k$, so we may apply the above construction. In this case, we may use the set $\{\Sym^{m_1}(F_1)\otimes\ldots\otimes\Sym^{m_n}(F_n), 0\le m_i\}$ as a set of representatives for the irreducible representations of $G$, where $F_i$ is the $G$-representation given by the $i$th projection 
$\SL_2^n\to \SL_2$ followed by the standard inclusion $\SL_2\hookrightarrow \GL_2$. Each such representation   comes with a canonical trivialization of $\det\phi$.

For $k$ a field of characteristic $\neq2$, and $G=N^n$, the tensor products of the representations $\rho_m$, $\rho_0^-$ of $N$ similarly give a complete set of representatives for  the irreducible representations of $N^n$, so for each such representation $\phi$  of $N^n$, we have a unique multi-index $\epsilon=(\epsilon_1,\ldots, \epsilon_n)\in \{0,1\}^n$ and a canonical isomorphism $\det\phi\cong (\rho_0^-)^{\otimes \epsilon_1}\circ p_1\otimes\ldots\otimes (\rho_0^-)^{\otimes \epsilon_n}\circ p_n$, where $p_i{\colon}N^n\to N$ is the $i$th projection.
\end{remark}

\section{Cohomology and Borel-Moore  homology}\label{sec:CohBMHom}
We recall from \cite{DJK} the general set-up for  twisted cohomology and twisted Borel-Moore homology and then specialize to the case of the Witt sheaf. We refer the reader to \cite{DJK}, \cite{Hoyois6} for details and proofs of the material in this section; a resum\'e of much of this material is also covered in \cite[\S2, \S3]{LR}. Our goal here is to assemble from the existing literature the structures and operations that we will need to use for the rest of the paper, giving the reader a ``user's guide'' to their application in the subsequent sections, without going into much detail about constructions or proofs.

We fix a noetherian separated base-scheme $B$ of finite Krull dimension. 

As in the previous sections, we let $\Sch/B$ denote the category of quasi-projective  $B$-schemes and $\Sm/B$ the full subcategory of smooth $B$-schemes; for $Z\in \Sch/B$ we let $\pi_Z{\colon}Z\to B$ denote the structure morphism. We have the motivic stable homotopy category $\SH(-)$ with its six functor formalism on 
$\Sch/B$ \cite{Hoyois6}. For $Y\in \Sch/B$, the unstable motivic category  $\sH(Y)$ is a certain localization of the (infinity) category $\Spc(Y)$ of {\em spaces over $Y$}, i.e., presheaves of simplicial sets on $Y$. We also have the (infinity) category $\Spc_\bullet(Y)$ of {\em pointed spaces over $Y$},  i.e., presheaves of pointed simplicial sets on $Y$ with corresponding localization $\sH_\bullet(Y)$ . The {\em motivic stable homotopy category over $Y$}, $\SH(Y)$, is a localization of an infinity category of $T$-spectra in $\Spc_\bullet(Y)$, where $T$ is the quotient presheaf $\A^1_Y/(\A^1_Y\setminus\{0_Y\})\in \Spc_\bullet(Y)$, and a $T$-spectrum in $\Spc_\bullet(Y)$ is a sequence $(E_0, E_1,\ldots)$ with the $E_n\in \Spc_\bullet(Y)$, together with bonding maps $\epsilon_n{\colon}\Sigma_TE_n\to E_{n+1}$, $n=0, 1,\ldots$.

The six-functor formalism consists of the data of two pairs of adjoint functors for each morphism $f{\colon}Y\to X$ in $\Sch/B$:
\[
f^*{\colon}\SH(X)\xymatrix{\ar@<3pt>[r]&\ar@<3pt>[l]}\SH(Y){\colon}f_*,\ 
f_!{\colon}\SH(Y)\xymatrix{\ar@<3pt>[r]&\ar@<3pt>[l]}\SH(X){\colon}f^!,
\]
together with internal Hom and tensor (smash) products, satisfying the properties outlined in \cite[Theorem 1.1. Theorem 6.18.]{Hoyois6}. For example, there is a natural transformation $\eta_f{\colon}f_!\to f_*$, which is an isomorphism if $f$ is projective (extended to $f$ proper in \cite{CD}).

Another of these properties is the {\em localization distinguished triangle} associated to a closed immersion $i{\colon}Z\to X$ in $\Sch/B$ with open complement $j{\colon}U\to X$:
\begin{equation}\label{eqn:LocDistTriangle}
j_!j^!\xrightarrow{\eta^!_!}\id_{\SH(X)}\xrightarrow{u^*_*}i_*i^*\xrightarrow{\delta}j_!j^![1]
\end{equation}
The map $u^*_*$ is the unit of the adjunction $i^*\dashv i_*$ and, as above, $\eta^!_!$ is the co-unit of the adjunction $j_!\dashv j^!$ and.

Following \cite[\S4.2]{DeligneDetCoh}, let $K(Z)$ denote the $K$-theory space of perfect complexes on $Z$ and let $\sK(Z)$ denote the  fundamental groupoid of $K(Z)$. For $p{\colon}V\to Z$ a vector bundle with zero-section $s{\colon}Z\to V$, we have the auto-equivalence of $\SH(Z)$ $\Sigma^V:=p_\#s_*$, with inverse $\Sigma^{-V}:=s^!p^*$. The association $V\mapsto \Sigma^V$ extends to a functor of groupoids
\[
\Sigma^{-}{\colon}\sK(Z)\to \Aut(\SH(Z)),\ v\mapsto  \Sigma^v,
\]
sending a locally free sheaf $\sV$ to the operator $\Sigma^{\V(\sV)}$, where $\V(\sV)\to Z$ is the vector bundle $\Spec_{\sO_Z}\Sym^*\sV$.  

In particular, we have canonical isomorphisms
\[
\Sigma^{v[1]}\cong \Sigma^{-v}\cong (\Sigma^v)^{-1}.
\] 
for $v\in D^\perf(Z)$, and  for each distinguished triangle in $D^\perf(Z)$ 
\[
v'\to v\to  v''\to v'[1]
\]
we have the canonical isomorphisms 
\[
\Sigma^{v}\cong \Sigma^{v'\oplus v''}\cong \Sigma^{v'}\circ \Sigma^{v''}
\cong \Sigma^{v''}\circ \Sigma^{v'}.
\]
For details, see \cite[\S2.1.1]{DJK}, \cite[Th\'eor\`eme 1.5.18]{Ayoub6} and \cite[Remark 2.4.5]{CD}.

For a pointed space $\sX$ over $Z$, we have the suspension functor $\Sigma_\sX$ on $\Spc_\bullet(Z)$, $\Sigma_\sX(\sY):=\sY\wedge \sX$, inducing a corresponding functor $\Sigma_\sX$ on $\SH(Z)$. With $S^1$ the constant presheaf on the pointed simplicial set $S^1$,  $\G_m$ the presheaf represented by the pointed scheme $(\A^1\setminus\{0\}, \{1\})$, and $\P^1$ the   presheaf  represented by the pointed scheme $(\P^1, \{\infty\})$, we have the corresponding suspension functors $\Sigma_{S^1}$,  $\Sigma_{\G_m}$ and $\Sigma_{\P^1}$. These all define invertible endofunctors on $\SH(Z)$ for every $Z\in \Sch/B$, and we have canonical isomorphisms
\[
\Sigma_{\P^1}\cong \Sigma_{S_1}\circ\Sigma_{\G_m}\cong \Sigma_{\G_m}\circ\Sigma_{S^1}.
\]
For $a, b\in \Z$, we write $\Sigma^{a,b}$ for $\Sigma_{S^1}^{a-b}\Sigma_{\G_m}^b$ and we have for $r\in \N$ the canonical isomorphisms
\[
\Sigma^{\sO^r_Z}\cong \Sigma^r_{\P^1}\cong\Sigma^{2r,r}.
\]
See \cite[\S2.1.5]{DJK} for a discussion of these suspension operators. 

We say a word about the functors $f_*, f^*$, $f_!$, $f^!$ and the adjunctions $f^*\dashv f_*$, $f_!, \dashv f^!$. 

For $f{\colon}Z\to W$ an arbitrary morphism in $\Sch/B$, the functor $f_*{\colon}\Spc(Z)\to \Spc(W)$ is induced by the pull-back functor
\[
f^{-1}{\colon}\Sm/W\to \Sm/Z,\quad f^{-1}(T\to W)=T\times_WZ\xrightarrow{p_2}Z,
\]
sending a (simplicial) presheaf $P$ on $\Sm/Z$ to $f_*(P)\in \Spc(W)$, defined by
\[
f_*(P)(T\to W)=P(T\times_WZ\xrightarrow{p_2}Z).
\]
The left adjoint $f^*$ of $f_*$ is the left Kan extension of the map sending a representable presheaf $h_{T\to W}\in \Spc(W)$, $h_{T\to W}(T'\to W):=\Hom_W(T',T)$, to the representable presheaf $h_{T\times_WZ\to Z}\in \Spc(Z)$.   

For $T\to W\in \Sm/W$,  usually write the representable presheaf $h_{T\to W}\in \Spc(W)$ simply as $T$, if the context makes the meaning clear. 

If $f{\colon}Z\to W$ is smooth, we have the functor 
\[
f\circ -{\colon}\Sm/Z\to \Sm/W,
\]
sending $g{\colon}S\to Z$ to $f\circ g{\colon}S\to W$, and $f^*$ then has the direct expression as the $()_*$ for this functor, that is
\[
f^*Q(S\xrightarrow{g}Z)=Q(S\xrightarrow{fg}W);
\]
one checks directly that this is left adjoint to $f_*$. 

As an example, for $\pi_X{\colon}X\to B$ in $\Sm/B$, $\pi_{X\#}(h_{\id{\colon}X\to X})$ is just $h_{\pi_X{\colon}X\to B}$. Denoting by $1_X$ the $T$-suspension spectrum of $X_+$ in $\SH(X)$, this says that $\pi_{X\#}(1_X)$ is just $T$-suspension spectrum $\Sigma^\infty_TX_+\in\SH(B)$. Similarly,   for $q{\colon}X\to Y$ a   morphism in $\Sm/B$, we have the induced morphism $q{\colon}h_{X\to B}\to h_{Y\to B}$ sending a $B$-morphism $g{\colon}S\to X$ to $q\circ g{\colon}S\to Y$, and the corresponding natural map $q:=\Sigma^\infty_Tq_+{\colon}\pi_{X\#}(1_X)\to \pi_{Y\#}(1_Y)$. In other words, we have a functor
\[
h_B:=\Sigma^\infty_T(-)_+{\colon}\Sm/B\to \SH(B).
\]

The exceptional functors are a bit more difficult to define.  The method of Hoyois \cite{Hoyois6} is, roughly speaking, to define $j_!=j_\#$ for an open immersion $j{\colon}U\to X$ and to define $p_!=p_*$ for a proper morphism $p{\colon}P\to X$. This is done on the stable level, i.e., for $\SH(-)$. Then one factors an arbitrary morphism $f{\colon}Y\to X$ as a composition $f=p\circ j$ for $j{\colon}Y\to P$ an open immersion and $p{\colon}P\to X$ a proper map, and defines $f_!:=p_*\circ j_\#$. One does this is a somewhat more sophisticated fashion, which shows that $f_!$ is independent of the choice of factorization, and one then needs to show that $(fg)_!=f_!\circ g_!$. The fact that $f_!$ admits the right adjoint $f^!$ is then a formal argument. 

From this construction, it follows that $f_!=f_*$ if $f$ is already proper. In fact, there is a natural transformation $\alpha_f{\colon}f_!\to f_*$ for arbitrary $f$ in $\Sch/B$, which is a natural isomorphism if $f$ is proper. For details, we refer the reader to \cite[\S6]{Hoyois6}.

The method of Ayoub \cite[Scholie 1.4.2]{Ayoub6} is somewhat different. He constructs the exceptional functors for closed immersions and for smooth morphisms separately, and then shows that, in factoring a given morphism $f$ as a composition $f=p\circ i$, with $p$ smooth and $i$ a closed immersion, the resulting composition $f_!:=p_!\circ i_!$ is independent of choices. In particular, if $f{\colon}Z\to W$ is a smooth morphism in $\Sch/B$, letting $\Omega_f=\Omega_{Z/W}$ be the locally free sheaf of relative differentials, one has by definition
\[
f^!=\Sigma^{\Omega_f}\circ f^*
\]
and taking adjoints gives the canonical isomorphism
\[
f_!=f_\#\circ\Sigma^{-\Omega_f}. 
\]
We call these the {\em purity isomorphisms}. In consequence,  we have a natural isomorphism $f_!f^!\cong f_\#f^*$. The purity isomorphisms are also established in \cite[Theorem 1.5]{Hoyois6}, but only for $f$ smooth and proper. 

Using the exceptional functors, one can also describe the functor $h_B$ as follows:  Given $q{\colon}X\to Y$ in $\Sm/B$, the morphism $q{\colon}\pi_{X\#}(1_X)\to \pi_{Y\#}(1_Y)$ arises from the isomorphisms $\pi_{X!}\pi_X^!\cong \pi_{X\#}\pi_X^*$, $\pi_{Y!}\pi_Y^!\cong \pi_{Y\#}\pi_Y^*$, $1_X\cong \pi_X^*1_B$, $1_Y\cong \pi_Y^*1_B$, and the natural transformation
\[
\pi_{X!}\pi_X^!\cong \pi_{Y!}q_!q^!\pi_Y^!\xrightarrow{\eta^!_!}\pi_{Y!}\pi_Y^!
\]
applied to $1_B$; here  $\eta^!_!$ is the counit of the adjunction $q_!\dashv q^!$.

\subsection{$\sE$-cohomology, Borel-Moore homology, products} 
For $Z\in \Sch/B$ and  $v\in \sK(Z)$, we have the twisted Borel-Moore motive
\[
Z/B(v)_\BM:=\pi_{Z!}(\Sigma^v1_Z)\in \SH(B)
\]
and for $\sE\in \SH(B)$ the twisted $\sE$-Borel-Moore homology
\[
\sE^\BM_{a,b}(Z/B,v):=\Hom_{\SH(B)}(\Sigma^{a,b}Z/B(v)_\BM, \sE).
\]
We write $Z/B_\BM$ for $Z/B(0)_\BM$, $\sE^\BM_{a,b}(Z/B)$ for $\sE^\BM_{a,b}(Z/B,0)$ and
$\sE^\BM(Z/B,v)$ for $\sE^\BM_{0,0}(Z/B,v)$. We often drop the $-/B$ from the notation if the base-scheme is understood.

Define the twisted $\sE$-cohomology by
\[
\sE^{a,b}(Z, v):=\Hom_{\SH(Z)}(1_Z, \Sigma^{a,b}\Sigma^v\pi_Z^*\sE)=\Hom_{\SH(B)}
(1_B, \Sigma^{a,b}\pi_{Z*}(\Sigma^v\pi_Z^*\sE)).
\]
If $p{\colon}W\to Z$ is a morphism in  $\Sch/B$, we will often write $W/B(v)_\BM$ for $W/B(p^*v)_\BM$, and similarly use the notation $\sE^\BM_{a,b}(W/B,v)$, $\sE^{a,b}(W, v)$, etc. For $\sE=1_B$, we write 
$H^\BM_{a,b}$ and $H^{a,b}$ for $(1_B)^\BM_{a,b}$ and $(1_B)^{a,b}$.

For $v\in \sK(Z)$, sending $p{\colon}W\to Z\in \Sch/Z$ to $\sE^{a,b}(W, v)$ defines a functor
\[
\sE^{a,b}(-, v){\colon}\Sch/Z^\op\to \Ab
\]
with pullback map $f^*{\colon}\sE^{a,b}(W, v)\to \sE^{a,b}(W', v)$ for a $Z$-morphism $f{\colon}W'\to W$ induced by the functor $f^*{\colon}\SH(W)\to \SH(W')$:
\begin{multline*}
\Hom_{\SH(W)}(1_W, \Sigma^{a,b}\Sigma^{p^*v}\pi_W^*\sE)\xrightarrow{f^*}
\Hom_{\SH(W')}(f^*1_W, \Sigma^{a,b}\Sigma^{f^*p^*v}f^*\pi_W^*\sE)\\=
\Hom_{\SH(W')}(1_{W'}, \Sigma^{a,b}\Sigma^{p^{\prime*}v}\pi_{W'}^*\sE).
\end{multline*}

If $\sE$ is a commutative ring spectrum in $\SH(B)$ (i.e. a commutative monoid object), the multiplication in $\sE$   induces associative external products
\[
\sE^{a,b}(Z, v)\times \sE^{a',b'}(Z', v')\to \sE^{a+a',b+b'}(Z\times_BZ', v+v');
\]
pulling back by the diagonal gives $\sE^{*,*}(Z)$ the structure of a bi-graded associative ring with unit and makes $\sE^{*,*}(Z, v)$ a bigraded $\sE^{*,*}(Z)$-module. The multiplication also
 induces the map (see \cite[2.1.10]{DJK})
\[
\pi_Z^!\sE\wedge \pi_Z^*\sE\to \pi_Z^!\sE.
\]
Via the adjunction isomorphism
\[
\sE^\BM_{a,b}(Z/B,v):=\Hom_{\SH(B)}(\Sigma^{a,b}Z/B(v)_\BM, \sE)=
\Hom_{\SH(Z)}(\Sigma^{a,b}\Sigma^v 1_Z, \pi_Z^!\sE)
\]
this multiplication map gives us  the associative cap product pairing
\[
\cap{\colon} \sE^\BM_{a,b}(Z/B,v)\times \sE^{c,d}(Z, w)\to \sE^\BM_{a-c,b-d}(Z/B,v-w)
\]
defined by
\begin{multline*}
 \sE^\BM_{a,b}(Z/B,v)\times \sE^{c,d}(Z, w)=\\
\Hom_{\SH(Z)}(\Sigma^{a,b}\Sigma^v 1_Z, \pi_Z^!\sE)\times\Hom_{\SH(Z)}(\Sigma^{-w}1_Z, \Sigma^{c,d}\pi_Z^*\sE)\\
\to \Hom_{\SH(Z)}(\Sigma^{a-c,b-d}\Sigma^{v-w}1_Z, \pi_Z^!\sE)\\
=\sE^\BM_{a-c,b-d}(Z/B,v-w)
\end{multline*}

\subsection{Purity isomorphisms} See \cite[Scholie 1.4.2]{Ayoub6},  \cite[Theorem 1.5]{Hoyois6}. For $X\in \Sm/B$, $v\in \sK(X)$, the adjunction $\pi_\#\dashv \pi_X^*$ gives the isomorphism 
\[
\sE^{a,b}(X,v):=\Hom_{\SH(B)}(\pi_{X\#}(\Sigma^{-v}1_X), \Sigma^{a,b}\sE).
\]
More generally, for $\sF\in \SH(X)$, we set $\sE^{a,b}(\sF,v):=\Hom_{\SH(B)}(\pi_{X\#}(\Sigma^{-v}\sF), \Sigma^{a,b}\sE)$. For example, if $i{\colon}Z\to X$ is a closed immersion in $\Sch/B$, we have the cohomology with support
\[
\sE^{a,b}_Z(X,v):= \sE^{a,b}(i_*(1_Z),v)=\Hom_{\SH(B)}(\pi_{X\#}(\Sigma^{-v}i_*1_Z), \Sigma^{a,b}\sE).
\]

The purity isomorphism $\pi_{X\#}\cong \pi_{X!}\circ\Sigma^{\Omega_{X/B}}$ gives the Poincar\'e duality isomorphism of cohomology with Borel-Moore homology
\begin{equation}\label{eqn:PBMurityIso}
\sE^{a,b}(X,v)\cong \Hom_{\SH(B)}(\pi_{X!}(\Sigma^{\Omega_{X/B}-v}1_X), \Sigma^{a,b}\sE)
=\sE^\BM_{-a,-b}(X/B, \Omega_{X/B}-v)
\end{equation}
Similarly if $i{\colon}Z\to X$ is a closed immersion in $\Sch/B$ with $X\in\Sm/B$, we have $i_*=i_!$, giving the isomorphism
\begin{multline}\label{eqn:BMPurityIsoSupp}
\sE^{a,b}_Z(X,v)\cong \Hom_{\SH(B)}(\pi_{X!}(\Sigma^{\Omega_{X/B}-v}i_!1_Z), \Sigma^{a,b}\sE)
\\\cong \Hom_{\SH(B)}(\pi_{Z!}(\Sigma^{i^*\Omega_{X/B}-i^*v}1_Z), \Sigma^{a,b}\sE)
 \\=\sE^\BM_{-a,-b}(Z/B, i^*\Omega_{X/B}-i^*v).
 \end{multline}
If in addition $Z\in \Sm/B$, then using the exact sequence
\[
0\to \sN_i\to i^*\Omega_{X/B}\to  \Omega_{Z/B}\to0
\]
with $\sN_i$ the conormal sheaf, \eqref{eqn:PBMurityIso} and \eqref{eqn:BMPurityIsoSupp} yield the purity isomorphism
\begin{equation}\label{eqn:PBMurityIso2}
\sE^{a,b}_Z(X,v)\cong \sE^{a,b}(Z, v-\sN_i).
\end{equation}

\subsection{Thom isomorphisms, oriented spectra} In this section, we assume that $B$ is an affine scheme, still noetherian and of finite Krull dimension. We refer the reader to  \cite[\S4.3]{BW} or  \cite[\S 3]{LR} for details on oriented and $\SL$-oriented spectra; oriented spectra are referred to as $\GL$-oriented spectra in {\it loc. cit.}. We also refer the reader to \cite{AnanSL} and \cite{PaninII, PaninI} for the basic properties of oriented and $\SL$-oriented cohomology theories. 

In the first treatments of this topic, an orientation for $\sE\in \SH(B)$ consists of the assignment of a {\em Thom class} $\th_V\in \sE^{2r,r}(\Th(V))$ for each vector bundle $V\to X$, $X\in \Sm/B$, where $r$ is the rank of $V$, satisfying appropriate axioms. The convention is to associate to $V$ the locally free coherent sheaf $\sV$ of sections of the dual bundle $V^\vee$, giving the isomorphism
\[
V\cong \V(\sV):= \Spec_{\sO_X}(\Sym^*\sV).
\]
With $\pi{\colon}X\to B$ the structure morphism, we have
\[
\Th(V)\cong \pi_\#(\Sigma^\sV1_X)
\]
in $\SH(B)$, giving the isomorphisms
\begin{multline*}
\sE^{2r,r}(\Th(V))=\Hom_{\SH(B)}(\pi_\#(\Sigma^\sV1_X), \Sigma^{2r,r}\sE)\\\cong 
\Hom_{\SH(X)}(1_X, \Sigma^{2r,r}\Sigma^{-\sV}\pi^*\sE)=\sE^{2r,r}(X, -\sV),
\end{multline*}
realizing $\th_V$ as an element in $\sE^{2r,r}(X, -\sV)$. We will use this convention in our discussion below.

Let $\sE$ be an oriented ring spectrum in $\SH(B)$.  The orientation gives rise to natural Thom classes $\th_v\in \sE^{2r,r}(X, -v)$ for $X\in \Sch/B$ and   $v\in \sK(X)$ of rank $r$, which give rise to Thom isomorphisms
\[
\cup\th_v{\colon}\sE^{*-2r, *-r}(X, w)\xrightarrow{\sim}\sE^{*,*}(X, w-v).
\]
Equivalently, we have the spectrum-level Thom isomorphism $\Sigma^{\sO_X^r-v}\pi^*\sE\cong \pi^*\sE$. 

Similarly, for $X\in \Sch/B$, and  $v\in \sK(X)$ of virtual rank $r$,  we have the  Thom isomorphism (depending on the choice of orientation) $\Sigma^{\sO_X^r-v}\pi_X^!\sE\cong  \pi_X^!\sE$,  giving the isomorphism
\[
\sE^\BM_{a-2r,b-r}(X/B,v+w)\cong \sE^\BM_{a,b}(X/B, w).
\]

If $\sE$ is an $\SL$-oriented ring spectrum,  we have similar Thom classes and Thom isomorphisms for pairs $(v,\rho)$ with $v\in \sK(X)$ and $\rho{\colon}\det v\xrightarrow{\sim}\sO_X$ a trivialization of $\det v$.  For $v\in \sK(X)$, $\det(v-\det v)=\det v\otimes(\det v)^{-1}$ has a canonical trivialization. Thus, we have canonical Thom classes for arbitrary $v\in \sK(X)$ (of virtual rank $r$), 
$\th_v\in \sE^{2r,r}(X,-v+\det v-\sO_X)$, giving  Thom isomorphisms
\[
\cup\th_v{\colon}\sE^{*-2r, *-r}(X, w)\xrightarrow{\sim} \sE^{*,*}(X, -v+w+\det v-\sO_X)
\]
induced by isomorphisms $\Sigma^{\sO^r_X-v}\pi^*\sE\cong \Sigma^{\sO_X-\det v}\pi^*\sE$.

For  $X\in \Sch/B$, $v\in \sK(X)$ we similarly have a canonical isomorphism
\[
\Sigma^{\sO^r_X-v}\pi_X^!\sE\cong \Sigma^{\sO_X-\det v}\pi_X^!\sE
\]
giving the isomorphism
\[
\sE^\BM_{a-2r,b-r}(X/B,v+w)\cong \sE^\BM_{a,b}(X/B, \det v-\sO_X+w).
\]

As the Thom classes and Thom  isomorphisms are often only defined for $X\in \Sm/B$ and only for $\sE$-cohomology (see e.g. \cite{Anan}) we briefly explain how these extend naturally to $X\in \Sch/B$ and to Borel-Moore homology. We discuss the case of $\SL$-orientations;  the oriented case is treated in exactly the same way.

One first reduces to the case of a twist by a vector bundle $V\to X$, $V=\V(\sV)$. Since we are working with quasi-projective $B$-schemes, we may apply a Jouanolou trick and reduce to the case of affine $X$; here is where we are using our assumption that $B$ is affine. Then $V$ is generated by global sections, and a choice of generating sections gives a morphism to a Grassmannian and an isomorphism of $V$ with the pullback of the tautological quotient vector bundle. Combining this with the given closed immersion in an open subscheme of $\P^n_B$ for some $n$, and applying a Jouanolou trick again, we see that there is a closed immersion
\[
i{\colon}X\to Y
\]
with $Y$ smooth and affine, a vector bundle $W$ on $Y$ and an isomorphism $\phi{\colon}V\to i^*W$. This gives us a class $\phi^*i^*\th_W\in \sE^{2r,r}(X, -\sV+ \det \sV-\sO_X)$. The axioms of Thom classes imply that locally on $Y$, $\th_W$ is just the image of the unit in $\sE^{0,0}(Y)$ by a suitable suspension isomorphism, so the same holds for $i^*\th_W$; a Mayer-Vietoris argument then shows that $i^*\th_W$ gives rise to Thom isomorphisms of the desired form.

Thus, we need only show that $i^*\th_W\in \sE^{2r,r}(X, -\sV+\det \sV-\sO_X)$ is independent of the choice of $i{\colon}X\to Y$, vector bundle $W\to Y$ and isomorphism $\phi$.   This all follows from

\begin{lemma} Let $i{\colon}X\to Y$ be a closed immersion of affine noetherian schemes, let $W_1, W_2$ be vector bundles on $Y$ and let $\bar\psi{\colon}i^*W_1\to i^*W_2$ be an isomorphism. Then there is an affine open neighborhood $j{\colon}U\to Y$ of $X$ in $Y$ and an isomorphism $\psi{\colon}j^*W_1\to j^*W_2$ with $i^*\psi=\bar\psi$
\end{lemma}

\begin{proof} Writing $Y=\Spec A$ and $X=\Spec A/I$ for some ideal $I\subset A$, $W_i$ is given by a finitely generated projective $A$-module $P_i$ and $\bar\psi$ is given by an isomorphism $\bar\psi{\colon}P_1/IP_1\to P_2/IP_2$. This gives us an element $\bar\theta\in P_1\otimes_AP_2^\vee/(IP_1\otimes_AP_2^\vee)$, which we may lift to $\theta\in P_1\otimes_AP_2^\vee$, giving us the map of $A$-modules $\tilde\psi{\colon}P_1\to P_2$. $\tilde\psi$ has a well-defined determinant $\det\tilde\psi\in \det(P_1)\otimes_A\det(P_2^\vee)$, with $i^*\det\tilde\psi$ a nowhere vanishing element of the invertible sheaf on $X$ corresponding to 
\[
\det(P_1)\otimes_A\det(P_2^\vee)/I\det(P_1)\otimes_A\det(P_2^\vee)=\det(P_1/IP_1)/\otimes_{A/I}\det(P_2^\vee/IP_2^\vee).
\]
The section $\det\tilde{\psi}$ of the invertible sheaf on $Y$ corresponding to $\det(P_1)\otimes_A\det(P_2^\vee)$ generates an ideal $(\det\tilde\psi)\subset A$, and since $i^*\det\tilde\psi$ is nowhere zero, $X$ is a closed subscheme of the affine open subscheme $U:=\Spec A[(\det\tilde\psi)^{-1}]$. Letting $j{\colon}U\to Y$ be the open immersion, we see that $\tilde\psi$ defines an isomorphism $\psi:=j^*\tilde\psi{\colon}j^*W_1\to j^*W_2$ lifting $\bar\psi$.
\end{proof}

Now suppose we have closed immersion $i_j{\colon}X\to Y_j$, vector bundles $W_j\to Y_j$ and isomorphisms $\phi_j{\colon}V\xrightarrow{\sim} i_j^*W_j$, $j=1,2$. Then we have the diagonal closed immersion $i:=(i_1, i_2){\colon}X\to Y_1\times_BY_2$, and the isomorphism $\bar\psi{\colon}\phi_1\circ\phi_2^{-1}{\colon}i^*p_2^*W_2\xrightarrow{\sim} i^*p_1^*W_1$. Lifting $\bar\psi$ to an isomorphism $\psi{\colon}j^*p_2^*W_2\xrightarrow{\sim} j^*p_1^*W_1$ on some open neighborhood $j{\colon}U\to Y$ of $X$,  and using the naturality of the Thom classes, we have $\psi^*\th_{W_1}=\th_{W_2}$, which implies that $\phi_1^*\th_{W_1}=\phi_2^*\th_{W_2}$. This gives us the well-defined Thom class in $\sE^{2r,r}(X, -\sV+\det \sV-\sO_X)$; the corresponding Thom isomorphisms follow. 

Another approach to constructing Thom isomorphisms for  $X\in \Sch/B$, $v\in \sK(X)$,  is given in \cite[\S4.3]{BW}.

Once we have Thom classes and Thom isomorphisms for $\sE$-cohomology, the Thom isomorphism for $\sE$-Borel-Moore homology follows by using the product action of 
 $\sE$-cohomology on $\sE$-Borel-Moore homology (cap product).

\begin{remark}\label{rem:LToLDual} In a similar vein, Ananyeskiy's result \cite[Lemma 4.1]{AnanSL}, that for $X\in \Sm/B$,  $L$ a line bundle on $X$ and $E$ a vector bundle on $X$, one has a canonical isomorphism $\Th_X(E\oplus L)\cong \Th_X(E\oplus L^{-1})$ in $\sH_\bullet(B)$, extends to the following statement: For $X\in \Sch/B$, with locally free sheaf $\sV$ and invertible sheaf $\sL$, one has a natural isomorphism of suspension functors on $\SH(X)$
\[
\Sigma^{\sV\oplus\sL}\cong  \Sigma^{\sV\oplus\sL^{-1}}.
\]
Using the splitting principle, this then induces corresponding isomorphisms on $\sE$-cohomology and $\sE$-Borel-Moore homology for arbitrary $v, w\in \sK(X)$
\begin{equation}\label{eqn:AAIso}
\sE^{*,*}(X, v+w)\cong \sE^{*,*}(X, v+w^\vee),\ \sE^\BM_{*,*}(X, v+w)\cong \sE^\BM_{*,*}(X, v+w^\vee).
\end{equation}
\end{remark}

\subsection{Functorialities}\label{subsec:Functorialities} Fix $\sE\in \SH(B)$. 
For $p{\colon}Z\to W$ a proper morphism in $\Sch/B$, and $v\in \sK(W)$, we have the {\em proper pull-back}
\[
p^*{\colon}W/B(v)_\BM\to Z/B(v)_\BM
\]
inducing the functorial {\em proper push-forward} on Borel-Moore homology
\[
p_*:=(p^*)^*{\colon}\sE^\BM_{a,b}(Z/B,v)\to \sE^\BM_{a,b}(W/B,v)
\]
The functorial proper pull-back is constructed as follows: We have the natural transformation $\alpha_p{\colon}p_!\to p_*$, which is an isomorphism since $p$ is proper. We also have the unit of adjunction $u^*_*{\colon}\id_{\SH(W)}\to p_*p^*$. This gives the natural transformation
\[
\pi_{W!}\xrightarrow{u^*_*}\pi_{W!}p_*p^*\xymatrix{\ar[r]^{\alpha^{-1}_p}&} 
\pi_{W!}p_!p^*=\pi_{Z!}p^*,
\]
which we apply to $\Sigma^v1_W$ to give 
\[
p^*{\colon}W/B(v)_\BM=\pi_{W!}\Sigma^v1_W\to
\pi_{Z!}p^*\Sigma^v1_W=\pi_{Z!}\Sigma^{p^*v}1_Z=Z/B(v)_\BM.
\]

For $q{\colon}P\to W$ a smooth morphism in $\Sch/B$, we have the relative K\"ahler differentials $\Omega_q:=\Omega_{P/W}$ and the functorial {\em smooth push-forward}
\[
q_!{\colon}P/B(v+\Omega_q)_\BM\to W/B(v)_\BM,
\]
inducing the functorial {\em smooth pullback} on Borel-Moore homology
\[
q^!:=(q_!)^*{\colon}\sE^\BM_{a,b}(W/B,v)\to \sE^\BM_{a,b}(P/B,v+\Omega_q).
\]
If $q{\colon}P\to W$ is a vector bundle on $W$ (or even a torsor for a vector bundle on $W$), then $q^!$ is an isomorphism. The smooth pushforward is constructed using the purity isomorphism $q_!\Sigma^{\Omega_q}\cong q_\#$ and the co-unit $\eta^*_\#{\colon}q_\#q^*\to \id_{\SH(W)}$ of adjunction: 
\begin{multline*}
P/B(v+\Omega_q)_\BM=\pi_{P!}\Sigma^{q^*v+\Omega_q}1_P=\pi_{W!}\Sigma^v q_!\Sigma^{\Omega_q}q^*1_W\\\cong
\pi_{W!}\Sigma^v q_\#q^*1_W\xrightarrow{\eta^*_\#}\pi_{W!}\Sigma^v 1_W=W/B(v)_\BM.
\end{multline*}

If $i{\colon}Z\to W$ is a regular immersion in $\Sch/B$, we have the Gysin pullback map
\[
i^!{\colon}\sE^\BM_{a,b}(W/B,v)\to \sE^\BM_{a,b}(Z/B,v-\sN_i)
\]
where $\sN_i$ is the conormal sheaf of $i$.   See \cite[\S2.4]{DJK} for the construction.  If $q{\colon}P\to W$ is a vector bundle and $0_P{\colon}W\to P$ is the zero-section, then  $q^!$ and $0_P^!$ are inverse isomorphisms.  

We call a morphism $f{\colon}Z\to W$ an {\em lci morphism} if $f$ admits a factorization as $f=q\circ i$, with $i{\colon}Z\to P$ a regular immersion and $q{\colon}P\to W$ a smooth morphism (both in $\Sch/B$). Since $Z\to B$ in $\Sch/B$ is quasi-projective over $B$, an lci morphism in $\Sch/B$ is a {\em smoothable} lci morphism in the sense of \cite{DJK}. Let $L_f$ be the relative cotangent complex of $f$. The lci pullback
\[
f^!{\colon}\sE^\BM_{a,b}(W/B,v)\to \sE^\BM_{a,b}(Z/B,v+L_f)
\]
is defined as $f^!:=i^!\circ q^!$, using the distinguished triangle
\[
i^*\Omega_q\to L_f\to \sN_i[1]\to i^*\Omega_p[1]
\]
to define the isomorphism $\sE^\BM_{a,b}(Z/B,v+L_f)\cong \sE^\BM_{a,b}(Z/B,v+\Omega_p-\sN_i)$, 
and is independent of the choice of factorization. Moreover, for composable lci morphisms $g{\colon}W\to U$, $f{\colon}Z\to W$, after using the distinguished triangle 
\[
f^*L_g\to L_{gf}\to L_f\to f^*L_g[1]
\]
to identify $\sE^\BM_{a,b}(Z/B,v+L_f+L_g)$ with $\sE^\BM_{a,b}(Z/B,v+L_{gf})$, we have $(gf)^!=f^!g^!$. See  \cite[Theorem 4.2.1]{DJK} for these facts. We will often unify the notation, writing $\sE^\BM_{a,b}(Z/B,v+L_q)$ for $\sE^\BM_{a,b}(Z/B,v+\Omega_q)$ in the case of a smooth morphism $q$ and $\sE^\BM_{a,b}(Z/B,v+L_i)$ for  $\sE^\BM_{a,b}(Z/B,v-\sN_i)$ in the case of a regular embedding $i$. Lci pullback is functorial for composition of lci morphisms

\subsection{$\sE$-cohomology with supports} Fix a commutative ring spectrum $\sE\in \SH(B)$. 
Let $j{\colon}U\to Z$ be an open immersion with closed complement $i{\colon}T\to Z$. The localization triangle
\[
j_!1_U\to 1_Z\to i_*1_{T}\to j_!1_U[1]
\]
identifies $ i_*1_{T}$ with the suspension spectrum of $Z/U$ in $\SH(Z)$, identifying the cohomology with support   
\[
\sE^{c,d}_T(Z, w):=\Hom_{\SH(Z)}(\Sigma^{v}i_{T*}1_T, \Sigma^{c,d}\pi_Z^*\sE).
\]
with $\sE^{c,d}(Z/U, w)$.  We let $i_{T*}{\colon}\sE^{c,d}_T(Z, w)\to \sE^{c,d}(Z, w)$ be the map (``forget supports'') induced by the co-unit $1_Z\to  i_*1_{T}=i_*i^*1_Z$. 

The localization triangle gives us the long exact  sequence for cohomology with support
\begin{equation}\label{eqn:SeqCohSupp}
\ldots\to \sE^{c,d}_T(Z, w)\xrightarrow{i_{T*}} \sE^{c,d}(Z, w)\xrightarrow{j^*} \sE^{c,d}(U, w)\xrightarrow{\delta} \sE^{c+1,d}_T(Z, w)\to\ldots\ .
\end{equation}

\subsection{Localization in Borel-Moore homology} Let $i{\colon}Z\to W$ be a closed immersion in $\Sch/B$ with open complement $j{\colon}U\to W$, and take $v\in \sK(W)$. Let $\pi_W{\colon}W\to B$ be the structure morphism. We have the {\em localization triangle}
\begin{equation}\label{eqn:BMLocTriangle}
U/B(j^*v)_\BM\xrightarrow{j_!}W/B(v)_\BM\xrightarrow{i^*}Z/B(i^*v)_\BM\xrightarrow{\delta}
U/B(j^*v)_\BM[1]
\end{equation}
formed by evaluating the localization triangle \eqref{eqn:LocDistTriangle}:
\[
j_!j^!\to \id_{\SH(W)}\to i_*i^*\to j_!j^![1]
\]
 on $\Sigma^v1_W$, applying $\pi_{W!}$ and  using the isomorphisms $j^!\Sigma^v1_W\cong j^*\Sigma^v1_W\cong \Sigma^{j^*v}1_U$, $i^*\Sigma^v1_W\cong \Sigma^{i^*v}1_Z$,  $\pi_{W!}j_!=\pi_{U!}$ and $\pi_{W!}i_*=\pi_{W!}i_!=\pi_{Z!}$. Mapping to $\Sigma^{-a,-b}\sE$ gives the long exact {\em localization sequence}
\begin{equation}\label{eqn:BMLocSeq}
\ldots\xrightarrow{j^!}\sE^\BM_{a+1,b}(U, v)\xrightarrow{\del} \sE^\BM_{a,b}(Z, v)\xrightarrow{i_*}
\sE^\BM_{a,b}(W, v)\xrightarrow{j^!}\sE^\BM_{a,b}(U, v)\xrightarrow{\del}\ldots
\end{equation}
If we fix a closed immersion $i_W{\colon}W\to X$ for some $X\in \Sm/B$ the various purity isomorphisms map this sequence isomorphically to the exact sequence in cohomology with supports \eqref{eqn:SeqCohSupp}:
\begin{multline}\label{mult:CohSuppLocSeq1}
\ldots\to\sE^{-a-1,-b}_U(V, \Omega_{X/B}-v)\xrightarrow{\delta} \sE^{-a,-b}_Z(X, \Omega_{X/B}-v)\to
\sE^{-a,-b}_W(X, \Omega_{X/B}-v)\\\xrightarrow{j_X^*}\sE^{-a,-b}_U(V, \Omega_{X/B}-v)\xrightarrow{\delta}\ldots
\end{multline}
with $j_X{\colon}V:=X\setminus i_W(Z)\to X$ the inclusion; if $W$ is smooth over $B$, we may take $i_W=\id_W$, giving an isomorphism of \eqref{eqn:BMLocSeq} with 
\begin{multline}\label{mult:CohSuppLocSeq2}
\ldots\to\sE^{-a-1,-b}(U, \Omega_{W/B}-v)\xrightarrow{\delta} \sE^{-a,-b}_Z(W, \Omega_{W/B}-v)\to
\sE^{-a,-b}(W, \Omega_{W/B}-v)\\\xrightarrow{j^*}\sE^{-a,-b}(U, \Omega_{W/B}-v)\xrightarrow{\delta}\ldots
\end{multline}

\begin{remark}\label{rem:Reduced} For $W\in \Sch/B$, applying localization to the closed immersion $i{\colon}W_\red\hookrightarrow W$ with empty open complement, we see that $i_*{\colon}\sE^\BM_{a,b}(W_\red, v)\to
\sE^\BM_{a,b}(W, v)$ is an isomorphism.
\end{remark}

\subsection{Euler class}
For $p{\colon}V\to Z$ a vector bundle on some $Z\in \Sch/B$, $V=\V(\sV)$, we have the Euler class $e(V)\in H^{0,0}(Z, \sV)$ defined as follows (see \cite[\S 3.1]{DJK} for details): Directly from the definitions, we have the isomorphism
\[
H^{a,b}(Z, v))\cong H_{-a,-b}^\BM(Z/Z,-v).
\]
Letting $s{\colon}Z\to V$ be the zero-section, we have the proper pushforward
\[
s_*{\colon}H^{0,0}(Z)=H_{0,0}^\BM(Z/Z)\to H_{0,0}^\BM(V/Z).
\]
As $s$ is a regular embedding with conormal sheaf $\sV$, we have the Gysin pullback
\[
s^!{\colon}H_{0,0}^\BM(V/Z)\to H_{0,0}^\BM(Z/Z, -\sV)\cong H^{0,0}(Z, \sV)
\]
Then one defines
\[
e(V):=s^!s_*(1_Z)\in H^{0,0}(Z, \sV).
\] 

For $\sE\in \SH(B)$ a commutative ring spectrum, applying the unit map $\epsilon_\sE{\colon}1_B\to \sE$ gives the Euler class $e^\sE(V)\in \sE(Z, \sV)$; this is just the image of the unit $1^\sE_Z\in \sE^{0,0}(Z)$ by the composition defined just as for $\sE=1_B$, 
\[
 \sE^{0,0}(Z)=\sE^\BM_{0,0}(Z/Z)\xrightarrow{s_*}\sE_{0,0}^\BM(V/Z)\xrightarrow{s^*} \sE_{0,0}^\BM(Z/Z, -\sV)\cong \sE^{0,0}(Z, \sV).
 \]
 If $\sE$ is oriented and $V$ has rank $r$, applying the Thom isomorphism gives us $e^\sE(V)\in \sE^{2r, r}(Z)$ and if $\sE$ is $\SL$-oriented,  the same gives us $e^\sE(V)\in \sE^{2r,r}(Z, \det \sV-\sO_Z)$.  We will nearly always use these ``Thom-twisted'' Euler classes when dealing with oriented or $\SL$-oriented $\sE$. We will sometimes write $e(V)$ for $e^\sE(V)$ if the spectrum $\sE$ is clear from the context.

 \subsection{Refined Gysin morphism} Let
 \[
\xymatrix{
X_0\ar[r]^-{\iota'}\ar[d]_{q}&X\ar[d]^{p}\\
Y_0\ar[r]^\iota&Y
}
\]
be a cartesian diagram in $\Sch/B$, with $\iota$ a regular immersion, and take $\sE\in \SH(B)$. We recall from \cite[Definition 4.2.5]{DJK} the {\em refined Gysin morphism}
\[
\iota^!_p{\colon}\sE^\BM(X,v)\to  \sE^\BM(X_0, v+L_\iota).
\]
This is actually defined more generally for $\iota$ a (smoothable) lci morphism, but we will not need this. Note that the refined Gysin morphism is just the Gysin morphism of \S\ref{subsec:Functorialities} in case $p$ is an identity map.
 
 \subsection{Compatibilities}\label{subsec:Relations} We list the various compatibilities of the operations of lci pullback, refined Gysin morphism and proper push-forward in the following proposition. We fix a  commutative monoid $\sE\in \SH(B)$. See \cite{DJK} for details; they write everything in terms of $\sE=1_B$, but these properties extend to arbitrary $\sE$ by applying the unit map $u{\colon}1_B\to \sE$ and the multiplication on $\sE$.  
 
\begin{proposition}\label{prop:GysinProperties} Take $Y\in \Sch/B$ with $v\in \sK(Y)$.\\[2pt]
 1 (compatibility with proper push-forward and smooth pullback). Let
\[
\xymatrix{
W_0\ar[r]^-{\iota''}\ar[d]_{q_2}\ar@/_20pt/[dd]_q&W\ar[d]^{p_2}\ar@/^20pt/[dd]^p\\
X_0\ar[r]^-{\iota'}\ar[d]_{q_1}&X\ar[d]^{p_1}\\
Y_0\ar[r]^\iota&Y
}
\]
be a commutative diagram in $\Sch/B$, with both squares cartesian.\\[2pt]
1a. Suppose $p_2$ is proper and  $\iota$ is a regular immersion. Then the diagram
\[
\xymatrix{
\sE^\BM(W,v)\ar[r]^-{\iota^!_p}\ar[d]^{p_{2*}}& \sE^\BM(W_0, v+L_\iota)\ar[d]^{q_{2*}}\\
\sE^\BM(X,v)\ar[r]^-{\iota^!_{p_1}}&\sE^\BM(X_0, v+L_\iota)
}
\]
commutes.\\[2pt]
1b. Suppose $p_2$ is smooth and  $\iota$ is a regular immersion. Then the diagram
\[
\xymatrix{
\sE^\BM(W,v+L_{p_2})\ar[r]^-{\iota^!_p}& \sE^\BM(W_0, v+L_\iota+L_{q_2})\\
\sE^\BM(X,v)\ar[u]^{p_{2}^!}\ar[r]^-{\iota^!_{p_1}}& \sE^\BM(X_0, v+L_\iota)\ar[u]^{q^!_{2}}
}
\]
commutes. Here we use the canonical isomorphism $L_{q_2}\cong i^{\prime\prime *}L_{p_2}$.\\[2pt]
1c. Suppose $p_2$ is an identity morphism,  $p_1$ is smooth and  $\iota$ is proper. Then $p_1^!\circ \iota_*=\iota'_*\circ q_1^!$.
2 (functoriality).  Let
\[
\xymatrix{
X_1\ar[r]^{\iota_1'}\ar[d]_{p_1}&X_0\ar[r]^{\iota_0'}\ar[d]_{p_0}&X\ar[d]^{p}\\
Y_1\ar[r]^{\iota_1}&Y_0\ar[r]^{\iota_0}&Y
}
\]
be a commutative diagram in $\Sch/B$, with both squares cartesian, and with $\iota_0$ and $\iota_1$ regular immersions. Then $\iota_0\circ\iota_1$ is a regular immersion and
\[
\iota_1^!\circ\iota_0^!=(\iota_0\circ\iota_1)^!{\colon}
\sE^\BM(X,v)\to \sE^\BM(X_1, v+L_{\iota_0\circ\iota_1}).
\]
Here we use the distinguished triangle
\[
  \iota_1^*L_{\iota_0}\to  L_{\iota_0\circ\iota_1}\to L_{\iota_1}\to\iota_1^*L_{\iota_0}[1]
\]
to identify $\sE^\BM(X_1, v+L_{\iota_1}+L_{\iota_0})$ with $\sE^\BM(X_1, v+L_{\iota_0\circ\iota_1})$.\\[2pt]
3 (excess intersection formula).  Let
\[
\xymatrix{
W_0\ar[r]^{\iota''}\ar[d]_{q_2}&W\ar[d]^{p_2}\\
X_0\ar[r]^{\iota'}\ar[d]_{q_1}&X\ar[d]^{p_1}\\
Y_0\ar[r]^\iota&Y
}
\]
be a commutative diagram in $\Sch/B$, with both squares cartesian, and with $\iota$ and $\iota'$ both regular immersions.  We have a natural surjection of locally free sheaves  on $X_0$, $q_1^*\sN_\iota\to \sN_{\iota'}$, with locally free kernel  $\sV$, giving the Euler class $e(q_2^*V)\in \sE(W_0, \sV)$, where $V=\V(\sV)$.  Then for $\alpha\in \sE^\BM(W, v)$, we have
\[
\iota^!(\alpha)=\iota^{\prime !}(\alpha)\cap e(q_2^*V)\in \sE^\BM(W_0, v+L_\iota)=
\sE^\BM(W_0, v+L_{\iota'}-\sV),
\]
where we use the exact sequence
\[
0\to \sV\to   q_1^*\sN_\iota\to \sN_{\iota'}\to 0
\]
to give the canonical isomorphism 
\[
\sE^\BM(W_0, v+L_\iota)\cong 
\sE^\BM(W_0, v+L_{\iota'}-\sV)).
\]
4. (Fundamental classes, Poincar\'e duality) Let $1^\sE\in \sE^{0,0}(B)=\sE^\BM(B/B)$ be the class of $u\circ \id_{1_B}$. Let $p_X{\colon}X\to B$ be an lci scheme over $B$, in $\Sch/B$ (an lci scheme for short). Define the {\em fundamental class} $[X]^\sE\in \sE^\BM(X/B, L_{X/B})$ by 
\[
[X]^\sE:=p_X^!(1^\sE).
\]
We often write $[X]$ for $[X]^\sE$ if $\sE$ is clear from the context.\\[2pt]
i. For $f{\colon}Y\to X$ an lci morphism of lci-schemes, we have $f^!([X]^\sE)=[Y]^\sE$. \\[2pt]
ii. For $X\in \Sm/B$, $L_{X/B}=\Omega_{X/B}$ and the cap product with $[X]^\sE$,
\[
[X]^\sE\cap-{\colon}\sE^{a,b}(X, v)\to \sE^\BM_{-a,-b}(X, \Omega_{X/B}-v),
\]
is an isomorphism, and is equal to  the purity isomorphism.\\[2pt]
5. (Localization) Let $i{\colon}Z\to X$ be a closed immersion in $\Sch/B$ with open complement $j{\colon}U\to X$, let $f{\colon}X'\to X$ be a morphism in $\Sch/B$, giving the diagram with all squares cartesian, defining the closed immersion $i'{\colon}Z'\to X'$ with open complement $j'{\colon}U'\to X'$,
\[
\xymatrix{
Z'\ar[r]^{i'}\ar[d]_{f_Z}&X'\ar[d]_f&U'\ar[l]_{j'}\ar[d]^{f_U}\\
Z\ar[r]^i&X&U\ar[l]_j
}
\]
Take $v\in \sK(X)$.\\[2pt]
i. Suppose $f$ is proper, Then the maps $f_*, f_{Z*}, f_{U*}$ define a map of distinguished localization triangles in $\SH(B)$
\[
\xymatrix{
Z'/B(v)_\BM\ar[r]^{i'_*}\ar[d]_{f_{Z*}}&X'/B(v)_\BM\ar[d]_{f'_*}\ar[r]^{j^{\prime!}}&U'/B(v)_\BM\ar[d]^{f_{U*}}\ar[r]^\del&Z'/B(v)_\BM[1]\ar[d]_{f_{Z*}[1]}\\
Z/B(v)_\BM\ar[r]^{i_*}&X'B(v)_\BM\ar[r]^{j^!}&U/B(v)_\BM\ar[r]^\del&Z/B(v)_\BM[1]
}
\]
ii. Suppose $f$ is smooth. Then the maps $f^!, f_Z^1, f_U^!$ define a map of distinguished localization triangles in $\SH(B)$
\[
\xymatrixcolsep{4.5pt}
\xymatrix{
Z/B(v)_\BM\ar[r]^{i_*}\ar[d]_{f_Z^!}&X/B(v)_\BM\ar[d]_{f^!}\ar[r]^{j^!}&U/B(v)_\BM\ar[d]^{f_U^!}\ar[r]^\del&Z/B(v)_\BM[1]\ar[d]_{f_Z^![1]}\\
Z'/B(v+L_f)_\BM\ar[r]^-{i'_*}&X'/B(v+L_f)_\BM\ar[r]^-{j^{\prime!}}&U'/B(v+L_f)_\BM\ar[r]^-\del&Z'/B(v+L_f)_\BM[1]
}
\]
\end{proposition}

\begin{proof} For (1)-(4), see \cite[Proposition 2.4.2, Example 2.4.3, Theorem 3.2.21, Theorem 4.2.1, Proposition 4.2.2]{DJK}. (5i) follows from \cite[Corollary 2.2.10]{DJK} and (5ii) follows from \cite[Proposition 2.4.7]{DJK}, both together with \cite[Proposition 2.3.3]{CD}.
\end{proof}

\begin{remark}[Relative pull-back]\label{rem:RelPullback} 1. Let 
\[
\xymatrix{
X_0\ar[r]^q\ar[d]^{f_0}&X_1\ar[d]^{f_1}\\
B_0\ar[r]^p&B_1
}
\]
be a Tor-independent cartesian square in $\Sch/B$, with $B_0$ and $B_1$ smooth over $B$. Take $v\in \sK(B_1)$. We define the relative pull-back
\[
p^!{\colon}\sE^\BM(X_1/B_1, v)\to \sE^\BM(X_0/B_0,v)
\]
as follows. Since $B_0$ and $B_1$ are smooth over $B$, the map $p$ is lci. Since the square is Tor-independent and cartesian, $q$ is an lci morphism and $L_q=f^*L_p$, giving   the lci pull-back $q^!{\colon}\sE^\BM(X_1/B_1, v)\to \sE^\BM(X_0/B_1,v+ L_p)$. 

The purity isomorphisms 
\[
p_{X_0!}=p_{B_1!}\circ (pf_{0})_!\cong p_{B_1\#}\circ\Sigma^{-L_{B_1/B}}\circ (pf_{0})_!
\]
and
\[
p_{X_0!}=p_{B_0!}\circ f_{0!}\cong p_{B_0\#}\circ\Sigma^{-L_{B_0/B}}\circ f_{0!}
\]
gives us the isomorphisms
\[
\sE^\BM(X_0/B_1,v+ L_p)\cong
 \sE^\BM(X_0/B,v+ L_p+L_{B_1/B})
\]
and
\[
\sE^\BM(X_0/B_0,v)\cong
 \sE^\BM(X_0/B,v + L_{B_0/B})
\]
The distinguished triangle
\[
p^*L_{B_1/B}\to L_{B_0/B}\to L_p\to p^*L_{B_1/B}
\]
gives the isomorphism
\[
\sE^\BM(X_0/B,v + L_p + L_{B_1/B})\cong \sE^\BM(X_0/B,v + L_{B_0/B}).
\]

Putting these together gives the isomorphism
\[
\vartheta_p{\colon}\sE^\BM(X_0/B_1, v+ L_p)\xrightarrow{\sim} \sE^\BM(X_0/B_0, v)
\]
and we define $p^!:=\vartheta_p\circ q^!$.\\[2pt]
2.  The  relative lci pull-back maps are functorial  with respect to adjacent cartesian squares
\[
\xymatrix{
\bullet\ar[r]^{q_2}\ar[d]&\bullet\ar[r]^{q_1}\ar[d]&\bullet\ar[d]\\
\bullet\ar[r]^{p_2}&\bullet\ar[r]^{p_1}&\bullet
}
\]
when defined. This follows from the  functoriality of pull-back for lci morphisms and the identities
\[
q_2^!\circ  \vartheta_{p_1}=\vartheta_{p_1}\circ q_2^!
\]
and
\[
\vartheta_{p_1p_2}=\vartheta_{p_2}\circ\vartheta_{p_1},
\]
this latter after using the canonical isomorphism $\sE^\BM(-/B, v+L_{p_1}+L_{p_2})\cong
\sE^\BM(-/B,v+L_{p_1p_2})$. The easy  proofs of these identities are left to the reader.
\end{remark}

\begin{ex}\label{ex:WittSheaf} We take $B=\Spec k$, $k$ a perfect field. We will be applying all the above constructions to cohomology/Borel-Moore homology for the sheaf of Witt rings $\sW$.  Recall Morel's theorem \cite[Theorem 5.2.6]{MorelICTP} identifying the heart of the {\em homotopy $t$-structure} on $\SH(k)$ with the abelian category of {\em homotopy modules}; for a homotopy module $\sM_*$ we write $\EM(\sM_*)\in \SH(k)$ for the corresponding object under this identification.

For the definition of the category of homotopy modules, we refer the reader to Morel \cite[\S5]{MorelICTP}; we note here that a homotopy module $\sM_*$ is a family $\{\sM_n\}_{n\in\Z}$ of Nisnevich sheaves of abelian groups on $\Sm/k$, satisfying additional properties, and together with isomorphisms $\epsilon_n{\colon}\sM_n\xrightarrow{\sim} \sHom(\G_m, \sM_{n+1})$; for details see \cite[Definition 5.2.4]{MorelICTP} and surrounding sections.

For $X\in \Sm/k$, the cohomology $\EM(\sM_*)^{p,q}(X)$ is given by the Nisnevich cohomology of $\sM_q$:
\[
\EM(\sM_*)^{p,q}(X)=H^{p-q}_\Nis(X, \sM_q).
\]

We have Nisnevich sheaf of Witt rings $\sW$ on $\Sm/k$. By the purity theorem \cite[Theorem A]{OP} of Ojanguren-Panin, $\sW$ is a so-called {\em unramified sheaf}; together with Morel's computation of $K^{MW}_*(F)[\eta^{-1}]$ for $F$ a field \cite[Lemma 6.3.9]{MorelICTP}, the homotopy module $\sW_*=\sK^{MW}_*[\eta^{-1}]$ has $\sW_n=\sW$ for each $n$, giving the associated Eilenberg-MacLane spectrum $\EM(\sW_*)\in \SH(k)$. $\EM(\sW_*)$ represents Witt-sheaf cohomology via canonical isomorphisms for $X\in \Sm/k$
and $a\in\Z$
\[
H^a_{\Nis}(X, \sW)\cong \EM^{a,0}(\sW_*)(X).
\]
$\EM(\sW_*)$ is an $\SL$-oriented theory \cite[\S3]{LevineEuler}, and since $\times\eta{\colon}\sW_n\to \sW_{n+1}$ is an isomorphism, we have the induced $\eta$-periodicity isomorphism 
\[
\EM^{a+r,r}(\sW_*)(X)\cong \EM^{a,0}(\sW_*)(X).
\]
for all $r\in \Z$. 

For $\sL$ an invertible sheaf on $X$,   there is an $\sL$ twisted version of the Grothendieck-Witt sheaf on $X$,  $\sGW(\sL)$, defined as the sheaf of (virtual) $\sL$-valued non-degenerate bilinear forms. This is a sheaf of modules for the sheaf of Grothendieck-Witt rings on $X$, $\sGW_X$,  and the $\sL$-twisted sheaf of Witt groups on $X$, $\sW(\sL)$, is defined by 
\[
\sW(\sL):=\sGW(\sL)\otimes_{\sGW_X}\sW_X
\]
Equivalently,  $\sW(\sL)= \sGW(\sL)/h\cdot \sGW(\sL)$, where $h$ is the hyperbolic form.   If $L\to X$ is the line bundle with sheaf of sections $\sL$, we also write $\sW(L)$ for $\sW(\sL)$.

 There is a canonical isomorphism
\[
H^a_{\Nis}(X, \sW(\sL))\cong \EM^{a,0}(\sW_*)(X, \sL-\sO_X) 
\]
and more generally
\[
H^a_{\Nis}(X, \sW(\det v))\cong \EM^{a,0}(\sW_*)(X, v), 
\]
for $v\in \sK(X)$ of virtual rank $0$.

We define the $\sL$-twisted $\sW$-Borel-Moore homology of $Z\in \Sch/k$ as
\[
H_a^\BM(Z, \sW(\sL)):=\EM(\sW_*)^\BM_{a, 0}(Z/B, \sO_Z-\sL)\cong
\EM(\sW_*)^\BM_{a+r, r}(Z/B, \sO_Z-\sL).
\]
Since $\EM(\sW_*)$ is $\SL$-oriented, we have, for $v\in \sK(Z)$ of virtual rank 0, canonical isomorphisms
\begin{multline*}
\EM(\sW_*)^\BM_{a, 0}(Z/B, v)\cong \EM(\sW_*)^\BM_{a, 0}(Z/B,  \det v-\sO_Z)\\\cong \EM(\sW_*)^\BM_{a, 0}(Z/B, \sO_Z-\det^{-1} v).
\end{multline*}
where $\det^{-1} v$ is the dual of $\det v$. For $i{\colon}Z\to X$ a closed immersion in $\Sch/k$ with $X$ smooth of dimension $d_X$ over $k$, the isomorphism \eqref{eqn:BMPurityIsoSupp} thus yields the purity isomorphism 
\[
H_{d_X-a}^\BM(Z/B, \sW(\sL))\cong H^a_Z(X, \sW(\omega_{X/k}\otimes \sL))
\]
with $H^*_Z(X,-)$ the usual cohomology with support. If $Z$ is also smooth over $k$ of codimension $c$ in $X$, then \eqref{eqn:PBMurityIso2} gives the isomorphism 
\[
H^a_Z(X, \sW(\sL))\cong H^{a-c}(Z, \sW(\det^{-1}\sN_i\otimes \sL)),
\]
with $\sN_i$ the conormal sheaf. 

By an abuse of notation, we denote the twisted $\EM(\sW_*)$-cohomology of some $Z\in \Sch/k$ and invertible sheaf $\sL$ on $Z$  as $H^a_{\Nis}(Z, \sW(\sL))$:
\[
H^a_{\Nis}(Z, \sW(\sL)):=\EM(\sW_*)^{a,0}(Z, \sL-\sO_Z)\cong \EM(\sW_*)^{a+b,b}(Z, \sL-\sO_Z),
\]
even though $\EM(\sW_*)^{a,0}(Z, \sL-\sO_Z)$ may not in general be given as the cohomology of some sheaf $\sW(\sL)$ on $Z$ for non-smooth $Z$. 
In case $Z$ is smooth over $k$, this does agrees with the twisted Witt-sheaf cohomology, as detailed above.

For a morphism $f{\colon}Y\to X$ in $\Sm/k$, the pullback map $f^*$ as defined above agrees with the one induced by the the pullback on the (twisted) Witt sheaf. If $f$ is proper of relative dimension $d$, applying the purity isomorphism gives the proper pushforward
\[
f_*{\colon}H^{a-d}(Y, \sW(f^*\sL\otimes \omega_f))\to H^a(X, \sW(\sL))
\]
where $\omega_f$ is the relative dualizing sheaf $\omega_{Y/k}\otimes f^*\omega_{X/k}^{-1}$. 

Cup products  are given by
\[
H^a(X, \sW(\sL)\times H^b(X, \sW(\sM))\to H^{a+b}(X, \sW(\sL\otimes\sM)),
\]
again, agreeing with the cup product on cohomology induced by product map on the sheaf level in case $X$ is smooth. Cap products are
\[
 H^\BM_b(X/B, \sW(\sM))\times H^a(X, \sW(\sL) \to H^\BM_{b-a}(X/B, \sW(\sL\otimes\sM)).
\]

For $V\to Z$ a rank $r$-vector bundle, $V=\V(\sV)$, we have the Euler class 
\begin{multline*}
e(V)\in \EM(\sW_*)^{0,0}(Z, \sV)\cong \EM(\sW_*)^{2r,r}(Z, \sV-\sO_Z^r)\\\cong
\EM(\sW_*)^{2r,r}(Z, \det\sV-\sO_Z)\cong H^r(Z, \sW(\det\sV))\\=
H^r(Z, \sW(\det^{-1}V)).
\end{multline*}
Here we write $\det^{-1}V$ to indicate the dual of the line bundle $\det V$, and  the dual arises because $\sV$ is the sheaf of sections of $V^\vee$. 

All the properties and compatibilities for cohomology and Borel-Moore homology in the general setting described above have their corresponding translation for the Witt-sheaf cohomology/Borel-Moore homology. In particular,    Remark~\ref{rem:LToLDual} gives us the natural isomorphisms
\[
H^a(X, \sW(\sL))\cong H^a(X, \sW(\sL\otimes\sM^{\otimes 2})),
\]
\[
H^\BM_a(X/B, \sW(\sL))\cong H^\BM_a(X/B, \sW(\sL\otimes\sM^{\otimes 2}));
\]
the isomorphism for cohomology agreeing with one induced by the canonical isomorphism of sheaves  $ \sW(\sL)\cong \sW(\sL\otimes\sM^{\otimes 2})$.
\end{ex}

\subsection{Rost-Schmid complexes}\label{subsec:RS} For later use, we recall the construction and basic properties of the Rost-Schmid complexes; we limit ourselves here to the case of the complexes with coefficients in the twisted Witt sheaves. We refer the reader to \cite[Chapter 5]{MorelA1} and to \cite{Feld1} for the details about the general construction. Feld \cite{Feld1} constructs the Rost-Schmid complex for a homotopy module $\sM=\{\sM_n\}_n$ (actually for a related object called a Milnor-Witt cycle module) and Morel works in a related ``unstable'' setting;  except for the proof of Lemma~\ref{lem:HtpyModBound}, we will only use the homotopy module $\{\sW_n=\sW\}_n$ associated to the sheaf of Witt groups. 
 
 Fix an $X\in \Sm/k$ and a line bundle $L$ on $X$. For $x\in X$ (not necessarily a closed point) we have the local ring $\sO_{X,x}$ with maximal ideal $\mathfrak{m}_x$ and residue field $k(x)$. Define the {\em orientation line} $\Lambda^X_x$ to be the one-dimensional $k(x)$-vector space $\det^{-1}(\mathfrak{m}_x/\mathfrak{m}_x^2)$, and we write $L_x$ for the  one-dimensional $k(x)$-vector space $i_x^*L$, where $i_x{\colon}x\to X$ is the inclusion. 
 
 For a field $F$ and a dimension one $F$-vector space $L$,  we have $\GW(F;L)$, the Grothendieck-Witt ring of $L$-valued non-degenerate quadratic forms. Product of quadratic forms makes $\GW(F;L)$ a $\GW(F)$-module. Letting $h$ denote as usual the hyperbolic form,  define  $W(F;L)$, the Witt ring of $L$-valued non-degenerate quadratic forms, by
 \[
 W(F;L):=\GW(F;L)/h\cdot \GW(F;L)=\GW(F;L)\otimes_{\GW(F)}W(F). 
 \]
We have the defined the sheaf version $\sW(\sL)$ in Example~\ref{ex:WittSheaf}, viewing  $L$ as an invertible sheaf $\sL$ on $\Spec F$, $W(F;L)$ is the same as $H^0(\Spec F, \sW(\sL))$.
 
 For $z\in X^{(n)}$ with closure $Z\subset X$ and for $y\in Z^{(1)}\subset X^{(n+1)}$, define the map
 \[
 \del_{z,y}{\colon}W(k(z); \Lambda^X_z\otimes_{k(z)}L_z)\to W(k(y); \Lambda^X_y\otimes_{k(y)}L_y)
 \]
(following \cite[\S 5.1]{MorelA1}) as follows. First suppose that $y$ is a regular point on $Z$, in other words, that the local ring $\sO_{Z,y}$ is a discrete valuation ring. Let $t\in \sO_{Z,y}$ be a uniformizing parameter. Since $Z$ is regular at $y$, the ideal  $\sI_{Z,y}$ of $Z$ in $\sO_{X,y}$ is cut out by a regular sequence and we have a canonical isomorphism
 \[
 \det^{-1}(\sI_{Z,y}/\sI_{Z,y}^2)\otimes_{\sO_{Z,y}}k(z)\cong \Lambda^X_z.
\]
 Moreover, let $\del/\del t$ be the basis of the $k(y)$ dual space of $(t)/(t^2)$, dual to the image of $t$. Then product with $\del/\del t$ defines an isomorphism
 \[
\del/\del t\wedge(-) {\colon}\det^{-1}(\sI_{Z,y}/\sI_{Z,y}^2)\otimes_{\sO_{Z,y}}k(y)\xrightarrow{\sim}\Lambda^X_y
\]
We have the additive map $\del_t{\colon}W(k(z))\to W(k(y))$ characterized by $\del_t(\<t\>)=1$, $\del_t(\<u\>)=0$ for $u\in \sO_{Z,y}^\times$. The map $\del_{z,y}$ is characterized by 
\[
\del_{z,y}(\alpha\otimes i_z^*\tau)=\del_t(\alpha)\otimes \del/\del t \wedge i_y^*\tau
\]
for $\tau$ a generating section of $\det^{-1}(\sI_{Z,y}/\sI_{Z,y}^2)$. 

If $y$ is not a regular point of $Z$, then one passes to the normalization, makes a similar construction to the above for each point $y'$ lying over $y$ and then takes a suitable pushforward to define $\del_{z,y}$. See \cite[\S 5.1]{MorelA1} for details.

 \begin{definition}[\hbox{\cite[Definition 5.7, Definition 5.11, Theorem 5.31]{MorelA1}}]\label{def:RS}  For  $X\in \Sm/k$ with a line bundle $L$ on $X$, the {\em Rost-Schmid complex} for $\sW(L)$ on $X$ is the complex $C^*_{RS}(X, \sW(L))$ with
 \[
 C^n_{RS}(X, \sW(L))=\oplus_{x\in X^{(n)}}W(k(x); \Lambda^X_x\otimes_{k(x)}L_x).
 \]
The differential $d_{RS}^n{\colon}C^n_{RS}(X, \sW(L))\to C^{n+1}_{RS}(X, \sW(L))$ is  given by
\[
d_{RS}^n(\sum_{z\in X^{(n)}}a_z):=\sum_{z\in X^{(n)}, y\in \bar{z}^{(1)}}\del_{z,y}(a_z). 
\]
One shows that this sum is finite and $d^{n+1}d^n=0$, giving a well-defined complex $(C^*_{RS}(X, \sW(L)), d^*_{RS})$.  See \cite[\S 5.1, \S 5.2]{MorelA1} for details. 
\end{definition}

\begin{Not}\label{Not:RSNotation}
For $X\in \Sm/k$, and  $Z\subset X$ a codimension $c$ subvariety with generic point $z\in Z$,   $\alpha\in W(k(z))$ and $\tau$ a generator of $\Lambda^X_z\otimes i_z^*L$, we have the corresponding $\Lambda^X_z\otimes i_z^*L$-valued quadratic form $\alpha\otimes \tau\in W(k(z); \Lambda^X_z\otimes i_z^*L)$. We write  $\alpha\otimes \tau\text{ on }Z$ for the element $\alpha\otimes\tau$ in the summand $W(k(z); \Lambda^X_z\otimes i_z^*L)$ of $C_{RS}^c(X, \sW(L))$ indexed by $z\in X^{(c)}$. If $W$ is a closed subset of $X$ containing $Z$, we also write $\alpha\otimes \tau\text{ on }Z$ for the corresponding element in the degree $c$ term of the Rost-Schmid complex with supports, $C_{RS, W}^c(X, \sW(L))\subset C_{RS}^c(Y, \sW(L))$.
\end{Not}

For $f{\colon}U\to X$ an \'etale map in $\Sm/k$, we have a well-defined pullback  of complexes $f^*{\colon}
C^*_{RS}(X, \sW(L))\to C^*_{RS}(Y, \sW(f^*L))$; if $f$ is the inclusion of a Zariski open subscheme $U$, the map $f^*{\colon}C^n_{RS}(X, \sW(L))\to C^n_{RS}(U, \sW(f^*L))$ is just the projection from $\oplus_{x\in X^{(n)}}W(k(x); L_x\otimes_{k(x)}\Lambda_x)$ to $\oplus_{x\in U^{(n)}}W(k(x); L_x\otimes_{k(x)}\Lambda_x)$. As this is clearly surjective, we have the identification of the homotopy fiber of the restriction map from $X$ to $U:=X\setminus Z$ with the subcomplex $C^*_{RS, Z}(X, \sW(L))$,  where 
\[
C^n_{RS, Z}(X, \sW(L))=\oplus_{x\in X^{(n)}\cap Z}W(k(x); \Lambda^X_x\otimes_{k(x)}L_x).
\]

Clearly the assignment $U\mapsto C^*_{RS}(U, \sW(L_{|U}))$ defines a presheaf of complexes on $X_\Nis$. We denote this presheaf by $\sC^*_{RS}(\sW(L))$, assuming the underlying scheme $X$ is clear from the context. We have the evident augmentation $\epsilon{\colon}\sW(L)\to \sC^*_{RS}(\sW(L))$.

\begin{theorem}[Morel \hbox{\cite[Theorem 5.41]{MorelA1}}]\label{thm:GerstenConjWittSheaf} For $X\in \Sm/k$ with line bundle $L$, the augmentation $\epsilon{\colon}\sW(L)\to \sC^*_{RS}(\sW(L))$ is a quasi-isomorphism for the Zariski topology. In other words,  for each $x\in X$, the complex $C^*_{RS}(\Spec\sO_{X,x}, \sW(L\otimes_{\sO_X}\sO_{X,x}))$ has
\[
H^n(C^*_{RS}(\Spec\sO_{X,x}, \sW(L\otimes_{\sO_X}\sO_{X,x})))=\begin{cases}0&\text{ for }n>0\\
\sW(L)_x&\text{ for }n=0
\end{cases}
\]
where the augmentation $\epsilon$ induces this identity for $n=0$. 
\end{theorem}

\begin{corollary}\label{cor:RSFacts} Take $X\in \Sm/k$ with line bundle $L$.\\[5pt]
1. For each $n\ge0$, there are natural isomorphisms
\[
H^n_\Nis(X, \sW(L))\cong H^n_\Zar(X, \sW(L))\cong H^n(C^*_{RS}(X, \sW(L))).
\]
In particular $H^n(X, \sW(L))=0$ for $n>\dim_kX$.\\[2pt]
2. Let $Z\subset X$ be a closed subset. Then for each $n\ge0$, we have a natural isomorphism
\[
H^n_Z(X, \sW(L))\cong H^n(C^*_{RS,Z}(X, \sW(L))).
\]
3. Let $Z\subset X$ be a closed subset and suppose $X$ has pure dimension $d$. Then for each $n\ge0$, we have a natural isomorphism
\[
H^\BM_n(Z, \sW(L))\cong H^{d-n}(C^*_{RS,Z}(X, \sW(L^{-1}\otimes\omega_{X/k}))),
\]
where $\omega_{X/k})$ is the canonical sheaf on $X$.
\end{corollary}

\begin{proof} For (1), the presheaves $(f{\colon}U\to X)\mapsto C^n_{RS}(U, \sW(f^*L))$ on $X_\Nis$ are sheaves  in the Nisnevich topology, and are acyclic in the Zariski and Nisnevich topologies, by \cite[Lemma 5.42]{MorelA1}. By Theorem~\ref{thm:GerstenConjWittSheaf}, 
\[
\sW(L)\xrightarrow{\epsilon}\sC^*_{RS}(\sW(L))
\]
is an acyclic resolution of $\sW(L)$ for both the Zariski and Nisnevich topologies, from which (1) follows immediately. 

 (2) follows from (1) and the fact that $C^*_{RS,Z}(X, \sW(L))$ represents the homotopy fiber of the restriction map $C^*_{RS}(X, \sW(L))\to C^*_{RS}(X\setminus Z, \sW(L))$, as discussed above. For (3), this follows from the isomorphism  \eqref{eqn:BMPurityIsoSupp}, together with the $\SL$-orientability and $\eta$-periodicity of $\EM(\sW_*)$, which together yield the identity
\[
H^\BM_n(Z, \sW(L))=H^{d-n}_Z(X, \sW(L^{-1}\otimes\omega_{X/k}))).
\]
\end{proof}

\begin{remark}\label{rem:CycModFunct} The functorialities of the Rost-Schmid cycle complexes over a perfect field has been discussed by Feld \cite{Feld1}. The theory has been extended to a general base-scheme by D\'eglise-Feld-Jin in \cite{DFJ}. As a special case, \cite[Proposition 3.4.3]{DFJ} implies that the isomorphisms of Corollary~\ref{cor:RSFacts}  are compatible with respect to pullback in Witt sheaf cohomology for morphisms $f{\colon}Y\to X$ in $\Sm/k$ and with respect to pushforward maps in Witt sheaf Borel-Moore homology for proper morphisms  in $\Sm/k$ and  pullback maps in Witt sheaf Borel-Moore homology for smooth morphisms  in $\Sm/k$.
\end{remark}

\section{Equivariant Borel-Moore Witt homology}\label{sec:EquivBMWitt}
In this section we recall  the construction of the $G$-equivariant $\sE$-cohomology and $G$-equivariant $\sE$-Borel-Moore homology functors for $\sE\in \SH(B)$ by Di Lorenzo and Mantovani  \cite[\S 1.2]{DiLorenzoMantovani}.  As described in {\it loc.\,cit.}, the theory is only for $\sE$ being the Eilenberg-MacLane spectrum of a homotopy module, and the theory we present here is in essence no more general  than this, as  $\sE$ needs a certain ``boundedness'' property for the resulting theory to be reasonable.  The general theory is currently being constructed by Mantovani \cite{Mantovani} and also appears in \cite{DAngelo}. The theory we give here is a twist of that of Di Lorenzo and Mantovani; this modification is made so that the theory will  compatible with the operation of forgetting the $G$-action on the Borel-Moore homology.

The construction of equivariant cohomology and Borel-Moore homology in the algebraic setting follows the method of Totaro \cite{Totaro}, Edidin-Graham \cite{EGEquivIntThy}, and Morel-Voevodsky \cite[\S 4]{MorelVoevodsky}.  Let $F_n\cong \A^n_B$ be the fundamental (left) representation of $\GL_n:=\GL_n/B$. For $j\ge0$, let $V_j=F_n^{n+j}$; we identify $V_j$ with the scheme of $n\times n+j$ matrices with the  usual left $\GL_n$-action. Let $W_j\subset V_j$ be the closed subscheme of matrices of rank $<n$ and let $E_j\GL_n:=V_j\setminus W_j$.   We have $\codim_{V_j}W_j=j+1$ and $\GL_n$ acts freely on $E_j\GL_n$ with quotient the Grassmannian $\Gr_B(n, n+j)$.   

Let $p_j{\colon}V_{j+1}=V_j\times F_n\to V_j$ be the projection to $V_j$, and let $i_j{\colon}V_j\to V_{j+1}=V_j\times F_n$ be the 0-section. The inclusion $E_j\GL_n\times F_n\hookrightarrow V_j\times F_n=V_{j+1}$ has image contained in $E_{j+1}\GL_n$; let $\eta_j{\colon}E_j\GL_n\times F_n\to E_{j+1}\GL_n$ be the corresponding open immersion, and  $p_j{\colon}E_j\GL_n\times F_n\to E_j\GL_n$, $s_j{\colon}E_j\GL_n\to E_j\GL_n\times F_n$ the projection and 0-section. Let $i_j{\colon}E_j\GL_n\to E_{j+1}\GL_n$ be the composition $\eta_j\circ s_j$.  We let $*_j\in E_j\GL_n$ be the image of the identify matrix in $E_0\GL_n=\GL_n$ under the composed inclusion $E_0\GL_n\to E_j\GL_n$.

We consider $E_j\GL_n$ as an object of $\Sm^{\GL_n}/B$. For $G$ a closed subgroup-scheme of $\GL_n$, smooth over $B$, we let $E_jG$ be the  object of $\Sm^G/B$ formed by restricting the  $\GL_n$-action on $E_j\GL_n$ to $G$. We will assume that the homogeneous space $G\backslash\GL_n$ exists as a smooth quasi-projective $B$-scheme and that $\GL_n\to G\backslash\GL_n$ is a $G$-torsor; by \cite[Proposition 7.17, Theorem 7.18 and its proof]{Milne}, this is the case if $B=\Spec k$, $k$ a field. 

 Following the proof of  \cite[Lemma 9]{EGEquivIntThy}  our assumptions on  $\GL_n\to G\backslash \GL_n$ implies that $G\backslash E_jG$ exists as a quasi-projective  $B$-scheme and the quotient map $q_{G,j}{\colon}E_jG\to B_jG$ is an \'etale $G$-torsor. Let $N_{G,j}:=G\backslash E_jG\times F_n$ with quotient map $q_{G,j}^0{\colon}E_jG\times F_n\to N_{G,j}$. The maps $i_j$, $s_j$, $\eta_j$ and $p_j$ induce corresponding regular embeddings $i_{G,j}{\colon}B_jG\to B_{j+1}G$, $s_{G,j}{\colon}B_jG\to N_{G,j}$, open immersion 
$\eta_{G,j}{\colon}N_{G,j}\to B_{j+1}G$,  and smooth morphism $p_{G,j}{\colon}N_{G,j}\to B_jG$, making $N_{G,j}$ the vector bundle over $B_jG$ corresponding via descent to the representation $G\hookrightarrow \GL_n$,   with zero-section $s_{G,j}$. We note that the vector bundle $N_{G,j}$ on $B_jG$ is canonically isomorphic to the normal bundle of $i_{G,j}$. Here is a diagram describing the construction.
\begin{equation}\label{eqn:QuotientDiag1}
\xymatrixrowsep{10pt}
\xymatrix{
E_jG\times F_n\ar[dd]_{p_j}\ar@{^(->}[ddrr]^{\eta_j}\ar[dddr]|(.61)\hole_(.67){q_{G,j}^0\hskip-7pt}\\\\
E_jG\ar[rr]^{i_j}\ar@/_10pt/[uu]_(.35){\hskip-1pt s_j}\ar[dddr]_{q_{G,j}}&&E_{j+1}G\ar[dddr]^{q_{G,j+1}}\\
&N_{G,j}\ar[dd]_(.35){p_{G,j}}\ar@{^(->}[ddrr]^{\eta_{G,j}}\\\\
&B_jG\ar@/_10pt/[uu]_(.4){s_{G,j}}\ar[rr]^{i_{G,j}}&&B_{j+1}G
}
\end{equation}

\begin{definition}
We let $\Sch^G_q/B$ be the full subcategory of $\Sch^G/B$ consisting of those $X$ for which all the fppf quotients  $G\backslash X\times_BE_jG$ are represented in $\Sch/B$; here $G$ acts diagonally on $X\times_BE_jG$. In case $B=\Spec k$, $k$ a field, we have $\Sch^G_q/B=\Sch^G/B$, by \cite[Lemma 9, Proposition 23]{EGEquivIntThy}. We similarly define $\Sm_q/B:=\Sm/B\cap \Sch^G_q/B$.
\end{definition}

Take $X\in \Sch_q^G/B$. We denote the quotient  $G\backslash X\times_BE_jG$ by  $X\times^GE_jG$. Let $N_{X,G, j}\subset X\times^GE_{j+1}G$ be the open subscheme $G\backslash X\times_B E_jG\times F_n$ with open immersion $\eta_{X,G,j}{\colon}N_{X, G,j}\to X\times^GE_{j+1}G$ induced by $\eta_j$.  The maps $i_j$, $p_j$ and $s_j$ induce maps 
\begin{gather*}
i_{X,G,j}{\colon}X\times^G E_jG\to X\times^GE_{j+1}G, \ p_{X,G,j}{\colon}N_{X, G,j}\to X\times^G E_jG,\\ 
s_{X,G,j}{\colon}X\times^GE_jG\to N_{X, G,j},
\end{gather*}
This gives $p_{X,G,j}{\colon}N_{X, G,j}\to X\times^GE_jG$ the structure of a vector bundle   with 0-section $s_{X,G,j}$ and we have $i_{X,G,j}=\eta_{X,G,j}\circ s_{X,G,j}$.

This gives us the diagram
\begin{equation}\label{eqn:QuotientDiag2}
\xymatrixrowsep{10pt}
\xymatrixcolsep{30pt}
\xymatrix{
X\times_BE_jG\times F_n\ar[dd]_{\id_X\times p_j}\ar@{^(->}[ddrr]^{\id_X\times \eta_j}\ar[dddr]|(.61)\hole_(.67){q_{X,G,j}^0\hskip-7pt}\\\\
X\times_BE_jG\ar[rr]^{\id_X\times i_j}\ar@/_10pt/[uu]_(.35){\hskip-1pt \id_X\times s_j}\ar[dddr]_{q_{X,G,j}}&&X\times_BE_{j+1}G\ar[dddr]^{q_{X,G,j+1}}\\
&N_{X,G,j}\ar[dd]_(.35){p_{X,G,j}}\ar@{^(->}[ddrr]^{\eta_{X,G,j}}\\\\
&X\times^GE_jG\ar@/_10pt/[uu]_(.4){s_{X,G,j}}\ar[rr]^{i_{X,G,j}}&&X\times^GE_{j+1}G
}
\end{equation}

The structure map $\pi_X{\colon}X\to B$ induces the projections $\pi_{X,j}{\colon}X\times_BE_jG\to E_jG$, $\pi^0_{X,j}{\colon}X\times_BE_jG\times F_n\to E_jG\times F_n$ and maps $\pi_{X,G,j}{\colon}X\times^GE_jG\to B_jG$ and $\pi^0_{X,G,j}{\colon}N_{X, G,j}\to N_{G,j}$. These give a map of diagram~\eqref{eqn:QuotientDiag2} to diagram~\eqref{eqn:QuotientDiag1}.   We write $EG$, $BG$, $X\times^GEG$, etc., for the corresponding Ind-objects $\{E_jG, i_j\}_j$, $\{B_jG, i_{G,j}\}_j$, $\{X\times^GE_jG, i_{X,G,j}\}_j$, etc.,. 

\begin{lemma}   $B_jG$ is smooth over $B$.
\end{lemma}

\begin{proof} Recall we are assuming that $G$ is smooth over $B$.
 For $G=\GL_n/B$, $B_jG$ is the Grassmannian $\Gr_B(n, n+j)$, which is smooth over $B$. For a general $G\subset \GL_n/B$, smooth over $B$,  we have the \'etale locally trivial fiber bundle
$B_jG\to \Gr_B(n,n+j)$ with fiber $\GL_n/G$, and \'etale locally trivial fiber bundle $\GL_n\to \GL_n/G$ with fiber $G$. Since $G$ is smooth over $B$, so is $ \GL_n/G$ and thus so is $B_jG$.
\end{proof}

\begin{remark} Morel-Voevodsky \cite[\S 4.3]{MorelVoevodsky} give a general version of the construction outlined above. Their construction depends on choices (including the ones we have used above), but it is shown in {\em loc.\,cit.}\! that  the resulting ind-object in the unstable motivic homotopy category $\sH(B)$ is independent of the choices.
\end{remark}

For $X\in \Sch^G/B$,  let $\sK^G(X)$ denote fundamental groupoid of  the $K$-theory space of perfect, $G$-linearized complexes on $X$; if $G$ is the trivial group, this is just the $K$-theory groupoid $\sK(X)$ mentioned in \S\ref{sec:CohBMHom}. Taking $v\in \sK^G(X)$  gives us for each $j$ the $G$-linearized perfect complex $p_1^*v$ on $X\times E_jG$, and by descent the perfect complex $v_j$ on $X\times^GE_jG$. Sending $v$ to $v_j$ defines the  map, natural in $X$,
\[
{(-)}_j{\colon}\sK^G(X)\to \sK(X\times^GE_jG).
\]
Similarly, for $V\to X$ a $G$-linearized vector bundle, we have the vector bundle $V_j\to X\times^GE_jG$ induced from the $G$-linearized vector bundle $p_1^*V\to X\times_BE_jG$ by descent. If $V=\V(\sV)$ for a $G$-linearized locally free sheaf $\sV$ on $X$, then $V_j=\V(\sV_j)$.

As $\sE$-cohomology $\sE^{a,b}(Y, v)$ is contravariantly functorial in $(Y,v)$, the ind-system $X\times^GEG:=\{X\times^GE_jG, i_{X,G,j}\}_j$ gives the pro-system $\{\sE^{a,b}(X\times^GE_jG, v_j), i_{X,G,j}^*\}_j$ for $v\in \sK^G(X)$. 
\begin{definition}\label{def:EquivCoh} We define
\[
\sE^{a,b}_G(X, v)=\lim_{\substack{\leftarrow\\j}}\sE^{a,b}(X\times^GE_jG, v_j).
\]
We write $\sE_G(X, v)$ for $\sE^{0,0}_G(X, v)$ and $\sE^{**}(BG, v)$ for $\sE^{**}_G(B, v)$.
 
For $B=\Spec k$,  $\sE=\EM(\sW_*)$, and $\sL$ a $G$-linearized invertible sheaf on $X$, we set
\[
H^a_G(X, \sW(\sL)):=\EM(\sW_*)^{a, 0}_G(X, \sL-\sO_X),
\]
with   $\sO_X$ having the canonical $G$-linearization. We write $H^a(BG, \sW(\sL))$ for the equivariant cohomology
$H^a_G(\Spec k, \sW(\sL))$
\end{definition}

\begin{remark}  1. The contravariant functoriality of $\sE$-cohomology induces a corresponding functoriality for 
$\sE^{a,b}_G(X, v)$ and hence for $H^*_G(X, \sW(\sL))$.\\[2pt]
2. The definition we give here is not the correct one for general $\sE$, due at the very least to the usual lack of exactness of the inverse limit. However, as we shall see below (Remark~\ref{rem:StabilityPurity}), for $X$ smooth, the inverse system used to define $H^a_G(X, \sW(\sL))$ is eventually constant, so this problem does not arise. In addition, the exactness shows up only in the localization sequence for a decomposition of a scheme into a closed subscheme and open complement, and this takes place in the setting of Borel-Moore homology. For Borel-Moore homology, we will use the Gysin maps to define an inverse system and define the equivariant Borel-Moore homology as the limit.  The fact that the corresponding inverse system is eventually constant is given by 
Lemma~\ref{lem:StrBoundImpBound} and Lemma~\ref{lem:HtpyModBound} below; see also \cite[\S 1.2]{DiLorenzoMantovani} for the general case of the Eilenberg-MacLane spectrum of a homotopy module. We give a proof of this fact here in Lemma~\ref{lem:HtpyModBound} below. This gives us the long exact localization sequence in the equivariant Borel-Moore homology for $\EM(\sW_*)$.
\end{remark}

\begin{proposition}\label{prop:EquivEulerClass}  Let $V\to X$ be a $G$-linearized vector bundle on $X\in \Sch^G_q/B$,  and let $\sV$ be the locally free sheaf of sections of $V^\vee$, with its canonical  $G$-linearization. Then the collection of Euler classes $\{e(V_j)\in \sE(X\times^GE_jG, \sV_j)\}_j$ yields a well-defined element $e_G(V)\in \sE_G(X, \sV)$. Moreover, the assignment $(V\to X)\mapsto e_G(V)\in \sE_G(X, \sV)$ is natural in the  $G$-linearized vector bundle $V\to X$: given a morphism  $f{\colon}Y\to X$ in $\Sch^G_q/B$, a  $G$-linearized vector bundle $W\to Y$ and an isomorphism of  $G$-linearized vector bundles $\tilde{f}{\colon}W\to f^*V$, the induced map $(f, \tilde{f})^*{\colon}\sE_G(X, \sV)\to \sE_G(Y, \sW)$ sends $e_G(V)$ to $e_G(W)$.
\end{proposition}

\begin{proof} For each $j$ we have the canonical isomorphisms $i_{X,G,j}^*V_{j+1}\cong V_j$,
$i_{X,G,j}^*\sV_{j+1}\cong \sV_j$, giving the identity
\[
i_{X,G,j}^*(e(V_{j+1}))=e(V_j)\in \sE(X\times^GE_jG, \sV_j),
\]
which proves the first assertion. For the second, the data $(f, \tilde{f})$ gives for each $j$ the identity
\[
(f\times^G\id, \tilde{f}\times^G\id)^*(e(V_j))=e(W_j)
\]
whence follows $(f, \tilde{f})^*(e_G(V))=e_G(W)$.
\end{proof} 

\begin{ex} Let $\rho{\colon}G\to \GL_N$ be a representation of $G$ over $B$, with total space $V(\rho)$. We consider $V(\rho)$ as a $G$-linearized vector bundle on $B$, giving the Euler class $e_G(\rho):=e_G(V(\rho))\in \sE_G(B, \sV(\rho))=\sE(BG, \sV(\rho))$. Here $\sV(\rho)$ is the $G$-linearized sheaf of sections of $V(\rho)^\vee$.
\end{ex}

\begin{definition}[Generic Euler class]\label{def:GenEulerClass}   Let $G$ be  group scheme over $k$, satisfying the hypotheses of Lemma~\ref{lem:LocalTrivGRep}. Let $V$ be a $G$-linearized vector bundle on some connected $Y\in \Sch^G/k$, with $G$-linearized locally free sheaf of section $\sV$. We suppose we are in Case 1, Case 2 or Case 3 of 
Construction~\ref{const:GenEuler} for $\sV$, so the generic representation class $[\sV^{gen}]$ is defined. 

Suppose $\sV^{gen}$ has even rank $2r$ (this is always the case in Case 2). Choose a representative representation $V^{gen}$ for $[\sV^{gen}]$, giving us   the Euler class $e_G(V^{gen})\in H^{2r}(BG, \sW(\sL))$, where $\sL=\det\sV^{gen}$.  Using the canonical isomorphism $\sW(\sL^{\otimes 2})\cong \sW$, we have $e_G(V^{gen})^2\in H^{4r}(BG, \sW)$ and we define $[e^{gen}_G(V)]\subset H^{4r}(BG, \sW)$ to be the subset
\[
[e^{gen}_G(V)]=\{u\cdot e_G(V^{gen})^2\mid u\in  W(k)^\times\}\subset H^{4r}(BG, \sW).
\]
\end{definition}
Note that $[e^{gen}_G(V)]$ is independent of the choice of representative $V^{gen}$ for $[\sV^{gen}]$ and depends only on the isomorphism class of $V$ as a $G$-linearized  vector bundle
on $Y$. Similarly, the localization
\[
H^*(BG, \sW)[[e^{gen}_G(V)]^{-1}]
\]
is equal to $H^*(BG, \sW)[1/e^{gen_G}(V)]$ for any $e^{gen}_G(V)\in [e^{gen}_G(V)]$, and we will often use $H^*(BG, \sW)[1/e^{gen}_G(V)]$ to denote $H^*(BG, \sW)[[e^{gen}_G(V)]^{-1}]$. Finally, if we have a trivialization of $\det V^{gen}$, then choosing one gives us  $e_G(V^{gen})\in H^{2r}(BG, \sW)$ and we have
\[
H^*(BG, \sW)[[e^{gen}_G(V)]^{-1}]=H^*(BG, \sW)[e_G(V^{gen})^{-1}].
\]

Having defined $G$-equivariant twisted cohomology, we now proceed to define the $G$-equivariant twisted Borel-Moore homology.

Fix $X\in \Sch^G_q/B$. For each $j$, we have the Tor-independent cartesian diagram
\[
\xymatrix{
X\times^GE_jG\ar[r]^{i_{X,G,j}}\ar[d]&X\times^GE_{j+1}G\ar[d]\\
B_jG\ar[r]^{i_{G,j}}&B_{j+1}G
}
\]
with $i_{G,j}$ a regular embedding.  Given $v\in \sK^G(X)$,  we have a canonical isomorphism $i_{X,G,j}^*v_{j+1}\cong v_j$, 
giving the pro-system $\{v_j\in \sK(X\times^GE_jG), i_{X,G,j}^*\}_j$. 
This gives us the relative pull-back map (see Remark~\ref{rem:RelPullback})
\[
i_{G,j}^!{\colon}\sE^\BM_{a,b}(X\times^GE_{j+1}G/B_{j+1}G, v_{j+1})\to 
\sE^\BM_{a,b}(X\times^GE_jG/B_jG, v_j)
\]
defining the pro-system $\{\sE^\BM_{a,b}(X\times^GE_jG/B_jG, v_j), i_{G,j}^!\}_j$. 

Similarly, we have the Tor-independent cartesian square
\[
\xymatrix{
X\times E_jG\ar[r]^{i_{X,j}}\ar[d]&X\times E_{j+1}G\ar[d]\\
E_jG\ar[r]^{i_j}&E_{j+1}G
}
\]
giving the pro-system $\{\sE^\BM_{a,b}(X\times E_jG/E_jG,  v_j), i_j^!\}_j$.

\begin{definition}\label{def:BMWIttHomology} For $X\in \Sch^G/B$ and $v\in \sK^G(X)$, define
\[
\sE^\BM_{G, a,b}(X/BG, v):=\lim_{\substack{\leftarrow\\ j}} \sE^\BM_{a,b}(X\times^GE_jG/B_jG, v_j).
\]
We write $\sE^\BM_G(X/BG, v)$ for $\sE^\BM_{G, 0,0}(X/BG, v)$, $\sE^\BM_{a,b}(BG, v)$ for $\sE^\BM_{G, a,b}(B/BG, v)$, and $H^\BM_{G,a,b}(X/BG, v)$ for $(1_B)^\BM_{G,a,b}(X/BG, v)$, etc. 

 For $B=\Spec k$, and $\sL$ a $G$-linearized invertible sheaf on $X$, define
\[
H^\BM_{G, a}(X, \sW(\sL))
:= \EM(\sW_*)^\BM_{a,0}(X/BG, \sO_X-\sL)
\]
and write $H^\BM_a(BG, \sW(\sL))$ for $H^\BM_{G, a}(\Spec k, \sW(\sL))$. If we need to explicitly indicate the base-field $k$, we write $H^\BM_{G, a}(X/k, \sW(\sL))$ for $H^\BM_{G, a}(X, \sW(\sL))$.
\end{definition}

\begin{lemma}\label{lem:Comparison} 1. The transverse cartesian square
\[
\xymatrix{
X\times E_jG\ar[d]\ar[r]^{q_{X,G,j}}& X\times^GE_jG\ar[d]\\
E_jG\ar[r]^{q_{G,j}}&B_jG
}
\]
induces the relative smooth pullback map
\[
q_{G,j}^!{\colon}\sE^\BM_{a,b}(X\times^GE_jG/B_jG, v_j)\to \sE^\BM_{a,b}(X\times E_jG/E_jG,  v)
\]
which gives rise to a map of pro-systems
\[
q_{G,*}^!{\colon}\{\sE^\BM_{a,b}(X\times^GE_jG/B_jG, v_j), i_{G,j}^!\}_j\to
\{\sE^\BM_{a,b}(X\times E_jG/E_jG,  v), i_j^!\}_j
\]
2. The transverse cartesian square
\[
\xymatrix{
X\times E_jG\ar[d]\ar[r]^-{p_1}& X\ar[d]\\
E_jG\ar[r]^-{p_{E_jG}}&B
}
\]
induces the relative smooth pullback map
\[
p_{E_jG}^!{\colon} \sE^\BM_{a,b}(X/B, v)\to  \sE^\BM_{a,b}(X\times E_jG/E_jG, v).
\]
giving a map of pro-systems
\[
p_{E_*G}^!{\colon}\{\sE^\BM_{a,b}(X/B, v)\}_j\to \{\sE^\BM_{a,b}(X\times E_jG/E_jG, v), i_j^!\}_j
\]
\end{lemma}

\begin{proof}  The maps  (1) and (2) are special cases of relative smooth pull-back (see Remark~\ref{rem:RelPullback}). The fact that the maps in (1) and (2) induce maps of pro-systems follows from the  functoriality of relative lci pull-back, as noted in Remark~\ref{rem:RelPullback}.
\end{proof}

\begin{proposition}\label{prop:EquivBMProperties}  Fix a commutative ring spectrum $\sE\in \SH(B)$. Let $f{\colon}X\to Y$ be  morphism in $\Sch_q^G/B$ and take $v\in\sK^G(Y)$. Consider the sequence of induced morphisms
\[
f_j:=f\times^G\id{\colon} X\times^GE_jG\to Y\times^GE_jG.
\]
\\[5pt]
1.  Suppose $f$ is an lci morphism. Then the $f_j$ are all lci morphisms and induce a map of pro-systems
\[
\{f_j^!\}_j{\colon} \{\sE^\BM_{a,b}(Y\times^GE_jG/B_jG, v_j)\}_j\to 
\{\sE^\BM_{a,b}(X\times^GE_jG/B_jG, v_j+L_{f, j})\}_j,
\]
giving rise to a well-defined map
\[
f^!{\colon} \sE^\BM_{G, a,b}(Y/BG, v)\to \sE^\BM_{G, a,b}(X/BG, v+L_f),
\]
2. Suppose $f$ is a proper morphism. Then the $f_j$ are all proper morphisms and induce a map of pro-systems
\[
\{f_{j*}\}_j{\colon} \{\sE^\BM_{a,b}(X\times^GE_jG/B_jG, v_j)\}_j\to 
\{\sE^\BM_{a,b}(Y\times^GE_jG/B_jG, v_j)\}_j,
\]
giving rise to a well-defined map
\[
f_*{\colon} \sE^\BM_{G, a,b}(X/BG, v)\to \sE^\BM_{G,\ a,b}(Y/BG, v),
\]
3. For general $f$, the maps $f_j$ induce a map of pro-systems
\[
f_j^*{\colon}\{\sE^{a,b}(Y\times^GE_jG, v_j)\}_j\to \{\sE^{a,b}(X\times^GE_jG, v_j)\}_j,
\]
giving rise to a well-defined map
\[
f^*{\colon} \sE^{a,b}_G(Y, v)\to \sE_G^{a,b}(X, v).
\]
4. Given   $v,w\in \sK^G(Y)$, the products
\begin{multline*}
\sE^{a,b}(Y\times^GE_jG, v_j)\times
\sE^{c,d}(Y\times^GE_jG, w_j)\xrightarrow{\cup} 
\sE^{a+c,b+d}(Y\times^GE_jG,v_j+ w_j)
\end{multline*}
and cap products
\begin{multline*}
\sE^\BM_{a,b}(Y\times^GE_jG/B_jG,  v_j)\times
\sE^{c,d}(Y\times^GE_jG, w_j)\\\xrightarrow{\cap} 
\sE^\BM_{a-c,b-d}(Y\times^GE_jG/B_jG,  v_j-w_j)
\end{multline*}
define maps of pro-systems, inducing a  product
\[
 \sE^{a,b}_G(Y, v)\times  \sE^{c,d}_G(Y, w)\xrightarrow{\cup}  \sE^{ a+c,b+d}_G(Y, v+w)
\]
and a cap product
\[
 \sE^\BM_{G, a,b}(Y/BG, v)\times  \sE^{c,d}_G(Y, w)\xrightarrow{\cap}  \sE^\BM_{G, a-c,b-d}(Y/BG, v-w).
\]
5.  (Fundamental classes, Poincar\'e duality) Let $1^\sE_G\in \sE^{0,0}_G(B)=\sE^\BM_G(B/BG)$ be the family of classes of $1^\sE_{B_jG}$. Let $p_X{\colon}X\to B$ be an lci scheme over $B$, in $\Sch^G_q/B$. Define the {\em fundamental class} $[X]^\sE_G\in \sE^\BM_G(X/BG, L_{X/B})$ by 
\[
[X]^\sE_G:=p_X^!(1^\sE_G).
\]
i. For $f{\colon}Y\to X$ an lci morphism of lci-schemes in $\Sch^G_q/B$, we have $f^!([X]^\sE_G)=[Y]^\sE_G$. \\[2pt]
ii. Take $v\in \sK^G(X)$. The cap product with $[X]^\sE_G$
\[
[X]^\sE_G\cap-{\colon}\sE^{a,b}_G(X, v)\to \sE^\BM_{G, -a,-b}(X, L_{X/B}-v)
\]
is an isomorphism for $X\in \Sm^G_q/B$.  
\end{proposition}

\begin{proof} (1) Note that each map $f\times\id{\colon}Y\times_BE_jG\to X\times_BE_jG$ is an lci morphism, and we have the cartesian diagram
\[
\xymatrix{
Y\times_BE_jG\ar[r]^{f\times\id}\ar[d]_{q_{Y,G,j}}& X\times_BE_jG\ar[d]^{q_{X,G,j}}\\
Y\times^GE_jG\ar[r]^{f_j}& X\times^GE_jG
}
\]
Since the map $q_{X,G,j}$ is faithfully flat, this implies that $f_j$ is an lci morphism, so $f_j^!$ is defined. 

By the functoriality of $(-)^!$, we have
\[
i_{G,j}^!\circ f_{j+1}^! =
f_j^!\circ i_{G, j}^!
\]
as maps  $\sE^\BM_{**}(X\times^GE_{j+1}G/B_{j+1}G,v_{j+1})\to \sE^\BM_{**}(Y\times^GE_jG/B_jG,v_j)$,
which proves (1).

For (2), we  note that the diagram
\[
\xymatrix{
Y\times_BE_jG\ar[d]^{f\times\id}\ar[r]_{\id\times i_{G,j}}& Y\times_BE_{j+1}G\ar[d]^{f\times \id}\\
X\times_BE_jG\ar[r]^{f\times\id}& X\times_BE_{j+1}G
}
\]
is cartesian and Tor-independent. Using faithful flatness again, this implies that the diagram
\[
\xymatrix{
Y\times^GE_jG\ar[d]^{f_j}\ar[r]_{i_{Y,G,j}}& Y\times^GE_{j+1}G\ar[d]^{f_{j+1}}\\
X\times^GE_jG\ar[r]^{i_{X,G,j}}& X\times^GE_{j+1}G}
\]
is also cartesian and Tor-independent. By the excess intersection formula \cite[Proposition 4.2.2]{DJK}, this says that we can identify the Gysin pullbacks $i_{Y,G,j}^!$ and $i_{X,G,j}^!$ with the refined Gysin pullback $i_{X, G,j}^!$, on $\sE^\BM_{**}(-/B, -)$. By Proposition~\ref{prop:GysinProperties},  this gives the identity
\[
i_{Y,G,j}^!\circ f_{j+1*}=f_{j*}\circ i^!_{X,G,j}.
\]
 
(3) follows immediately from the functoriality of $f\mapsto f^*$.

For (4), for $i{\colon}Z\to W$ a regular immersion in $\Sch/B$, and $\alpha, \beta$ in $\sE$-cohomology and $\gamma$ in $\sE$-Borel-Moore homology, all on $W$, and  we have the identity in $\sE$-cohomology
\[
i^*(\alpha)\cdot i^*(\beta)=i^*(\alpha\cdot \beta)
\]
and the identity in $\sE$-Borel-Moore homology
\[
i^!(\gamma)\cap i^*(\beta)=i^!(\gamma\cap\beta).
\]
This proves (4).

(5) follows by applying Proposition~\ref{prop:GysinProperties}(4) to each of the $B_jG$-schemes $X\times^GE_jG$.
\end{proof}

\begin{remark}  Having the operations of $f^!$ for lci $f$ and $f_*$ for proper $f$ on the equivariant Borel-Moore homology $\sE^\BM_{G, **}(-/-,-)$, and products in cohomology and cap products in Borel-Moore homology, all the relations  discussed in \S\ref{subsec:Relations} in the non-equivariant case pass to the equivariant case. Indeed, we have seen that the various properties of a morphism $f$, such as being proper or lci, pass to the induced morphism on each   each term in the pro-system defining the equivariant theory. The same holds for properties such as ``$f$ is smooth'', ``$f$ is \'etale'', etc., and also for a square of morphisms to commute, to be cartesian or to be cartesian and Tor-independent. Thus, each relation   discussed in \S\ref{subsec:Relations} in the non-equivariant case can be checked on the individual terms in the relevant pro-system, which then becomes a question in the the non-equivariant case.

For example, the products in equivariant cohomology make $\sE^{**}_G(X,*)$ a tri-graded ring (the third grading being the twist), and make $\sE^\BM_{G, **}(X/BG,*)$ a tri-graded $\sE^{-*,-*}_G(X,-*)$-module. We have a projection formula 
\[
f_*(\alpha\cap f^*\beta)=f_*(\alpha)\cap \beta
\]
for $f$ a proper map, $\alpha$ in Borel-Moore homology and $\beta$ in cohomology, and we have the identity  $f^!g_*=g'_*f^{\prime!}$ for each Tor-independent cartesian square
\[
\xymatrix{
Z\ar[r]^{g'}\ar[d]_{f'}&W\ar[d]^f\\
X\ar[r]^g&Y
}
\]
\end{remark}

\begin{remark}\label{rem:TrivialAction}
Suppose the $G$-action on $X$ is trivial. Then we have the canonical isomorphism 
$X\times^G E_jG\cong X\times_BB_jG$, giving us the cartesian square
\[
\xymatrix{
X\times^G E_jG\ar[r]\ar[d]^{\pi_{X,G,j}}&X\ar[d]^{p_X}\\
B_jG\ar[r]^{p_{B_jG}}&B
}
\]
with $p_{B_jG}$ smooth. The relative pullback maps $p_{B_jG}^!$ thus gives the map of pro-systems
\[
p_{B_*G}^!{\colon}\{\sE^\BM_{a,b}(X/B, v)\}_j\to \{\sE^\BM_{a,b}(X\times^G E_jG/B_jG, v_j)\}_j
\]
giving the map
\[
p_{BG, X}^!{\colon}\sE^\BM_{a,b}(X/B, v)\to \sE^\BM_{G, a,b}(X/BG, v).
\]
\end{remark}

\begin{definition}\label{def:Bounded}
1. Let $\sE$ be a commutative ring spectrum in $\SH(B)$. We call $\sE$ {\em bounded} if for each
 fixed $X\in \Sch_q^G/B$, integers $(a,b)$,  and $v\in \sK^G(X)$, the 
 pro-systems
 $\{\sE^\BM_{a,b}(X\times^GE_jG/B_jG, v_j), i_{G,j}^!\}_j$ and $\{\sE^\BM_{a,b}(X\times E_jG/E_jG, v), i_{j}^!\}_j$  are  eventually constant.\\[2pt]
2. We say that $\sE\in \SH(B)$ is {\em strongly bounded} if, given integers $a,b, d$,  there is a constant $C$ such that,
given  $Y\to B$ in $\Sm/B$, and a dense open immersion $j{\colon}U\to V$ in $\Sch/Y$ with closed complement $W$, such that   $W$ has codimension $\ge C$ in $V$ and $V$ has relative dimension $\le d$ over $Y$, the restriction map 
 \[
 j^*{\colon}\sE^\BM_{a,b}(V/Y, v)\to \sE^\BM_{a,b}(U/Y, v)
 \]
 is an isomorphism for all $v\in \sK(V)$ of virtual rank $0$.\\[2pt]
 3. We say that $\sE\in \SH(B)$ is {\em uniformly strongly bounded} if $\sE$ is strongly bounded,  the constant $C$ depends only on $d$ and $a-b$, $C=C(d, a-b)$, and $C(d, a-b)\ge C(d, a-b+1)$
 \end{definition}
 
\begin{remark} We have already remarked that our Definition~\ref{def:EquivCoh} for equivariant cohomology should be considered as provisional. Definition~\ref{def:BMWIttHomology} is also not in general the correct definition for similar reasons.  A better definition would be to consider the  pro-object of (derived) mapping spaces in the infinity categories $\mathbf{SH}(B_jG)$, 
\[
\{\Map_{\mathbf{SH}(B_jG)}(\pi_{X,G,j!}\Sigma^{a,b}\Sigma^v1_{X\times^GE_jG},\sE)\}_j, 
\]
define
\[
\Map^G_{a,b}(X/BG, v, \sE):=\lim_j \Map_{\mathbf{SH}(B_jG)}(\pi_{X,G,j!}\Sigma^{a,b}\Sigma^v1_{X\times^GE_jG},\sE)
\]
(in the infinity category of spectra)
and then set 
\[
\sE^\BM_{G, a,b}(X/BG,v):=\pi_0\Map^G_{a,b}(X/BG, v, \sE). 
\]
Equivariant $\sE$-cohomology would be similarly defined via the pro-system 
\[
\{\Map_{\mathbf{SH}(X\times^GE_jG)}(1_{X\times^GE_jG},\Sigma^{a,b}\Sigma^v p_{X\times^GE_jG}^*\sE)\}_j, 
\]
giving the mapping spectrum
\[
\Map_G^{a,b}(X, v,\sE):=\lim_j \Map_{\mathbf{SH}(X\times^GE_jG)}(1_{X\times^GE_jG},\Sigma^{a,b}\Sigma^v p_{X\times^GE_jG}^*\sE),
\]
and one would set $\sE^{a,b}_G(X, v):=\pi_0\Map_G^{a,b}(X, v,\sE)$.

In fact, this approach has been successfully carried out in the Ph.D. thesis of Alessandro D'Angelo \cite{DAngelo}, in which he extends many of the results of this paper to arbitrary $\SL$-oriented $\sE$ for which the algebraic Hopf map $\eta$ acts invertibly, what D'Angelo calls an 
$\SL[\eta^{-1}]$-oriented theory.

For $\sE$ bounded, the refined definition of equivariant Borel-Moore homology agrees with the naive one (see Lemma~\ref{lem:BoundContr} below). In particular, this is the case for $B=\Spec k$ and $\sE=\EM(\sM_*)$, $\sM_*$ a homotopy module.

Using the Poincar\'e duality isomorphism of Proposition~\ref{prop:EquivBMProperties}, we also see that the refined definition of equivariant $\sE$-cohomology agrees with the naive one for $\sE$ bounded and $X$ smooth over $B$.
\end{remark}

\begin{lemma}\label{lem:BoundContr} Suppose $\sE$ is bounded. \\[2pt]
1.  The map of pro-systems
 \[
p_{E_*G}^!{\colon}\{\sE^\BM_{a,b}(X/B, v)\}_j\to \{\sE^\BM_{a,b}(X\times E_jG/E_jG, v), i_j^!\}_j
\]
of Lemma~\ref{lem:Comparison} induces an isomorphism
\[
p_{E_*G}^*{\colon}\sE^\BM_{a,b}(X/B, v)\to \lim_j\sE^\BM_{a,b}(X\times E_jG/E_jG, v)
\]
2. The natural map   
\[
\pi_0\Map^G_{a,b}(X/BG, v, \sE)\to \sE^\BM_{G, a,b}(X/BG, v)
\]
is an isomorphism.
\end{lemma}

\begin{proof}   We may assume that $v$ has virtual rank zero; replacing $\sE$ with $\Sigma^{-a,-b}\sE$, we may assume $a=b=0$.

For (1), by  \cite[\S 4, Proposition 2.3]{MorelVoevodsky}, the  projection $p_1{\colon}X\times EG\to X$ becomes an isomorphism in  $\sH(X)$, and thus
\[
p_1^*{\colon}\sE^\BM_{0,0}(X/B, v)=\Hom_{\SH(X)}(1_X, \Sigma^{-v}p_X^!\sE)\to
\Hom_{\SH(X)}(\Sigma^\infty_{\P^1}X\times EG_+, \Sigma^{-v}p_X^!\sE)
\]
is an isomorphism.  

We have
\begin{align*}
\Hom_{\SH(X)}(\Sigma^\infty_{\P^1}X\times EG_+, \Sigma^{-v}p_X^!\sE)&\cong
\pi_0\holim_j\Map_{\mathbf{SH}(X)}(\Sigma^\infty_{\P^1}X\times E_jG_+, \Sigma^{-v}p_X^!\sE)\\
&\cong \pi_0\holim_j\Map_{\mathbf{SH}(X)}(p_{1\#}(1_{X\times E_jG}), \Sigma^{-v}p_X^!\sE)\\
&\cong \pi_0\holim_j\Map_{\mathbf{SH}(X)}(1_{X\times E_jG}, \Sigma^{-v}p_1^*p_X^!\sE)\\
&\cong \pi_0\holim_j\Map_{\mathbf{SH}(X)}(\Sigma^v1_{X\times E_jG}, p_2^!p_{E_jG}^*\sE)\\
&\cong \pi_0\holim_j\Map_{\mathbf{SH}(X)}(p_{2!}(\Sigma^v1_{X\times E_jG}), p_{E_jG}^*\sE)
\end{align*}
But $\pi_a\Map_{\mathbf{SH}(X)}(p_{2!}(\Sigma^v1_{X\times E_jG}), p_{E_jG}^*\sE)
=\sE^\BM_{a,0}(X\times E_jG/E_jG, v)$, so the assumption that $\sE$ is bounded implies that
the inverse system 
\[
\{\pi_1\Map_{\mathbf{SH}(X)}(p_{2!}(\Sigma^v1_{X\times E_jG}), p_{E_jG}^*\sE)\}_j
\]
is eventually constant, so 
\[
R^1\lim_j\pi_1\Map_{\mathbf{SH}(X)}(p_{2!}(\Sigma^v1_{X\times E_jG}), p_{E_jG}^*\sE)=0.
\]
From the Milnor sequence
\begin{multline*}
0\to R^1\lim_j\pi_1\Map_{\mathbf{SH}(X)}(p_{2!}(\Sigma^v1_{X\times E_jG}), p_{E_jG}^*\sE)\\\to  \pi_0\holim_j\Map_{\mathbf{SH}(X)}(p_{2!}(\Sigma^v1_{X\times E_jG}), p_{E_jG}^*\sE)\\\to
 \pi_0\lim_j\pi_0\Map_{\mathbf{SH}(X)}(p_{2!}(\Sigma^v1_{X\times E_jG}), p_{E_jG}^*\sE)\to0
 \end{multline*}
 we see that  the natural map
\begin{multline*}
\pi_0\holim_j\Map_{\mathbf{SH}(X)}(p_{2!}(\Sigma^v1_{X\times E_jG}), p_{E_jG}^*\sE)\\\to
\lim_j\pi_0\Map_{\mathbf{SH}(X)}(p_{2!}(\Sigma^v1_{X\times E_jG}), p_{E_jG}^*\sE)\\=
\lim_j\sE^\BM_{0,0}(X\times E_jG/E_jG, v)
\end{multline*}
is an isomorphism. Thus
\[
\sE^\BM_{0,0}(X/B, \sV)\to \lim_j\sE^\BM_{0,0}(X\times E_jG/E_jG, v)
\]
is an isomorphism.

The proof of (2) is the same: the assumption that $\sE$ is bounded implies that the natural map
\begin{multline*}
\pi_0\Map^G_{a,b}(X/BG, v, \sE)\to
\lim_j\pi_0\Map_{\mathbf{SH}(B_jG)}(\pi_{X,G,j!}(\Sigma^{a,b}\Sigma^v1_{X\times^G E_jG}), p_{B_jG}^*\sE)\\=
\lim_j\sE^\BM_{a,b}(X\times^GE_jG/B_jG, v)=\sE^\BM_{G, a,b}(X/BG, v)
\end{multline*}
is an isomorphism.
\end{proof}
 
\begin{definition}\label{def:ForgetG}  Take $X\in \Sch^G/B$, $v\in \sK^G(X)$ and let $\sE\in \SH(B)$ be a commutative ring spectrum.\\[5pt]
1. Let $*_{X,j}{\colon}X\to X\times E_jG$ be the inclusion $*_{X,j}(x)=(x,*_j)$, where $*_j$ is the base-point in $E_jG$ and let $\rho_{X,G,j}$ be the composition
\[
X\xrightarrow{*_{X,j}}X\times E_jG\to X\times^GE_jG,
\]
giving the map of pro-systems $\rho_{X,G}:=\{\rho_{X,G,j}\}{\colon}X\to \{X\times^GE_jG\}_j$. Define the ``forget $G$'' map on $\sE$-cohomology as
\[
\rho_{X,G}^*{\colon} \sE^{*,*}_G(X, v)\to \sE^{*,*}(X, v).
\]
2. Suppose $\sE$ is bounded. \\[2pt]
Define the ``forget $G$'' map
\[
\rho^\BM_{X,G}{\colon}\sE^\BM_{G, a,b}(X/BG, v)\to \sE^\BM_{a,b}(X/B, v),
\]
by the zig-zag diagram
\[
\sE^\BM_{G, a,b}(X/BG, v)
\xrightarrow{\lim\pi_{G,*}^!}  \lim_j\sE^\BM_{a,b}(X\times E_jG/E_jG, v)
\xleftarrow{\sim}\sE^\BM_{a,b}(X/B, v).
\]
\end{definition}

\begin{lemma}\label{lem:StrBoundImpBound} If $\sE$ is a strongly bounded commutative ring spectrum, then $\sE$ is bounded.
 \end{lemma}
 
\begin{proof} We check the boundedness condition on a fixed $X\in \Sch^G/B$ and for $\sE^\BM_{a,b}(-, v)$ for given $a, b$ and $v\in \sK^G(X)$. We may assume that $X$ is connected, so $v$ has a constant virtual rank $r$ on $X$; replacing $v$ with $v-\sO_X^r$ and changing $(a, b)$ to $(a+2r, b+r)$, we may assume from the start that $v$ has virtual rank zero.  Take an integer $d$ with $\dim_BX\le d$. 
 
 Recall that we have represented $G$ as a closed subgroup-scheme of $\GL_n/B$, and that $E_jG$ is an open subscheme of $\A^{n(n+j)}_B$ with closed complement $W_j$ of codimension $j+1$. Moreover,   we have the open immersion $\eta_j{\colon}E_jG\times F_n\to E_{j+1}G$, with closed complement $ W_j\times F_n\setminus W_{j+1}$ of codimension $j+1$. 
 
 $G$ acts freely on $E_jG$ and on $E_{j+1}G$, so taking the quotients, we have induced open immersion $\eta_{X,G,j}{\colon}X\times^G(E_jG\times F_n)\to X\times^GE_{j+1}G$ with closed complement $G\backslash (W_j\times F_n\setminus W_{j+1}) $ also of codimension $j+1$. On the other hand, 
$p_{X,G,j}{\colon}X\times^G(E_jG\times F_n)\to X\times^G E_jG$ is the  vector bundle $N_{X,G,j}\to  X\times^G E_jG$,  so the pullback map
\[
p_{X,G,j}^!{\colon}\sE^\BM_{a,b}(X\times^G(E_jG\times F_n)/B_{j+1}G, v)\to
\sE^\BM_{a,b}(X\times^G E_jG/B_jG, v)
\]
is an isomorphism. 

 Suppose $\sE$ is strongly bounded. Note that $X\times^G E_{j+1}G$ has dimension $\le d$ over $B_{j+1}G$, independent of $j$. Take the integer $C$ as in Definition~\ref{def:Bounded} for $a,b,  d$. Then for $j+1\ge C$, the restriction map
 \[
\eta_{X,G,j}^*{\colon}\sE^\BM_{a,b}(X\times^G E_{j+1}G/B_{j+1}G, v)\to 
 \sE^\BM_{a,b}(X\times^G(E_jG\times F)/B_{j+1}G, v) 
 \]
 is an isomorphism, hence $\sE$ is bounded.
 \end{proof}
 
 Recall our discussion of Morel's category of homotopy modules in Example~\ref{ex:WittSheaf}.    We will be mainly interested in the homotopy module $\sW_*=\{\sW\}_n$ computing Witt-sheaf cohomology, but in the interest of formulating a more general result for equivariant Borel-Moore homology, we briefly allow ourselves the added generality. 
 
\begin{lemma}\label{lem:HtpyModBound} Take $B=\Spec k$, $k$ a perfect field and let $\sM_*$ be a homotopy module. Then $\sE:=\EM(\sM_*)$ is uniformly strongly bounded.
\end{lemma} 

\begin{proof} We recall that the homotopy module $\sM_*$ is in particular a sheaf on $\Sm/k$ of graded modules for the sheaf of graded rings $\sK^{MW}_*$.   For $X\in\Sm/k$, and $\sL$ an invertible sheaf, we thus may form the twisted version $\sM_n(\sL)$ as the sheaf on $X_\Nis$ defined by
\[
\sM_n(\sL):=\sM_{n|X_\Nis}\otimes_{\sK^{MW}_{0|X}}\sK^{MW}_0(\sL)=
\sM_{n|X_\Nis}\otimes_{\sG\sW_X}\sG\sW(\sL).
\]
For $Z\subset X$ a closed subset, we have the Rost-Schmid complex for $\sM_*(\sL)$ on $X$ with supports in $Z$, which computes the (Nisnevich) cohomology   with supports, $H^p_Z(X, \sM_q(\sL))$.  We have discussed the main case of interest for us, the homotopy module $\{\sW_n=\sW\}_n$ that computes Witt sheaf cohomology in \S\ref{subsec:RS}; for the proof of the Lemma, we briefly sketch the construction for a general homotopy module.  We refer the reader to \cite{Feld, Feld1} for details; we will use the results of \cite{Feld} to view the homotopy module $\sM_*$ as a  Milnor-Witt cycle module, and then the Rost-Schmid complex is Feld's Milnor-Witt cycle complex \cite[\S 5.2]{Feld1}. We retain the notation used in \S\ref{subsec:RS}.

The Rost-Schmid complex  for $\sM_q(\sL)$ with supports in $Z$ is the complex
\begin{multline*}
C^*_{RS,Z}(X, \sM_q(\sL)):= \oplus_{x\in X^{(0)}\cap Z}\sM_q(k(x);\sL_x)\xrightarrow{\del} \ldots\\\to 
 \oplus_{x\in X^{(j)}\cap Z}\sM_{q-j}(k(x); \Lambda^X_x\otimes_{k(x)}\sL_x)\to \ldots
 \end{multline*}
 with the ``general'' term written here in degree $j$. The boundary map has a similar form to the explicit one we described for $\sM_q=\sW$ and we will not need its definition here. This computes the cohomology with support as
  \[
 H^p_Z(X, \sM_q(\sL))\cong H^p(C^*_{RS,Z}(X, \sM_q(\sL))). 
 \]
 
Moreover, the Eilenberg-MacLane spectrum $\EM(\sM_*)$ associated to the homotopy module $\sM_*$ represents the cohomology with support via a canonical isomorphism
\[
 H^p_Z(X, \sM_q(\sL))\cong \EM(\sM_*)^{p+q,q}(i_{Z*}1_Z, \sL\setminus\sO_X)).
\]
Applying the above description of cohomology with support to a rank $r$ vector bundle $V\to X$ with 0-section $0_V$,  and comparing with the line bundle $\det V\to X$ with 0-section $0_{\det V}$, this gives us the canonical isomorphisms
\begin{align*}
\EM(\sM_*)^{p+2r,q+r}(\Th(V))&\cong H^{p+r-q}_{0_V}(V, \sM_{q+r})\\
&\cong H^{p-q}(X, \sM_q(\det V))\\
&\cong H^{p-q+1}_{0_{\det V}}(\det V, \sM_{q+1})\\
&\cong \EM(\sM_*)^{p+2,q+1}(\Th(\det V)).
\end{align*}
This shows in particular that $\sE:=\EM(\sM_*)$ is $\SL$-oriented.
 
Now take $Y\in \Sm/k$, $j{\colon}U\to V$ a dense open immersion in $\Sch/Y$ with closed complement $W$, and with $\dim_YV\le d$. Take $v\in \sK(V)$ of virtual rank zero. We may assume that $Y$ is integral, let $d_Y=\dim_kY$. Choose a closed immersion $V\subset X$ in $\Sch/Y$ with $X\in \Sm/k$ of dimension $d_X$ over $k$. Let  $\sL=\det(v+\Omega_{Y/k})$. We have
\begin{multline*}
\sE^\BM_{a,b}(V/Y, v)=\sE^\BM_{a,b}(V/k, v+\Omega_{Y/k})=
\sE^\BM_{a+2d_Y, b+d_Y}(V/k, v+\Omega_{Y/k}-\sO_V^{d_Y})
\\=\sE^{2d_X-a-2d_Y, d_X-b-d_Y}_V(X, v+\Omega_{Y/k}-\sO_V^{d_Y})\\
=\sE^{2d_X-a-2d_Y, d_X-b-d_Y}_V(X, \sL-\sO_V) =H^{d_X-d_Y+b-a}_V(X, \sM_{d_X-b-d_Y}(\sL)).
\end{multline*}
Similarly, $U$ is a closed subscheme of $X\setminus W$,
\[
\sE^\BM_{a,b}(U/Y, \sV)=H^{d_X-d_Y+b-a}_U(X\setminus W, \sM_{d_X-b-d_Y}(\sL))
\]
and the restriction from $V$ to $U$ is induced by the restriction from $X$ to $X\setminus W$.

Computing $H^{d_X-d_Y+b-a}_V(X, \sM_{d_X-b-d_Y}(\sL))$ by using the Rost-Schmid complex on $X$ with supports in $V$, we see that
 $H^{d_X-d_Y+b-a}_V(X, \sM_{d_X-b-d_Y}(\sL))$ only involves the points of codimension $d_X-d_Y+b-a, d_X-d_Y+b-a+1$ and $d_X-d_Y+b-a-1$  on $X$. If $W$ has codimension $\ge c$ on $V$, then since $\dim_kV=\dim_YV+d_Y\le d+d_Y$, then $W$ has dimension $\le d+d_Y-c$ over $k$, hence has codimension $\ge c+ d_X-d_Y-d$ on $X$. We compute 
 $H^{d_X-d_Y+b-a}_U(X\setminus W, \sM_{d_X-b-d_Y}(\sL))$ similarly.

 Thus if 
\[
c\ge b-a+d+2
\]
removing $W$ from $X$ will induce an isomorphism
\[
H^{d_X-d_Y+b-a}_V(X, \sM_{d_X-b-d_Y}(\sL))\to 
H^{d_X-d_Y+b-a}_U(X\setminus W, \sM_{d_X-b-d_Y}(\sL))
\]
Taking $C(d, a-b):=d-(a-b)+2$, we see that  $\sE$ is uniformly strongly bounded.
\end{proof}

\begin{remark}
These last two results, Lemma~\ref{lem:StrBoundImpBound} and Lemma~\ref{lem:HtpyModBound}, in somewhat different form, are due to Di Lorenzo and Mantovani   \cite[Lemma 1.2.8, 1.2.9]{DiLorenzoMantovani}. The definition of equivariant Borel-Moore homology via the limit mapping spectrum is from \cite{Mantovani}, where the deeper properties of this theory are discussed. This construction also appears in \cite{DAngelo}, where the equivariant cohomology and Borel-Moore homology are defined using the motivic stable homotopy theory of the associated quotient stacks, and the identification with the theory given as the limit mapping spectrum is found in \cite[Proposition 1.4.25, Corollary 1.4.32]{DAngelo}. We will only be looking at the case $\sE=\EM(\sW_*)$, which recovers Witt-sheaf cohomology, so the naive version defined as the limit will suffice for our purposes. 
\end{remark}

\begin{remark}\label{rem:StabilityPurity} Let $\sE\in \SH(B)$ be bounded.\\[5pt]
1. Take $X\in \Sch^G_q/B$ and  let $\sL$  be a $G$-linearized invertible sheaf on $X$. 
By Lemma~\ref{lem:StrBoundImpBound} and Lemma~\ref{lem:HtpyModBound}, we have the canonical isomorphism
\[
H^\BM_{G, a}(X,\sW(\sL))\cong H^\BM_{a}(X\times^GE_jG/B_jG, \sW(\sL_j))
\]
for all $j\gg 0$.\\[2pt]
2.  For $X\in \Sch^G_q/B$,  let $\sL$ and $\sM$ be $G$-linearized invertible sheaves on $X$. 
From (1), the cap product maps 
\begin{multline*}
H^\BM_{a}(X\times^GE_jG/B_jG, \sW(\sM_j))\times H^n(X\times^GE_jG, \sW(\sL_j))  \\\xrightarrow{\cap}
H^\BM_{a-n}(X\times^GE_jG/B_jG, \sW(\sL_j\otimes \sM_j))
\end{multline*}
give rise a well-defined associative cap products
\[
\cap{\colon}H^\BM_{G, a}(X, \sW(\sM))\times H^n_G(X, \sW(\sL)) \to
H^\BM_{G, a-n}(X, \sW(\sL\otimes \sM)).
\]
Using the pullback $H^{-*}(BG, \sW)\to H^{-*}_G(X, \sW)$, this makes $H^\BM_{G, *}(X, \sW(\sM))$ into a graded $H^{-*}(BG, \sW)$-module.
\end{remark}         

\begin{proposition}\label{prop:BMLocalization}   For $X\in \Sch^G_q/B$, let $i{\colon}Z\to X$ be a closed $G$-stable subscheme with open complement $j{\colon}V\to X$, with $Z$ and $V$ in  $\Sch^G_q/B$. Suppose $\sE$ is bounded and take $v\in \sK^G(X)$. \\[5pt]
1. The family of localization sequences \eqref{eqn:BMLocSeq} for the closed immersions $Z\times^GE_jG\to  X\times^GE_jG$ with open complements $V\times^GE_jG\to X\times^GE_jG$ induces a long exact sequence of $\sE^{**}(BG)$-modules
\begin{multline*}
\ldots\to \sE^\BM_{G,a, b}(Z/BG,v)\xrightarrow{i_*}\sE^\BM_{G,a,b}(X/BG,v)\xrightarrow{j^*}
\sE^\BM_{G,a,b}(V/BG,v)\\\xrightarrow{\delta}\sE^\BM_{G, a-1,b}(Z/BG,v)\to\ldots
\end{multline*}

If we take $\sE=\EM(\sW_*)$, then for  a $G$-linearized invertible sheaf  $\sL$ on $X$, we have
 the  long exact sequence of  $H^*(BG,\sW)$-modules
\begin{multline*}
\ldots\to H^\BM_{G,a}(Z,\sW(i^*\sL))\xrightarrow{i_*}H^\BM_{G,a}(X,\sW(\sL))\xrightarrow{j^*}
H^\BM_{G,a}(V,\sW(j^*\sL))\\\xrightarrow{\delta}H^\BM_{G, a-1}(Z,\sW(i^*\sL))\to\ldots
\end{multline*}
2. Let $f{\colon}X'\to X$ be a morphism in  $\Sch^G_q/B$, giving the cartesian diagram in $\Sch^G/B$
\[
\xymatrix{
Z'\ar[r]^{i'}\ar[d]^{f_Z}&X'\ar[d]^f&U'\ar[l]_{j'}\ar[d]^{f_U}\\
Z\ar[r]^i&X&U\ar[l]_j
}
\]
Suppose that $Z'$ and $U'$ are in $\Sch^G_q/B$. If $f$ is smooth, then localization sequences in (1) are natural with respect to smooth pullback $(f_Z^!, f^!, f_U^!)$ (after twisting appropriately) and if $f$ is proper, then localization sequences in (1) are natural with respect to proper pushforward $(f_{Z*}, f_*, f_{U*})$.
\end{proposition}

\begin{proof} For (1), the assertion for $\sE=\EM(\sW_*)$ follows from the general case of bounded $\sE$ by Lemma~\ref{lem:StrBoundImpBound} and Lemma~\ref{lem:HtpyModBound}. To prove the general case, we have the long exact localization sequence in  $\sE^\BM_a(-/B_\ell G, v_\ell)$ associated to the closed immersion $Z\times^GE_\ell G\to  X\times^GE_\ell G$ with open complement $V\times^GE_\ell G\to X\times^GE_\ell G$ for each  $j$. For $?=X, Z, V$, we factor the map $i^!_{?,G,\ell}$ as $s_{?,G,\ell}^!\circ\eta_{?,G,\ell}^!=(p_{?,G,\ell}^!)^{-1}\circ\eta_{?,G,\ell}^!$. Since both $\eta_{?,G,\ell}$ and $p_{?,G,\ell}$ are smooth, Proposition~\ref{prop:GysinProperties}(5ii) says that the maps $\eta_{?,G,\ell}^!$ and the isomorphisms $p_{?,G,\ell}^!$ define maps of the relevant localization distinguished triangles, hence so do the maps $i^!_{?,G,\ell}$. Thus the system $\{i^!_{?,G,\ell}, ?=Z,X,V\}_\ell$ define a map of pro-systems of long exact localization sequences. Since $\sE$ is bounded, this gives  the desired long exact  localization  sequence in any given finite range of values of $a$ (for fixed $b$) by taking $\ell$ sufficiently large.

The assertion (2) follows from Proposition~\ref{prop:GysinProperties}(5) applied to the diagram
\[
\xymatrix{
Z'\times^GE_\ell G\ar[r]^{i'\times\id}\ar[d]^{f_Z\times\id}&X'\times^GE_\ell G\ar[d]^{f\times\id}&U'\times^GE_\ell G\ar[l]_{j'\times\id}\ar[d]^{f_U\times\id}\\
Z\times^GE_\ell G\ar[r]^{i\times\id}&X\times^GE_\ell G&U\times^GE_\ell G\ar[l]_{j\times \id}
}
\]
again for $\ell$ sufficiently large.
\end{proof}

\section{Witt-sheaf cohomology of $BN$ and $\BSL_2$}\label{sec:WittCohBG} In this section we work over a fixed perfect base field $k$ of characteristic $\neq2$, giving us the base-scheme $B:=\Spec k$.   We let $F= k^2$ be the fundamental 2-dimensional representation of $\SL_2\subset\GL_2$, defined via right matrix multiplication on row vectors, let $T_1\subset \SL_2$ be the diagonal torus, and let $N\supset T_1$ be the normalizer of $T_1$ in $\SL_2$. We recall that $N$ has the unique non-trivial character $\rho_0^-{\colon}N\to \G_m$, where $\rho_0^-$ restricts to the trivial character of $T_1$, and $\rho_0^-(\sigma)=-1$. Thus $\Pic(BN)=\Z/2$, where the non-trivial element $\gamma\in \Pic(BN)$ arises from the representation $\rho_0^-$. 

Recall from \cite[Proposition 5.5]{LevineBG} the computation of the Witt-sheaf cohomology of $\BSL_2$ and $BN$.
\begin{theorem} \label{thm:Presentation} Let $E_2\to \BSL_2$ be the rank 2 bundle associated to the representation $F$, and let $e\in H^2(\BSL_2,\sW)$ be the corresponding Euler class $e_{\SL_2}(E_2)$. \\[5pt]
1. $H^*(\BSL_2,\sW)$ is the polynomial algebra
\[
H^*(\BSL_2,\sW)=W(k)[e].
\]
2. The map $p^*{\colon}H^*(\BSL_2,\sW)\to H^*(BN,\sW)$ induced by the inclusion $N\hookrightarrow \SL_2$ is injective; we write $e\in  H^*(BN,\sW)$ for $p^*e$. There is a canonical element $x\in H^0(BN, \sW)$ with
\[
H^*(BN,\sW)=W(k)[x,e]/((1+x)\cdot e, x^2-1)
\]
3. There is a rank two bundle $\sT$ on $BN$ with determinant $\gamma$, giving the Euler class  $e(\sT)\in H^2(BN,\sW(\gamma))$,  such that  $H^{*\ge2}(BN,\sW(\gamma))$ a free $W(k)[e]$-module with generator $e(\sT)$, and with $x\cdot e(\sT)=-e(\sT)$. There is an additional element $y\in H^{0}(BN,\sW(\gamma))$ such that $H^{0}(BN,\sW(\gamma))=W(k)\cdot y$ a free $W(k)$-module with generator $y$.\end{theorem}

\begin{remark}\label{rem:Correction} There is an error in the statement of \cite[Proposition 5.5]{LevineBG} regarding the structure of 
$H^{*}(BN,\sW(\gamma))$ as a $W(k)[e]$-module, as given in Theorem~\ref{thm:Presentation}(3): the existence of the element $y$ was not mentioned. This was corrected in \cite[Proposition 2.2.14]{DAngelo} to the statement given above.  

The result \cite[Proposition 2.2.14]{DAngelo} mentioned above is for an arbitrary $\SL[\eta^{-1}]$-oriented theory; for the case of Witt sheaf cohomology, one can give an explicit construction of a generator $y\in H^{0}(BN,\sW(\gamma))$, and a description of the structure of $H^*(BN, \sW(\gamma))$ as a $H^*(\BSL_2, \sW)$-module.  We now proceed to describe $y$, and to refine Theorem~\ref{thm:Presentation} to give a complete description of the $\Z\times\Z/2$-graded $W(k)$-algebra $H^*(BN,\sW)\times H^*(BN,\sW(\gamma))$. Details will be given in \S\ref{sec:SpectSeqEulerClass} below; see Theorem~\ref{thm:CohBN} for the end result. 

We also give a similar description of $H^*_N(z, \sW)$, where $z=\Spec k(\sqrt{a})$ and $N$ acts on $z$ through the $\Z/2$-quotient $N/T_1$, acting by conjugation on $k(\sqrt{a})$. This is done in \S\ref{sec:TwistedQuadrics}.
\end{remark}

As in \cite[\S2]{LevineBG}, we represent $BN$ as $(N\backslash \SL_2)\times^{\SL_2}E\SL_2$. We have the closed immersion $\sq{\colon}\P(F)\to \P(\Sym^2F)$ induced by the squaring map $\sq{\colon}F\to \Sym^2F$; let $Y:=\P(\Sym^2F)\setminus \P(F)$. We have shown in \cite[\S5, pg. 276]{LevineBG} that there is an $\SL_2$-equivariant isomorphism $N\backslash \SL_2\cong Y$. The bundle $\sT$ of Theorem~\ref{thm:Presentation} is the one induced via these isomorphisms from the tangent bundle $T_Y$ with its canonical $\SL_2$-linearization, and the Euler class $e(\sT)$ is thus the $\SL_2$-equivariant Euler class of $T_Y$, which naturally lives in $H^2_{\SL_2}(Y, \sW(\omega_{Y/k}))$. Finally, via the isomorphism $BN\cong Y\times^{\SL_2}E\SL_2$, we have a canonical isomorphism of $\omega_{Y/k}$ with $\gamma$ (Lemma~\ref{lem:SignOmega}). Thus, we have  canonical isomorphisms
\[
\phi_i{\colon}H^i_{\SL_2}(Y, \sW)\xrightarrow{\sim} H^i(BN, \sW),\ 
\phi_i'{\colon}H^i_{\SL_2}(Y, \sW(\omega_{Y/k}))\xrightarrow{\sim} H^i(BN, \sW(\gamma)).
\]
We will give details of these identifications in \S\ref{sec:SpectSeqEulerClass}.

As a refinement of Theorem~\ref{thm:Presentation}, we give a complete description of the $\Z\times\Z/2$-graded $W(k)$-algebra $H^*(BN, \sW)\times H^*(BN, \sW(\gamma))$. We first introduce some notation.

\begin{Not}\label{not:QuadraticForms}
Let $e_0, e_1$ be the standard basis of $F=k^2$. Let $X_0, X_1$ be the dual basis, giving the basis of linear forms on $\P(F)$. We have the corresponding basis $e_0^2, e_0e_1, e_1^2$ of $\Sym^2F$ and the dual basis $T_0, T_1, T_2$ of linear forms on $\P(\Sym^2F)$. With these  choices, we have the corresponding isomorphisms $\P(F)\cong \P^1_k$, $\P(\Sym^2F)\cong \P^2_k$.  

The squaring map $\sq{\colon}F\to  \Sym^2F$ is given explicitly as
\[
sq(x_0e_0+x_1e_1)=(x_0e_0+x_1e_1)^2=x_0^2e_0^2+2x_0x_1e_0e_1+x_1^2e_1^2.
\]
Thus, with respect to our coordinates $T_0, T_1, T_2$ for $\P(\Sym^2F)$ and $X_0, X_1$ for $\P(F)$, we have
\[
sq^*(T_0)=X_0^2,\ sq^*(T_1)=2X_0X_1,\ sq^*(T_2)=X_1^*,
\]
which shows that $\sq{\colon}\P(F)\to \P(\Sym^2F)$ is a closed immersion. Letting $D:=sq(\P(F))$ be the image curve, we see that $D$ has defining equation $Q:=T_1^2-4T_0T_2$. 

The section $Q$ of $\sO_{\P(\Sym^2F)}(2)$ is nowhere vanishing on $Y$ and thus defines the rank one quadratic form on $Y$, 
\begin{equation}\label{eqn:QDef}
\<Q\>{\colon}\sO_Y(-1)\to \sO_Y,
\end{equation}
sending a local section $z$ of $\sO_Y(-1)$ to $Q\cdot z^2$. One checks that $Q$ is an $\SL_2$-invariant section of $\sO_{\P(\Sym^2F)}(2)$, hence  $\<Q\>$ is an $\SL_2$-invariant form,   and thus defines a class $[\<Q\>]\in H^0_{\SL_2}(Y, \sW)$. We then define $x:=\phi_0([\<Q\>])\in H^0(BN, \sW)$. 

To define $y$, the sheaf $\sO_{\P(\Sym^2F)}(1)$ has a canonical $\SL_2$-linearization, inducing an $\SL_2$-linearization on $\sO_Y(m)$ for all $m\in \Z$.  Let $\sV$ be the locally free $\SL_2$-linearized sheaf $\sO_Y(-2)\otimes F$, with its $\SL_2$-linearization induced from that of $\sO_Y(-2)$ and the right $\SL_2$-action on $F$. We have the $\SL_2$-invariant global generator $\Omega_2$ of $\omega_{\P(\Sym^2F)/k}(3)\cong \sO_{\P(\Sym^2F)}$ given by
\[
\Omega_2=T_0dT_1dT_2-T_1dT_0T_2+T_2dT_0dT_1.
\]

Writing a local section $v$ of $\sV$ as
\[
x=x_0\cdot e_0+x_1\cdot e_1
\]
with the $x_i$ local sections of $\sO_Y(-2)$, we define the $\omega_{Y/k}$-valued quadratic form $\tilde{Q}$ on $\sV$ by
\begin{equation}\label{eqn:TildeQ}
\tilde{Q}(x)=2(T_0x_0^2+T_1x_0x_1+T_2x_1^2)\otimes  \Omega_2;
\end{equation}
this is the form associated to the matrix
\[
\tilde{B}:=\begin{pmatrix}2T_0&T_1\\T_1&2T_2\end{pmatrix}
\]
after a suitable twist. Since $\tilde{B}$ has determinant $4T_0T_2-T_1^2$, $\tilde{Q}$ is non-degenerate. One checks that with respect to the given $\SL_2$-linearizations on $\sV$ and $\omega_{Y/k}$, $\tilde{Q}$ is an $\SL_2$-invariant form, hence gives rise to an element $[\tilde{Q}]\in H^0_{\SL_2}(Y, \sW(\omega_{Y/k}))$.  Our element $y$ is the corresponding element  $\phi_0'([\tilde{Q}])\in H^0(BN, \sW(\gamma))$. 
\end{Not}

\begin{theorem}\label{thm:Presentation2} Let  $e\in H^2(BN, \sW)$ and $\tilde{e}:=e(\sT)\in H^2(BN, \sW(\gamma))$ be as in Theorem~\ref{thm:Presentation}.  Let $x=\phi_0([\<Q\>])\in H^0(BN, \sW)$ and let  $y=\phi_0'([\tilde{Q}])\in H^0(BN, \sW(\gamma))$. Then 
\[
H^*(BN,\sW)\times H^*(BN,\sW(\gamma))=W(k)[e, \tilde{e}, x,y]/I,
\]
where $I$ is the ideal with generators $x^2-1$, $(1+x)e$, $(1+x)\tilde{e}$, $\tilde{e}^2-4e^2$, $(1+x)y$, $y^2-2(1-x)$,  $ye-\tilde{e}$, and $y\tilde{e}-4e$.
\end{theorem}

The proof of Theorem~\ref{thm:Presentation2} will be given in \S\ref{sec:SpectSeqEulerClass}.

\begin{remark} The relation $\tilde{e}^2=4e^2$ follows from the  relations $(1+x)e=0$, $y^2=2(1-x)$ and  $ye=\tilde{e}$.  
\end{remark}

\begin{corollary}\label{cor:Presentation} The map
\[
W(k)[e, e^{-1}]=H^*(\BSL_2,\sW)[e^{-1}]\to H^*(BN,\sW)[e^{-1}]
\]
is an isomorphism of $\Z$-graded $W(k)$-algebras; the element $x\in H^0(BN, \sW)$ of Theorem~\ref{thm:Presentation} is the image of $-1\in W(k)$. The $Z\times \Z/2$-graded $W(k)$-algebra  $(H^*(BN,\sW)\times H^{*}(BN,\sW(\gamma))[e^{-1}]$ is presented as
\[
(H^*(BN,\sW)\times H^{*}(BN,\sW(\gamma))[e^{-1}]=W(k)[e, e^{-1}, y]/(y^2-4)
\]
with $y\in H^0(BN,\sW(\gamma))$. The Euler class $e(\sT)$  satisfies
\[
e(\sT)=ye\in H^{*}(BN,\sW(\gamma))[e^{-1}].
\]
In particular
$H^{*}(BN,\sW(\gamma))[e^{-1}]$ is a free $H^*(BN,\sW)[e^{-1}]=W(k)[e, e^{-1}]$-module with generator $e(\sT)$ in degree 2, or with generator $y$ in degree 0.
\end{corollary}

\begin{proof} Recalling that $H^*(\BSL_2,\sW)=W(k)[e]$,  this all follows directly from the presentation of $H^*(BN,\sW)\times H^{*}(BN,\sW(\gamma))$ as $W(k)$-algebra, given in Theorem~\ref{thm:Presentation2}. 
\end{proof}

We finish this section by deriving some elementary but useful consequences of Theorem~\ref{thm:Presentation} and Corollary~\ref{cor:Presentation}.

 Since $\EM(\sW_*)$ is bounded, we have the ``forget $G$'' maps for $X\in \Sch^G/B$ with $G$-linearized invertible sheaf $\sL$ (see Definition~\ref{def:ForgetG}): 
\[
H^*_G(X, \sW(\sL))\to H^*(X, \sW(\sL)),\ 
H^\BM_{G,a}(X, \sW(\sL))\to H^\BM_{a}(X, \sW(\sL)), 
\]

Suppose $X\in \Sch^G/B$ has the trivial $G$-action. Let $\sL$ be an  invertible sheaf  on $X$ with the trivial $G$-linearization. This gives us the pull-back map of Remark~\ref{rem:TrivialAction},
\[
p_{BG,X}^!{\colon}H^\BM_a(X, \sW(\sL))\to H^\BM_{G,a}(X, \sW(\sL)).
\]
 
\begin{corollary}\label{cor:Trivial} We take $G=\SL_2$ or $G=N$. Suppose $X\in \Sch^G/k$ has the trivial $G$-action and let $\sL$ be an invertible sheaf on $X$,   with trivial $G$-linearization.  Let $\chi{\colon}G\to \G_m$ be a character, and let $\sL(\chi)$ denote the corresponding $G$-linearized invertible sheaf on $X$.  Then $p_{BG,X}^!$ and the cap product define an isomorphism
\[
\cap\circ p_{BG,X}^!{\colon}H^{-*}(BG,\sW(\chi))\otimes_{W(k)}H^\BM_*(X/B, \sW(\sL))\to 
H^\BM_{G,*}(X/B, \sW(\sL(\chi))).
\]
Moreover, the products in $H^*(-,\sW)$ define isomorphisms
\[
H^*(BG^n, \sW)\cong H^*(BG, \sW)^{\otimes_{W(k)} n}
\]
for $G=\SL_2, N$.  
\end{corollary}

\begin{proof} The second assertion follows from the first using induction, Lemma~\ref{lem:StrBoundImpBound} and Lemma~\ref{lem:HtpyModBound}, and the purity isomorphism.

For the first assertion, we note that $H^n(BG,\sW(\chi))$ is a free finite rank $W(k)$-module for each $n$; this follows from Theorem~\ref{thm:Presentation}. We have the Quillen-type spectral sequence
\[
E^1_{p,q}=\oplus_{x\in X_{(p)}}H^\BM_{p+q}(k(x)/k, \sW(\sL\otimes k(x)))\Rightarrow
H^\BM_{p+q}(X/k, \sW(\sL))
\]
and a corresponding sequence for $H^\BM_G(-/k, \sW(\sL(\chi)))$
\[
E^1_{G, p,q}=\oplus_{x\in X_{(p)}}H^\BM_{G, p+q}(k(x)/k, \sW(\sL(\chi)\otimes k(x)))\Rightarrow
H^\BM_{G,p+q}(X/k, \sW(\sL(\chi)).
\]
Since $H^*(BG, \sW(\chi))$ is a flat $W(k)$-module, the first of these gives rise to a spectral sequence
\begin{multline*}
{\mathbb E}^{1}_{p,q}=\oplus_nH^n(BG, \sW(\chi)) \otimes_{W(k)}\oplus_{x\in X_{(p)}}H^\BM_{p+q+n}(k(x)/k, \sW(\sL\otimes k(x)))\\\Rightarrow
\oplus_n H^n(BG, \sW(\chi))\otimes_{W(k)}H^\BM_{p+q}(X/k, \sW(\sL)),
\end{multline*}
 with differentials of the form $\id\otimes d_r^{p,q+n}$. The  maps $\cap\circ p_{BG,Z}^!$ for $Z$ a subscheme of $X$  induce a map of spectral sequences
\[
\cap\circ p_{BG,-}^!{\colon}{\mathbb E}^{*}_{*,*}\to E^{*}_{G, *,*}
\]

For $x\in  X_{(p)}$, we have 
\[
H^\BM_{p+q}(k(x)/k, \sW(\sL\otimes k(x)))\cong \begin{cases}0&\text{ for }q\neq0,\\
W(k(x), \sL\otimes k(x))&\text{ for }q=0,
\end{cases}
\]
and it follows from Proposition~\ref{thm:Presentation}  and the purity isomorphism that cap product induces isomorphisms
\begin{multline*}
H^n(BG, \sW(\chi))\otimes_{W(k)}H^\BM_p(k(x)/k, \sW(\sL\otimes k(x)))\\\xrightarrow{\sim}
H^n(BG_{k(x)}, \sW(\sL(\chi)))\cong H^\BM_{G, p-n}(k(x)/k, \sW(\sL(\chi))).
\end{multline*}
Thus $\cap\circ p_(BG, -)^!{\colon}{\mathbb E}^{*}_{*,*}\to E^{*}_{G, *,*}$ is an isomorphism on the $E^1$-terms and since the spectral sequences are strongly convergent, this proves the first assertion.
\end{proof}

We have the group scheme $\SL_2^n$  over $k$, with $i$th factor $\SL_2^{(i)}$. We let   $F_i$ denote the representation of $\SL_2^n$ induced by $F$ via  the $i$th projection. For $m_1,\ldots, m_n$ non-negative integers, let $\Sym^{m_1,\ldots, m_n}$ denote the $\SL_2^n$-representation $\Sym^{m_1}(F_1)\otimes\ldots\otimes\Sym^{m_n}(F_n)$; in characteristic zero, these are exactly the irreducible representations of $\SL_2^n$.

For a rank $r$ $k$-representation $\phi$ of a linear algebraic group-scheme $G$, we have the corresponding $G$-linearized vector bundle $V(\phi)$ on $\Spec k$, and the associated Euler class $e_G(V(\phi))\in H^r(BG, \sW(\det^{-1} \phi))$; we often write $e_G(\phi)$ for $e_G(V(\phi))$.

\begin{proposition} \label{prop:NonZeroEuler Class} 1. The Euler class $e_{\SL_2^n}(\Sym^{m_1,\ldots, m_n})$ is zero if and only if all the $m_i$ are even.\\[5pt]
2. If  at least one $m_i$ is odd, then the localization map
\[
H^*(\BSL_2^n,\sW)\to H^*(\BSL_2^n, \sW)[e_{\SL_2^n}(\Sym^{m_1,\ldots, m_n})^{-1}]
\]
is injective.
\end{proposition}

\begin{proof}  Recall that $\Sym^m(F)$ has rank $m+1$, so $\Sym^{m_1,\ldots, m_n}$ has odd rank if and only if all the $m_i$ are even. For a vector bundle $V\to Y$ over some smooth $k$-scheme $Y$, $e(V)=0$ in $H^*(Y, \sW(\det^{-1}V))$ if $V$ has odd rank (\cite[Proposition 3.4]{LevineEuler}), so $e_{\SL_2^n}(\Sym^{m_1,\ldots, m_n})=0$ if all the $m_i$ are even.

For the remaining assertions,  we may assume that $m_1$ is odd. We have the map $Bi{\colon}\BSL_2\to \BSL_2^n$ induced by the inclusion $i{\colon}\SL_2\to \SL_2^n$ as the first factor. Since $i^*\Sym^{m_1,\ldots, m_n}=\Sym^{m_1}(F)^r$ with $r=\prod_{j=2}^n(m_j+1)$, we have $Bi^*e_{\SL_2^n}(\Sym^{m_1,\ldots, m_n})=e_{\SL_2}(\Sym^{m_1})^r$. By \cite[Theorem 8.1]{LevineBG} for $\Char k>m_1$, extended to arbitrary characteristic in \cite{BW}, 
\[
e_{\SL_2}(\Sym^{m_1})=m_1!!\cdot e^{m_1+1}\in H^*(\BSL_2, \sW)=W(k)[e], 
\]
where
$m_1!!:=m_1(m_1-2)\cdots  3\cdot 1$. 

Since $H^*(\BSL_2^n, \sW)=W(k)[e_1,\ldots, e_n]$,  the localization map 
\[
H^*(\BSL_2^n,\sW)\to H^*(\BSL_2^n, \sW)[e(\Sym^{m_1,\ldots, m_n})^{-1}]
\]
is injective if and only if $e_{\SL_2^n}(\Sym^{m_1,\ldots, m_n})$ acts as a non-zero divisor on $W(k)$. We can detect this after passing to $W(k)[e_1,\ldots, e_n]/(e_2,\ldots, e_n)=W(k)[e]=H^*(\BSL_2,\sW)$, which reduces us to the case $n=1$. Since $e_{\SL_2}(\Sym^{m_1})=m_1!!\cdot e^{m_1+1}$, we need to show that $m_1!!$ is a non-zero divisor on $W(k)$. This follows from the fact that $m_1!!$ is an odd integer and all the $\Z$-torsion in $W(k)$ is 2-primary. In particular, $e_{\SL_2^n}(\Sym^{m_1,\ldots, m_n})\neq0$, finishing the proof of (1).
\end{proof}

\begin{corollary}\label{cor:NonZeroEulerClass} Suppose $k$ has characteristic zero. Let $\phi$ be a representation of $\SL_2^n$ such that $e_{\SL_2^n}(\phi)\in H^*(\BSL_2^n, \sW)$ is non-zero. Then the localization map
\[
H^*(\BSL_2^n, \sW)\to H^*(\BSL_2^n, \sW)[e_{\SL_2^n}(\phi)^{-1}]
\]
is injective.
\end{corollary}

\begin{proof} Since the characteristic is zero, $\BSL_2^n$ is linearly reductive and thus $\phi$ is a direct sum of the irreducible representations $\Sym^{m_1,\ldots, m_n}$, so 
\[
V(\phi)\cong \oplus_iV(\Sym^{m_1^{(i)},\ldots, m_n^{(i)}}).
\]
Thus $e_{\SL_2^n}(\phi)=\prod_ie(\Sym^{m_1^{(i)},\ldots, m_n^{(i)}})$,  so 
\[
e_{\SL_2^n}(\phi)\neq0\Rightarrow e_{\SL_2^n}(\Sym^{m_1^{(i)},\ldots, m_n^{(i)}})\neq0\ \forall i.
\]
The result then follows from Proposition~\ref{prop:NonZeroEuler Class}.
\end{proof}

\section{Localization for $\SL_2^n$-actions}\label{sec:SLLOc} In this section we prove our  localization theorem for $\SL_2^n$. We assume our base-field $k$ has characteristic $\neq2$.

\begin{lemma}\label{lem:P1Vanishing} Let $\SL_2$ act on $\P^1$ by the standard action, giving the $N$-action on $\P^1$ by restriction. Then
\[
H^*_{\SL_2}(\P^1,\sW)[e^{-1}]=0=H^*_N(\P^1,\sW)[e^{-1}].
\]
\end{lemma}

\begin{proof}  Let $G=\SL_2$ and let $\sG\subset G\times\P^1$ be the isotropy subgroup-scheme of the diagonal section in $\P^1\times \P^1$; the fiber over $p\in \P^1$ is thus the isotropy subgroup-scheme $G_p$. Let $L\subset F\times \P^1$ be the rank one sub-bundle with fiber over $p\in \P^1$ the 1-dimension subspace of $F_{k(p )}$, which is $\sG_p$-invariant. This gives us the rank 1 $G$-sub-bundle $\sL\subset  EG\times \P^1\times F$, giving us the exact sequence of bundles on $G\backslash EG \times \P^1$
\[
0\to G\backslash\sL\to G\backslash \sF\to G\backslash(\sF/\sL)\to 0
\]
where $\sF=EG \times \P^1\times F$. Since $ G\backslash\sL$  has rank one, we have $e(G\backslash\sL)=0$ in $H^2(G\backslash EG\times \P^1, \sW(\det^{-1}(G\backslash\sL))$, and thus
\[
e(G\backslash \sF)=e(G\backslash\sL)\cdot e(G\backslash (\sF/\sL))=0.
\]
Since $e(G\backslash \sF)=\pi^*e$, where $\pi{\colon}G\backslash EG \times \P^1\to BG$ is the canonical projection, we have $H^*_{\SL_2}(\P^1,\sW)[e^{-1}]=0$. 

The same proof shows  that $H^*_N(\P^1,\sW)[e^{-1}]=0$.
\end{proof}

Let $p_i{\colon}\SL_2^n\to \SL_2$ be the projection on the $i$th factor.
\begin{lemma}\label{lem:MaxParabolic} Take $G=\SL_2^n$, let $L\supset k$ be a field extension and let $P\subset G\times_kL$ be a maximal parabolic subgroup-scheme.  Then there is a Borel subgroup-scheme $B\subset \SL_2/L$ such that $P=\pi_i^{-1}(B)$ for some $i$. Moreover $B$ is conjugate to the standard Borel subgroup-scheme $B_0$ of upper triangular matrices in $\SL_2/L$ by an element $g\in \SL_2(L)$.  
\end{lemma}

\begin{proof} Replacing $k$ with $L$, we may assume $L=k$.  Let $\SL_2^{(i)}\subset \SL_2^n$ be the $i$th factor in $\SL_2^n$. Let $\pi{\colon}\SL_2^n\to \SL_2^n/P$ be the quotient map.
We claim that $\pi(\SL_2^{(i)})=\SL_2^{(i)}/(P\cap \SL_2^{(i)})$ is closed in $\SL_2^n/P$.  If not, take a closed point $x$ in the closure but not in $\pi(\SL_2^{(i)})$. We pass to $k(x)$, and changing notation, may assume the $k(x)=k$. The isotropy group of $x$ in $\SL_2^n$  is conjugate to $P$  but since $\SL_2^{(i)}$ is normal in $\SL_2^n$, the isotropy group of $x$ in $\SL_2^{(i)}$ is conjugate in $\SL_2^{(i)}$ to $\SL_2^{(i)}\cap P$, and thus the orbit $\SL_2^{(i)}\cdot x$ has the same dimension as  $\pi(\SL_2^{(i)})$. This is impossible, since each component of the $\SL_2^{(i)}$-stable closed subscheme $\overline{\pi(\SL_2^{(i)})}\smallsetminus \pi(\SL_2^{(i)})$ has strictly smaller dimension than that of $\pi(\SL_2^{(i)})$. This implies that $\pi(\SL_2^{(i)})$ is proper over $k$.

Thus if  $P\cap \SL_2^{(i)}$ is a proper subgroup-scheme of $\SL_2^{(i)}$, then $P\cap \SL_2^{(i)}$
 is a Borel subgroup-scheme $B$ of $\SL_2^{(i)}$. By \cite[Theorem 20.9]{BorelLAG}, $B$ is conjugate to the  standard Borel subgroup-scheme
 of upper triangular matrices in $\SL_2^{(i)}$ by an element of $\SL_2^{(i)}(k)$. 
 
 Since $P$ is a maximal parabolic subgroup-scheme, there is exactly one index $i$ for which 
 $P\cap \SL_2^{(i)}$ is not the whole group $\SL_2^{(i)}$, which proves the lemma.
\end{proof}

\begin{definition} Take indices $i\neq j$, $1\le i,j\le n$, and let $G_{i,j}\subset \SL_2^n$ be the closed subgroup-scheme with points $(g_1,\ldots, g_n)$, $g_i=g_j$. 

Let $K\supset k$ be a field extension. We call a closed subgroup-scheme $G\subset \SL_2^n/K$ a {\em diagonal} subgroup-scheme if there are indices $i\neq j$ such that $G$ is $K$-conjugate to $G_{i,j}\times_kK$. 
\end{definition}

Let $p_i{\colon}\BSL_2^n\to \BSL_2$ be the projection on the $i$th-factor and let $e_i=p_i^*(e)\in H^2(\BSL_2^n, \sW)$. We have the K\"unneth formula (Corollary~\ref{cor:Trivial})
\[
H^*(\BSL_2^n, \sW)\cong H^*(\BSL_2,\sW)^{\otimes_{W(k)}n}\cong W(k)[e_1,\ldots, e_n],
\]
and the isomorphism of modules over 
$H^*(\BSL_2^n, \sW)\cong  W(k)[e_1,\ldots, e_n]$,
\[
H_{\BSL_2^n}^*((\P^1)^n, \sW)\cong H_{\SL_2}^*(\P^1, \sW)^{\otimes_{W(k)}n}.
\]

\begin{proposition}\label{prop:DiagonalVanishing}  For $G=\SL_2^n$, $e_i-e_j\in H^2(BG, \sW)$ maps to zero in $H^2_G(G/G_{i,j},\sW)$.
\end{proposition}

\begin{proof} For notational simplicity, we take $i=1, j=2$. Sending $(h_1,\ldots, h_n)\in G$ to $h_2h_1^{-1}\in \SL_2$ gives an isomorphism of 
the homogeneous space $G/G_{1,2}$ with $\SL_2$, with $G$ acting on $\SL_2$ by 
\[
(g_1,\ldots, g_n)\cdot g= g_2gg_1^{-1}
\]
This shows that $G/G_{1,2}\times^GEG\cong \BSL_2^{n-1}$, with the canonical map $G/G_{1,2}\times^GEG\xrightarrow{\pi} BG=\BSL_2^n$ being the $1,2$-diagonal map $\pi=(p_1, p_1)\times p_2\times\ldots\times p_{n-1}$. Clearly $\pi^*(e_1)=\pi^*(e_2)$, whence the result.
\end{proof} 

We recall the notion of an orbit $O$ for a $G$-action on a scheme $X$, and the field of $G$-invariant functions on $O$, $k_O$, from Definition~\ref{def:Orbit}.
\begin{definition}\label{def:localizing} Take $X\in \Sch^{\SL_2^n}/k$ We call the $\SL_2^n$-action on $X$ {\em localizing} if for each orbit $O\subset X$  of positive dimension, there is a closed subgroup-scheme $P_O\subset \SL_2^n/k_O$, with $P_O$ either a maximal parabolic or a diagonal subgroup-scheme, and a morphism $O\to \SL_2^n/P_O$ in $\Sch^{\SL_2^n/k_O}/k_O$.
\end{definition}

\begin{proposition}\label{prop:BasicSplitVanishing} Take $X\in \Sch^{\SL_2^n}/k$.  Let $e_*:=\prod_{i=1}e_i\cdot \prod_{1\le i<j\le n}e_i-e_j\in H^*(BG,\sW)$ and suppose that \\[5pt]
i. the  $\SL_2^n$-action on $X$ is localizing\\
ii. $X^{\SL_2^n}=\0$.
\\[5pt]
Let $\sL$ be a $\SL_2^n$-linearized invertible sheaf on $X$. Then $H^\BM_{\SL_2^n\,*}(X, \sW(\sL))[1/e_*]=0$.
\end{proposition}

\begin{proof}   If $k$ is not perfect, and of characteristic $p>2$,  taking a base-change to the perfect closure $k^{perf}$ of $k$ induces an isomorphism 
\[
H^\BM_{\SL_2^n,*}(-/k, \sW(-))\to H^\BM_{\SL_2^n,*}(-/k^{perf}, \sW(-)), 
\]
so if the result holds over $k^{perf}$, then the result holds for $k$. Thus we may assume that $k$ is perfect.

We first show that for each orbit $O\subset X$, we have $H^\BM_{\SL_2^n/k_O\,*}(O/k_O, \sW)[1/e_*]=0$. 

By assumption, there is a morphism $f{\colon}O\to (\SL_2^n/k_O)/P$ in $\Sch^{\SL_2^n/k_O}/k_O$ for $P
\subset \SL_2^n/k_O$ a closed subgroup-scheme that is either a maximal parabolic or a diagonal subgroup-scheme. Suppose that $P$ is a maximal parabolic. By Lemma~\ref{lem:MaxParabolic},  as $k_O$ scheme with $ \SL_2^n/k_O$-action, $(\SL_2^n/k_O)/P$ is isomorphic $\P^1_{k_O}$, with the action given by the projection on the $i$th factor followed by the standard action of $\SL_2/k_O$ on $\P^1_{k_O}$. 

It follows from Lemma~\ref{lem:P1Vanishing} that, under the pullback map for $O\to \Spec k$,   $e_i\in H^2(\BSL^n_2, \sW)$ is sent to zero  in 
$H^2_{\SL_2^n/k_O}((\SL_2^n/k_O)/P, \sW)$. Thus  $e_*$ maps to zero in 
$H^*_{\SL_2^n/k_O}(O, \sW)$ and hence $H^*_{\SL_2^n/k_O}(O, \sW)[1/e_*]=0$.  

Similarly, since $e_i-e_j$ maps to zero in $H^2_{\SL_2^n}(\SL_2^n/G_{ij}, \sW)$ by Proposition~\ref{prop:DiagonalVanishing}, we see that
$H^*_{\SL_2^n/k_O}(O, \sW)[1/e_*]=0$ if there is an equivariant morphism $O\to (\SL_2^n/G_{ij})_{k_O}$. 

As $H^\BM_{\SL_2^n/k_O, *}(O/k_O, \sW)[1/e_*]$ is a module over $H^*_{\SL_2^n/k_O}(O, \sW)[1/e_*]$, this shows that $H^\BM_{\SL_2^n/k_O, *}(O/k_O, \sW)[1/e_*]=0$ for each orbit $O\subset X$.

Following Proposition~\ref{prop:Quotient}, we may assume that $X$ has a decomposition $X=\amalg_{i=1}^r X_i$ such that the quotients $\SL_2^n\backslash X_i$ all exist. Using the localization sequence  for $H^\BM_{\SL_2^n, *}(-, \sW(\sL))$ (Proposition~\ref{prop:BMLocalization}), we reduce to the case in which the quotient $\pi{\colon}X\to Y:=\SL_2^n\backslash X$ exists.

Using localization for  $H^\BM_{\SL_2^n, *}(-, \sW(\sL))$ again, and passing to the colimit in the usual way, we form the strongly convergent Quillen-type spectral sequence
\[
E^1_{p,q}=\oplus_{y\in Y_{(p)}}H^\BM_{\SL_2^n, p+q}(\pi^{-1}(y)/k, \sW(\sL))\Rightarrow 
H^\BM_{\SL_2^n, p+q}(X, \sW(\sL)).
\]
Here we define $H^\BM_{\SL_2^n, p+q}(\pi^{-1}(y)/k, \sW(\sL))$ as the colimit of 
$H^\BM_{\SL_2^n, p+q}(\pi^{-1}(U)/k, \sW(\sL))$ as $U$ runs over neighborhoods of $y$ in $Y$.
Since $k$ is perfect, we have 
\[
H^\BM_{\SL_2^n, p+q}(\pi^{-1}(y)/k, \sW(\sL))\cong
H^\BM_{\SL_2^n/k(y), q}(\pi^{-1}(y)/k(y), \sW(\sL)),
\]
rewriting our spectral sequence as
\[
E^1_{p,q}=\oplus_{y\in Y_{(p)}}H^\BM_{\SL_2^n/k(y), q}(\pi^{-1}(y)/k(y), \sW(\sL))\Rightarrow 
H^\BM_{\SL_2^n, p+q}(X, \sW(\sL)).
\]

As inverting $e_*$ is exact, this gives the 
spectral sequence
\[
E^1_{p,q}=\oplus_{y\in Y_{(p)}}H^\BM_{\SL_2^n/k(y), q}(\pi^{-1}(y)/k(y), \sW(\sL))[1/e_*]\Rightarrow 
H^\BM_{\SL_2^n, p+q}(X, \sW(\sL))[1/e_*].
\]
Since each $\pi^{-1}(y)$ is an $\SL_2^n$-orbit $O$ with quotient $\Spec k_O=\Spec k(y)$, all the $E_1$-terms are zero,  hence $H^\BM_{\SL_2^n\,,*}(X, \sW(\sL))[1/e_*]=0$, as claimed.
 \end{proof}

\begin{theorem}[Atiyah-Bott localization for $\SL_2^n$]\label{thm:MainSLLoc} Let $k$ be a field of characteristic $\neq2$.   Take $X\in \Sch^{\SL_2^n}/k$, let $\sL$ be a $\SL_2^n$-linearized invertible sheaf on $X$, and suppose that the $\SL_2^n$-action on $X$ is localizing.   Let $e_*\in W(k)[e_1,\ldots, e_n]$ be the polynomial
\[
e_*=\prod_{i=1}^ne_i\cdot\prod_{1\le j<i\le n}e_i-e_j
\]
Then  the inclusion $i{\colon}X^{\SL_2^n}\to X$ induces an isomorphism
\[
i_*{\colon}H^{\BM}_*(X^{\SL_2^n}, \sW(i^*\sL))\otimes_{W(k)}W(k)[e_1,\ldots, e_n, 1/e_*]\to H^{\BM}_{\SL_2^n, *}(X, \sW(\sL))[ 1/e_*]
\]
\end{theorem}

\begin{proof} Let $U=X\setminus X^{\SL_2^n}$.  By Proposition~\ref{prop:BasicSplitVanishing}, we have $H^{\BM}_{\SL_2^n, *}(X, \sW(\sL))[ 1/e_*]=0$.  Using the localization sequence (Proposition~\ref{prop:BMLocalization}) for 
$H^{\BM}_{\SL_2^n, *}$ with respect to the closed immersion $i{\colon}X^{\SL_2^n}\to X$, and inverting $e_*$, it follows that 
\[
i_*{\colon} H^{\BM}_{\SL_2^n, *}(X^{\SL_2^n}, \sW(i^*\sL))[ 1/e_*]
\to  H^{\BM}_{\SL_2^n, *}(X, \sW(\sL))[ 1/e_*]
\]
is an isomorphism.  

Since the $\SL_2^n$-action on $X^{\SL_2^n}$ is trivial, we have the K\"unneth isomorphism (Corollary~\ref{cor:Trivial}, note that $\SL_2^n$ has only the trivial character)
\[
H^{\BM}_{\SL_2^n, *}(X^{\SL_2^n}, \sW(i^*\sL))\cong H^\BM_*(X^{\SL_2^n}, \sW(i^*\sL))\otimes_{W(k)}H^{-*}(\BSL_2^n, \sW).
\]
Also,
\[
H^*(\BSL_2^n, \sW)=W(k)[e_1,\ldots, e_n]
\]
again by Theorem~\ref{thm:Presentation} and  Corollary~\ref{cor:Trivial}. This completes the proof. 
 \end{proof}
 
\section{Homogeneous spaces and localization for $N^n$}\label{sec:NnLoc}
We continue to assume that our base-field $k$ has characteristic $\neq2$.
 
We consider the subgroup-scheme $N^n\subset \SL_2^n$, the normalizer of the diagonal torus $T_n=T_1^n\subset \SL_2^n$.  We write $\iota(t)\in N$ for the matrix $\begin{pmatrix}t&0\\0&t^{-1}\end{pmatrix}$, $\sigma$ is the matrix $\begin{pmatrix}0&1\\-1&0\end{pmatrix}$, and for $\lambda\in k^\times$, we write $\lambda\cdot\sigma$ for $\iota(\lambda)\cdot \sigma$. We let $\iota_j{\colon}N\to N^n$ be the inclusion as the $j$th factor; we write $\iota_j(t)$ for $\iota_j(\iota(t))$. We set $\sigma_j:=\iota_j(\sigma)$. We use the standard coordinates $t_1,\ldots, t_n$ on $T_n$, that is 
\[
t_i(\iota_j(t))=\begin{cases} t&\text{ for }i=j\\1&\text{ otherwise.}\end{cases}
\]

For $I\subset \{1,\ldots, n\}$, let $\sigma_I=\prod_{i\in I}\sigma_i\in N^n$ and for $\lambda_*=(\lambda_1,\ldots, \lambda_n)\in T_n$, let $\lambda_{*I}=\prod_{i\in I}\iota_i(\lambda_i)$. For $\lambda\in k^\times$, we let $\lambda_I=\prod_{i\in I}\iota_i(\lambda)$. We write $T_I$ for the subtorus $\prod_{i\in I}\iota_i(T_1)$ of $T_n$. 

For $m_*:=(m_1,\ldots, m_n)$ an $n$-tuple of integers $m_i\ge 1$, and $\lambda_*\in T_n(k)$, let $\Lambda_I(\lambda_*, m_*)$ be the closed subscheme of $T_n$ with ideal generated by elements $t_i^{m_i}-\lambda_i$, $i\in I$, and $t_i-1$, $i\not\in I$. Let $\mu_{m_*,I}=\Lambda_I(1, m_*)=\prod_{i\in I}\iota_i(\mu_{m_i})$, so $\Lambda_I(\lambda_*, m_*)$ is a $\mu_{m_*,I}$-torsor.  

 Let $\pi{\colon}N^n\to N^n/T_n=(\Z/2)^n$ be the projection. We identify $(\Z/2)^n$ with the set of subsets $I\subset \{1,\ldots, n\}$, sending $I$ to the tuple $\sum_{i\in I}e_i$, with $e_i$ the $i$th standard basis vector in $(\Z/2)^n$. For a subgroup-scheme $H$ of $N^n$, we write $\pi(H)$ for $\pi(H(\bar{k}))$ and set $\bar\sigma_i:=\pi(\sigma_i)$.

\begin{definition} Let $F\supset k$ be a field, let $\chi{\colon}T_n\to \G_m$ be a non-trivial even character (i.e., $\chi(-\id)=+1$), and take $\lambda\in F^\times$. Define the subgroup-scheme $\Lambda(\chi, \lambda)\subset N^n_F$ by
\[
\Lambda(\chi, \lambda):=\chi^{-1}(1)_F\amalg \sigma_{1,\ldots, 1}\cdot  \chi^{-1}(\lambda)
\]
Note that  $\pm \id$ is in  $\Lambda(\chi, \lambda)$ since $\chi(-\id)=1$. 
\end{definition}

\begin{remark} $\Lambda(\chi, \lambda)$ contains $\chi^{-1}(1)$ as a normal subgroup-scheme with quotient group-scheme $\Lambda(\chi, \lambda)/\chi^{-1}(1)\cong \Z/2$, with generator the image $\bar\sigma$ of $\sigma_{1,\ldots, 1}\cdot  \chi^{-1}(\lambda)$. Since $\chi$ is non-trivial, $\chi^{-1}(1)_F$ and $\Lambda(\chi, \lambda)$ both have dimension one over $F$.
\end{remark}

\begin{construction}\label{con:HomogType}
Let $\chi{\colon}T_n\to \G_m$ be a non-trivial character and let $F\supset k$ be a field We construct four types of $N^n_F$ homogeneous spaces $X$ of dimension one over $F$; let $F_X$ denote the $N^n_F$-invariant subring of $\sO_X(X)$. \\[5pt]
Type a. $X=(N^n/\chi^{-1}(1))_F$.\\[2pt]
Type b.  $\chi$ is even and $X=N^n_F/\Lambda(\chi,\lambda)$ for some  $\lambda\in F^\times$.
We can also describe $X$ as a two-step quotient: starting with the quotient $(N^n/\chi^{-1}(1))_F$ of type (a), the quotient group $\Lambda(\chi, \lambda)/\chi^{-1}(1)=\<\bar\sigma\>\cong \Z/2$ acts and $X=(N^n/\chi^{-1}(1))_F/\<\bar\sigma\>$.

A bit more explicitly, we identify $T_N/\chi^{-1}(1)$ with $\G_m$ via the character $\chi{\colon}T_n\to \G_m$, so as a scheme, we have
\[
(N^n/\chi^{-1}(1))_F=\amalg_{I\in (\Z/2)^n}\sigma_I\cdot \G_{m\, F}.
\]
 For $I=(\epsilon_1,\ldots, \epsilon_n)\in (\Z/2)^n$, let $(-1)^I=(\iota_1((-1)^{\epsilon_1}\id),\ldots, \iota_n((-1)^{\epsilon_n}\id))\in T_n$ and let $I^c=(1-\epsilon_1,\ldots, 1-\epsilon_n)$. The (right) action $*$ of $\bar\sigma$ on $\amalg_{I\in (\Z/2)^n}\sigma_I\cdot \G_{m\, F}$ is  then given by
\[
(\sigma_I\cdot x)*\bar\sigma=\sigma_{I^c}\chi((-1)^I)\cdot(\lambda/x).
\]

Indeed, for $t\in T_n$, we have $t\cdot \sigma_{(1,\ldots, 1)} =\sigma_{(1,\ldots, 1)}\cdot (1/t)$, for $I\in (\Z/2)^n$, we have $\sigma_I\cdot \sigma_{(1,\ldots, 1)}=\sigma_{I^c} \cdot (-1)^I$, and in $\G_m:=T_n/\chi^{-1}(1)$, the subscheme $\chi^{-1}(\lambda)$ maps to $\lambda\in \G_m$. 
\\[2pt]
Type c$\pm$. Let $F'\supset F$ be a degree two extension.\\[2pt]
Type c+. Suppose $\chi$ is even, and choose $\lambda\in F^\times$.  Let $\tau$ denote the conjugation of $F'$ over $F$, which we extend to an action on $F'[x,x^{-1}]$ by $\tau(x)=x$. As above, the character $\chi$ defines an isomorphism of $T_n$-homogeneous spaces $T_n/\chi^{-1}(1)_{F'}\cong \G_{m\, F'}$, with $T_n$ acting by $\chi$ on $\G_{m\, F'}$, 

To describe the induced $N^n$-action on the quotient scheme,
\[
(N^n/\chi^{-1}(1))_{F'}=\amalg_{I\in (\Z/2)^n}\sigma_I\cdot \G_{m\, F'}.
\]
 we use the following notation. For $I=(i_1,\ldots, i_n)\in (\Z/2)^n$, and $t=(t_1,\ldots, t_n)\in T_n$, let
 \[
 t^{(-1)^I}=(t_1^{(-1)^{i_1}},\ldots, t_n^{(-1)^{i_n}})
 \]
 Recall that $(-1)^I=((-1)^{i_1},\ldots,(-1)^{i_n})\in T_n$. 
 Given a second index $J$, we have $I+J, IJ\in (\Z/2)^n$ using the component-wise addition and multiplication. Then $t\in T_n$ acts on the component $\sigma_I\cdot \G_{m\, F'}$   by 
\[
t*\sigma_I\cdot x=\sigma_I \cdot \chi(t^{(-1)^I})\cdot x
\]
and $\sigma_J$ acts by 
\[
\sigma_J*\sigma_I\cdot x=\sigma_{I+J}\cdot \chi((-1)^{IJ})\cdot x.
\]

We define an action of $\<\bar\sigma\>=\Z/2$ on $\amalg_{I\in (\Z/2)^n}\sigma_I\cdot \G_{m\, F'}$ as a twisted version of the action described in type (b). Namely, on the component $\sigma_I\cdot \G_{m\, F'}$, the map
\[
\bar\sigma^*{\colon}\sigma_{I^c}\cdot F'[x, x^{-1}]\to \sigma_I\cdot F'[x, x^{-1}]
\]
is the map 
\[
\bar\sigma^*(\sigma_{I^c}\cdot f(x))=\sigma_{I}\cdot f^\tau(\chi((-1)^I)\cdot \lambda/x).
\]
Since $\lambda$ is in $F^\times$, $(-1)^I\cdot (-1)^{I^c}=-\id$ and $\chi(-\id)=1$, we see that $(\bar\sigma^*)^2=\id$, giving an action of $\<\bar\sigma\>$ on $\amalg_{I\in (\Z/2)^n}\sigma_I\cdot \G_{m\, F'}$ as an $F$-scheme.  One can similarly check that the $N^n_F$-action on $\amalg_{I\in (\Z/2)^n}\sigma_I\cdot \G_{m\, F'}$ commutes with the $\bar\sigma$-action, making the quotient
\[
X:=(\amalg_{I\in (\Z/2)^n}\sigma_I\cdot \G_{m\, F'})/\<\bar\sigma\>
\]
an $N^n_F$-homogeneous space. 
\\[2pt]
Type c-.  Suppose $\chi$ is odd: $\chi(-\id)=-1$. Take $\lambda\in F^{\prime \times}$ such that $\tau(\lambda)=-\lambda$; if we choose an $a\in F^\times$ with $F'=F(\sqrt{a})$, then $\lambda=\lambda_0\cdot \sqrt{a}$ with $\lambda_0\in F^\times$. 

We identify $N^n/\chi^{-1}(1))_{F'}$ with $\amalg_{I\in (\Z/2)^n}\sigma_I\cdot \G_{m\, F'}$ as above, and define
\[
\bar\sigma^*{\colon}\sigma_{I^c}\cdot F'[x, x^{-1}]\to \sigma_I\cdot F'[x, x^{-1}]
\]
by the same formula as in case (c+). The fact that $\chi(-\id)=-1$ and that $\tau(\lambda)=-\lambda$ yield that $(\bar\sigma^*)^2=\id$ in this case as well: we just check on $f(x)=x$, to find
\begin{multline*}
(\bar\sigma^*)^2(\sigma_{I^c}\cdot x)=\bar\sigma^*(\sigma_I\cdot \chi((-1)^I)\lambda/x)=
\sigma_{I^c}\cdot \chi((-1)^I)\chi((-1)^{I^c})\frac{\tau(\lambda)}{\lambda}x\\
=\sigma_{I^c}\cdot \chi(-\id)\cdot (-1)\cdot x=\sigma_{I^c}\cdot x.
\end{multline*}
One checks again that the right $\bar\sigma$-action commutes with the left $N^n_F$-action, 
giving as in case (c+) the $N^n$-homogeneous space 
\[
X:=[N^n/\chi^{-1}(1))_{F'}]/\<\bar\sigma\>. 
\]
\end{construction}

\begin{lemma}\label{lem:HomogSpaceStructure} Let $X$ be a homogeneous space for $N^n$, smooth and of dimension one over $\Spec k_X:=N^n\backslash X$ and let $Y=T_n\backslash X$, with induced structure of a $N^n/T_n=(\Z/2)^n$-homogeneous space over $k_X$. Suppose that $T_1$ acts non-trivially on $X$ via each inclusion $\iota_j$. Then \\[5pt]
1. As $Y$-scheme, $X\cong \G_{m Y}$. \\[2pt]
2.  Fix a closed point $y_0$ of $Y$, let  $A\subset (\Z/2)^n$ be isotropy group of $y_0$ and let $B\subset A$ be the kernel of the action map $A\to \Aut_{k_X}(k_X(y_0))$. Then there are three cases \\[2pt]
Case a. $B=A=\{0\}$ \\[2pt]
Case b. $A=\<(1,\ldots, 1)\>=B$\\[2pt]
Case c. $A=\<(1,\ldots, 1)\>$, $B=\{0\}$.
\\[2pt]
3. We refer to Construction~\ref{con:HomogType} for the three types of homogeneous spaces; in all cases, we have $F=k_X$. In Case (a)  $X$ is of type (a), in Case (b) $X$ is of type  (b), and in both cases (a) and (b) we have  $k_X(y_0)=k_X$. In Case $(c )$,   $X$ is of type (c+) if $\chi$ is even and is of type type (c-) if $\chi$ is odd; in both cases  $k_X(y_0)$ is the degree degree two extension $F'$ of $F$. 

In particular, Construction~\ref{con:HomogType} gives a complete classification of the $N^n$-homogeneous spaces satisfying the hypotheses listed in the statement of the Lemma. 
\end{lemma}

\begin{proof} By construction, $T_n$ acts on $X_y$ for each point $y\in Y$, with $y=T_n\backslash X_y$. Thus the image of $T_n$ in $\Aut_y(X_y)$ is a quotient torus acting freely on $X_y$, hence is isomorphic to $\G_{m\, y}$ as group-scheme over $y$. This defines the character $\chi_y{\colon}T_n\to \G_m$ via which $T_n$ maps to $\G_{m\, y}\subset \Aut_y(X_y)$.  By Hilbert's theorem 90, $X_y\cong \G_{m\,y}$ as $\G_{m\, y}$-torsor, with $T_n$ acting via $\chi_y$. This proves (1).

For (2), let $\chi=\chi_{y_0}$; we refer to Construction~\ref{con:HomogType} for the notation. For $I\in (\Z/2)^n=\{0, 1\}^n$, let $\chi^{(I)}$ be the character
 \[
 \chi^{(I)}(t)=\chi(t^{(-1)^I}).
 \]

Take $I\in A$,  that is, $\sigma_I(y_0)=y_0$. Let $t$ be a point of $T_n$ and let $x$ be a point of $X_{y_0}$.  Let $*$ denote the $N^n$-action on $X$ and $\cdot$ the usual multiplication on $X_{y_0}=\G_{m\, y_0}$.  Since $t\sigma_I=\sigma_I t^{(-1)^I}$, we have
 \[
 \chi(t)\cdot (\sigma_I*x)= (t\sigma_I)*x=(\sigma_It^{(-1)^I})*x=\sigma_I*(\chi^{(I)}(t)\cdot x).
\]

 Our assumption that $\iota_j(T_1)$ acts non-trivially on $X$ for each $j$ implies that 
 $\iota_j(T_1)$ acts non-trivially on $X_{y_0}=\G_{m\, y_0}$ for each $j$. Indeed, we take the point $1\in X_{y_0}=\G_{m\, y_0}$ and suppose $\iota_j(T_1)$ acts trivially on $X_{y_0}$. Take $I\in (\Z/2)^n$ with $y=\sigma_I(y_0)$. Then $\iota_j(T_1)*(\sigma_I*1)=\sigma_I*1$, hence $\sigma_I*(\iota_j(T_1)*1)=\sigma_I*1$, whence $\iota_j(T_1)*1=1$. Since $t*\sigma_I*1=\sigma_I*1$ if and only if $t$ acts trivially on $X_y$, we see that $\iota_j(T_1)$ acts non-trivially on $X_{y_0}$, as claimed. This implies that $\chi(t_1,\ldots, t_n)=\prod_{i=1}^nt_i^{m_i}$ with no $m_i=0$.

Suppose $I=(\epsilon_1,\ldots, \epsilon_n)$ is not in $\< (1,\ldots, 1)\>$, so there is an $\epsilon_i=0$ and an $\epsilon_j=1$.  Then for $x\in \G_{m,y_0}$ and $s, t\in \G_m$, we have
 \[
t^{m_i}\cdot \sigma_I*x=\chi(\iota_i(t))\cdot \sigma_I*x=\sigma_I*(t^{m_i}\cdot x)
\]
and similarly
\[
s^{m_j}\cdot \sigma_I*x= \chi(\iota_j(s))\cdot \sigma_I*x = \sigma_I*(s^{-m_j}\cdot x).
 \]
Since both $m_i$ and $m_j$ are non-zero, this is impossible: take $t=u^{m_j}$ and $s=u^{-m_i}$ for $u\in k^\times$ suitably general and let $x=\sigma_I^{-1}*1$.\footnote{If $k$ is a finite field, we can base-extend to a transcendental extension $k(z)$, and it suffices to show that for the extended action we have $A\subset \< (1,\ldots, 1)\>$} Then we have
 \begin{multline*}
 u^{m_im_j}=u^{m_im_j}\cdot   \sigma_I*x=\chi(\iota_i(u^{m_j}))\cdot \sigma_I*x
= \sigma_I*(u^{m_im_j}\cdot x)\\
= \chi(\iota_j(u^{-m_i}))\cdot \sigma_I*x=u^{-m_im_j}\cdot \sigma_I*x= u^{-m_im_j},
\end{multline*}
which is a contradiction.  Thus $A\subset \< (1,\ldots, 1)\>$. Since $Y$ is a homogeneous space for the commutative group $N^n/T_n\cong (\Z/2)^n$, it follows that $A$ is the isotropy subgroup $(\Z/2)^n_y\subset (\Z/2)^n$ for each point $y\in Y$.

Before we proceed further, we consider a finite extension of fields $F\subset K$, the scheme $\G_{m\, K}=\Spec K[x, x^{-1}]$,  an automorphism $\tau{\colon}K\to K$ over $F$ and an automorphism $\bar\sigma{\colon}\G_{m\, K}\to \G_{m\, K}$ as $F$-scheme, such that $\bar\sigma$ is an automorphism over $\tau$; we may extend $\tau$ to an $F$ automorphism, also denoted $\tau$,  of $K[x, x^{-1}]$ by setting $\tau(x)=x$.  The automorphism $\bar\sigma$ is given by an $F$-algebra homomorphism
\[
\bar\sigma^*{\colon}K[x, x^{-1}]\to K[x, x^{-1}]
\]
such that $\tau^{-1}\circ\bar\sigma^*$ is a $K$-algebra automorphism. Thus, there is an $\epsilon\in \{\pm1\}$ and a $\lambda\in K^\times$ such that $\tau^{-1}\circ\bar\sigma^*(x)=\tau(\lambda)\cdot x^\epsilon$, that is, $\bar\sigma^*(x)=\lambda\cdot x^\epsilon$, and this uniquely determines $\bar\sigma^*$. 
 
If $A=\{0\}$, then $Y=(N^n/T_n)_{k_X}=\amalg_I\sigma_I\cdot \Spec k_X$ and thus $X$ has the $k_X$-rational point $1\in X_{y_0}=\G_{m\, y_0}$. The isotropy group of $1\in X_{y_0}$ is clearly $\chi^{-1}(1)$, giving the isomorphism  of $N^n$-homogeneous spaces.
\[
X\cong N^n/\chi^{-1}(1).
\]

If $A=B= \<(1,\ldots, 1)\>$, then $k_X(y_0)=k_X$ and
 $Y=\Spec k_X\times (\Z/2)^n/\< (1,\ldots, 1)\>$.  Choosing a set of representatives $\{I_1,\ldots, I_{2^{n-1}}\}\subset (\Z/2)^n$ for $(\Z/2)^n/\,(1,\ldots, 1)\>$, with $I_1=(0,\ldots, 0)$, we use $\sigma_{I_j}*-$ to identify the fiber $X_{y_0}$ with $X_{\sigma_{I_j}*y_0}$, where we take $y_0$ corresponding to $I_1$.
  
 $\sigma_{1,\ldots,1}$ acts as a $k_X$-automorphism of $X_{y_0}=\G_{m, k_X}$. Since $\sigma_{1,\ldots,1}(t_1,\ldots, t_n)=(t_1^{-1},\ldots, t_n^{-1})\sigma_{1,\ldots, 1}$, the $k_X$-automorphism 
\[
\sigma_{1,\ldots,1}^*{\colon}k_X[x,x^{-1}]\to k_X[x,x^{-1}]
\]
must be of the form   $\sigma^*_{1,\ldots,1}(x)=\lambda/x$ for some $\lambda\in k_X^\times$. In particular, $-\id=\sigma_{1,\ldots, 1}^2$ acts as the identity on $X_{y_0}$, so $\chi(-\id)=1$.  By transport of structure via $\sigma_{I_j}*-$, we have the same action of $\sigma_{1,\ldots, 1}$ on each fiber $X_y$, $y\in Y$.

 Writing
 \[
 N^n/\chi^{-1}(1)=\amalg_I\sigma_I\cdot T_n/\chi^{-1}(1),
 \]
 taking $y_0$ corresponding to $I_1$ and using $1\in X_{y_0}=\G_{m y_0}$ to define the action map $N^n\to X$, we see that $\Lambda(\chi, \lambda)$ acts trivially, giving the map of $N^n$-homogeneous spaces
 \[
 N^n/\Lambda(\chi, \lambda)\to X
 \]
 which is easily seen to be an isomorphism.
 
 Now suppose $A= \< (1,\ldots, 1)\>$, $B=\{0\}$. We use the set of representatives $\{I_1,\ldots, I_{2^{n-1}}\}\subset (\Z/2)^n$ for $(\Z/2)^n/\,(1,\ldots, 1)\>$ as above, making $y$ isomorphic to  $\Spec k_X(y_0)$, by writing $y=\sigma_{I_j}*y_0$, and giving the isomorphism $\sigma_{I_j}^*{\colon}k_X(y)\xrightarrow{\sim} k_X(y_0)$. Similarly, this gives the isomorphism of $k_X(y_0)$-schemes
\[
X/\chi^{-1}(1)=\G_{m\, Y}\cong \amalg_j \sigma_{I_j}\cdot \G_{m, k_X(y_0)}.
\]
Since $B=0$,  $\sigma_{1,\ldots, 1}$ acts by a non-trivial involution $\tau$ on $k_X(y_0)$ and thus $k_X(y_0)$ is a degree two Galois extension of $k_X$ (we have assumed that $\Char k\neq2$). As in the previous case, $\sigma_{1,\ldots,1}$ acts on $X_{y_0}=\G_{m, k_X(y_0)}$, but this time as an automorphism over $k_X$,  acting by $\tau$ on $k_X(y_0)$. As above, we have 
 \[
  \sigma_{1,\ldots, 1}*(t^{m_i}\cdot x)=
  \sigma_{1,\ldots, 1}*(\iota_i(t)*x) =t^{-m_i}\cdot  \sigma_{1,\ldots, 1}*x,
  \]
  so the $k_X$-automorphism $\sigma^*_{1,\ldots, 1}$ of $k_X(y_0)[x,x^{-1}]$ induced by $\sigma_{1,\ldots, 1}$  satisfies
  \[
  \sigma^*_{1,\ldots, 1}(x)=\lambda/x
  \]
  for some $\lambda\in k_X(y_0)^\times$, with  $\sigma^*_{1,\ldots, 1}$ acting
   by $\tau$ on $k_X(y_0)$. 
   
For $y=\sigma_{I_j}\cdot y_0$, we identify the fiber $X_y$ with $\Spec k_X(y_0)[x,x^{-1}]$ via the multiplication isomorphism $\sigma_{I_j}*-{\colon}X_{y_0}\to X_y$. By transport of structure, the action of $\sigma_{1,\ldots, 1}$ on $X_y$ is given, via this isomorphism, by $\sigma^*_{1,\ldots, 1}(x)=\lambda/x$, and by $\tau$ on $k_X(y_0)$. 
 
 Suppose $\chi(-\id)=+1$. Then $\sigma_{1,\ldots, 1}^2=-\id$ acts as the identity on $X_y$, so
 $ \tau(\lambda)=\lambda$ and thus $\lambda$ is in $k_X^\times$. Similarly, if $\chi(-\id)=-1$, then $\sigma_{1,\ldots, 1}^2(x)=-x$, so $\tau(\lambda)=-\lambda$, and $\lambda=\lambda_0\cdot\sqrt{a}$, where $k_X(y)=k_X(y_0)=k_X(\sqrt{a})$, and $\lambda_0, a$ are in $k_X^\times$.

 In both cases, this action defines an involution $\bar\sigma$ of $N^n/\chi^{-1}(1)_{k_X(\sqrt{a})}=\amalg_j X_{\sigma_{I_j}*y_0}$ over $k_X$, defined as the disjoint union of the involutions $\sigma_{1,\ldots, 1}$ on $X_y$ described above.   As in cases (a), (b), the action map $N^n/\chi^{-1}(1)\to X$ defined by the point $1\in \G_{m k_X(y_0)}=X_{y_0}$ descends to an isomorphism of  $(N^n/\chi^{-1}(1)_{k_X(\sqrt{a})})/\<\bar\sigma\>$ with $X$, showing we are in case (c+) if $\chi(-\id)=+1$ and case (c-) if  $\chi(-\id)=-1$.
 \end{proof}

\begin{remark}\label{rem:Dim1Quotient} If $X$ is an arbitrary homogeneous space for $N^n$ of positive dimension,  then $X$ admits a morphism $X\to \bar{X}$ of $N^n$-homogeneous spaces with $\dim \bar{X}=1$ and with $N^n\backslash X=N^n\backslash \bar{X}$. Indeed, $\iota_j(T_1)$ will act nontrivially on $X$ for some $j$ and since $\iota_j(N)$ is a normal subgroup-scheme of $N^n$, we have the $N^n$-homogeneous space $\iota_j(N)\backslash X$ of dimension $\dim X-1$ and with the same quotient by $N^n$ as for $X$.

Similarly, if $\iota_j(T_1)$ acts trivially on $X$, we can take the quotient $X\to \bar X:= \iota_j(N)\backslash X$ with $\dim X=\dim\bar X$. Thus, there is an $N^r$-homogeneous space $X'$ and  an equivariant map $X\to X'$ over a projection $N^n\to N^r$,  with $X'$ of dimension 1, and with $\iota_j(T_1)$ acting non-trivially on $X'$ for $i=1,\ldots, r$. Moreover, the induced map on the quotients $N^n\backslash X\to N^r\backslash X'$ is an isomorphism.
\end{remark}

Recall that for a $G$-scheme $X\in \Sch^G/k$, the structure map $p_X{\colon}X\to \Spec k$ defines the pullback $p_X^*{\colon}H^*(BG, \sW)\to H^*_G(X, \sW)$. 

\begin{lemma}\label{lem:NnVanishing} We have the class $e\in H^2(BN, \sW)$ described in Theorem~\ref{thm:Presentation}; for $i=1,\ldots,n$, we let $e_i\in H^2(BN^n, \sW)$ be the class $\pi_i^*(e)$, where $\pi_i{\colon}BN^n\to BN$ is the map induced by the $i$th projection $N^n\to N$. Let $X$, $\chi$ be as in Lemma~\ref{lem:HomogSpaceStructure}, write $\chi(t_1,\ldots, t_n)=\prod_{i=1}^nt_i^{m_i}$. 
\\[5pt]
In Case (a), $e_i$ goes to zero in $H^2_{N^n}(X, \sW)$ for all $i$.\\[2pt]
 In Case (b), the Euler class $e(\otimes_i\pi_i^*\tilde{\sO}(m_i))\in H^*(BN, \sW)$ goes to zero in $H^*_{N^n}(X, \sW)$.\\[2pt]
  In Case (c+)  (i.e., $\sum_im_i$ even), the class $e(\otimes_i\pi_i^*\tilde{\sO}(m_i))\in H^*(BN, \sW)$ goes to zero in $H^*_{N^n}(X, \sW)$.\\[2pt]
   In Case (c-) (i.e., $\sum_im_i$ odd), choose $j$ with $m_j$ odd. Then  $e((\otimes_{i\neq j}\pi_i^*\tilde{\sO}(m_i))\otimes \tilde{\sO}(2m_j))\in H^*(BN, \sW)$ goes to zero in $H^*_{N^n}(X, \sW)$.
\end{lemma}

\begin{proof} We use the fact that the $H^*(-, \sW(\det))$-valued  Euler characteristic is multiplicative in short exact sequences and vanishes on odd rank bundles.

We write $X$ as a right $N^n$-homogeneous space. If we are in Case (a), then 
\[
X\times^{N^n}E\GL_2\cong \ker\chi\backslash N^n_{k_X}\times^{N^n}E\GL_2\cong B\ker\chi_{k_X}
\]
The pullback of $\pi_i^*\tilde{O}(1)$ to $B\ker\chi$ is a direct sum $\sO(\bar\pi_i)\oplus \sO(-\bar\pi_i)$, where $\bar\pi_i{\colon}\ker\chi\to \G_m$ is the restriction of the projection. Thus $\pi_i^*(e_i)=0$.

In Case (b)  there is a $\Z/2$-action on $N^n/\chi^{-1}(1)_{k_X}$ with $X$ the quotient. This gives the identity
\[
X\times^{N^n}E\GL_2\cong (\Spec k_X)^{\times \Z/2}B\chi^{-1}(1).
\]
and the $\Z/2$-acting via $\bar\sigma_{1,\ldots, 1}$ (the image of $\sigma_{1,\ldots, 1}$ in $N^n/T_n$).  Since $\chi(t_1,\ldots, t_n)=\prod_it_i^{m_i}$, the tensor product $\otimes \rho_{m_i}\circ\pi_i$ contains $\rho_{\chi}\oplus\rho_{\chi^{-1}}$ a summand. This pulls back to $B\chi^{-1}(1)$ to the trivial rank two bundle, with $\bar\sigma_{1,\ldots, 1}$ acting as 
 with 
\[
\chi^{-1}(\lambda^{-1})\bar\sigma_{1,\ldots,1}=\begin{pmatrix}0&\lambda^{-1}\\(-1)^{\sum_im_i}\lambda&0\end{pmatrix}
=\begin{pmatrix}0&\lambda^{-1}\\\lambda&0\end{pmatrix}
\]
since $\sum_im_i$ is even. This  has the two invariant subspaces spanned by the $\pm1$ eigenvectors
\[
\begin{pmatrix}1\\\lambda\end{pmatrix},\ \begin{pmatrix}-1\\\lambda\end{pmatrix}
\]
Thus $\otimes_i\tilde{\sO}(m_i)$ pulls back to a bundle which contains two rank 1 summands, hence $e(\otimes_i\tilde{\sO}(m_i))$ goes to zero.

In the cases (c$\pm$), we have
\[
X\times^{N^n}E\GL_2\cong (\Spec k_X(\sqrt{a}))^{\times \Z/2}B\chi^{-1}(1).
\]
with $\Z/2$ acting on $k_X(\sqrt{a})$ by conjugation over $k_X$. In case  (c+), we have essentially the same computation as  is the same as case (b), noting that in case  $\sum_im_i$ is even, $\lambda$ is in $k_X$, so the eigenvectors we chose define subspaces that are invariant under the conjugation action of $k(\sqrt{a})$ over $k_X$. In case (c-),  $\sum_im_i$ is odd, so there is an odd $m_j$. We replace $X$ with the quotient $X':=\iota_j(\{\pm1\})\backslash X$, giving the character $\chi'(t_1,\ldots, t_m)=\prod_{i\neq j}t_i^{m_i}\cdot t_j^{2m_j}$ with $\sum_{i\neq j}m_i+2m_j$ even, and we are back in case (c+) with $e(\otimes_{i\neq j}\pi_i^*\tilde{\sO}(m_i)\otimes \tilde{\sO}(2m_j))$ going to zero in $H^*_{N^n}(X', \sW)$. Pulling back further to $H^*_{N^n}(X, \sW)$, this class still goes to zero. 
\end{proof}

We have the canonical map $\Z[e_1,\ldots, e_n]\to W(k)[e_1,\ldots, e_n]=H^*(\BSL_2^n,\sW)$. For  $Y\in \Sch^{N^n}/k$, this gives $H^*_{N^n}(Y, \sW)$ the structure of a  $\Z[e_1,\ldots, e_n]$-module.

 \begin{proposition}\label{prop:NnVanishing} Let $X$ be a homogeneous space for $N^n$ of positive dimension. Then there is a non-zero homogeneous polynomial $P\in\Z[e_1,\ldots, e_n]$ of positive degree such that  $H^*_{N^n}(X, \sW)[P^{-1}]=0$. 
 \end{proposition}
 
 \begin{proof} By Remark~\ref{rem:Dim1Quotient},  there is a projection $\phi{\colon} N^n\to N^r$, an $N^r$-homogeneous space $\bar{X}$ of dimension one such that $\iota_j(T_1)$ acts non-trivially on $\bar{X}$ for $i=1,\ldots, r$,   and   a $\phi$-equivariant map $X\to \bar{X}$.  Changing notation, we may assume that $X$ has dimension one and  $\iota_j(T_1)$ acts non-trivially on $\bar{X}$ for $i=1,\ldots, n$.  
  
  Referring to Lemma~\ref{lem:NnVanishing}, in Case (a), we can take $P=e_i$ for any $i$ between 1 and $n$. In Cases (b), (c +), we have the Euler class $e(\tilde{\sO}(m_1,\ldots, m_n))$ and in Case (c+), we have the  Euler class $e(\tilde{\sO}(m_1,\ldots, 2m_j, \ldots, m_n))$ that go to zero in $H^*_{N^n}(X, \sW)$  note that  each $m_j$ is strictly positive. 
  
  For the tensor product of two rank two bundles $V_1, V_2$, we have 
\[
e(V_1\otimes V_2)= e(V_1)^2-e(V_2)^2;\ p_1(V_1\otimes V_2)=2e(V_1)^2+2e(V_2)^2
\]
(see \cite[Proposition 9.1]{LevineBG}). Thus, if after some pullback we can write
\[
V_1\otimes V_2 \cong W_1\oplus W_2
\]
for rank 2 bundles $W_1, W_2$, 
we have 
\[
e(W_1)e(W_2)= e(W_1\oplus W_2)=e(V_1)^2-e(V_2)^2
\]
and
\[
e(W_1)^2+e(W_2)^2= p_1(W_1\oplus W_2)= 2e(V_1)^2+2e(V_2)^2.
\]

For a commutative ring $R$, the subring of  symmetric polynomials $R[s,t]^{\Z/2}\subset R[s,t]$ is the polynomial ring $R[st, (s+t)]$. Thus, for $F(s,t)\in R[s,t]^{\Z/2}$,  there is a polynomial $G(x,y)\in R[x,y]$ with 
\[
F(e(W_1)^2, e(W_2)^2)=G((e(V_1)^2-e(V_2)^2)^2, 2e(V_1)^2+2e(V_2)^2).
\]

Using this and Ananyevskiy's $\SL_2$-splitting principle \cite[Theorem 6]{Anan}, an induction in $r\ge2$ shows that there is a non-zero homogeneous polynomial $Q\in\Z[x_1,\ldots, x_r]$  of degree $2^{r-1}$, such that, for for $V_1, \ldots, V_r$ rank two bundles on some $X\in \Sm/k$, 
\[
e(V_1\otimes\ldots\otimes V_r)=Q(e(V_1)^2,\ldots, e(V_r)^2)
\]
Indeed, for $r=2$, we just use the formula for $e(V_1\otimes V_2)$ written above.  For $r>2$, we have  Ananyevskiy's special linear  flag bundle $\SF(2,4;V_1\otimes V_2)\to \BSL_2^r$ \cite[Definition 25]{Anan}; let $X\to  \SF(2,4;V_1\otimes V_2)$ be the affine space bundle of splittings of the universal flag, giving the map $\pi{\colon}X\to \BSL_2^r$. By construction, there are rank 2 bundles $W_1, W_2$ on $X$ with trivialized determinants, together with an isomorphism $\pi^*(V_1\otimes V_2)\cong W_1\oplus W_2$, compatible with the respective trivializations of the determinant. By  Ananyevskiy's splitting principle \cite[Theorem 6]{Anan},  
\[
\pi^*{\colon}H^*(\BSL_2^r, \sW)\to H^*(X, \sW)
\]
is injective. We have
\[
\pi^*(V_1\otimes \ldots\otimes V_r)\cong (W_1\oplus W_2)\otimes V_3\otimes\ldots V_r
\]
and hence
\[
\pi^*(e(V_1\otimes \ldots V_r))=e(W_1\otimes V_3\otimes\ldots V_r)\cdot e(W_2\otimes V_3\otimes\ldots V_r).
\]
Applying the induction hypothesis to  $W_i\otimes V_3\otimes\ldots V_r$, $i=1,2$, and using the identities discussed above, gives the result.

In addition,  we have $e(\tilde{\sO}(m))=\pm m\cdot e(\tilde{\sO}(1))$ if $m$ is odd, and $e(\tilde{\sO}(m))^2=m^2e(\tilde{\sO}(1))^2$ if $m$ is even \cite[Theorem 7.1]{LevineBG}. Since $e(\tilde{\sO}(1))$ is the image of the canonical generator $e\in H^2(\BSL_2,\sW)$, it follows that for $n\ge2$, $e(\tilde{\sO}(m_1,\ldots, m_n))\in H^{2^n}(BN^n, \sW)$ is the image of a non-zero homogeneous polynomial $P\in \Z[e_1,\ldots, e_n]$.

If $n=1$, then if $m$ is odd, we have $e(\tilde{\sO}(m))=\pm m\cdot e$, and we use $P=m\cdot e$. If $m$ is even then $e(\tilde{\sO}(m))^2=m^2\cdot e^2$ and we use $P=m\cdot e$.
  \end{proof}
  
  \begin{remark}\label{rem:type} Looking at the proof of Proposition~\ref{prop:NnVanishing}, we have associated to each positive  $N^n$-homogeneous space $X$ a finite list of characters of $T_n$ (depending on the quotients we took to reduce the dimension down to one), and for each character three different action-types, (a), (b), (c). 
The polynomial $P$ depends only on the character $\chi$ and the action-type.  The polynomial $Q:=e_1\cdot P$ only depends on $\chi$.

If $k$ admits a real embedding, then as the composition $\Z\to W(k)\to W(\R)=\Z$ is the identity, the polynomial $P$ is non-zero in $H^{*\ge2}(\BSL_2^n,\sW)$. 
  \end{remark}
 
 \begin{theorem}[Atiyah-Bott localization for $N^n$]\label{thm:NnLocalization}\label{thm:ABLocNr}  Take $X$  in $\Sch^{N^n}/k$, and let $i{\colon}X^{T_n}\hookrightarrow 
 X$ be the inclusion. Let $\sL$ be an invertible sheaf on $X$ with an $N^n$-linearization. Then there is a homogeneous polynomial   $P\in\Z[e_1,\ldots, e_n]$ of positive degree such that 
 \[
 i_*{\colon}H^\BM_{*, N^n}(X^{T_n}/k, \sW(i^*\sL))[1/P]\to H^\BM_{*, N^n}(X/k,\sW(\sL))[1/P] 
 \]
 is an isomorphism. If the base-field $k$ admits a real embedding, then $P$ is non-zero in $H^*(\BSL_2^n,\sW)$. 
 \end{theorem}
 
 \begin{proof} As in the proof of Proposition~\ref{prop:BasicSplitVanishing}, we may assume that  $k$ is perfect.  By Remark~\ref{rem:Reduced}, we may assume that $X$ is reduced.
 
 Using  the localization sequence in equivariant Borel-Moore homology (Proposition~\ref{prop:BMLocalization}), we can replace $X$ with $X\setminus X^{T_1}$ and can assume that $X^{T_n}=\0$, and which reduces us to showing that, if $X$ has no 0-dimensional orbits, then we can find a $P$ as in the statement, such that
 $H^\BM_{*, N^n}(X,\sW(\sL))[1/P]=0$.  
 
By Proposition~\ref{prop:Quotient}, there is a finite decomposition of $X$ as $X=\amalg_\alpha X_\alpha$ such that for each $\alpha$ the quotient scheme $X_\alpha\to N^n\backslash X_\alpha$ exists as a   quasi-projective $k$-scheme. Using the localization sequence again, we reduce to the case in which the quotient $Z:=N^n\backslash X$ exists as a  quasi-projective $k$-scheme. Localizing again, we may assume that $Z$ is integral and that $\pi{\colon}X\to Z$ is equi-dimensional. Since we are assuming that $X^{T_n}=\0$, $\pi$ must have strictly positive relative dimension. Using localization once more, we may replace $Z$ with any dense open subscheme. 
 Since $N^n$ is smooth over $k$, it follows from Proposition~\ref{prop:Quotient} that we may assume that $\pi{\colon}X\to Z$ is smooth.

We may replace $N^n$ with the quotient of $N^n$ by the factors $\iota_j(N)$ for which  $\iota_j(T_1)$ acts trivially; changing notation, we may assume that $\iota_j(T_1)$ acts non-trivially for each $j$. Again passing to a dense open subscheme of $Z$ and changing notation, we may assume that for each point $z$ of $Z$, $\iota_j(T_1)$ acts non-trivially on the fiber $X_z$ for each $j$. 

Let $Z_{(q)}$ be the set of dimension $q$ points of $Z$.   Since $k$ is perfect, the closure $\bar{z}$ of a point $z\in Z_{(b)}$ has a dense open subscheme $U(z)\subset \bar{z}\subset Z$ that is smooth of dimension $b$ over $k$; we may also take $U(z)$ such that $\Omega_{U(z)/k}\cong \sO_{U(z)}^b$. Since $X_z$ is not of finite type over $k$, we define $H^\BM_{N^n, *}(X_z/k, \sW(\sL\otimes \sO_{X_z}))$  as the colimit over models,
 \[
 H^\BM_{N^n, *}(X_z/k, \sW(\sL\otimes \sO_{X_z}))=\colim_{U(z)\subset \bar{z}}H^\BM_{N^n, *}(\pi^{-1}(U(z))/k, \sW(\sL\otimes \sO_{U(z)})),
 \]
 with the open subschemes $U(z)$ as above.  

A  filtration by closed subschemes
\[
\0=Z_{-1}\subset  Z_0\subset\ldots \subset Z_b\subset\ldots\subset Z
\]
with $Z_b$ of pure dimension $b$ over $k$ gives a corresponding $N^n$-stable filtration of $X$
\[
\0=X_{-1}\subset   X_0\subset\ldots \subset X_b\subset\ldots\subset X
\]
with $X_b=\pi^{-1}(Z_b)$. This gives a corresponding collection of localization sequences for the open-closed decompositions
\[
X_{b-1}\hookrightarrow X_b\hookleftarrow X_b\setminus X_{b-1},
\]
which yields the Quillen-type spectral sequence of homological type
\[
 E^1_{a,b}=H^\BM_{N^n, a+b}(X_b\setminus X_{b-1}/k, \sW(\sL\otimes \sO_{X_b\setminus X_{b-1}})\Rightarrow H^\BM_{N^n, a+b}(X/k,\sW(\sL)).
 \]

Passing to the colimit over the filtrations on $Z$ gives us  the Leray spectral sequence of homological type
 \[
 E^1_{a,b}=\oplus_{z\in Z_{b}}H^\BM_{N^n, a+b}(X_z/k, \sW(\sL\otimes \sO_{X_z}))\Rightarrow H^\BM_{N^n, a+b}(X/k,\sW(\sL)).
 \]
 Looking at the description of $H^\BM_*(-, \sW(-))$ via the Rost-Schmid complex in \S\ref{subsec:RS}, we find a canonical isomorphism
 \[
 H^\BM_{N^n, a}(X_z/k(z), \sW(\sL\otimes \sO_{X_z}))\cong
 H^\BM_{N^n, a+b}(X_z/k, \sW(\sL\otimes \sO_{X_z}))
 \]
 for $z\in Z_{(b)}$.  This rewrites our spectral sequence as
\[
 E^1_{a,b}=\oplus_{z\in Z_{b}}H^\BM_{N^n, a}(X_z/k(z), \sW(\sL\otimes \sO_{X_z}))\Rightarrow H^\BM_{N^n, a+b}(Y,\sW(\sL)).
 \]

Referring to  Remark~\ref{rem:type}, we have  the discrete invariant $\chi_z$ giving the character of the action of $T_n$ on $X_z$; shrinking $Z$ again if necessary, we may assume that $\chi_z=\chi_{z'}$ for all $z, z'\in Z$ and let $\chi$ denote this common character.
 
Since $X\to Z$ is smooth,  $X_z$ is smooth over $k(z)$ for each $z\in Z$. Each $X_z$ is an $N^n$-homogeneous space over $k(z)$ of positive dimension, so by Proposition~\ref{prop:NnVanishing}, there is a non-zero polynomial $P_z\in \Z[e_1,\ldots, e_n]$ such that  $H^*_{N^n}(X_z, \sW)[1/P_z]=0$.  Following Remark~\ref{rem:type}, there is  non-zero polynomial $P\in \Z[e_1,\ldots, e_n]$  that depends only on the character $\chi$ and is divisible by $P_z$ for all $z\in Z$,  thus 
\[
H^*_{N^n}(X_z, \sW)[1/P]=0\ \forall z\in Z. 
\]

 As the Borel-Moore homology  $H^\BM_{N^n, *}(X_z/k(z), \sW(\sL\otimes \sO_{X_z}))$ is an  
 $H^*_{N^n}(X_z, \sW)$ module, we have
 \begin{multline*}
 H^*_{N^n}(X_z, \sW)[1/P]=0,\ \forall z\in Z\\\Rightarrow H^\BM_{N^n, *}(X_z/k(z), \sW(\sL\otimes \sO_{X_z}))[1/P]=0,\ \forall z\in Z\\
 \Rightarrow  H^\BM_{*, N^n}(X,\sW(\sL))[1/P]=0,
 \end{multline*}
 the last implication following from the spectral sequence.
 \end{proof}

\section{Homogeneous spaces and localization for $N$}\label{sec:NLoc}

 Theorem~\ref{thm:NnLocalization}  reduces the computation of $H^\BM_{N^n, *}(Y,\sW)$, after suitable localization,  to the case of a scheme with only 0-dimensional orbits. For general $n$, the situation is still quite complicated, so we will restrict our discussion to the case of the group $N$. We will see at the end of this section that this case will be sufficient for many applications. We continue to assume that $k$ is a field of characteristic $\neq2$, and work in $\Sch^N/k$ throughout this section.

 \begin{remark}\label{rem:NLocDetails}  We apply Lemma~\ref{lem:NnVanishing} in case $n=1$.
Let $X$ be an $N$-homogeneous space of dimension one, with quotient $N\backslash X=\Spec k_X$ and associated character $\chi_m{\colon}T_1\to \G_m$, $\chi_m(t)=t^m$ as given by Lemma~\ref{lem:HomogSpaceStructure}. We may assume that $m>0$. 
  
  We consider the four types of homogeneous spaces:\\[5pt]
 Type a. In this case $X\cong (N/\mu_m)_{k_X}$ and $e:=e(\tilde{\sO}(1))$ goes to zero in $H^*_N(X, \sW)$, so we invert $P=e$ to kill $H^*_N(X, \sW)$.\\[2pt]
 Type b. In this case, $m$ is even and there is a $\lambda\in k_X^\times$ with  $X\cong N_{k_X}/\Lambda(\chi,\lambda^{-1})$. $e(\tilde{\sO}(m))$ goes to zero in $H^*_N(X, \sW(\det\tilde{\sO}(2))$, giving us the element $e(\tilde{\sO}(m))^2=m^2e^2$ in $H^*(\BSL_2,\sW)$.  If we take $P=m\cdot e$, then $H^*_N(X, \sW)[P^{-1}]=0$.\\[2pt]
 Type c$\pm$. In this case $X\cong [N_{k_X(\sqrt{a})}/\mu_m]/{\Z/2}$, for a certain $\Z/2$ action on $N_{k_X(\sqrt{a})}/\mu_m$ as $k_X$-scheme, with $\Z/2$ acting on $k_X(\sqrt{a})$ as  conjugation  over $k_X$. In type c$+$, $m$ is even and we need to invert $P=m\cdot e$ to kill $H^*_N(X, \sW)$. In type c$-$, $m$ is odd and we need to invert $2m\cdot e$ to kill $H^*_N(X, \sW)$.\\[5pt]
 In summary, we invert $e$ in case (a), $m\cdot e$ in  cases (b) and (c$+$) and $2m\cdot e$ in case (c$-$). Thus, in all cases except for (a), we are inverting 2, which unfortunately is often a zero-divisor on $W(k)$. 
 
Thus, the polynomial $P$ given by Theorem~\ref{thm:NnLocalization} in case of $n=1$ is of the form $M\cdot e$, where $M$ is the least common multiple of the coefficients $m$ discussed above, as one ranges over the (finitely many) dimension one orbit types of the $N$-action on $Y$.
 
If $k_X$ already contains a square root of -1, one can remove the factors of 2 from  $m$ if $m$ is even in the above analysis.
 \end{remark}

 We let $\bar\sigma$ denote the image of $\sigma\in N$ in $N/T_1$, so $N/T_1=\<\bar\sigma\>\cong \Z/2$.

 We recall Definition~\ref{def:StrictIntro} from the introduction. 
 
 \begin{definition}\label{def:Strict} 1. Take $X\in \Sch^N/k$.  We let $X^{T_1}_N$ denote the union of the integral components  $Z\subset X^{T_1}$ such that  $\bar\sigma\cdot Z=Z$,  and let $X^{T_1}_\Ind$ be the union of the integral components $Z\subset X^{T_1}$ such that $\bar\sigma\cdot Z\cap Z=\0$. 
 \\[5pt]
2. We call the $N$-action {\em semi-strict} if  $X^{T_1}_\red=X^{T_1}_N\cup X^{T_1}_\Ind$.
\\[5pt]
3. If the $N$-action on $X$ is semi-strict, we say the $N$-action on $X$ is {\em strict} if
$X^{T_1}_N\cap X^{T_1}_\Ind=\0$ and we can decompose $X^{T_1}_N$ as a disjoint union of two $N$-stable  closed  subschemes 
\[
X^{T_1}_N=X^N\amalg X^{T_1}_\fr 
\]
where the $N/T_1$-action on $X^{T_1}_\fr$ is free. 
 \end{definition}

\begin{remark} 1. If  $X^{T_1}_\red$ is a union of closed points of $X$, then the $N$-action is strict.\\[2pt]
2.  If  $X^{T_1}_\red$ is a disjoint union of its integral components, then the $N$-action is semi-strict, and the $N$-action is strict if and only if each integral  component of $X^N$ is an integral  component of $X^{T_1}$. In case $X^{T_1}_\red$ is a disjoint union of its integral components, the action thus fails to be strict exactly when there is an integral component $C$ of $X^{T_1}$ that is $\bar\sigma$ stable, but with $C^{\bar\sigma}$ a proper subscheme of $C$.

In case $X$ is smooth over $k$, then $X^{T_1}$ is also smooth over $k$ and  is thus a disjoint union of its integral components.
\end{remark}

\begin{lemma}\label{lem:Rank1LocalIso} Take $Y\in\Sch^N/k$ and let $\sL$ be an $N$-linearized invertible sheaf on $Y$. Suppose that the $T_1$ action on $Y$ is trivial and that the quotient map $q{\colon}Y\to \bar{Y}:=N\backslash Y$ is an \'etale degree two cover. Then each $\bar{y}\in \bar{Y}$ admits an open neighborhood $V$ such that, letting $U=q^{-1}(V)$ with inclusion $j{\colon}U\to Y$, $j^*\sL$ is isomorphic to $\sO_U$ as $N$-linearized sheaf, where we give $\sO_U$ the canonical $N$-linearization induced by the $N$-action on $U$. 
\end{lemma}

\begin{proof}  This follows from Lemma~\ref{lem:NLocalTriv} with $\sV=\sL$. Indeed, since $\sL$ has rank one, we have an $N$-stable open neighborhood $j{\colon}U\to Y$ of $q^{-1}(\bar{y})$ and an isomorphism of $N$-linearized sheaves $j^*\sL\cong\sO_U\otimes^{\tau,\sigma}_kV(\rho_0)=\sO_U$. 
\end{proof}

We recall from Theorem~\ref{thm:Presentation} that $H^*(\BSL_2,\sW)=W(k)[e]$, with $e$ the Euler class of the tautological rank 2 bundle on $\BSL_2$,  and that
\[
H^*(BN,\sW)=W(k)[x,e]/((1+x)\cdot e, x^2-1)
\]
for a certain $x\in H^0(BN, \sW)$. We also recall that the group scheme $N$ has two characters, the trivial character $\rho_0$ and the character $\rho_0^-$, with $\rho_0^-(\sigma)=-1$, and $\rho_0^-$ restricted to $T_1$ the trivial character. 

For a character $\chi$ of some group-scheme $G$ over $k$, and $F$ an extension field of $k$, we let $F(\chi)$ denote the representation of $G$ on $F$ with character $\chi$. The character $\chi$ determines an invertible sheaf $\sL_\chi$ on $BG$, and we write $\sW(\chi)$ for the sheaf $\sW(\sL_\chi)$ on $BG$.

\begin{remark} In what follows we will be using some computations from \S\ref{sec:TwistedQuadrics}, but these are obtained independently of the results in this section, so there is no circularity.
\end{remark}

The following result gives information on $H^\BM_{N, *}(X^{T_1}, \sW(\sL))$.

\begin{lemma}\label{lem:Coh0DimlNOrbit} Take $X\in \Sch^N/k$ and let $\sL$ be an $N$-linearized invertible sheaf on $X$. \\[5pt]
1. Let $i_0{\colon}X^N\to X$ be inclusion. Then for each connected component $X^N_i$ of $X^N$ there is an invertible sheaf $\sL_i$ on $X^N_i$ with trivial $N$-action and  a character $\chi_i$ of $N$ with $i_0^*\sL_{|X^N_i}\cong \sL_i\otimes_kk(\chi_i)$ as $N$-linearized invertible sheaf. Moreover
\[
H^\BM_{N, *}(X^N, \sW(i_0^*\sL))\cong \oplus_i H^\BM_*(X^N_i, \sW(\sL_i))\otimes_{W(k)} H^{-*}(BN, \sW(\chi_i)).
\]
2. Let  $i_\Ind{\colon}X^{T_1}_\Ind\to X$ be the inclusion.  Then
\[
H^\BM_{N,*}(X^{T_1}_\Ind, \sW(i_\Ind^*\sL))[e^{-1}]=0.
\]
\\[2pt]
3. Let $i_Y{\colon}Y\to X^{T_1}_N$ be a locally closed $N$-stable subscheme. Suppose that  the $N/T_1$-action on $Y$ is free; let 
\[
q_\fr{\colon}Y\to \bar{Y}:=N\backslash Y
\]
be the quotient map. Take $\bar{z}\in \bar{Y}$ and let $z=q_\fr^{-1}(\bar{z})$, with inclusion $i_z{\colon}z\to X$, giving the degree two \'etale extension $k(\bar{z})\subset \sO_z$.   Then  the isomorphism of Lemma~\ref{lem:Rank1LocalIso} induces the isomorphism of $H^*(BN,\sW)$ modules
\[
H^*_N(z,\sW(i_z^*\sL))\cong H^*_N(z,\sW).
\]
Moreover, writing $\sO_z=k(\bar{z})[X]/T^2-a$ for suitable $a\in k(\bar{z})^\times$, then as algebra over $H^*(BN_{k(\bar{z})}, \sW)=W(k(\bar{z}))[x, e]/(x^2-1, (x+1)e)$, we have
\[
H^*_N(z,\sW)=H^*(BN_{k(\bar{z})}, \sW)[y]/(x-\<a\>, I_a\cdot y, y^2-2(\<1\>-\<a\>),I_a\cdot e).
\]
where $I_a\subset W(k(\bar{z}))$ is the ideal $\im(\Tr_{\sO_z/k(\bar{z})})\subset W(k(\bar{z}))$.
Finally, letting $\bar{W}(k(\bar{z})):=W(k(\bar{z}))/\im(\Tr_{\sO_z/k(\bar{z})})$, we have
\[
H^*_N(z,\sW)[e^{-1}]=\bar{W}(k(\bar{z}))[e, e^{-1}].
\]
\end{lemma}

\begin{proof} For (1), we may assume that $X^N$ is connected.  As $i_0^*\sL$ is an invertible sheaf on $X^N$ with an action of $N$ over the trivial action on $X^N$, $N$ acts by a character $\chi$ on $i_0^*\sL$, and we have the isomorphism of $N$-linearized sheaves $i_0^*\sL=\sL_0\otimes_kk(\chi)$, with $\sL_0$ having the trivial $N$-action. The isomorphism
\[
H^\BM_{N, *}(X^N, \sW(i_0^*\sL))\cong 
H^\BM_*(X^N, \sW(\sL_i))\otimes_{W(k)} H^{-*}(BN, \sW(\chi))
\]
follows from Corollary~\ref{cor:Trivial}.

Let $C_1,\ldots, C_{2r}$ be the irreducible components of $X^{T_1}_\Ind$ with $\bar\sigma(C_{2i-1})=C_{2i}$. We prove (2) by induction on $r$. If  $r=1$, then  
$X^{T_1}_\Ind=C_1\amalg C_2$ with $C_2=\bar\sigma\cdot C_1$. Thus as scheme with $N$-action, $X^{T_1}\cong (N/T_1)\times_kC_1$, with the trivial action on $C_1$, so
\[
X^{T_1}_\Ind\times^NEN\cong C_1\times^{T_1}EN\cong C_1\times_k BT_1.
\]
Thus there is a character $\chi$ of $T_1$,  an invertible sheaf $\sL_0$ on  $C_1$ and an isomorphism 
\[
H^\BM_{N,*}(X^{T_1}_\Ind, \sW(i_\Ind^*\sL))\cong
H^\BM_{ *}(C_1, \sW(\sL_0))\otimes_{W(k)}H^{-*}(BT_1, \sW(\chi)).
\]
As $H^{-*}(BT_1, \sW(\chi))$ is a $H^{-*}(BT_1, \sW)$-module and $H^{-*}(BT_1, \sW)[e^{-1}]=0$, this proves (2) in case $r=1$. 

In general, write ${X^{T_1}_\Ind}=C\cup  C'$ with  $C=C_1\amalg C_2$ and $C'=C_3\cup\ldots\cup C_{2r}$. By induction, $H^\BM_{N,*}(C, \sW(\sL))[1/e]=H^\BM_{N,*}(C', \sW(\sL))[1/e]=0$. Moreover, 
$C\cap C'=(C_1\cap C')\amalg \bar\sigma\cdot (C_1\cap C')$, so the argument for the case $r=1$ shows that $H^\BM_{N,*}(C\cap C', \sW(\sL))[1/e]=0$. The localization sequence for 
$H^\BM_{N,*}(-, \sW(\sL))$ then shows that $H^\BM_{N,*}(X^{T_1}_\Ind, \sW(i_\Ind^*\sL))[1/e]=0$.

For (3),  Lemma~\ref{lem:Rank1LocalIso} gives us the isomorphism $H^*_N(z,\sW(i_z^*\sL))\cong H^*_N(z,\sW)$.  Suppose first that $z$ is an integral $k$-scheme, that is, $a$ is not a square in $k(\bar{z})^\times$. Then we use Theorem~\ref{thm:Nontrivial0DimlCase} to compute $H^*_N(z, \sW)$ and $H^*_N(z, \sW)[e^{-1}]$. If $a$ is a square in $k(\bar{z})^\times$, then we use Remark~\ref{rem:Nontrivial0DimlCase}. 

\end{proof}

\begin{theorem}[Atiyah-Bott localization for an $N$-action]\label{thm:AtiyahBottLocalizationN}  Take $X\in \Sch^N/k$  and let $i{\colon}X^{T_1}_N \hookrightarrow X$ be the inclusion. Let $\sL$ be an invertible sheaf on $X$ with an $N$-linearization. Suppose the $N$-action on $X$ is semi-strict. Then there is an $M\in \Z\setminus\{0\}$ such that
 \[
 i_*{\colon}H^\BM_{N,*}(X^{T_1}_N, \sW(i^*\sL))[1/Me]\to H^\BM_{N,*}(X,\sW(\sL))[1/Me] 
 \]
 is an isomorphism.   
\end{theorem}

\begin{proof} Since the action is semi-strict, we have the closed immersion $X^{T_1}_\Ind \cup X^{T_1}_N\hookrightarrow X$ with open complement $X\setminus X^{T_1}$. By Lemma~\ref{lem:Coh0DimlNOrbit}, 
$H^\BM_{N,*}(X^{T_1}_\Ind , \sW(i^*\sL))[1/e]=0$. If we apply Lemma~\ref{lem:Coh0DimlNOrbit} to 
the scheme $X^{T_1}_N\in \Sch^N/k$, we see that $H^\BM_{N,*}(X^{T_1}_\Ind\cap X^{T_1}_N , \sW(i^*\sL))[1/e]=0$, so by localization, 
\[
H^\BM_{N,*}(X^{T_1}\setminus X^{T_1}_N , \sW(i^*\sL))[1/e]= H^\BM_{N,*}(X^{T_1}_\Ind\setminus X^{T_1}_N , \sW(i^*\sL))[1/e]=0
\]
 as well, and thus the inclusion $X^{T_1}_N\to X^{T_1}$ induces an isomorphism 
\[
H^\BM_{N,*}(X^{T_1}_N, \sW(\sL))[1/e]\to H^\BM_{N,*}(X^{T_1}, \sW(\sL))[1/e]. 
\]
The result then follows from Theorem~\ref{thm:NnLocalization}.
\end{proof}

 \begin{remark} For $n>1$, it is not in general the case that for $N^n\cdot y\subset Y$ a 0-dimensional orbit consisting of more than one point, that $H^*_{N^n}(N^n\cdot y, \sW)=0$ after inverting some non-zero $P$, so one cannot in general localize $H^\BM_{*, N^n}(Y, \sW)$ to the fixed points.  We have not made a computation of the $N^n$-cohomology of a 0-dimension orbit in all cases. 
 \end{remark}

We conclude this section by showing how to reduce an $N^n$-action to an action of $N$. 
 
 \begin{proposition} Given odd integers $a_1,\ldots, a_n$, there is a unique homomorphism of group schemes $\rho_{a_1,\ldots, a_n}{\colon}N\to N^n$ with 
 \[
 \rho_{a_1,\ldots, a_n}(\iota(t))=\prod_{i=1}^n\iota_i(t^{a_i}),\ \rho(\sigma)=\prod_{i=1}^n\iota_i(\sigma)
 \]
 \end{proposition}
 
 \begin{proof} Uniqueness is clear; for existence, we note that $\rho_{a_1,\ldots,a_n}$ restricted to $T_1$ is a homomorphism $T_1\to T_n$. We have
 \begin{align*}
  \rho_{a_1,\ldots, a_n}(\iota(t))\circ\sigma)&=\prod_{i=1}^n\iota_i(\iota(t)^{a_i}\cdot \sigma) \\
  &=\prod_{i=1}^n\iota_i(\sigma\cdot\iota(t)^{-a_i})\\
  &=\rho_{a_1,\ldots, a_n}(\sigma\cdot \iota(t)^{-1}) 
  \end{align*}
 and
 \begin{align*}
\rho_{a_1,\ldots, a_n}(\sigma^2)&=\rho_{a_1,\ldots, a_n}(\prod_{i=1}^n\iota_i(-\id_N))\\
&=
\prod_{i=1}^n\iota_i((-1)^{a_i}\id_N)\\
&=-\id_{N^n}=(\prod_{i=1}^n\iota_i(\sigma))^2\\
&= \rho_{a_1,\ldots, a_n}(\sigma)^2.
 \end{align*}
The relations defining $N$ are thus respected by $\rho_{a_1,\ldots, a_n}$, so $\rho_{a_1,\ldots, a_n}$ yields a well-defined homomorphism.
 \end{proof}
 
 \begin{remark} Suppose we have a $k$-scheme $Y$ with an $N^n$ action. By a general choice of positive odd integers $a_1,\ldots, a_n$, the $T_1$-action on $Y$ via the co-character $t\mapsto (t^{a_1},\ldots, t^{a_n})$ will have the same fixed point set as for $T_n$; this follows from the fact that for $t_\gen\in T_n$ a geometric generic point over $k$, $Y^{t_\gen}=Y^{T_n}$, and so the 
co-characters $t\mapsto (t^{a_1},\ldots, t^{a_n})$ for which $Y^{T_1}\neq Y^{T_n}$ all have image in a fixed proper closed subscheme of $T_n$.

Thus, if $Y$ has an $N^n$-action with 0-dimensional $T_n$-fixed point locus, for a general choice of positive odd integers $a_1,\ldots, a_n$, the $N$-action on $Y$ induced by $\rho_{a_1,\ldots, a_n}$ will also have 0-dimensional $T_1$-fixed point locus, and the $N$-action with thus be strict.
\end{remark}

\section{A Bott residue theorem for Borel-Moore Witt homology}\label{sec:BottRes}
We consider a codimension $r$ regular embedding   $i{\colon}Y\to X$ in $\Sch^G/k$. Let $\sN_i$ be the conormal sheaf of $i$, giving the normal bundle $N_i\to Y$, $N_i:=\V(\sN_i)$.  $N_i$ carries a canonical $G$-linearization, so we have the Euler class $e_G(N_i)\in H^r_G(Y, \sW(\det^{-1}N_i))$. 

\begin{lemma}\label{lem:PushPull} Let $\sL$ be an invertible sheaf on $X$. The map 
\[
i^!i_*{\colon}H^\BM_{G, *}(Y, \sW(i^*\sL))\to H^\BM_{G, *-r}(Y, \sW(i^*\sL\otimes\det^{-1}N_i)) 
\]
is cap product with $e_G(N_i)\in H^*_G(Y, \sW(\det^{-1}N_i))$.
\end{lemma}

\begin{proof} We have the corresponding regular embedding  $i_j{\colon}Y\times^G E_jG\to X\times^G E_jG$. Then the sheaf $p_1^*\sN_i$ on $Y\times E_jG$ is the conormal sheaf of $i\times\id{\colon}
Y\times E_jG\to X\times E_jG$ and the conormal  sheaf $\sN_{i_j}$  of $i_j$    is the sheaf  constructed  from $p_1^*\sN_i$ by descent via the $G$-action. Thus the result in the $G$-equivariant setting follows from the case $G=\{\id\}$.

In this case, the result is a consequence of the self-intersection formula of D\'eglise-Jin-Khan \cite[Corollary 4.2.3]{DJK}, itself a special case of the excess intersection formula \cite[Proposition 4.2.2]{DJK}; see also the statement in Proposition~\ref{prop:GysinProperties}.
\end{proof}

\begin{remark}[Invertible elements in Witt-cohomology] We briefly return to the general setting. Let $G\subset \GL_n$ be a group-scheme over $k$, take $X\in \Sch^G/k$ and let $\sL$ be a $G$-linearized invertible sheaf on $X$. Products in equivariant Witt cohomology gives $H^{2*}_G(X, \sW):=\oplus_mH^{2m}_G(X, \sW)$ the structure of a commutative, graded ring, and makes $H^{2*}_G(X, \sW(\sL))$ a graded $H^{2*}_G(X, \sW)$-module. Using the canonical isomorphism $H^*_G(X, \sW(\sL^{\otimes 2}))\cong H^*_G(X, \sW)$ (see Example~\ref{ex:WittSheaf}), we also have a commutative graded product 
\[
H^{2*}_G(X, \sW(\sL))\times H^{2*}_G(X, \sW(\sL))\to H^{2*}_G(X, \sW)
\]

For $s\in H^{2m}_G(X, \sW)$, we say a homogeneous element $x\in H^{2*}_G(X, \sW(\sL))[s^{-1}]$   is invertible if there is a homogeneous element $y\in H^{2*}_G(X, \sW(\sL)[s^{-1}]$ with $xy=1\in H^{2*}(X, \sW)[s^{-1}]$. Clearly, a homogeneous element $x\in H^{2*}_G(X, \sW(\sL)[s^{-1}]$ is invertible if and only if $x^2\in H^{2*}_G(X, \sW)[s^{-1}]$ is invertible in the commutative ring $H^{2*}_G(X, \sW)[s^{-1}]$, in the usual sense.
\end{remark}

\begin{lemma}\label{lem:Nilpotence} We take $G=\SL_2^n$ or $G=N$. Let $V$ be  a $G$-linearized vector bundle of rank $2r$ on some connected $Y\in \Sch^G/k$; we suppose that the assumptions of Construction~\ref{const:GenEuler}, Case 1, Case 2 or Case 3, hold, giving us the generic Euler class $[e^{gen}_G(V)]\subset H^{4r}(BG, \sW)$. Choose an element $e^{gen}_G(V)\in [e^{gen}_G(V)]$.   Then $e_G(V)\in H^{2r}_G(Y, \sW(\det^{-1}V))$ is invertible in $H^{2*}_G(Y, \sW(\det^{-1}V))[e^{gen}_G(V)^{-1}]$. 
\end{lemma}

\begin{proof} 
Since $e_G(V)$ is invertible in 
$H^{2*}_G(Y, \sW(\det^{-1}\sV))[e^{gen}_G(V)^{-1}]$ if and only if $e_G(V\oplus V)=e_G(V)^2$ is invertible in $H^{2*}_G(Y, \sW)[e^{gen}_G(V\oplus V)]$, we may assume that $\det V\cong \sL^{\otimes 2}$ for some $G$-linearized invertible sheaf $\sL$, and $\det V^{gen}\cong k(\chi)^{\otimes 2}$ for some character $\chi$ of $G$. It thus suffices to show that $e_G(V)$ is invertible in $H^{2*}_G(Y, \sW)[e_G(V^{gen})^{-1}]$, for some choice of representative $G$-representation $V^{gen}$ for the isomorphism  class $[\sV^{gen}]$, in case $\det V$ and $\det V^{gen}$ are trivial. 

It follows from Construction~\ref{const:GenEuler} that, for each $G$-stable trivializing open subscheme $j_{U_1}{\colon} U_1\hookrightarrow Y$, we have $j^*_{U_1}e_G(V)-e_G(V^{gen})=0$ in $H^*_G(U_1, \sW)$,  so $j^*_{U_1}e_G(V)/e_G(V^{gen})=1$ in 
$H^*_G(U_1, \sW)[e_G(V^{gen})^{-1}]$. Let $i{\colon}Z_1\to Y$ be the closed complement of $U_1$. The long exact sequence for cohomology with support gives us a (non-unique) element $\alpha_1\in H^*_{G, Z_1}(Y, \sW)[e_G(V^{gen})^{-1}]$ with $i_*(\alpha_1)=e_G(V)/e_G(V^{gen})-1$. Given a second   $G$-stable trivializing open subscheme $j_{U_2}{\colon}U_2\to Y$, we have
\[
\alpha_1\cdot (j^*_{U_2}[e_G(V)/e_G(V^{gen})-1])=0\in H^*_{G, Z\cap U_2}(U_2, \sW)[e_G(V^{gen})^{-1}]
\]
so letting $Z_2=Y\setminus U_2$,  there is an element $\alpha_2 \in H^*_{G, Z_1\cap Z_2}(Y, \sW)[e_G(V^{gen})^{-1}]$  mapping to $\alpha_1\cdot [e_G(V)/e_G(V^{gen})-1]$ in $H^*_{G, Z_1}(Y, \sW)[e_G(V^{gen})^{-1}]$, so $\alpha_2$ maps to $[e_G(V)/e_G(V^{gen})-1]^2$ in $H^*_G(Y, \sW)[e_G(V^{gen})^{-1}]$. Taking $U_1,\ldots, U_r$ a Zariski open cover of $Y$ by $G$-stable trivializing open subschemes (which exists by Remark~\ref{rem:G-Triv}), we see that 
$[e_G(V)/e_G(V^{gen})-1]^r=0$ in $H^*_G(Y, \sW)[e^{gen}_G(V)^{-1}]$, so $e_G(V)$ is invertible in 
$H^*_G(Y, \sW)[e_G(V^{gen})^{-1}]$.
\end{proof}

\begin{remark}\label{rem:FixedPointSpecialCases}  If $\sV$ has odd rank, then $[e^{gen}_G(V)]=0$, so there is nothing to prove.
\end{remark}

For $X\in \Sch^N/k$,   we have defined the closed subscheme $X^{T_1}_N$ of $X$ (see Definition~\ref{def:Strict}). In the statement of the next result,  for $X\in \Sch^{\SL_2^n}/k$, we define $|X|^{\SL_2^n}:=X^{\SL_2^n}$ and $|X|^N:=X^{T_1}_N$.

\begin{theorem}\label{thm:BottResidue} We take $G=\SL_2^n$ or $G=N$. Take $X\in \Sch^G/k$ and let $\sL$ be a $G$-linearized invertible sheaf on $X$. We assume that  for each connected component $i_j{\colon}|X|^G_j\to X$ of $|X|^G$, $i_j$ is a regular embedding. Let $N_{i_j}$ be the normal bundle of $i_j$.

If $G=N$, we assume that for each $j$,  the hypotheses of   Case 1, Case 2, or Case 3 of
Construction~\ref{const:GenEuler} for $V=N_{i_j}$ holds. We  let $e_G^{gen}(N_i)$ be the product of the Euler classes $e^{gen}_G(N_{i_j})$.  We also assume that the $N$-action on $X$ is semi-strict.

 If $G=\SL_2^n$, we suppose   that the $G$-action on $X\setminus |X|^G$ is  localizing (Definition~\ref{def:localizing}). We also assume that either $k$ has characteristic zero, or $\dim_k |X|^G=0$. 
 
In case $G=\SL_2^n$, let  
\[
P=\prod_{i=1}^ne_i\dot\prod_{1\le i<j\le n}e_i-e_j\in H^{*\ge2}(\BSL_2^n, \sW)
\]
and in case $G=N$, take $P=p\cdot M\cdot e\in H^{*\ge2}(\BSL_2, \sW)$, with $M>0$ the integer defined in  Remark~\ref{rem:NLocDetails} and $p$ the exponential characteristic.  

Decomposing $H^\BM_{G*}(|X|^G,\sW(i^*\sL))$ according to the connected components $|X|^G_j$ gives the isomorphism 
\[
H^\BM_{G*}(|X|^G,\sW(i^*\sL))[P^{-1}e^{gen}(N_i)^{-1}]\cong
\prod_jH^\BM_{G*}(|X|_j^G,\sW(i_j^*\sL))[P^{-1}e^{gen}(N_{i_j})^{-1}]
\]
Then the inverse of the isomorphism (Theorem~\ref{thm:MainSLLoc} for $G=\SL_2^n$, Theorem~\ref{thm:AtiyahBottLocalizationN}  for $G=N$)
\[
i_*{\colon}H^\BM_{G*}(|X|^G,\sW(i^*\sL))[P^{-1}e^{gen}(N_i)^{-1}]\xrightarrow{\sim} H^\BM_{G*}(X,\sW(\sL)[P^{-1}e^{gen}_G(N_i)^{-1}] 
\]
is the map 
\[
x\mapsto \prod_ji_j^!(x)\cap e_G(N_{i_j})^{-1}
\]
\end{theorem}

\begin{remark}\label{rem:BottResidueSmoothCase}
If $X$ is in  $\Sm^N/k$ , then  $X^{T_1}$ is smooth, the $N$-action is semi-strict and each   $i_j{\colon}(X^{T_1}_N)_j\to X$ of $X^{T_1}_N$ is a regular embedding. Moreover, $N_{i_j}=N_{i_j}^\mov$ for each $j$, so we are in Case 2 of 
Construction~\ref{const:GenEuler}, and we may apply Theorem~\ref{thm:BottResidue}.
\end{remark}

\begin{proof}[Proof of   Theorem~\ref{thm:BottResidue}]  By Lemma~\ref{lem:Nilpotence}, $e_G(N_{i_j})$ is invertible in the localization 
\[
H^{*}_G(|X|_j^G,\sW(\det^{-1}N_{i_j}))[P^{-1}e_G^{gen}(N_{i_j})^{-1}]
\]
 for each $j$. Moreover the map $i_j^!$ is of the form
\begin{multline*}
H^\BM_{G,*}(X,\sW(\sL))[P^{-1}e^{gen}_G(N_i)^{-1}] \\
\xrightarrow{i_j^!}
H^\BM_{G,*}(|X|_j^G,\sW(i_j^*\sL\otimes\det N^{-1}_{i_j}))[P^{-1}e_G^{gen}(N_{i_j})^{-1}]
\end{multline*}
and we have the cap product
\begin{multline*}
H^\BM_{G,*}(|X|_j^G,\sW(i_j^*\sL\otimes\det N^{-1}_{i_j}))[P^{-1}e_G^{gen}(N_{i_j})^{-1}]\\
\times
H^{-*}_G(|X|_j^G,\sW(\det^{-1}N_{ij}))[P^{-1}e^{gen}_G(N_{i_j})^{-1}]\\
\xrightarrow{\cap}
H^\BM_{G,*}(|X|_j^G,\sW(i_j^*\sL))[P^{-1}e^{gen}_G(N_{i_j})^{-1}]
\end{multline*}
 so the formula in the statement of the Theorem makes sense.  

Given $x\in H^\BM_{G,*}(X,\sW(\sL)[P^{-1}e^{gen}(N_i)^{-1}]$, it follows from Theorem~\ref{thm:MainSLLoc}  (for $G=\SL_2^n$) or Theorem~\ref{thm:AtiyahBottLocalizationN} (for $G=N$), that there are 
\[
y_j\in H^\BM_{G,*}(|X|_j^G,\sW(i_j^*\sL))[P^{-1}e_G^{gen}(N_{i_j})^{-1}]
\]
with
\[
x=\sum_ji_{j*}(y_j)
\]
But then by  Lemma~\ref{lem:PushPull}, 
\[
i_j^!(x)=i_j^!i_{j*}(y_j)= y_j\cap e_G(N_{i_j}).
\]
so $y_j=i_j^!(x)\cap e_G(N_{i_j})^{-1}$.
\end{proof} 

\begin{remark} Take $G, X$ as in Theorem~\ref{thm:BottResidue}; if $G=\SL_2^n$ we suppose the $G$-action is localizing. Suppose that $|X|^G$ is a finite set of closed points of $X$,  and that $\pi{\colon}X\to \Spec k$ is smooth and proper over $k$. Then all the hypotheses of Theorem~\ref{thm:BottResidue} are satisfied. Suppose in addition that we have a $G$-linearized vector bundle $V$ of rank $d:=\dim_kX$ on $X$ and a $G$-equivariant isomorphism $\rho{\colon} \det^{-1} V\xrightarrow{\sim} \omega_{X/k}\otimes \sM^{\otimes 2}$ for some $G$-linearized invertible sheaf $\sM$. Then we have the $G$-equivariant Euler class
\[
e_G(V)\in H^d_G(X, \sW(\det^{-1}V))\xymatrix{\ar[r]^\rho_\sim&} H^d_G(X, \sW(\omega_X))
=H^\BM_{G,0}(X,\sW),
\]
and via $\pi_*$, we have the equivariant degree
\[
\pi_*^\rho(e_G(V))\in H_{G,0}^\BM(\Spec k,\sW)=H^0(BG, \sW),
\]
We similarly have the usual Euler class $e(V)\in H^d(X, \sW(\det^{-1}V))=H^\BM_{0}(X,\sW)$ and the degree,
\[
\pi_*^\rho(e(V))\in H^0(\Spec k, \sW)=W(k).
\]

Using our version of the Bott residue theorem, Theorem~\ref{thm:BottResidue}, we can compute
$\pi_*^\rho(e_G(V))$ in a suitable localization of  $H^0(BG, \sW)$ by taking the appropriate trace map applied to $e_G(i^*V)/e(N_i))$, $i{\colon}|X|^G\to X$ the inclusion. In the case $G=N$, for the points $x\in X^{T_1}_\fr\subset X^{T_1}_N$, we will need to use the computation of the push-forward given in Corollary~\ref{cor:Nontrivial0DimlCasePushforward}.
\end{remark}

\begin{ex}[Some trivial examples] 1. Let $\SL_2^n$ act on $\P^{2n-1}$ by 
\begin{multline*}
(g_1,\ldots, g_n)([x_0,x_1,\ldots, x_{2n-2}, x_{2n-1}])\\=
[(g_1\cdot \begin{pmatrix}x_0\\ x_1\end{pmatrix})^t, (g_2\cdot\begin{pmatrix}x_2\\ x_3\end{pmatrix})^t,\ldots, (g_n\cdot \begin{pmatrix}x_{2n-2}\\ x_{2n-1}\end{pmatrix})^t]
\end{multline*}
or on $\P^{2n}$ by the same formula with trivial action on the last coordinate $x_{2n}$. For 
$X:=\P^{2n-1}$, we have $X^G=\0$ and for $X=\P^{2n}$ with have $X^G=\{[0,\ldots, 0,1]\}$, with normal bundle the fundamental representation $F_1\oplus\ldots\oplus  F_n$. We give the tangent bundles $T_{\P^{2n-1}}$ and $T_{\P^{2n}}$ the induced $G$-linearization, use the canonical orientation $\det^{-1}T_X=\omega_X$, and apply Theorem~\ref{thm:BottResidue} to compute
\[
\pi_*(e_G(T_{\P^{2n-1}}))=0\in H^0(BG, \sW)[e^{-1}]
\]
\[
\pi_*(e_G(T_{\P^{2n}}))=1\in H^0(BG, \sW)[e^{-1}]
\]
Since the localization map $H^0(BG, \sW)\to H^0(BG, \sW)[e^{-1}]$ is injective, this implies
\[
\pi_*(e_G(T_{\P^{2n-1}}))=0\in H^0(BG, \sW)
\]
\[
\pi_*(e_G(T_{\P^{2n}}))=1\in H^0(BG, \sW)
\]
and restricting to the fiber over the base-point of $BG$ recovers the known results
\[
\pi_*(e(T_{\P^{2n-1}}))=0\in W(k)
\]
\[
\pi_*(e(T_{\P^{2n}}))=1\in W(k).
\]
2. We let $G$ act on $\P^{2n-1}$ or $\P^{2n}$ as in (1) and consider the induced action on a Grassmann variety $\Gr(m, 2n)$ or $\Gr(m, 2n+1)$. Let   $e_{j-1}$ denotes the $j$th standard basis vector of $k^{2n}$ or $k^{2n+1}$ and let $F_j$ denote the span of $e_{2j-2}, e_{2j-1}$, considered as a representation of $G$ via the $j$th factor $\SL_2^{(j)}$.  It is clear that the $G$-invariant linear subspaces of $\P^{2n-1}$ are exactly those  of the form $\P(F_{i_1}\oplus\ldots F_{i_r})\cong \P^{2r-1}$; for $\P^{2n}$, we have all these together with those of the form 
$\P(F_{i_1}\oplus\ldots F_{i_r}\oplus k\cdot e_{2n}$, where $e_j$ denotes the $j$th standard basis vector of $k^{2n}$. 

Thus $\Gr(m, 2n)^G=\0$ for $m$  odd. Letting $x_{i_1,\ldots, i_r}=F_{i_1}\oplus\ldots \oplus F_{i_r}\in \Gr(2r, 2n)^G$, the normal bundle $N_{x_{i_1,\ldots, i_r}}$ is
\[
N_{x_{i_1,\ldots, i_r}}=(F_{i_1}\oplus\ldots\oplus F_{i_r})^\vee\otimes (F_{j_1}\oplus\ldots\oplus F_{j_{n-r}})
\]
where $\{j_1,\ldots, j_{n-r}\}$ is the complement of $\{i_1,\ldots, i_r\}$ in $\{1,\ldots, n\}$. For $i\neq j$,  
$F_i^\vee \otimes F_j\cong \Sym^{\epsilon_1,\ldots, \epsilon_n}$ with $\epsilon_\ell=1$ for $\ell=i,j$, $\epsilon_\ell=0$ otherwise. Thus, the Euler class $e(N_{x_{i_1,\ldots, i_r}})$ is non-zero, giving
\[
\pi_*(e(T_{\Gr(2r, 2n)}))=\binom{n}{r}\cdot 1\in W(k)[e_1,\ldots, e_n]
\]
\[
\pi_*(e(T_{\Gr(2r+1, 2n)}))=0\in W(k)[e_1,\ldots, e_n]
\]

For $Gr(2r, 2n+1)$, we have the same description of the $G$-invariant even-dimensional subspaces. The normal bundle is
\[
N_{x_{i_1,\ldots, i_r}}=(F_{i_1}\oplus\ldots\oplus F_{i_r})^\vee\otimes (F_{j_1}\oplus\ldots\oplus F_{j_{n-r}}\oplus k\cdot e_{2n})
\]
giving terms of the form $\Sym^{\epsilon_1,\ldots, \epsilon_n}$ with $0\le \epsilon_j\le 1$  and with $1\le \sum_\ell\epsilon_\ell 2$. Thus  $e(N_{x_{i_1,\ldots, i_r}})\neq0$ and
\[
\pi_*(e(T_{\Gr(2r, 2n+1)}))=\binom{n}{r}\cdot 1\in W(k)[e_1,\ldots, e_n]
\]
In the case of odd dimensional subspaces, the $G$-invariant ones are all of the form 
$F_{i_1}\oplus\ldots\oplus F_{i_r})\oplus k\cdot e_{2n}$, giving a normal bundle of the form 
\[
N_{x_{i_1,\ldots, i_r}}=(F_{i_1}\oplus\ldots\oplus F_{i_r}\oplus k\cdot e_{2n})^\vee\otimes (F_{j_1}\oplus\ldots\oplus F_{j_{n-r}})
\]
Thus $e(N_{x_{i_1,\ldots, i_r}})\neq0$ and 
\[
\pi_*(e(T_{\Gr(2r+1, 2n+1)}))=\binom{n}{r}\cdot 1\in W(k)[e_1,\ldots, e_n]
\]
\end{ex}

\section{Equivariant cohomology of twisted quadrics}\label{sec:TwistedQuadrics}

We turn to the study of the case of a 0-dimensional orbit for $N$. Necessarily the torus $T_1$ acts trivially, giving an action of $N/T_1=\<\bar\sigma\>\cong \Z/2$.  If we fix the invariant function field to be $k$, there are three possibilities. One is  $\Spec k$ with trivial action, another is $\Spec k\amalg \Spec k$ with $\bar\sigma$ exchanging the two points. Both of these have  equivariant cohomology that is easy to compute: the first just gives us $H^*(BN, \sW)$ or $H^*(BN, \sW(\gamma))$ (see Theorem~\ref{thm:Presentation}), and the second is $H^*(BT_1, \sW)=W(k)$, supported in degree 0, or $H^*(BT_1, \sW(\sO(1)))=0$.  

In this section, we consider remaining case of an $N$-action on a dimension zero,  integral $k$-scheme $z$ with residue field $k(z)$ a degree two extension of $k$, and with $\bar\sigma$ acting as Galois conjugation. Since the characteristic is assumed to be different from 2, we can write $k(z)=k(\sqrt{a})$ for some $a\in k^\times\setminus k^{\times2}$,  with $\sigma^*(\sqrt{a})=-\sqrt{a}$; to be precise, we fix a separable closure $k^\sep$ of $k$ and fix a choice of square root $\sqrt{a}\in k^\sep$. In addition we fix an isomorphism 
\begin{equation}\label{eqn:SqrtIso}
\theta_x{\colon}k[x]/(x^2-a)\xrightarrow{\sim} k(\sqrt{a})
\end{equation}
by sending $x$ to $\sqrt{a}$.

We let $X_a:=\Spec k(\sqrt{a})\times^N\SL_2$. Thus
\begin{equation}\label{eqn:ChangeGroup}
H^*_N(z, \sW)\cong H^*_{\SL_2}(X_a, \sW),\quad H^*_N(z, \sW)\cong H^*_{\SL_2}(X_a, \sW).
\end{equation}
One can use Hilbert's theorem 90 to show that $\Pic_N(z)=\Pic(X_a)=\{0\}$, so there is no need to consider twisted coefficients.

We retain the Notations~\ref{not:QuadraticForms}.  We let $F$ be the (right) representation  of $\SL_2$ corresponding to the inclusion $\SL_2\subset \GL_2$, giving  (right) $\SL_2$-action on the projective line $\P(F)$ and on the projective plane $\P(\Sym^2F)$. Giving $F$ the standard basis $e_0, e_1$, we have the dual basis $X_0, X_1$ of linear forms on $\P(F)$ and the basis $T_0, T_1, T_2$ of linear forms on $\P(\Sym^2F)$ dual to $e_0^2, e_0e_1, e_1^2$. In particular, we have the $\SL_2$-scheme $Y:=\P(\Sym^2F)\setminus\P(F)$, the non-degenerate quadratic form $\<Q\>$ on $\sO_Y(-1)$ defined by $Q(T_0,T_1, T_2):=T_1^2-4T_0T_2$, the embedding $\sq{\colon}\P(F)\hookrightarrow \P(\Sym^2F)$ with image curve $D$ defined by $Q=0$, and the $\omega_Y$-valued non-degenerate quadratic form $\tilde{Q}$ on $\sV:=\sO_Y(-2)\otimes F$ given by the symmetric matrix 
\[
\tilde{B}:=\begin{pmatrix}2T_0&T_1\\T_1&2T_2\end{pmatrix}
\]
with  discriminant $-Q$. The forms $Q$ and $\tilde{Q}$ are $\SL_2$-invariant, giving rise to elements $[\<Q\>]\in H^0_{\SL_2}(Y, \sW)$ and  $[\tilde{Q}]\in H^0_{\SL_2}(Y, \sW(\omega_Y))$, respectively. 

We have the 2-1 map $\bar\pi{\colon}\P^1\times\P^1\to \P^2$ defined by 
\[
\bar\pi([x_0,x_1], [y_0, y_1])=[x_0y_0, x_0y_1+x_1y_0, x_1y_1]. 
\]
Via the (modified) Segre embedding $i{\colon}\P^1\times\P^1\to \P^3$, 
\[
i([x_0,x_1], [y_0, y_1])=[x_0y_0, x_0y_1+x_1y_0, x_1y_1, x_0y_1-x_1y_0], 
\]
$i(\P^1\times\P^1)$ is the quadric defined by $T_3^2-T_1^2+4T_0T_2$, the map $\bar\pi$ is the restriction of the projection $\pi{\colon}\P^3\setminus\{(0,0,0,1)\}\to \P^2$, $[T_0, T_1, T_2, T_3]\mapsto [T_0, T_1, T_2]$, and $\bar\pi$ defines a 2-1 cover of $\P^2$, ramified over $D$.

We embed $\P^2:=\P(\Sym^2F)$ in $\P^3:=\P(\Sym^2F\oplus k)$ by $[t_0, t_1, t_2]\mapsto [t_0, t_1, t_2,0]$. 
and let $\SL_2$ act on $\P^3$ via its action on $\P^2$, with trivial action on the last coordinate $T_3$. Via this embedding, we consider $D$ as a curve in $\P^3$, stable under the $\SL_2$-action. 

We let $\bar{X}_a$ denote the twisted form of $\P^1\times\P^1$ defined by the $\Z/2$-action on 
$\Spec k(\sqrt{a})\times_k  \P^1\times\P^1$ with generator acting by conjugation on $k(\sqrt{a})$ and by exchange of factors on $\P^1\times\P^1$; in other words, $\bar{X}_a$ is the restriction of scalars of $\P^1_{k(\sqrt{a})}$ over $k$. The diagonal (right) action of $\SL_2$ on $\P^1\times\P^1$ gives an $\SL_2$-action on $\bar{X}_a$. The extension of scalars of $\bar\pi$ to  
\[
\bar\pi_{k(\sqrt{a})}{\colon}\Spec k(\sqrt{a})\times_k  \P^1\times\P^1\to \Spec k(\sqrt{a})\times_k  \P^2
\]
defines by descent the finite morphism $\bar\pi_a{\colon}\bar{X}_a\to \P^2$ over $k$. 

\begin{lemma}\label{lem:Identification} 1. The $k$-scheme $\bar{X}_a$ is $\SL_2$-equivariantly isomorphic to the quadric $V(Q_a)$ in $\P^3$ defined by $Q_a:=T^2_3-a(T_1^2-4T_0T_2)$, with $\bar\pi_a$ equal to the restriction of $\pi$ to $V(Q_a)$. \\[2pt]
2. The $k$-scheme $X_a$ is $\SL_2$-equivariantly isomorphic to the complement of $D\subset\P^2\subset \P^3$ in $V(Q_a)$. \end{lemma}

\begin{proof} (1) Let $k(\sqrt{a})[\P^1\times\P^1]$ denote the bi-homogeneous coordinate ring of 
$(\P^1\times\P^1)_{k(\sqrt{a})}$. The   invariants of bi-degree $\{(n,n)\}_{n\ge0}$ of the $\Z/2$-action on $k(\sqrt{a})[\P^1\times\P^1]$  is the $k$-subalgebra of $k(\sqrt{a})[\P^1\times\P^1]$ generated by $t_0:=x_0y_0$, $t_2:= x_1y_1$, $t_1:=x_0y_1+x_1y_0$ and $t_3:=\sqrt{a}(x_0y_1-x_1y_0)$. This defines the $\SL_2$-equivariant embedding $\bar{i}{\colon}\bar{X}_a\hookrightarrow\P^3$ with $\bar{i}^*(T_i)=t_i$, and with image the quadric defined by $Q_a$.  One easily sees that $\bar\pi_a{\colon}\bar{X}_a\to \P^2$ goes over to the restriction of $\pi$ to $V(Q_a)$. 

 Let $T\subset \GL_2$, $T_1\subset \SL_2$ be the diagonal tori. Then  $T_1\backslash \SL_2\cong T\backslash\GL_2\cong \P^1\times\P^1\setminus\Delta$ by  the map
\[
\begin{pmatrix}x_0&x_1\\y_0&y_1\end{pmatrix}\mapsto ([x_0,x_1],[y_0,y_1])
\]
The action by $\sigma\in \SL_2$ becomes the exchange of factors. Thus
$X_a$ is isomorphic to the quotient of $\Spec k(\sqrt{a})\times_k  \P^1\times\P^1\setminus\Delta$ by the $\Z/2$-action used in the definition of $\bar{X}_a$. Sending 
$([x_0,x_1],[y_0,y_1])$ to $[t_0, t_1, t_2, t_3]$ sends $\Delta$ to $D\subset \P^2\subset \P^3$ and gives the desired isomorphisms of $\bar{X}_a$ with $V(Q_a)$ and of  $X_a$ with $V(Q_a)\setminus D$.
\end{proof} 

We henceforth identify $\bar{X}_a$ with $V(Q_a)$ and $X_a$ with $V(Q_a)\setminus D$.

As a quadric in $\P^3$ with rational point $p:=[1,2,1,0]$, $\bar{X}_a$ is birationally isomorphic to $\P^2$ via projection from $p$, $\bar{\mathfrak{p}}{\colon}\bar{X}_a\dashrightarrow\P^2$. Using the explicit rational map $\mathfrak{p}{\colon}\P^3\dashrightarrow\P^2$, $[t_0,\ldots, t_3]\mapsto [t_0+t_2-t_1, t_0-t_2,  t_3]$,  $\bar{\mathfrak{p}}$ is the restriction of $\mathfrak{p}$ to $\bar{X}_a$. Via $\bar{\mathfrak{p}}$,  the point $p$ gets blown up to the line $\ell:=V(T_0)$ in $\P^2$, and the curve $D=\bar{X}_a\cap (T_3=0)$ gets sent by $\bar{\mathfrak{p}}$ isomorphically to the line $\ell':=V(T_2)$ in $\P^2$.

The (projective) tangent plane $T_p\bar{X}_a$ to $\bar{X}_a$ at $p$ is the plane $V(T_0+T_2-T_1)$; let $C_p:= T_p\bar{X}_a\cap \bar{X}_a$. Over $k(\sqrt{a})$,  $C_p$ is the union of the two lines contained in $\bar{X}_a$ and containing $p$; the two lines are conjugate over $k$ and $C_p$ is irreducible as $k$-scheme. $C_p$ can also be defined by equations
\[
T_3^2-a(T_2-T_0)^2=0=T_0+T_2-T_1. 
\]
Via $\bar{\mathfrak{p}}$, $C_p$ get blown down to the point $y\in \P^2$ defined by $T_0=0$, $T_2^2=aT_1^2$; viewed in the other direction, we construct $\bar{X}_a$ from $\P^2$ by first blowing up $y$, and then blowing down the proper transform $E \subset\Bl_y\P^2$ of $\ell$ to the point $p\in \bar{X}_a$. 

We may thus factor $\bar{\mathfrak{p}}$ as $\phi\circ\psi^{-1}$, where $\phi{\colon}\Bl_y\P^2\to \P^2$ and $\psi{\colon}\Bl_p\bar{X}_a\to \bar{X}_a$ are the blow-down maps, together with the identification $\Bl_y\P^2=\Bl_p\bar{X}_a$. 


 The normalization $C_p^N$ of $C_p$ is isomorphic to $\P^1_{k(y)}$; indeed, $C_p^N$ is the exceptional divisor of $\Bl_y\P^2\to \P^2$, which is isomorphic to $\P^1_{k(y)}$, the isomorphism depending on a choice of local parameters $(x_1, x_2)$ for $\sO_{\P^2, y}$. Since $y$ is the closed subscheme $V(T_0, T_2^2-aT_1^2)$ of $\P^2$, we have the isomorphism \eqref{eqn:SqrtIso} $\theta_0:=\theta_{T_2/T_1}{\colon}k(y)\to k(\sqrt{a})$, giving us an isomorphism $\theta_p{\colon}\P^1_{k(\sqrt{a})}\xrightarrow{\sim}C_p^N$. 
 
 We summarize the geometric picture in the following diagram.
\[
\xymatrixcolsep{4pt}
\xymatrixrowsep{8pt}
\xymatrix{
&&&&0_{k(y)}\ar@{^(->}[drr]\ar@{^(->}[dll]\\
&&\P^1_{k(y)}\ar@{}[rr]|-\subset\ar[ddll]\ar[ddrrrrrrr]|(.36)\hole|(.82)\hole&&  \Bl_y\P^2=\Bl_p\bar{X}_a\ar[dd]^(.6)\phi\ar[ddrrrrrrr]^(.77)\psi&&E\ar[ddllll]|(.24)\hole|(.43)\hole|(.52)\hole\ar@{}[ll]|-\supset\ar@{}[dd]|-{\hole\hole}\ar[ddr]|(.5)\hole\\\\
y\ar@{}[rr]|-\in&&\ell\ar@{}[rr]|-\subset&&\P^2&\ell'\ar@{}[l]|-\supset&&p\ar@{}[rr]|-\in&&C_p\ar@{}[rr]|(.45)\subset &&
\bar{X}_a\ar@{-->}@/^20pt/[lllllll]^{\bar{\mathfrak{p}}}&&\ar@{}[ll]|-\supset D\ar@/^20pt/[llllllll]|(.6)\hole|(.65)\hole\\\\}
\]
Note that $y$ is not on $\ell'$, since $k(y)\supsetneq k$ and $\ell\cap \ell'$ is a single $k$-rational point.

For an arbitrary point $q\in D\subset \bar{X}_a$, let $T_q\bar{X}_a\subset \P^3$ be the (projective) tangent plane to $\bar{X}_a$ at $q$ and let $C_q=T_q\bar{X}_a\cap \bar{X}_a$.  Since $C_p$ is defined above as $T_p\bar{X}_a\cap  \bar{X}_a$, our notation in case $q=p$ is unambiguous.

Take $q\in D$, $q\neq p$. Then $p$ is not in $C_q$, so the birational map $\bar{\mathfrak{p}}{\colon}\bar{X}_a\dashrightarrow \P^2$ is defined on all of $C_q$; it follows from the fact that $C_q$ and $C_p$ intersect transversely on $\bar{X}_a$ that $\bar{\mathfrak{p}}{\colon}C_q\to \bar{\mathfrak{p}}(C_q)$ is in fact an isomorphism. The point $q$ is sent to a point $q':=\bar{\mathfrak{p}}(q)$ on the line $\ell'=V(T_2)$ and $\bar{\mathfrak{p}}(C_q)$ is the cone $C_{q'}(y)$ over $y$ with vertex $q'$. For  $q\in D(k)\setminus\{p\}$, we henceforth identify $C_q$ with the cone $C_{q'}(y)$ via $\bar{\mathfrak{p}}$.

Thus, for $q\neq p$ a $k$-point of $D$, the normalization $C_q^N$ of $C_q$ is isomorphic to the normalization $C_{q'}(y)^N$. As $C_{q'}(y)$ is the cone in $\P^2_k$ over $y$ with vertex $q'\in \P^2(k)$, $C_{q'}(y)^N$ is isomorphic to $\P^1_{k(y)}$; using the isomorphism $\theta_0{\colon}k(y)\xrightarrow{\sim}k(\sqrt{a})$ as above, this gives an isomorphism $\theta_q{\colon}\P^1_{k(\sqrt{a})}\xrightarrow{\sim} C_q^N$.
 
As a preliminary, we recall the well-known computation of the Witt-sheaf cohomology for projective spaces (see, e.g., \cite[Theorem 11.7]{FaselPBF} for the cohomology of the Milnor-Witt sheaves $\sK^{MW}_j$, noting that $\sW=\sK^{MW}_j$ for $j<0$).

\begin{lemma}\label{lem:PnCoh} Take $U\in\Sm/k$. For $n>0$, we have
\[
H^i(\P^n\times U, \sW)=\begin{cases} H^0(U, \sW)&\text{ for }i=0\\0&\text{ for }0<i<n\\0&\text{ for }i=n\text{ even}\\
H^0(U, \sW)&\text{ for }i=n\text{ odd.}
\end{cases}
\]
and
\[
H^i(\P^n\times U, \sW(\sO(1)))=\begin{cases}0&\text{ for }0\le i<n\\0&\text{ for }i=n\text{ odd}\\
H^0(U, \sW)&\text{ for }i=n\text{ even.}
\end{cases}
\]
The identity $H^0(U, \sW)=H^0(\P^n\times U, \sW)$ is induced by pull-back from $U$, and in cases $H^n(\P^n\times U, \sW)=H^0(U, \sW)$ or $H^n(\P^n\times U, \sW(\sO(1)))=H^0(U, \sW)$, the identity is induced by push-forward for the inclusion $p\times U\hookrightarrow \P^n\times U$ for  a $k$-point $p$.
\end{lemma}

\begin{proof}[Sketch of proof] Use the localization sequence for $H^*(-,\sW(\sO(m)))$ for the closed immersion $i{\colon}\P^{n-1}\times U\to \P^n\times U$ with open complement $j{\colon}\A^n\times U\to \P^n\times U$, and show that the boundary map
\[
H^0(U, \sW)=H^0(\A^n\times U, \sW)\to H^0(\P^{n-1}\times U, \sW(\sO(m-1)))
\]
is an isomorphism for $m$ odd and the zero map for $m$ even and $n=1$; for this, it suffices to make the computation after replacing $U$ with $\Spec k(U)$.  This proves the result for $n=1$ and the result in general follows by induction.
\end{proof}

Next, we will compute the Witt sheaf cohomology of $\bar{X}_a$, but first we need a discussion on orientations. 

\begin{definition} Let $Y$ be a smooth $k$-scheme. Recall that an {\em orientation} on $Y$ is a pair $(\rho, \sL)$ with $\sL$ an invertible sheaf on $Y$ and $\rho{\colon}\omega_{Y/k}\xrightarrow{\sim} \sL^{\otimes 2}$ an isomorphism.  Given an orientation $(\rho,\sL)$ on $Y$, we call a local section $\eta$ of $\omega_Y$ {\em oriented} if $\rho(\eta)=\lambda^2$ for some local section $\lambda$ of $\sL$.

More generally, given a morphism $f{\colon}Y\to X$ in $\Sm/k$, we have the relative dualizing sheaf $\omega_f:=\omega_{Y/k}\otimes f^*\omega_{X/k}^{-1}$. An orientation for $f$  is a 
pair $(\rho, \sL)$ with $\sL$ an invertible sheaf on $Y$ and an isomorphism $\rho{\colon}\omega_f\xrightarrow{\sim}\sL^{\otimes 2}$. We define the notion of a local oriented section of $\omega_f$ as for $\omega_{Y/k}$.
\end{definition}

Orientations arise when one needs to compare pullbacks (defined via cohomology) and proper pushforwards (defined using Borel-Moore homology) as they give a way of ``untwisting'' the twists arising in the comparison isomorphism. We concentrate on the case of Witt-sheaf cohomology. An orientation $o:=(\rho, \sL)$ on some $Y\in \Sm/k$ gives an isomorphism
\[
\tau_o{\colon}H^*(Y, \sW(\omega_Y))\xrightarrow{\sim} H^*(Y, \sW)
\]
as the composition  $H^*(Y, \sW(\omega_Y))\xrightarrow{\rho_*}H^*(Y, \sW(\sL^{\otimes 2}))\xrightarrow{can} H^*(Y, \sW)$. Given orientations $o:=(\rho_X, \sL_X)$ for $X$ and $o':=(\rho_Y, \sL_Y)$ for $Y$, and a proper morphism $f{\colon}Y\to X$ of relative dimension $d$, we  define the pushforward map $f_*^{o,o'}{\colon}H^*(Y, \sW)\to H^{*-d}(X, \sW)$ as the unique map making
\[
\xymatrix{
H^*(Y, \sW(\omega_Y))\ar[d]^{f_*}\ar[r]^-{\tau_{o'}}_-\sim&H^*(Y, \sW)\ar[d]^{f_*^{o,o'}}\\
H^{*-d}(X, \omega_X)\ar[r]^-{\tau_o}_-\sim&H^{*-d}(X, \sW)
}
\]
commute. We will often drop the various decorations if the context makes the meaning clear.

For our purposes, we will need orientations for $\bar{X}_a$, $D$, $\P^1$ and $\P^1_{k(\sqrt{a})}$, which we now fix. 

Give $\P^n$ the standard homogeneous coordinates $T_0,\ldots, T_n$. We have the generating global section $\Omega_n$ of $\omega_{\P^n/k}(n+1)$, 
\begin{equation}\label{eqn:OrOmega}
\Omega_n:=\sum_{i=0}^n(-1)^idT_idT_0\wedge\ldots\wedge\widehat{dT_i}\wedge\ldots\wedge dT_n
\end{equation}
which we declare to be an oriented global generator for $\omega_{\P^n/k}(n+1)$.

\begin{lemma}\label{lem:OriGen} Let $t_i=T_i/T_0$. $i=1,2,\ldots$.
\begin{enumerate}
\item  For $n$ odd, $dt_1\cdots dt_n$ is an oriented generator of $\omega_{\P^n}$ over  $T_0\neq0$.
\item For $n=1$, let $t=t_1$.   Then $dt$ is an oriented generator of $\omega_{\P^1}$ over $T_0\neq$ and $-d(t^{-1})$ is an oriented generator over $T_1\neq 0$. The same holds for $\P^1_{k(\sqrt{a})}$ by taking the base-extension.
\item We identify $D\subset\P^3$ with $\P^1$ via the isomorphism $\phi{\colon}\P^1\to D$, $\phi([x_0, x_1])=[x_0^2, 2x_0x_1, x_1^2,0])$. Then the oriented local generator $dt$ of $\omega_{\P^1}$  induces the oriented generator $(1/2)dt_1$ (restricted to $D$) of $\omega_D$ over $T_0\neq0$.
\item Identify $\bar{X}_a$ with $V(Q_a)\subset \P^3$ as above. Then $(-1/4a)dt_1dt_3$ (restricted to $\bar{X}_a$)  is an oriented generator of $\omega_{\bar{X}_a}$ over $T_0\neq0$.
\end{enumerate}
\end{lemma}

\begin{proof} Since $\Omega_n$ is an oriented global generator for $\omega_{\P^n/k}(n+1)$, the generator $\Omega_n/T_0^{n+1}$ is a generator for  $\omega_{\P^n/k}$ over $T_0\neq0$ and is an oriented generator if $n$ is odd. Globally, for $n$ odd, we have the orientation $\omega_{\P^n/k}\cong \sO_{\P^n}(-(n+1)/2)^{\otimes 2}$ by sending a local section $\lambda$ of $\sO_{\P^n}(-(n+1)/2)$ to $\Omega_n\cdot \lambda^2$. 
 
For  $\P^1$, this gives us the oriented generators $dt_1=dt=\Omega_1/T_0^2$ on $T_0\neq0$ and $-dt^{-1}=\Omega_1/T_1^2$ on $T_1\neq0$. By pullback, this gives us an orientation and local oriented sections of $\omega_{\P^1_{k(\sqrt{a})}}$, and via the isomorphism $\phi{\colon}\P^1\to D$, we have the same for $D$. Since $\phi^*(t_1)=2t$, we have shown (1) and (2).

For (3),  we use the residue isomorphism to induce an orientation for $\bar{X}_a$ from that of $\P^3$. 

Recall that for $i_W{\colon}W\hookrightarrow Y$ a smooth codimension one closed subscheme of some $Y\in \Sm/k$, we have the ideal sheaf $\sI_W$, a locally principal sheaf of ideals with dual $\sO_Y(W)$, where a section of $\sO_Y(W)$ is locally of the form $g/f$ with $g$ a local section of $\sO_Y$ and $f$ a local section of $\sI_W$.   This gives the invertible sheaf $\omega_Y(W):=\omega_Y\otimes_{\sO_Y}\sO_Y(W)$, where a section is locally of the form $\eta/f$, with $\eta$ a local section of  $\omega_Y$ and $f$ a local section of $\sI_W$. 

The residue map is defined by writing a local section $\eta/f$ of $\omega_Y(W)$ locally as $(df/f)\wedge \tau$ and then defining $\res_W((df/f)\wedge \tau):=i_W^*\tau$. One then shows that this gives a well-defined isomorphism 
\[
\res_W{\colon}i_W^*(\omega_Y(W))\xrightarrow{\sim} \omega_W.
\]

If we now have an orientation $\rho{\colon}\omega_Y\to \sL^{\otimes 2}$ and an isomorphism $\sO_Y(W)\cong \sM^{\otimes 2}$, this gives an isomorphism $i_W^*((\omega_Y(W)))\cong i_W^*(\sL\otimes \sM)^{\otimes 2}$, which in turn gives the isomorphism $\rho_W{\colon}\omega_W\to 
i_W^*(\sL\otimes \sM)^{\otimes 2}$ using the residue map. 

For $\bar{X}_a=V(Q_a)\subset \P^3$, as $Q_a$ is a section of $\sO_{\P^3}(2)$, multiplication by $Q_a$ defines an isomorphism  $\sO_{\P^3}(\bar{X}_a)\cong \sO_{\P^3}(2)\cong \sO_{\P^3}(1)^{\otimes 2}$, and thus the residue map and our orientation for $\P^3$ gives us the isomorphism $\rho_{\bar{X}_a}{\colon}\omega_{\bar{X}_a}\xrightarrow{\sim}\sO_{\bar{X}_a}(2-4)\cong 
\sO_{\bar{X}_a}(-1)^{\otimes 2}$. 

An explicit oriented section of  $\omega_{\bar{X}_a}$ on $T_0\neq0$ can be constructed as follows. We have our oriented section $dt_1dt_2dt_3$ of $\omega_{\P^3}$ on $T_0\neq0$, which gives us the oriented section   $(T_0^2/Q_a)dt_1dt_2dt_3=dt_1dt_2dt_3/q_a$ of $\omega_{\P^3}(\bar{X}_a)$, where $q_a=t_3^2-a(t_1^2-4t_2)=Q_a/T_0^2$. Using the relation $dq_a=2t_3dt_3-2at_1dt_1+4adt_2$, we have 
\[
dt_1dt_2dt_3/q_a=(dq_a/q_a)\wedge (-1/4a)dt_1dt_3
\]
so $\res_{\bar{X}_a}(dt_1dt_2dt_3/q_a)=(-1/4a)dt_1dt_3$ is an oriented generator of $\omega_{\bar{X}_a}$ over $T_0\neq0$.
\end{proof}

\begin{remark} \label{rem:SLInvOrientation} The isomorphisms $\omega_{\P^n/k}(n+1)\cong \sO_{\P^n}$ using the generator $\Omega_n$ are $\SL_{n+1}$-equivariant, since the chosen global section $\Omega_n$ is $\SL_{n+1}$-invariant. This implies that our orientations for $\bar{X}_a$, $D$, $\P^1$ and $\P^1_{k(\sqrt{a})}$ are all $\SL_2$-equivariant.
\end{remark}

\begin{remark} \label{rem:RSMultilinearAlgebra}
An additional remark concerning the translation to the pushforward in the Rost-Schmid complex is in order. If we have a codimension $c$ closed  immersion $i{\colon}Z\to Y$ in $\Sm/k$, we have pushforward maps $i_{*,n}{\colon}C^n_{RS}(Z, \sW(\omega_Z))\to C^{n+c}_{RS}(Y, \sW(\omega_Y))$, giving  $i_{*}{\colon}C^*_{RS}(Z, \sW(\omega_Z))\to C^{*}_{RS}(Y, \sW(\omega_Y))[c]$ with induced map on cohomology $i_*{\colon}H^n(Z, \sW(\omega_Z))\to H^{n+c}(Y, \sW(\omega_Y))$ the same as the pushforward given by the six-functor formalism (see \cite[Proposition 3.4.3]{DFJ}  and Remark~\ref{rem:CycModFunct}). 

To make  $i_{*,n}$ explicit,  if we have $z\in Z^{(n)}$, then via $i$, $z$ is also in $Y^{(n+c)}$. Let $i_z^Z{\colon}z\to Z$, $i_z^Y{\colon}z\to Y$ be the inclusions. We need to define a map
\[
i_{*z}{\colon}W(k(z),\Lambda^Z_z\otimes i_z^{Z*}\omega_Z)\to W(k(z),\Lambda^Y_z\otimes i_z^{Y*}\omega_Y),
\]
and then 
\[
i_{*,n}{\colon}\oplus_{z\in Z^{(n)}}W(k(z),\Lambda^Z_z\otimes i_z^{Z*}\omega_Z)\to 
\oplus_{y\in Y^{(n+c)}}W(k(y),\Lambda^Y_y\otimes i_y^{Y*}\omega_Y)
\]
is $\prod_{z\in Z^{(n)}}i_{*z}$. 

The map $i_{*z}$ arises from an isomorphism $\rho^{Y/Z}_z{\colon}\Lambda^Z_z\otimes i_z^{Z*}\omega_Z\xrightarrow{\sim}
\Lambda^Y_z\otimes i_z^{Y*}\omega_Y$ by just applying  $\rho^{Y/Z}_z$ to the twists. To define $\rho^{Y/Z}_z$, let $\mathfrak{m}_z^Z\subset \sO_{Z,z}$ and $\mathfrak{m}_z^Y\subset \sO_{Y,z}$ be the respective maximal ideals. We have the short exact sequences 
\[
0\to \mathfrak{m}^Z_z/(\mathfrak{m}_z^Z)^2\to i^{Z*}_z\Omega_{Z/k}\to \Omega_{z/k}\to 0
\]
\[
0\to \mathfrak{m}^Y_z/(\mathfrak{m}_z^Y)^2\to i^{Y*}_z\Omega_{Y/k}\to \Omega_{z/k}\to 0
\]
and taking determinants induce the isomorphisms 
\[
\Lambda^Y_z\otimes i^{Y*}_z\omega_{Y/k}\xrightarrow{\sim} \omega_{z/k}
\]
\[
\Lambda^Z_z\otimes i^{Z*}_z\omega_{Z/k}\xrightarrow{\sim} \omega_{z/k},
\]
giving the isomorphism $\rho^{Y/Z}_z$. To be perfectly explicit, if we have a short exact sequence of finitely generated projective $A$ modules for some commutative ring $A$
\[
0\to P'\to P\to P''\to0
\]
we define the isomorphism $\det_AP'\otimes \det_AP''\to \det_AP$ by choosing a splitting $s{\colon}P''\to P$ and sending $(q_1\wedge\ldots\wedge q_n)\otimes (t_1\wedge\ldots\wedge t_m)$ to $q_1\wedge\ldots\wedge q_n\wedge s(t_1)\wedge\ldots\wedge s(t_m)$. This is independent of the choice of splitting, so extends to an isomorphism $\det\sV'\otimes_{\sO_Y}\det\sV''\xrightarrow{\sim} \det\sV$ for each exact sequence $0\to \sV'\to \sV\to \sV''\to0$ of locally free coherent $\sO_Y$-modules.  Similarly we have the isomorphism $\det^{-1}\sV\cong
\det\sV^\vee$ defined by the perfect pairing $\det\sV^\vee\otimes \det\sV\to \sO_Y$ given on local sections $v_1,\ldots, v_n$ of $\sV$, $l_1,\ldots, l_n$ of $\sV^\vee$, $n=\rnk\sV=\rnk\sV^\vee$ by 
\[
\<l_1\wedge\ldots\wedge l_n, v_1\wedge\ldots\wedge v_n\>=\sum_{\sigma\in S_n}\prod_{i=1}^n l_i(v_{\sigma(i)}).
\]

Now if we have parameters $t_1,\ldots, t_{m+n}$ for  $\mathfrak{m}^Y_z$ such that the images $\bar{t}_1,\ldots,  \bar{t}_n$ of  $t_1,\ldots, t_n$ are parameters for $\mathfrak{m}^Z_z$ and  $\bar{t}_{n+1}=\cdots=\bar{t}_{n+m}=0$ in $\sO_{Z,z}$, and we have  $\tilde{s}_{m+n+1},\ldots, \tilde{s}_{m+n+r}\in k(z)$ such that $d\tilde{s}_{m+n+1},\ldots, d\tilde{s}_{m+n+r}$ is a $k$-basis of $\Omega_{k(z)/k}$,  then $\del/\del t_1\wedge\ldots\wedge \del/\del t_{m+n}$ is a $k(z)$ basis for $\Lambda^Y_z$ dual to the basis $dt_1\wedge\ldots\wedge dt_{m+n}$ of $\sI^Y_z/(\sI^Y_z)^2$, $\del/\del \bar{t}_1\wedge\ldots\wedge \del/\del \bar{t}_n$ is a $k(z)$ basis for $\Lambda^Z_z$ dual to the basis $d\bar{t}_1\wedge\ldots\wedge d\bar{t}_n$ of $\sI^Z_z/(\sI^Z_z)^2$, and if we lift $\tilde{s}_j$ to $\bar{s}_j\in \sO_{Z,z}$,  and then further lift to $s_j\in \sO_{Y,z}$, then $d\bar{t}_1\wedge\ldots\wedge d\bar{t}_n\wedge d\bar{s}_{m+n+1}\wedge\ldots\wedge  d\bar{s}_{m+n+r}$ is a $k(z)$-basis of $i_z^{Z*}\omega_{Z/k}$, $dt_1\wedge\ldots\wedge dt_{m+n}\wedge ds_{m+n+1}\wedge\ldots  ds_{m+n+r}$ is a $k(z)$-basis of $i_z^{Y*}\omega_{Y/k}$ and $d\tilde{s}_{m+n+1}\wedge\ldots\wedge d\tilde{s}_{m+n+r}$ is a $k(z)$-basis for $\omega_{k(z)/k}$. Then our isomorphisms $\Lambda^Y_z\otimes i^{Y*}_z\omega_{Y/k}\to \omega_{k(z)/k}\to  \Lambda^Z_z\otimes i^{Z*}_z\omega_{Z/k}$ are given by
\begin{multline*}
\del/\del t_1\wedge\ldots\wedge \del/\del t_{m+n}\otimes dt_1\wedge\ldots\wedge dt_{m+n}\wedge ds_{m+n+1}\wedge\ldots  \wedge ds_{m+n+r}\\ \mapsto 
d\tilde{s}_{m+n+1}\wedge\ldots\wedge d\tilde{s}_{m+n+r}\\\mapsto
\del/\del \bar{t}_1\wedge\ldots\wedge \del/\del \bar{t}_n\otimes d\bar{t}_1\wedge\ldots\wedge d\bar{t}_n\wedge d\bar{s}_{m+n+1}\wedge\ldots  \wedge d\bar{s}_{m+n+r}.
\end{multline*}
\end{remark}

For the rest of this section, we fix the point $q$ to be $[1,0,0,0]$ and we fix the morphism $p_q{\colon}\P^1_{k(\sqrt{a})}\to \bar{X}_a$ to be the map $p_q([x_0,x_1])=[x_0+x_1, x_1, 0, \sqrt{a}x_1]$. The map $p_q$ has image $C_q$ and gives a birational morphism $\P^1_{k(\sqrt{a})}\to C_q$, which in turn fixes our choice of isomorphism $\theta_q{\colon}\P^1_{k(\sqrt{a})}\to C_q^N$.  Note that $p_q^*(\sO_{\bar{X}_a}(1))\cong \sO_{\P^1_{k(\sqrt{a})}}(1)$.

\begin{lemma} \label{lem:OrientedExtension} Let $g$ be the rational function $t_2/t_3$ on $\bar{X}_a$, giving a section of the ideal sheaf $\sI_{C_q}$ over $T_0T_3\neq0$. Then $(1/2)p_q^*(\del/\del g)$ is an oriented global generator of $\omega_{p_q}$. 
\end{lemma}

\begin{proof} We use standard homogeneous coordinates $x_0, x_1$ on $\P^1_{k(\sqrt{a})}$. We have the local oriented generator $(-1/4a)dt_1dt_3$ for $\omega_{\bar{X}_a}$ over $T_0\neq0$ given by Lemma~\ref{lem:OriGen}.  Let $s_i=T_i/T_1$, We have
\begin{align*}
(-1/2)ds_0[d(s_2/s_3)-t_2d(1/t_3)]&=(-1/2)(d(t_1^{-1})(d(t_2/t_3)-t_2d(1/t_3))\\&=(1/2t_3)t_1^{-2}dt_1dt_2=(-1/4a)t_1^{-2}dt_1dt_3.
\end{align*}
so $(-1/2)ds_0[d(s_2/s_3)-t_2d(1/t_3))$ is an oriented generator of $\omega_{\bar{X}_a}$ over $T_0T_1T_3\neq0$. 
Since $p_q^*(t_2)=0$,   $(-1/2)p_q^*ds_0d(s_2/s_3)$ is an oriented generator of $p_q^*\omega_{\bar{X}_a}$ over $x_0x_1\neq0$. Setting $t=x_1/x_0$, 
$-dt^{-1}=t^{-2}dt$ is an oriented generator of $\omega_{\P^1_{k(\sqrt{a})}}$ over $x_1\neq 0$ (again by Lemma~\ref{lem:OriGen}), $p_q^*ds_0=dt^{-1}$, and thus
\begin{align*}
(-dt^{-1})\otimes (p_q^*(-1/2)ds_0d(s_2/s_3))^{-1}&=
(-dt^{-1})\otimes(dt^{-1}p_q^*(-1/2)d(s_2/s_3))^{-1}\\&=(1/2)p_q^*\del/\del(s_2/s_3)\\&=(1/2)p_q^*\del/\del(t_2/t_3)=(1/2)p_q^*\del/\del g
\end{align*}
is an oriented generator of $\omega_{p_q}$ over $x_0x_1\neq0$.

Note that $g=s_2/s_3$ is a regular function on $\bar{X}_a\setminus (T_1T_3=0)$ and generates $\sI_{C_q}$ on that open subset. Moreover $V(T_1T_3)\cap C_q=\{q\}$. Since $\rho_q{\colon}\P^1_{k(\sqrt{a})}\to C_q$ is an isomorphism over $C_q\setminus\{q\}$, it follows that $p_q^*(\del/\del g)$ generates $\omega_{p_q}$ over $x_0\neq0$. But since $p_q^*\sO_{\bar{X}_a}(1)\cong \sO_{\P^1_{k(\sqrt{a})}}(1)$, we have
\[
\omega_{p_q}=\omega_{\P^1_{k(\sqrt{a})}}\otimes p_q^*\omega_{\bar{X}_a}^{-1}\cong
\sO_{\P^1_{k(\sqrt{a})}}(-2)\otimes p_q^*\sO_{\bar{X}_a}(2)\cong \sO_{\P^1_{k(\sqrt{a})}}.
\]
Thus, a section of $\omega_{p_q}$ that generates over $x_0\neq0$ is a global generator, completing the proof.
\end{proof}
 
 \begin{lemma}\label{lem:OrientedSectionD} $(-2a)\cdot i_D^*(\del/\del t_3)$ is a local oriented section of $\omega_{i_D}$ over $D\setminus (T_0=0)$.
 \end{lemma}
 
 \begin{proof} We work on the affine open subsets of $\P^3$, $\bar{X}_a$ and $D$ given by $T_0\neq0$, $\Omega_3/T_0^4=dt_1dt_2dt_3$ is an oriented generating section of $\omega_{\P^3}$ over $T_0\neq0$,   $(-1/4a)dt_1dt_3$ restricts to  an oriented generating section of $\omega_{\bar{X}_a}$ over $\bar{X}_a\setminus V(T_0)$ and $dt_1/2$ restricts to an oriented generating section of $\omega_D$, again over $D\setminus V(T_0)$. This gives us the oriented local generator $(dt_1/2)\otimes i_D^*((-1/4a)dt_1dt_3)^{-1}$ for $\omega_D\otimes i_D^*\omega_{\bar{X}_a}^{-1}$.

The short exact
\[
0\to \sI_D/\sI_D^2\to i_D^*\Omega_{\bar{X}_a}\to \Omega_D\to 0
\]
gives us the canonical isomorphism
\[
\alpha{\colon} \omega_D\otimes i_D^*\omega_{\bar{X}_a}^{-1}\xrightarrow{\sim} (\sI_D/\sI_D^2)^\vee=\omega_{i_D},
\]
with
\[
\alpha(dt_1\otimes (dt_1dt_3)^{-1})=i_D^*(\del/\del t_3).
\]
This yields $(-2a)i_D^*(\del/\del t_3)$ as an oriented local section of $ \omega_{i_D}$, as claimed.
\end{proof}

With all these details now being settled, we are ready to compute.

\begin{lemma}\label{lem:WittCohbarXa} 1. Pullback by the projection $\pi_{\bar{X}_a}{\colon}\bar{X}_a\to \Spec k$ defines an isomorphism  $H^0(\bar{X}_a,\sW)\cong W(k)$.\\[2pt]
2. The pushforward  
\[
p_{q*}{\colon}H^0(\P^1_{k(\sqrt{a})}, \sW(\omega_{\P^1_{k(\sqrt{a})}}))\to H^1(\bar{X}_a, \sW(\omega_{\bar{X}_a}))
\]
is an isomorphism, which, via our choice of orientations,  induces an isomorphism $\tilde{p}_{q*}{\colon}W(k(\sqrt{a}))=H^0(\P^1_{k(\sqrt{a})},\sW)\xrightarrow{\sim} H^1(\bar{X}_a, \sW)$. Here 
 the equality is the pullback $\pi_{\P^1_{k(\sqrt{a})}}^*{\colon}W(k(\sqrt{a}))\xrightarrow{\sim} H^0(\P^1_{k(\sqrt{a})},\sW)$ (see Lemma~\ref{lem:PnCoh}).
\\[2pt]
3. Recall the point $p=[1,2,1,0]\in \bar{X}_a$. The inclusions $i_p{\colon}p\to \bar{X}_a$, $i^D_p{\colon}p\to D$ and $i_D{\colon}D\to \bar{X}_a$ induce a commutative diagram of isomorphisms  
\[
\xymatrix{
W(k)=H^0(p, \sW)\ar[r]^{i_{p*}}\ar[dr]_{i^D_{p*}}&H^2(\bar{X}_a, \sW(\omega_{\bar{X}_a}))\\
&H^1(D, \sW(\omega_D))\ar[u]_{i_D*}
}
\]
which by our choice of orientations for $\bar{X}_a$ and $D$ induces
 a commutative diagram of isomorphisms  
\[
\xymatrix{
W(k)=H^0(p, \sW)\ar[r]^{i_{p*}}\ar[dr]_{i^D_{p*}}&H^2(\bar{X}_a, \sW)\\
&H^1(D, \sW)\ar[u]_{i_D*}
}
\]
4. For $i\ge3$, $H^i(\bar{X}_a, \sW)=0$.
\end{lemma}

\begin{proof} The assertion (4) follows from the fact that the Rost-Schmid complex for $\sW$ on $\bar{X}_a$, whose $i$th cohomology computes $H^i(\bar{X}_a, \sW)$, is supported in degrees $[0,2]$. 

Suppose we have a smooth surface $Y$ over $k$, a closed point $w\in Y$ and let $\tilde{Y}\to Y$ be the blow-up of $Y$ at $w$, with exceptional divisor $E\cong \P^1_w$. Let $V=Y\setminus\{w\}$. We have the exact localization sequences
\begin{equation}\label{eqn:LocSeq0}
\ldots\to H^{i-2}(w, \sW)\xrightarrow{i_{w*}} H^i(Y, \sW(\omega_Y))\to H^i(V, \sW(\omega_V))\xrightarrow{\del}H^{i-1}(w, \sW)\to\ldots
\end{equation}
and
\begin{multline}\label{mult:LocSeq1}
\ldots\to H^{i-1}(\P^1_w, \sW(\omega_{\P^1_w})))\xrightarrow{i_{\P^1_w*}} H^i(\tilde{Y}, \sW(\omega_{\tilde{Y}}))\\\to H^i(V, \sW(\omega_V))\xrightarrow{\del}H^{i}(\P^1_w, \sW(\omega_{\P^1_w}))\to\ldots
\end{multline}
We twist \eqref{eqn:LocSeq0} by $\omega_Y^{-1}$ to give 
\begin{equation}\label{eqn:LocSeq3}
\ldots\to H^{i-2}(w, \sW(i_w^*\omega_Y^{-1}))\xrightarrow{i_{w*}} H^i(Y, \sW)\to H^i(V, \sW)\xrightarrow{\del}H^{i-1}(w, \sW(i_w^*\omega_Y^{-1}))\to\ldots
\end{equation}
We twist \eqref{mult:LocSeq1} by $\omega_{\tilde{Y}}^{-1}$ and use our canonical isomorphism 
\[
\sO_{\P^1_w}(-1)\cong N_{\P^1_w/\tilde{Y}}\cong \omega_{\P^1_w}\otimes i_{\P^1_w}^*\omega_{\tilde{Y}}^{-1}
\]
to give the exact localization sequence 
\begin{multline}\label{mult:LocSeq2}
\ldots\to H^{i-1}(\P^1_w, \sW(\sO_{\P^1_w}(-1))))\xrightarrow{i_{\P^1_w*}} H^i(\tilde{Y}, \sW)\\\to H^i(V, \sW)\xrightarrow{\del}H^{i}(\P^1_w, \sW(\sO_{\P^1_w}))\to\ldots
\end{multline}
 
Since $H^j(\P^1_w,\sW(\sO(-1)))=0$ for all $j$ (Lemma~\ref{lem:PnCoh}), the sequence \eqref{mult:LocSeq2} shows that $ H^i(\tilde{Y}, \sW)\cong H^i(V, \sW)$ for all $i$.
From \eqref{eqn:LocSeq0}, we have $H^0(Y, \sW)\cong H^0(V, \sW)$, so the pullback $H^0(Y, \sW)\to H^0(\tilde{Y}, \sW)$ is an isomorphism. 

 We recall (Lemma~\ref{lem:PnCoh} again) that  $H^i(\P^2, \sW)=0$ for $i>0$ and that the pullback  $\pi_{\P^2}^*{\colon}W(k)=H^0(\Spec k, \sW)\to H^0(\P^2, \sW)$ is an isomorphism. The comments in the previous paragraphs thus prove (1). 
 
 For (3), we take $Y=\P^2$, $w=y$ in  \eqref{eqn:LocSeq3}. Since    $H^2(\P^2, \sW)=0$,  we see that $H^2(\P^2\setminus \{y\},\sW)=0$, and then \eqref{mult:LocSeq2} shows that $H^2(\Bl_y\P^2, \sW)=0$. We now take $Y=\bar{X}_a$, $w=p$. Since $\Bl_p\bar{X}_a$ is isomorphic to $\Bl_y\P^2$, we see by a similar argument that  $H^2(\bar{X}_a\setminus \{p\},\sW)=0$. Since $\bar{X}_a$ is oriented, so is $\bar{X}_a\setminus\{p\}$, so $H^2(\bar{X}_a\setminus \{p\},\sW(\omega_{\bar{X}_a\setminus \{p\}}))=0$ as well .
 
We have the push-forward map  $i_{p*}{\colon}W(k)=H^0(p, \sW)\to H^2(\bar{X}_a,\sW(\omega_{\bar{X}_a})$, and since $H^2(\bar{X}_a\setminus \{p\},\sW)=0$,  the  localization sequence \eqref{eqn:LocSeq0} shows that $i_{p*}$ is surjective. We also have the pushforward $\pi_{\bar{X}_a*}{\colon}H^2(\bar{X}_a,\sW(\omega_{\bar{X}_a}))\to W(k)$ with $\pi_{\bar{X}_a*}\circ i_{p*}=\id$, so $i_{p*}$ is an isomorphism, proving (3). 
 
 For (2), our  map $p_q$ and the birational map $\bar{\mathfrak{p}}{\colon}\bar{X}_a\dashrightarrow \P^2$  induces a birational morphism $\rho_q:=\bar{\mathfrak{p}}\circ p_q{\colon} 
 \P^1_{k(\sqrt{a})}\to C_{q'}(y)$, which is an isomorphism over  $C_{q'}(y)\setminus\{q'\}$.  
 Writing  $\infty_{k(\sqrt{a})}\in  \P^1_{k(\sqrt{a})}$ for $[0,1]$ and $0_{k(\sqrt{a})}$ for $[1,0]$,  we have $\rho_q^{-1}(q')=0_{k(\sqrt{a})}$ and $\rho_q^{-1}(y)=\infty_{k(\sqrt{a})}$.

Letting $\A^1_{k(\sqrt{a})}\subset  \P^1_{k(\sqrt{a})}$ be the complement of $\infty_{k(\sqrt{a})}$,  the inclusion $C_{q'}(y)\subset\P^2$ gives  the finite morphism $\rho_0{\colon}\A^1_{k(\sqrt{a})}=\P^1_{k(\sqrt{a})}\setminus\{\infty_{k(\sqrt{a})}\}\to \P^2\setminus\{y\}$. Via the standard homogeneous coordinates $X_0, X_1$ on $\P^1_{k(\sqrt{a})}$,   $\A^1_{k(\sqrt{a})}$ is the open subscheme $X_0\neq 0$, with coordinate $t:=X_1/X_0$, and with $\rho_q^{-1}(q')=0_{k(\sqrt{a})}$ the origin $t=0$.

 Let $\rho{\colon}\P^1_{k(\sqrt{a})}\to \P^2$ be the composition of $\rho_q$ with the inclusion $C_{q'}(y)\hookrightarrow\P^2$. Let $j{\colon}\P^2\setminus\{y\}\to \P^2$, $\bar{j}{\colon}\A^1_{k(\sqrt{a})}\to \P^1_{k(\sqrt{a})}$ be the inclusions. We have the relative dualizing sheaves $\omega_\rho:=\omega_{\P^1_{k(\sqrt{a})}/{k(\sqrt{a})}}\otimes \rho^*\omega_{\P^2/k}^{-1}$ and $\omega_{\rho_0}:=\omega_{\A^1_{k(\sqrt{a})}/{k(\sqrt{a})}}\otimes \rho_0^*\omega_{\P^2/k}^{-1}$, with $\omega_{\rho_0}=\bar{j}^*\omega_{\rho}$, giving us the commutative diagram 
 \[
 \xymatrix{
 H^i(\P^1_{k(\sqrt{a})}, \sW(\omega_\rho))\ar[r]^{\bar{j}^*}\ar[d]^{\rho_*}&
 H^i(\A^1_{k(\sqrt{a})}, \sW(\omega_{\rho_0}))\ar[d]^{\rho_{0*}}\\
 H^i(\P^2, \sW)\ar[r]^{j^*}& H^i(\P^2\setminus\{y\}, \sW)
 }
 \]

Let $i_\infty{\colon}\infty_{k(\sqrt{a})}\to \P^1_{k(\sqrt{a})}$, $i_y{\colon}y\to \P^2$ be the inclusions,  let $N_{\infty}$, $N_y$ be the corresponding normal bundles, and let $\bar\rho{\colon}\infty_{k(\sqrt{a})}\to y$ be the map induced by $\rho$. Since the map $\rho$ is unramified, we have the canonical isomorphism
 \[
 i_\infty^*\omega_\rho\otimes N_\infty\cong \bar\rho^*\det N_y.
 \]
 We choose an isomorphism $k(y)\cong \det N_y$,  inducing by the above isomorphism an isomorphism $i_\infty^*\omega_\rho\otimes N_\infty\cong k(\infty_{k(\sqrt{a})})$.
 This gives us the commutative diagram
\begin{equation}\label{eqn:ComLocDiag}
\xymatrix{
H^0(\A^1_{k(\sqrt{a})},\sW(\omega_{\rho_0})))\ar[d]^{\rho_{0*}}\ar[r]^-\del&H^0(\infty_{k(\sqrt{a})},\sW)\ar[d]_\wr^{\bar\rho_*}\ar[r]&H^1(\P^1_{k(\sqrt{a})}, \sW(\omega_\rho))\ar[d]^{\rho_*}\\
H^1(\P^2\setminus\{y\}, \sW)\ar[r]^-\del&H^0(y, \sW)\ar[r]&H^2(\P^2, \sW)
}
\end{equation}
with rows the respective localization sequences; note that the isomorphism $\bar\rho_*$ is the canonical one induced by the isomorphism $\bar\rho{\colon} \infty_{k(\sqrt{a})}\to y$, and is independent of any choices.

Since $\rho^*\sO_{\P^2}(1)\cong \sO_{\P^1_{k(\sqrt{a})}}(1)$, we have an isomorphism $\phi_0{\colon}\omega_\rho\xrightarrow{\sim}\sO_{\P^1_{k(\sqrt{a})}}(1)$, defined as the composition
\[
\omega_\rho=\omega_{\P^1_{k(\sqrt{a})}/{k(\sqrt{a})}}\otimes \rho^*\omega_{\P^2/k}^{-1}\cong
 \sO_{\P^1_{k(\sqrt{a})}}(-2)\otimes  \sO_{\P^1_{k(\sqrt{a})}}(3)= \sO_{\P^1_{k(\sqrt{a})}}(1).
 \]
Choosing an isomorphism $\bar{j}^*\sO_{\P^1_{k(\sqrt{a})}}(1)\cong \sO_{\A^1_{k(\sqrt{a})}}$ induces an isomorphism $\phi_1{\colon}\omega_{\rho_0}\xrightarrow{\sim}\sO_{\A^1_{k(\sqrt{a})}}$, giving rise to a commutative diagram of isomorphisms
\[
\xymatrix{
\omega_{\rho_0}\ar[r]^-\sim\ar[d]^{\phi_1}_\wr&\bar{j}^*\omega_\rho\ar[d]^{\bar{j}^*(\phi_0)}_\wr\\
\sO_{\A^1_{k(\sqrt{a})}}\ar[r]^-\sim&\bar{j}^*\sO_{\P^1_{k(\sqrt{a})}}(1)
}
\]
This gives an isomorphism of the top row in \eqref{eqn:ComLocDiag} to the corresponding part of the localization sequence
\begin{multline*}
\ldots \to H^0(\P^1_{k(\sqrt{a})},\sW(\sO(1)))\xrightarrow{\bar{j}^*}
H^0(\A^1_{k(\sqrt{a})},\sW)\xrightarrow{\del} H^0(\infty_{k(\sqrt{a})},\sW)\\
\xrightarrow{i_{\infty*}} H^1(\P^1_{k(\sqrt{a})},\sW(\sO(1)))\to \ldots
\end{multline*}
Since 
\[
H^0(\P^1_{k(\sqrt{a})},\sW(\sO(1)))=0=H^1(\P^1_{k(\sqrt{a})},\sW(\sO(1)))
\]
the boundary map $\del$ in the  top row of \eqref{eqn:ComLocDiag} an isomorphism. Similarly, the vanishing of $H^i(\P^2, \sW)$, $i=1,2$, implies that the boundary map $\del$ in the bottom row of \eqref{eqn:ComLocDiag} an isomorphism as well.  Thus the map
\[
\rho'_{0*}{\colon}H^0(\A^1_{k(\sqrt{a})},\sW)\to H^1(\P^2\setminus\{y\}, \sW)
\]
induced by $\rho_{0*}$ and $\phi_1$ is an isomorphism. 

Next, we look on the blow-up $\Bl_y\P^2\to \P^2$.  As we have seen above, the localization sequence \eqref{mult:LocSeq1} and the vanishing of $H^*(\P^1_{k(y)}, \sW(\sO(-1)))$ shows that the restriction map $H^1(\Bl_y\P^2, \sW)\to H^1(\P^2\setminus\{y\}, \sW)$ is an isomorphism. 

 The map $\rho{\colon}\P^1_{k(\sqrt{a})}\to \P^2$ lifts uniquely to a map $\tilde\rho{\colon}\P^1_{k(\sqrt{a})}\to \Bl_y\P^2$. Letting $\phi{\colon} \Bl_y\P^2\to \P^2$ be the blow-down map,  with exceptional divisor $\phi^{-1}(y)$ we have
 \[
 \omega_{\Bl_y\P^2}\cong \phi^*\omega_{\P^2}\otimes \sO_{\Bl_y\P^2}(\phi^{-1}(y)).
 \]
 Since $\omega_{\P^2}=\sO_{\P^2}(-3)$, $\tilde\rho^*\sO_{\P^2}(-3)=  \sO_{\P^1_{k(\sqrt{a})}}(-3)$ and $\tilde\rho^{-1}(\phi^{-1}(y))=\infty_{k(y)}$, we 
 have an isomorphism 
 \begin{equation}\label{eqn:phi2}
\phi_2{\colon}\tilde\rho^*(\omega_{\Bl_y\P^2})\xrightarrow{\sim} \sO_{\P^1_{k(\sqrt{a})}}(-3)\otimes \sO_{\P^1_{k(\sqrt{a})}}(\infty_{k(y)})\cong \sO_{\P^1_{k(\sqrt{a})}}(-2)\cong \omega_{\P^1_{k(\sqrt{a})}}, 
\end{equation}
inducing an isomorphism $\omega_{\tilde\rho}\cong \sO_{\P^1_{k(\sqrt{a})}}$, and giving us the push-forward map
\[
\tilde\rho_*{\colon}W(k(\sqrt{a}))=H^0(\P^1_{k(\sqrt{a})}, \sW)\to H^1(\Bl_y\P^2, \sW).
\]
Since the restriction of $\tilde\rho_*$  over $\P^2\setminus{y}\subset \Bl_y\P^2$ is our isomorphism $\rho'_{0*}$, we see that $\tilde\rho_*$ is an isomorphism as well, depending on the choice of  $\phi_2$. 

Finally, we map down to $\bar{X}_a$ by the blow-down map $\psi{\colon} \Bl_y\P^2=\Bl_p\bar{X}_a\to \bar{X}_a$. As the exceptional divisor $E=\P^1_p$ of $\psi$ is the proper transform of the line $\ell$, we have $\tilde\rho^{-1}(E)=\0$, so $\tilde\rho$ factors through the inclusion $\tilde{\jmath}{\colon}\Bl_y\P^2\setminus E\to \Bl_y\P^2$. The map  $p_q{\colon}\P^1_{k(\sqrt{a})}\to \bar{X}_a$ is just the map $\psi\circ\tilde\rho$; since $\psi$ is an isomorphism in a neighborhood of the image of $\tilde\rho$, the isomorphism $\omega_{\tilde\rho}\cong \sO_{\P^1_{k(\sqrt{a})}}$ defined above induces an isomorphism $\omega_{p_q}\cong \sO_{\P^1_{k(\sqrt{a})}}$, giving us  a well-defined push-forward map 
\[
p_{q*}{\colon}W(k(\sqrt{a}))=H^0(\P^1_{k(\sqrt{a})},\sW)\to H^1(\bar{X}_a, \sW).
\]
Using the localization sequence \eqref{mult:LocSeq1}, applied to $\Bl_p\bar{X}_a=\Bl_y\P^2$, the restriction map
\[
H^1(\Bl_y\P^2, \sW)=H^1(\Bl_p\bar{X}_a, \sW)\xrightarrow{\tilde{\jmath}^*} H^1(\Bl_p\bar{X}_a\setminus E, \sW)=
H^1(\bar{X}_a\setminus \{p\}, \sW)
\]
is an isomorphism. Let $j{\colon}\bar{X}_a\setminus\{p\}\to \bar{X}_a$ be the inclusion. Since $\tilde{\rho}_*$ is an isomorphism and $\tilde{\rho}$ factors through $\tilde{\jmath}$,  the fact that $\tilde{\jmath}^*$ is an isomorphism  implies that 
$p_{q*}$ defines an isomorphism after composing with $j^*{\colon}H^1(\bar{X}_a, \sW)\to 
H^1(\bar{X}_a\setminus\{p\}, \sW)$. The proof will be finished once we show that  $j^*$ an isomorphism.

For this, we have the exact sequence
\begin{multline*}
0\to H^1(\bar{X}_a, \sW)\xrightarrow{j^*}
H^1(\bar{X}_a\setminus\{p\}, \sW)\xrightarrow{\del}H^0(p, \sW)\xrightarrow{i_{p*}}
H^2(\bar{X}_a, \sW)\to \ldots
\end{multline*}
We have already shown in (3) that $i_{p*}$ is an isomorphism, so $j^*$ is an isomorphism, as desired.

\end{proof}

\begin{remark} \label{rem:SL2Inv} Suppose we have an $\SL_2$-action on some $Y\in \Sm/k$ and an $\SL_2$-linearized invertible sheaf $\sL$ on $Y$. Then the induced action of $\SL_2(k)$ on $H^*(Y, \sW(\sL))$ is trivial.  Indeed, the map $\SL_2\to \A^2\setminus\{0\}$ sending a matrix to its first column represents $\SL_2$ as an $\A^1$-bundle over $\A^2\setminus\{0\}$. For each $x\in \A^2\setminus\{0\}(k)$, we can connect $x$ to $(1,0)$ by a chain of  (at most) two maps $\A^1_k\to \A^2\setminus\{0\}$, thus we can connect each $g\in \SL_2(k)$ to the identity by a chain of at most  two maps $\A^1_k\to \SL_2$.

On the other hand, given a map $f{\colon}\A^1_k\to \SL_2$, we have the   map $F{\colon}\A^1\times  Y\to Y$ induced by the $\SL_2$-action on $Y$. Given an element $\alpha\in H^n(Y, \sW(\sL))$, the $\A^1$-homotopy invariance of $H^n(-, \sW(-))$ implies that there is an element $\beta\in H^n(Y, \sW(\sL))$ with 
\[
F^*(\alpha)=p_2^*(\beta)
\]
which then shows that $(F\circ i_0)^*(\alpha)=\beta=(F\circ i_1)^*(\alpha)$, where $i_0, i_1$ are the 0- and 1-sections, respectively, hence the $\SL_2$-action on $H^n(Y, \sW(\sL))$ is trivial, as claimed.

In particular, the isomorphisms described in Lemma~\ref{lem:WittCohbarXa} are all invariant under the $\SL_2$-action on $\bar{X}_a$.
\end{remark}

Next, we  describe the push-forward and pullback maps in Witt cohomology for the map $i{\colon}D\to \bar{X}_a$. Let $\iota{\colon}W(k)\to W(k(\sqrt{a}))$ be the base-change map and let $\Tr=\Tr_{k(\sqrt{a})/k}{\colon}W(k(\sqrt{a}))\to W(k)$ be the trace map for the finite separable extension $k\subset k(\sqrt{a})$.


\begin{lemma} \label{lem:iDComputation} With respect to our chosen orientations and the isomorphisms given by Lemma~\ref{lem:WittCohbarXa}, the following holds.\\[5pt]
1.  The map $W(k)\cong H^0(D, \sW)\xrightarrow{i_{D*}} H^1(\bar{X}_a,\sW)\cong W(k(\sqrt{a}))$ sends $\alpha$ to $-\iota(\alpha)$.\\[2pt]
2. The map $W(k)\cong H^1(D, \sW)\xrightarrow{i_{D*}} H^2(\bar{X}_a,\sW)\cong W(k)$ is the identity map.\\[2pt]
3. The map $W(k(\sqrt{a}))\cong H^1(\bar{X}_a,\sW)\xrightarrow{i_D^*} H^1(D, \sW)\cong W(k)$ send   $\alpha$ to $\<2\>\cdot\Tr_{k(\sqrt{a})/k}(\<\sqrt{a}\>\cdot\alpha)$.
\end{lemma}

\begin{proof} We use the Rost-Schmid complex to prove (1).  Letting $g$ be the rational function $t_2/t_3$ on $\bar{X}_a$, we use $\<g\>\cdot\alpha\in W(k(\bar{X}_a))$ on $\bar{X}_a$ to give the relation $\del_{RS} (\<g\>\alpha)=0$ in $H^1(\bar{X}_a,\sW)$. Explicitly, 
\[
\del_{RS}(\<g\>\alpha)=p_{q*}(\iota(\alpha)\otimes p_q^*\del/\del g \text{ on }\P^1_{k(\sqrt{a})})+(\alpha\otimes i_D^*\del/\del g\text{ on }D)
\]
Since $d(1/g)=-1/g^2\cdot dg$ and $\<-1\>=-1$ in $W(k)$, , this gives 
\[
p_{q*}[\iota(\alpha)\otimes p_q^*\del/\del g \text{ on }\P^1_{k(\sqrt{a})}]=[\alpha\otimes i_D^*\del/\del (1/g)\text{ on }D] \text{ in } 
H^1(\bar{X}_a,\sW).
\]

We note that $i_D^*(t_2)=x_1^2$, and $t_3=0$ on $D$, so 
\[
 i_D^*\del/\del (1/g)=  i_D^*\del/\del (t_3/t_2)= x_1^2i_D^*\del/\del t_3
 \]
Lemma~\ref{lem:OrientedSectionD} says that $i_{D*}(\alpha)$ is represented by $\<-2a\>\alpha\otimes  i_D^*\del/\del t_3$ on $D$, so 
 \[
 i_{D*}(\alpha)=\<-2a\>[\alpha\otimes i_D^*\del/\del (1/g)\text{ on }D]
 \]
 in $H^1(\bar{X}_a, \sW)$. 

By Lemma~\ref{lem:OrientedExtension}, $(1/2)p_q^*\del/\del g$ is an oriented section of $\omega_{p_q}$. Thus
\[
p_{q*}[\iota(\alpha)\otimes p_q^*(\del/\del g) \text{ on }\P^1_{k(\sqrt{a})}]=\<2\>\tilde{p}_{q*}\iota(\alpha)
\]
in $H^1(\bar{X}_a,\sW)$. Thus  
\[
i_{D*}(\alpha)=\<-a\>\tilde{p}_{q*}(\iota(\alpha)) =-\<a\>\cdot\tilde{p}_{q*}(\iota(\alpha)) 
\]
in $H^1(\bar{X}_a,\sW)\cong W(k(\sqrt{a}))$. Since  we use $\tilde{p}_{q*}{\colon} W(k(\sqrt{a}))=H^0(\P^1_{k(\sqrt{a})}, \sW)\to 
H^1(\bar{X}_a,\sW)$  to identify $H^1(\bar{X}_a,\sW)$ with $W(k(\sqrt{a}))$, and $\<a\>=1$ in $W(k(\sqrt{a}))$,   this completes the proof of (1).

For (2), we use our chosen orientations to define the push-forward maps. This gives us the commutative diagram
\[
\xymatrix{
H^0(p, \sW)\ar[r]^{i^D_{p*}}\ar[dr]^{{i^{Do}_{p*}}}\ar@/^20pt/[rr]^{i_{p*}}\ar@/_40pt/[drr]_{i^o_{p*}}&H^1(D, \sW(\omega_D))\ar[d]^{\tau_{o_D}}_\wr\ar[r]^{i_{D*}}&H^2(\bar{X}_a, \sW(\omega_{\bar{X}_a}))\ar[d]_\wr^{\tau_{o_{\bar{X}_a}}}\\
&H^1(D, \sW)\ar[r]^{i^o_{D*}}&H^2(\bar{X}_a, \sW)
}
\]
where $i^D_p{\colon}p\to D$, $i_p{\colon}p\to \bar{X}_a$ are the inclusions and we write $i^{Do}_{p*}$, $i^o_{p*}$ and $i_{D*}^o$ for the maps induced by the orientations. In particular, we have a commutative diagram of isomorphisms 
\[
\xymatrix{
H^0(p, \sW)\ar[r]^{{i^{Do}_{p*}}}\ar[dr]_{i^o_{p*}}&H^1(D, \sW)\ar[d]^{i^o_{D*}}\\
&H^2(\bar{X}_a, \sW)
}
\]
and we use $i^{Do}_{p*}$ and $i^o_{p*}$ to give our identifications $H^0(p, \sW)\cong H^1(D, \sW)$ and
$H^0(p, \sW)\cong H^2(\bar{X}_a, \sW)$ in Lemma~\ref{lem:WittCohbarXa}. This proves (2).  

For (3), consider the  cartesian diagram 
\[
\xymatrix{
0_{k(\sqrt{a})}\ar[r]^{i_0}\ar[d]_{p_{q0}}&\P^1_{k(\sqrt{a})}\ar[d]_{p_q}\\
D\ar[r]_{i_D}&\bar{X}_a
}
\]
defining the morphism $p_{q0}$ as the projection with image $q\in D(k)$. Since this cartesian square is transverse,  we have the canonical isomorphism $\beta{\colon}\omega_{p_{q0}}\xrightarrow{\sim}i_0^*(\omega_{p_q})$, constructed as follows. The diagram of $k(0_{k(\sqrt{a})}$ vector spaces
\[
\xymatrix{
0\ar[r]\ar[d]&i_0^*T_{\P^1_{k(\sqrt{a})}, 0_{k(\sqrt{a})}}\ar[d]^{dp_q}\\
p_{q0}^*T_{D, q}\ar[r]^-{di_D}&(p_q\circ i_0)^*T_{\bar{X}_a,q}
}
\]
is a cartesian square of injective linear maps. We then have the canonical isomorphisms
\[
i_0^*\omega_{p_q}=\coker(dp_q),\ \omega_{p_{q0}}=p_{q0}^*T_{D, q}
\]
and the composition
\[
\omega_{p_{q0}}=p_{q0}^*T_{D, q}\xymatrix{\ar@{^(->}[r]^{di_D}&}
(p_q\circ i_0)^*T_{\bar{X}_a,q}\twoheadrightarrow \coker(dp_q)=i_0^*\omega_{p_q}
\]
is our canonical isomorphism $\beta{\colon}\omega_{p_{q0}}\xrightarrow{\sim}i_0^*\omega_{p_q}$. Since $p_q$ has the explicit form (in the affine coordinates $t:=x_1/x_0$, $t_i=T_i/T_0$) 
\[
p_q(t)=(t+1, 0, \sqrt{a}t)=(t_1,t_2,t_3),
\]
we see that $\im(dp_q)$ is the $k(\sqrt{a})$-subspace of $(p_q\circ i_0)^*T_{\bar{X}_a,q}$ spanned by 
$(p_q\circ i_0)^*(\del/\del t_1+\sqrt{a}\del/\del t_3)$.

 Using $\beta$ to identify $i_0^*\omega_{p_q}$ with $\omega_{p_{q0}}$ gives us the commutative diagram 
\begin{equation}\label{eqn:PushPullDiag}
\xymatrix{
H^0(0_{k(\sqrt{a})}, \sW(\omega_{p_{q0}}))\ar[d]^{p_{q0*}}&\ar[l]_-{\beta^{-1}\circ i_0^*}\ar[d]^{p_{q*}}H^0(\P^1_{k(\sqrt{a})}, \sW(\omega_{p_q}))\\
H^1(D, \sW)&\ar[l]_-{i_D^*}H^{1}(\bar{X}_a, \sW)
}
\end{equation}
where the pushforward maps are the canonical ones. This follows from \cite[Proposition 10.7]{Feld1} and Remark~\ref{rem:CycModFunct}.

We factor $p_{q0}{\colon} 0_{k(\sqrt{a})}\to D$ as $\pi{\colon}0_{k(\sqrt{a})}\to q$ composed with the inclusion $i_q{\colon}q\to D$. This factors $p_{q0*}$ as
\[
H^0(0_{k(\sqrt{a})}, \sW(\omega_{p_{q0}}))\xrightarrow{\pi_*}H^0(q, \sW(\omega_{i_q}))\xrightarrow{i_{q*}}H^1(D, \sW)
\]
where $\pi_*{\colon}H^0(0_{k(\sqrt{a})}, \sW(\omega_{p_{q0}}))\to H^0(q, \sW(\omega_{i_q}))$ is the twist of $\pi_*{\colon}H^0(0_{k(\sqrt{a})}, \sW(\omega_{\pi}))\to H^0(q, \sW)$ using the canonical isomorphisms
\[
\omega_{p_{q0}}=p_{q0}^*\omega_D^{-1}\cong \pi^*i_q^*\omega_D^{-1}=\pi^*\omega_{i_q}
\]
and (noting that $\pi$ is \'etale)
\[
\omega_\pi\cong k(0_{k(\sqrt{a})}).
\]
We trivialize $\omega_{i_q}$ using our local oriented section $2\del/\del t_1$ of $\omega_D^{-1}=T_D$ near $q$, so we trivialize $\omega_{i_q}$ using  $2i_q^*\del/\del t_1$. This gives us our chosen isomorphism $\tilde{i}_{q*}{\colon}W(k)\xrightarrow{\sim}H^1(D, \sW)$ and rewrites $p_{q0*}$ as the composition
\[
H^0(0_{k(\sqrt{a})}, \sW(\omega_\pi))\xrightarrow{\pi_*}H^0(q, \sW)=W(k)\xymatrix{\ar[r]^-{\tilde{i}_{q*}}_-\sim&}H^1(D, \sW)
\]
Our canonical isomorphism $\omega_\pi\cong k(0_{k(\sqrt{a})})$ identifies $H^0(0_{k(\sqrt{a})}, \sW(\omega_\pi))$ with $W(k(\sqrt{a}))$ and rewrites $\pi_*$ as the trace map $\Tr_{k(\sqrt{a})/k}{\colon}W(k(\sqrt{a}))\to W(k)$. This gives the commutative diagram
\[
\xymatrix{
W(k(\sqrt{a}))\ar[d]^{\Tr_{k(\sqrt{a})/k}}\ar[r]^-{\tau^o}_-\sim&H^0(0_{k(\sqrt{a})}, \sW(\omega_{p_{q0}}))\ar[ddl]^{p_{q0*}}\\
W(k)\ar[d]^{\tilde{i}_{q*}}_\wr\\
H^1(D, \sW)
}
\]
where $\tau^o{\colon}W(k(\sqrt{a}))\to H^0(0_{k(\sqrt{a})}, \sW(\omega_{p_{q0}}))$ is the isomorphism induced by our choice of generating section $2p_{q0}^*\del/\del t_1$.

For the right-hand side of \eqref{eqn:PushPullDiag}, we  use our oriented section   $p_q^*((-1/4a)\del/\del t_3)$ for $\omega_{p_q}$ to give the isomorphism $\tau^{o'}{\colon}H^0(\P^1_{k(\sqrt{a})}, \sW)\xrightarrow{\sim}
H^0(\P^1_{k(\sqrt{a})}, \sW(\omega_{p_q}))$. We may also use $\beta^{-1}(i_{0_{k(\sqrt{a})}}^*(p_q^*((-1/4a)\del/\del t_3))$ to trivialize $\omega_{p_{q0}}$, instead of our previous choice $2p_{q0}^*\del/\del t_1$. This gives the commutative diagram
\[
\xymatrix{
H^0(0_{k(\sqrt{a})}, \sW)\ar[d]_\wr^{\tau^{o'}}\\
H^0(0_{k(\sqrt{a})}, \sW(\omega_{p_{q0}}))\ar[d]^{p_{q0*}}&\ar[l]_-{\tilde{i}_0^*}\ar[d]^{\tilde{p}_{q*}}H^0(\P^1_{k(\sqrt{a})}, \sW)\ar[ul]_{i_{0_{k(\sqrt{a})}}^*}\\
H^{1}(D, \sW)&\ar[l]_-{i_D^*}H^{1}(\bar{X}_a, \sW)
}
\]
where $\tilde{i}_0^*$ is the composition  
\[
H^0(\P^1_{k(\sqrt{a})}, \sW)\xrightarrow{\tau^{o'}}H^0(\P^1_{k(\sqrt{a})}, \sW(\omega_{p_q}))\xrightarrow{\beta^{-1}i_0^*}H^0(0_{k(\sqrt{a})}, \sW(\omega_{p_{q0}})).
\]
The composition 
\[
W(k(\sqrt{a}))\xrightarrow{\pi_{\P^1_{k(\sqrt{a})}}^*}H^0(\P^1_{k(\sqrt{a})}, \sW)\xrightarrow{\tilde{p}_{q*}}
H^{1}(\bar{X}_a, \sW)
\]
 is what we have used to define our chosen isomorphism $\alpha{\colon}W(k(\sqrt{a}))\xrightarrow{\sim}
H^{1}(\bar{X}_a, \sW)$. 

We thus have two generating sections of $\omega_{p_{q0}}$: $s_1:=p_{q0}^*(2\del/\del t_1)$ and $s_2:=\beta^{-1}(i_0^*p_q^*((-1/4a)\del/\del (t_2/t_3))$. Since $\omega_{p_{q0}}$ is a one-dimensional  vector space over $k(\sqrt{a})$, there is a unique $\lambda\in k(\sqrt{a})^\times$ with $s_2=\lambda\cdot s_1$. Thus, the composition of isomorphisms
\[
W(k(\sqrt{a}))\xymatrix{\ar[r]^{s_1^{-1}}_\sim&}H^0(0_{k(\sqrt{a})}, \sW(\omega_{p_{q0}}))
\xymatrix{\ar[r]^{s_2}_\sim&}W(k(\sqrt{a}))
\]
is the isomorphism $\times\<\lambda\>{\colon}W(k(\sqrt{a}))\xrightarrow{\sim} W(k(\sqrt{a}))$. 

Putting this together, we have the commutative diagram
\[
\xymatrix{
&W(k(\sqrt{a}))\ar[dl]^{\times\<\lambda\>}&\ar[l]_-\sim^-{i_0^*}H^0(\P^1_{k(\sqrt{a})}, \sW)\ar[dd]_\wr^{\tilde{p}_{q*}}&W(k(\sqrt{a}))\ar[l]^-{\pi_{\P^1_{k(\sqrt{a})}}^*}_-\sim\ar[ddl]^\alpha\ar@/_20pt/[ll]_{\id}\\
W(k(\sqrt{a}))\ar[d]^{\Tr_{k(\sqrt{a})}}\\
W(k)\ar[r]_-\sim^-{\tilde{i}_{q*}}&H^1(D, \sW)&H^1(\bar{X}_a, \sW)\ar[l]_{i_D^*}
}
\]
with $\tilde{i}_{q*}$ and $\tilde{p}_{q*}$ our chosen isomorphisms.  Thus, the composition $W(k(\sqrt{a}))\cong H^1(\bar{X}_a, \sW)\xrightarrow{i_D^*} H^1(D, \sW)=W(k)$  sends $\alpha\in W(k(\sqrt{a}))$ to  $\Tr_{k(\sqrt{a})/k}(\lambda\cdot \alpha)\in W(k)$. It remains to compute $\lambda$.

We give $(p_q\circ i_0)^*T_q(\bar{X}_a)$ the basis $(p_q\circ i_0)^*\del/\del t_1$, $(p_q\circ i_0)^*\del/\del t_3$. In this basis  $s_1=(2,0)$, $\im(dp_q)$ is the $k(\sqrt{a})$ subspace of $k(\sqrt{a})^2$ spanned by $(1, \sqrt{a})$, and $\beta(s_2)$ is the image of $(0, -1/4a)$ under the quotient map 
\[
k(\sqrt{a})^2\twoheadrightarrow k(\sqrt{a})^2/\<(1, \sqrt{a})\>
\]

Since $(2,0)-2(1,\sqrt{a})=(0,-2\sqrt{a})$ we see that $\beta(s_1)=8a\sqrt{a}\cdot \beta(s_2)$, so modulo squares in $k(\sqrt{a})$, we have $\lambda=2\sqrt{a}$. As $\Tr_{k(\sqrt{a})/k}$ is $W(k)$-linear, we have $\Tr_{k(\sqrt{a})/k}(\<2\sqrt{a}\>\alpha)= \<2\>\cdot\Tr_{k(\sqrt{a})/k}(\<\sqrt{a}\>\alpha)$, proving (2).
\end{proof}

We let 
\[
\Tr^{\<\sqrt{a}\>}{\colon}W(k(\sqrt{a}))\to W(k)
\]
denote the modified trace map $\alpha\mapsto \Tr_{k(\sqrt{a})/k}(\<\sqrt{a}\>\cdot \alpha)$ and write $\Tr$ for the usual trace map $\Tr_{k(\sqrt{a})/k}$.

We now use the Leray spectral sequence to compute the $\SL_2$-equivariant cohomology of $\bar{X}_a$.  Recall from Theorem~\ref{thm:Presentation} the computation $H^*(\BSL_2,\sW)=W(k)[e]$, where $e$ is the $\SL_2$-equivariant Euler class $e_{\SL_2}(E_2)$ of the tautological rank two vector bundle $E_2\to \BSL_2$.

\begin{proposition}\label{prop:EquivWittCohBarXa} Letting  $\pi_{\bar{X}_a}{\colon}\bar{X}_a\to \Spec k$ be the structure map,  $H^*_{\SL_2}(\bar{X}_a,\sW)$ has the following description as $H^*(\BSL_2,\sW)$-module.\\[5pt]
1. The  map $W(k)=H^0(\BSL_2,\sW)\xrightarrow{\pi_{\bar{X}_a}^*} H^0_{\SL_2}(\bar{X}_a,\sW)$ is an isomorphism.\\[2pt]
2. $H^{2}_{\SL_2}(\bar{X}_a,\sW)$ consists of two cyclic $W(k)$-summands, with respective generators the equivariant Euler classes $\pi_{\bar{X}_a}^*e$ and $e_{\SL_2}(T_{\bar{X}_a})$. They both have annihilator  $\Tr(W(k(\sqrt{a}))\subset W(k)$, giving the isomorphism
\[
H^{2}_{\SL_2}(\bar{X}_a,\sW)\cong (W(k)/\im(\Tr))\cdot \pi_{\bar{X}_a}^*e\oplus (W(k)/\im(\Tr))\cdot e_{\SL_2}(T_{\bar{X}_a}).
\]
In particular,  the natural map
\[
H^2(\BSL_2, \sW)=W(k)\cdot e\xrightarrow{\pi_{\bar{X}_a}^*} H^{2}_{\SL_2}(\bar{X}_a,\sW)
\]
has kernel   $\im(\Tr)\cdot e\subset W(k)\cdot e$. \\[2pt]
3. $H^1_{\SL_2}(\bar{X}_a,\sW)\cong \im(\iota{\colon}W(k)\to W(k(\sqrt{a}))$.\\[2pt]
4. $H^{2n+1}_{\SL_2}(\bar{X}_a,\sW)=0$ for $n\ge1$ and the multiplication map
\[
\times e^n{\colon}H^{2}_{\SL_2}(\bar{X}_a,\sW)\to H^{2n+2}_{\SL_2}(\bar{X}_a,\sW)
\]
is an isomorphism for all $n\ge1$. \\[5pt]
In addition, the kernel of $\Tr^{\<\sqrt{a}\>}$ is equal to the image of the base-change map $\iota{\colon}W(k)\to W(k(\sqrt{a}))$, and the kernel of $\iota$ is the ideal $(1-\<a\>)\subset W(k)$.   The image of $\Tr$ is the kernel of the map $W(k)\to W(k)$ given by multiplication by $1-\<a\>$.
\end{proposition}

\begin{proof} 
Recall from \S\ref{sec:EquivBMWitt} that we use the schemes $U_{2,n}$ of $2\times n$ matrices of rank $2$ to define the $\SL_2$-equivariant cohomology.  Let $\pi_n{\colon} \SL_2\backslash\bar{X}_a\times U_{2,n}\to 
\SL_2\backslash U_{2,n}$ be the projection. This is a Zariski locally trivial fiber bundle with fiber $\bar{X}_a$, giving for each $n$   the spectral sequence
\[
E_{2,n}^{p,q}=H^p(\SL_2\backslash U_{2,n}, R^q\pi_{n*}\sW)\Rightarrow H^{p+q}(\SL_2\backslash\bar{X}_a\times U_{2,n},\sW).
\]
Taking $n\gg 0$ gives us  the spectral sequence
\[
E_2^{p,q}=H^p(\BSL_2, R^q\pi_{{\bar{X}_a}*}\sW)\Rightarrow H^{p+q}_{\SL_2}(\bar{X}_a,\sW).
\]
This is justified by the results in \S\ref{sec:EquivBMWitt}, especially Lemma~\ref{lem:StrBoundImpBound} and Lemma~\ref{lem:HtpyModBound},  to show that the transition maps from $U_{2,n+1}$ to $U_{2,n}$ induce isomorphism
 \[
 H^p_{\SL_2}(U_{2,n+1}, R^q\pi_{n+1*}\sW)\to H^p_{\SL_2}(U_{2,n}, R^q\pi_{n*}\sW)
 \]
 and
 \[
 H^{p+q}(\SL_2\backslash\bar{X}_a\times U_{2,n+1},\sW)\to
 H^{p+q}(\SL_2\backslash\bar{X}_a\times U_{2,n},\sW).
\]
 for $n\gg 0$ (depending on $p$ and $p+q$).  
  
 Let $p{\colon}U_{2,n}\times_kk(\sqrt{a})\to U_{2,n}$ be the projection and  let $\sW_{k(\sqrt{a})}=p_*\sW$; we use the same notation for the corresponding sheaves on $\BSL_2$. With our choices, the isomorphisms (Lemma~\ref{lem:WittCohbarXa})
\[
H^0(\bar{X}_a, \sW)\cong W(k),\ H^1(\bar{X}_a, \sW)\cong W(k(\sqrt{a})), H^2(\bar{X}_a,\sW)\cong W(k)
\]
are all $\SL_2$-invariant (Remark~\ref{rem:SL2Inv}), so we have canonical   isomorphisms
\[
R^0\pi_{{\bar{X}_a}*}(\sW)\cong \sW\cong R^2\pi_{{\bar{X}_a}*}(\sW),\ R^1\pi_{{\bar{X}_a}*}(\sW)\cong  \sW_{k(\sqrt{a})}
\]
This gives the isomorphism
\[
E_2^{p,q}=\begin{cases} H^p(\BSL_2, \sW)&\text{ for }q=0,2\\ H^p(\BSL_2, \sW_{k(\sqrt{a})})&\text{ for }q=1\\0&\text{ otherwise.}\end{cases}
\]

Recall from Theorem~\ref{thm:Presentation} that $H^p(\BSL_2,\sW)=W(k)e^n$ for $p=2n$, and is 0 for $p$ odd. Also, by Gersten's conjecture for the Witt sheaves (more generally for homotopy modules, see e.g \cite[Theorem 5.41]{MorelA1}), we have $R^np_*\sW=0$ for $n>0$, so we have the isomorphism $H^p(\BSL_2, \sW_{k(\sqrt{a})})=H^p(\BSL_2\times_kk(\sqrt{a}), \sW)$. This gives the isomorphism
\[
E_2^{p,q}(\bar{X}_a)=\begin{cases} W(k)\cdot e^n&\text{ for }q=0,2, p=2n\\ W(k(\sqrt{a}))\cdot e^n&\text{ for }q=1, p=2n,\\0&\text{ otherwise.}\end{cases}
\]
In particular, the sequence degenerates at $E_3$. Using the multiplicative structure on the spectral sequence, we reduce to computing the differentials 
\[
d_2^{0,2}(\bar{X}_a){\colon}W(k)=H^0(\BSL_2, H^2(\bar{X}_a,\sW))\to H^2(\BSL_2, H^1(\bar{X}_a, \sW))=W(k(\sqrt{a}))
\]
and
\[
d_2^{0,1}(\bar{X}_a){\colon}W(k(\sqrt{a}))=H^0(\BSL_2, H^1(\bar{X}_a,\sW))\to H^2(\BSL_2, H^0(\bar{X}_a, \sW))=W(k).
\]

To facilitate the computation of the differentials, we compare with the Leray spectral sequence for the $\SL_2$-equivariant cohomology of $\P(F)$:
\[
E_2^{p,q}(\P(F))=H^p(\BSL_2, R^q\pi_*(\sW))\Rightarrow H^{p+q}_{\SL_2}(\P(F),\sW)
\]
where now $\pi$ is the family of projections $\P(F)\times U_{2,n}\to U_{2,n}$. Recall \cite[Lemma  5]{LevineBG} that 
\[
H^n_{\SL_2}(\P(F),\sW)=\begin{cases} 0&\text{ for }n>0,\\ W(k)&\text{ for }n=0.
\end{cases}
\]
 Computing as above, we find
\[
E_2^{p,q}(\P(F))=\begin{cases} W(k)\cdot e^n&\text{ for }q=0,1,  p=2n\\0&\text{ else.}\end{cases}
\]
The fact that $e$ maps to zero in $H^2_{\SL_2}(\P(F), \sW)$ implies that the differential 
\[
d_2^{0,1}(\P(F)){\colon}W(k)=H^0(\BSL_2, H^1(\P(F), \sW))\to H^2(\BSL_2, H^0(\P(F), \sW))=W(k)\cdot e
\]
is a $W(k)$-linear isomorphism, hence is given by multiplication by a unit in $W(k)$. 

We have the $\SL_2$-equivariant inclusion $i_D{\colon}D=\P(F)\to \bar{X}_a$ and the $\SL_2$-equivariant orientation of $\SL_2$-linearized bundles $\omega_{i_D}\cong \sO_{\P(F)}(-1)^{\otimes 2}$. This induces a map of spectral sequences, 
\[
i_{D*}{\colon}E_*^{*,*}(\P(F))\to E_*^{*,*+1}(\bar{X}_a)
\]
with the map on the $E_2$-terms  induced from  
$i_{D*}{\colon}H^q(\P(F), \sW)\to H^{q+1}(\bar{X}_a,\sW)$. Since $i_{D*}{\colon}H^1(D, \sW)\to H^2(\bar{X}_a, \sW)$ is an isomorphism (Lemma~\ref{lem:WittCohbarXa}), and 
$i_{D*}{\colon}H^0(D, \sW)\to H^1(\bar{X}_a, \sW)$  is the base-change map $\iota{\colon}W(k)\to W(k(\sqrt{a})$ (see Lemma~\ref{lem:iDComputation}), the differential 
$d_2^{0,1}(\P(F))$ determines $d_2^{0,2}(\bar{X}_a)$ explicitly by
\[
d_2^{0,2}(\bar{X}_a)(\alpha)=\iota(d_2^{0,1}(D)(\bar\alpha))\in W(k(\sqrt{a}))
\]
for $\alpha\in W(k)=H^0(\BSL_2, H^2(\bar{X}_a,\sW))$, where $\bar\alpha=i_{D*}^{-1}(\alpha)$. 
As $d_2^{0,1}(\P(F))$ is an isomorphism, this gives
\[
\ker d^{0,2}_2(\bar{X}_a)=\ker \iota,\ \im\, d^{0,2}_2(\bar{X}_a)=\im\, \iota
\]

Similarly, $i_D^*$ induces a map of spectral sequences
\[
i_D^*{\colon}E_*^{*,*}(\bar{X}_a) \to E_*^{*,*}(\P(F)).
\]
Since $i_D^*{\colon}W(k)=H^0(\bar{X}_a, \sW)\to H^0(D,\sW)=W(k)$ is minus the identity map and $i_D^*{\colon}W(k(\sqrt{a}))=H^1(\bar{X}_a, \sW)\to H^1(D, \sW)=W(k)$ is the map $\<2\>\cdot\Tr^{\<\sqrt{a}\>}$ (Lemma~\ref{lem:iDComputation}), we have
\[
d^{0,1}_2(\bar{X}_a)(\alpha)=d_2^{0,1}(\P(F))(-\<2\>\Tr^{\<\sqrt{a}\>}(\alpha))\in W(k)
\]
for $\alpha\in W(k(\sqrt{a}))=H^0(\BSL_2, H^1(\bar{X}_a,\sW))$. Thus, as $-\<2\>$ is a unit, 
\[
\ker d_2^{0,1}(\bar{X}_a)=\ker \Tr^{\<\sqrt{a}\>},\ \im\, d_2^{0,1}(\bar{X}_a)=\im\, \Tr^{\<\sqrt{a}\>};
\]
the second identity following since $d_2^{0,1}(\P(F))$ is multiplication by a unit in $W(k)$ and $\Tr^{\<\sqrt{a}\>}$ is $W(k)$-linear. Finally, we note that $\im(\Tr^{\<\sqrt{a}\>})=\im(\Tr)$, since $\<\sqrt{a}\>$ is a unit in $W(k(\sqrt{a}))$.

Multiplying by $e^n$ gives us corresponding expressions for the kernel and image of $d^{2n, q}$, $q=1,2$, $n\ge1$.  We conclude by invoking  Lam's exact  triangle
\[
\xymatrix{
W(k)\ar[rr]^{\iota}&&W(k(\sqrt{a}))\ar[dl]^{\Tr^{\<\sqrt{a}\>}}\\
&W(k)\ar[ul]^{\times(1-\<a\>)},
}
\]
which is exact at all three points  (see e.g. \cite[34.12, Prop. 34.1]{EKM}), 
and recalling that   $\im \Tr^{\<\sqrt{a}\>}=\im\Tr$. 
\end{proof}

We now use our computation of the equivariant cohomology 
\[
H^*_{\SL_2}(\P(\Sym^2F)\setminus \P(F),\sW),\ H^*_{\SL_2}(\P(\Sym^2F)\setminus \P(F),\sW(\omega_{\P(\Sym^2F)\setminus \P(F)}))
\]
 given by Theorem~\ref{thm:CohBN} below; the computations in \S\ref{sec:SpectSeqEulerClass} are independent of this section, avoiding a possible circularity. Recall the projection $\bar{\pi}_a{\colon}\bar{X}_a\to \P(\Sym^2F)=\P^2$ defined at the beginning of this section,  let $j{\colon}X_a\to \bar{X}_a$ be the inclusion and   let $\pi_a{\colon}X_a\to \P(\Sym^2F)\setminus \P(F)$ be the projection induced by $\bar{\pi}_a$. Since $\pi_a$ is \'etale, we have the canonical isomorphism $\omega_{X_a}\cong \pi_a^*\omega_{\P(\Sym^2F)\setminus \P(F)}$
 
 We have the $\SL_2$-invariant element $\<Q\>\in H^0(\P(\Sym^2F)\setminus\P(F), \sW)$ defined by the polynomial $Q:=T_1^2-4T_0T_2$, giving the class $[\<Q\>]\in H^0_{\SL_2}(\P(\Sym^2F)\setminus\P(F), \sW)$, and the 
 $\SL_2$-invariant, $\omega_{\P(\Sym^2F)\setminus\P(F)}$-valued quadratic form $\tilde{Q}$ \eqref{eqn:TildeQ},  giving the class $[\tilde{Q}]\in H^0_{\SL_2}(\P(\Sym^2F)\setminus\P(F), \sW(\omega_{\P(\Sym^2F)\setminus \P(F)}))$.
 
\begin{Not}\label{not:QuadFormNotation}
We set $[\<Q\>]_a:=\pi_a^*[\<Q\>]\in H^0_{\SL_2}(X_a, \sW)$, and 
$[\tilde{Q}]_a:=\pi_a^*[\tilde{Q}]\in H^0_{\SL_2}(X_a, \sW(\omega_{X_a}))$. 
\end{Not}

\begin{lemma} Orienting $X_a$ by the restriction of our chosen orientation for $\bar{X}_a$, the section $2p^*\Omega_2/T_3^3$ is an $\SL_2$-invariant, global oriented generator for $\omega_{X_a}$, giving an $\SL_2$-equivariant orientation $\omega_{X_a}\cong \sO_{X_a}$. 
\end{lemma}

\begin{proof} Since $\Omega_2$ is an $\SL_2$-invariant global generator for $\omega_{\P(\Sym^2F)}(3)$ and $\pi_a$ is \'etale, $\pi_a^*\Omega_2$ is an $\SL_2$-invariant global generator for $\omega_{X_a}(3)$. Since $T_3$ is an invariant global generator for $\sO_{X_a}(1)$, it follows that $2\pi_a^*\Omega_2/T_3^3$ is an $\SL_2$-invariant global generator for $\omega_{X_a}$; we now check the orientation condition. It suffices to check on some Zariski dense open subscheme of $X_a$.

Set $t_i=T_i/T_0$. Since $\Omega_2/T_0^3=dt_1dt_2$, we have
\[
 \frac{2\pi_a^*\Omega_2}{T_3^3} =\frac{2}{t_3^3}dt_1dt_2.
\]
Using the relation $(-1/4a)dt_1dt_3=(1/2t_3)dt_1dt_2$ on $X_a$, we see that 
\[
\frac{2p^*\Omega_2}{T_3^3}=(4/t_3^2)(1/2t_3)dt_1dt_2=(4/t_3^2)(-1/4a)dt_1dt_3
\]
on $X_a$, so (Lemma~\ref{lem:OriGen}) $2\pi_a^*\Omega_2/T_3^3$ is an oriented global generator of $\omega_{X_a}$. 
\end{proof}
Using $2\pi_a^*\Omega_2/T_3^3$ as global oriented generator, we henceforth consider $[\tilde{Q}]_a$ as an element of $H^0_{\SL_2}(X_a, \sW)$.
  
\begin{proposition}\label{prop:EquivWittCohXaComp}  1. $j^*{\colon}H^p_{\SL_2}(\bar{X}_a, \sW)\to H^p_{\SL_2}(X_a, \sW)$ is an isomorphism for $p\ge2$.\\[2pt]
2. $H^1_{\SL_2}(X_a, \sW)=0$\\[2pt]
3. The localization sequence for $j$ gives rise to an  exact sequence
\[
0\to H^0_{\SL_2}(\bar{X}_a,\sW)\xrightarrow{j^*} H^0_{\SL_2}(X_a,\sW)\xrightarrow{\del} 
H^0_{\SL_2}(D, \sW)\xrightarrow{i_{D*}}
H^1_{\SL_2}(\bar{X}_a, \sW)\to 0.
\]
4. Via our chosen isomorphism $W(k)\cong H^0_{\SL_2}(D, \sW)$, the image of $\del$ is the ideal $(1-\<a\>)W(k)$ and $[\tilde{Q}]_a\in H^0_{\SL_2}(X_a,\sW)$ generates a $W(k)$-module summand of $H^0_{\SL_2}(X_a,\sW)$ which determines a $W(k)$-splitting of  the surjection $\del{\colon} H^0_{\SL_2}(X_a,\sW)\to (1-\<a\>)W(k)$.  \\[2pt]
5.   $[\<Q\>]_a=\<a\>$ and $[\<Q\>]_a\cdot [\tilde{Q}]_a=-[\tilde{Q}]_a$.
\end{proposition}

\begin{proof} (1), (2) and the exact sequence of (3) except for the surjectivity of $i_{D*}$  follow from the localization sequence (Proposition~\ref{prop:BMLocalization})
\[
\ldots\to H^{p-1}_{\SL_2}(D, \sW)\xrightarrow{i_{D*}} H^p_{\SL_2}(\bar{X}_a, \sW)\xrightarrow{j^*}
H^p_{\SL_2}(X_a, \sW)\xrightarrow{\del}H^p_{\SL_2}(D, \sW)\to\ldots
\]
and the computation \cite[Lemma  5]{LevineBG}
\[
H^p_{\SL_2}(D, \sW)=H^p_{\SL_2}(\P(F), \sW)=\begin{cases}0&\text{ for }p>0\\W(k)&\text{ for }p=0.\end{cases}
\]
The surjectivity of $i_{D*}$ follows from Lemma~\ref{lem:iDComputation}(1) and Lemma~\ref{lem:WittCohbarXa}(3). These also imply that, via our identification of $H^0_{\SL_2}(D, \sW)$ with $W(k)$,  the image of $\del$ is the kernel of the base-extension map $\iota{\colon}W(k)\to W(k(\sqrt{a}))$, which by Lam's exact triangle is the ideal $(\<1\>-\<a\>)W(k)$. 

To complete the proof of (4), it thus suffices to show\\[5pt]
i. $\del([\tilde{Q}]_a)=u\cdot (\<1\>-\<a\>)$ for some unit $u\in W(k)$\\[2pt]
ii. For $\alpha\in W(k)$, if $\alpha\cdot (\<1\>-\<a\>)=0$ in $W(k)$, then $\alpha\cdot [\tilde{Q}]_a=0$.
\\[5pt]
For (i), recall the definition \eqref{eqn:TildeQ} of $\tilde{Q}$ as the $\omega_{\P(\Sym^2F)}$-valued quadratic form 
\[
\tilde{Q}(x_0e_0+x_1e_1):=2(T_0x_0^2+T_1x_0x_1+T_2x_1^2)\otimes \Omega_2
\]
where the $x_i$  are sections of $\sO_{\P(\Sym^2F)}(-2)$ and $e_0, e_1$ are the standard basis vectors in $F=k^2$ (viewed as right $\SL_2$-representation). If we fix a linear form $L=L(T_0, T_1, T_2)$, we can trivialize $\sO_{\P(\Sym^2F)}(-2)$ on the open subset $U_L$ defined by $L\neq0$, and then we have the induced form $\tilde{Q}_L$ on $\sO_{U_L}^2$ defined by
\[
\tilde{Q}_L(x_0,x_1):=2(t_0x_0^2+t_1x_0x_1+t_2x_1^2)\otimes \Omega_2/L^3
\]
where $t_i=T_i/L$.

Taking $L=T_0$, we  diagonalize $\tilde{Q}_{T_0}$ to give the isometric form
\[
\tilde{q}=[\<1\>+\<-Q/T_0^2\>]\otimes \frac{2\Omega_2}{T_0^3};\quad Q=T_1^2-4T_0T_2.
\]

We pull back to $X_a$, use $2\pi_a^*\Omega_2/T_3^3$ to trivialize $\omega_{X_a}$ and note the relation $aQ/T_0^2=(T_3/T_0)^2$ on $X_a$, giving
\[
[\pi_a^*\tilde{q}]=[(\<1\>+\<-a\>)\<t_3\>]\in H^0(X_a, \sW),\ t_i:=T_i/T_0.
\]
For our localization sequence, recall we have the oriented section $(-2a)\del/\del t_3$ for the normal bundle of $D$ in $\bar{X}_a\setminus\{T_0=0\}$ (Lemma~\ref{lem:OrientedSectionD}). Thus
\[
\del_D([\tilde{Q}]_a)=\<-2a\>\del_{t_3}((\<1\>+\<-a\>)\<t_3\>)=\<-2a\>(\<1\>+\<-a\>)=\<2\>(1-\<a\>).
\]
which proves (i). 

If $\alpha\cdot (\<a\>-\<1\>)=0$, then  by our computation above, we have
\[
\alpha\cdot [\tilde{Q}]_a=[\alpha\cdot(\<1\>-\<a\>)\<t_3\>]=0,
\]
proving (ii).

Finally, (5) follows from the relations on $X_a$
\[
(T_1^2-4T_0T_2)/T_3^2=a^{-1},\ \<a\>\cdot (\<1\>-\<a\>)\<t_3\>= (\<a\>-\<1\>)\<t_3\>.
\]
\end{proof}

\begin{definition} For $a\in k^\times$, let $I_a\subset W(k)$ denote the kernel of the multiplication map $W(k)\to W(k)$, $x\mapsto x\cdot(1-\<a\>)$. 
\end{definition}
By Lam's exact triangle, $I_a$ is also the image of the trace map 
\[
\Tr_{k(\sqrt{a})/k}{\colon}W(k(\sqrt{a}))\to W(k). 
\]

\begin{proposition}\label{XaWittCoh} We retain the Notations~\ref{not:QuadFormNotation}. Let $\bar{W}(k):=W(k)/I_a$. Then as $H^*(\BSL_2,\sW)=W(k)[e]$-algebra, 
\[
H^*_{\SL_2}(X_a,\sW)\cong W(k)[e, y]/(y^2-2(\<1\>-\<a\>), I_a\cdot y, I_a\cdot e),
\]
by the map sending $y$ to $[\tilde{Q}]_a$. Moreover, $[\tilde{Q}]_a\cdot e=e_{\SL_2}(T_{X_a})$, $[\<Q\>]_a=\<a\>$, $[\<Q\>]_a\cdot [\tilde{Q}]_a=-[\tilde{Q}]_a$, $[\tilde{Q}]_a^2=2\cdot(\<1\>- \<a\>)$ and 
$[\tilde{Q}]_a^2\cdot e=4e$.  
\end{proposition}

\begin{proof} Write $Y$ for $\P(\Sym^2F)\setminus \P(F)$. Since $\pi_a$ is \'etale, we have $\pi_a^*T_Y\cong T_{X_a}$ and hence $e_{\SL_2}(T_{X_a})=\pi_a^*(e_{\SL_2}(T_Y))$. With Lemma~\ref{lem:EdgeHom}(3) and Proposition~\ref{prop:MultiplicativeStructure}(1), this gives the relations 
\[
[\tilde{Q}]_a^2=2\cdot(\<1\>- [\<Q\>]_a),\ [\tilde{Q}]_a\cdot e=e_{\SL_2}(T_{X_a}).
\]
Moreover, since $X_a$ is the affine variety defined by $T_3^2=aQ$, $T_3\neq0$, we have
\[
[\<Q\>]_a=\<a\>\in H^0(X_a, \sW)
\]
which gives us the relation $[\tilde{Q}]_a^2=2\cdot(\<1\>- \<a\>)$. 

By Proposition~\ref{prop:EquivWittCohXaComp} and Proposition~\ref{prop:EquivWittCohBarXa}, the degree $\ge 2$ $W(k)$-submodule of $H^*_{\SL_2}(X_a,\sW)$ is the free 
$(W(k)/I_a)[e]$-module with basis $e, e_{\SL_2}(T_{X_a})$. The degree 0 $W(k)$-submodule of $H^*_{\SL_2}(X_a,\sW)$ is $W(k)\cdot 1\oplus W(k)\cdot [\tilde{Q}]_a$, with $W(k)\cdot [\tilde{Q}]_a\cong (1-\<a\>)W(k)$. By the above relations,   $H^0_{\SL_2}(X_a,\sW)$ is thus isomorphic as $W(k)$-algebra to $W(k)[y]/(I_a\cdot y, y^2-2(\<1\>-\<a\>)$. Thus the $W(k)[e]$-algebra
\[
W(k)[e, y]/(I_a\cdot y, y^2-2(\<1\>-\<a\>),I_a\cdot e)
\]
is isomorphic to $H^*_{\SL_2}(X_a,\sW)$ by sending $y$ to $[\tilde{Q}]_a$.

Since $\<1\>+\<a\>= \Tr_{k(\sqrt{a}/k}(\<2\>)$,  and $I_a=\im\Tr_{k(\sqrt{a})/k}$, we have $\<a\>e=-e$, so 
\[
[\tilde{Q}]_a^2e=2\cdot(\<1\>- \<a\>)\cdot e=4e.
\]

The remaining relations follow from Proposition~\ref{prop:EquivWittCohXaComp}.
\end{proof}

We achieve our main goal of the section, computing the $N$-equivariant Witt cohomology in the case of a non-trivial action on an integral, 0-dimensional scheme, with the following result.

\begin{theorem}\label{thm:Nontrivial0DimlCase} Let $z$ be an integral finite $k$-scheme with $k(z)$ a degree two extension of $k$ and with $N$ acting on $z$ by letting $\sigma$ act as conjugation over $k$.  Then as an $H^*(BN,\sW)=W(k)[e,x]/((1+x)e, x^2-1)$-algebra, we have
\[
H^*_{N}(z,\sW)\cong W(k)[e, x][y]/(x-\<a\>,y^2-2(\<1\>-\<a\>),I_a\cdot y, I_a\cdot e).
\]
After inverting $e$, we have 
\[
H^*_{N}(z,\sW)[e^{-1}]=\bar{W}(k)[e, y, e^{-1}]/(y^2-2(\<1\>-\<a\>)),
\]
where $\bar{W}(k)=W(k)/I_a$, and the map $H^*_{N}(z,\sW)\to H^*_{N}(z,\sW)[1/e]$ sends $x$ to $-1$. 
\end{theorem}

\begin{proof} As a direct consequence of Proposition~\ref{XaWittCoh} and \eqref{eqn:ChangeGroup}, we have 
\[
H^*_{N}(z,\sW)\cong W(k)[e][y]/(y^2-2(\<1\>-\<a\>),I_a\cdot y, I_a\cdot e)
\]
as algebra over $H^*(\BSL_2, \sW)=W(k)[e]$. The $\SL_2$-invariant section $Q:=T_1^2-4T_0T_2$ of $\sO_{\P^2}(2)$ gives rise to the section $x\in H^0_{\SL_2}(\P^2\setminus D,\sW)\cong H^0(BN, \sW)$ used in the presentation of $H^0(BN, \sW)$ given in Theorem~\ref{thm:Presentation} (see \cite[Proposition 5.5]{LevineBG}). 
Since $X_a\to \P^2\setminus D$ is the double cover defined by $T_3^2=aQ$, the element $x$ maps to $\<a\>\in H^0_{\SL_2}(X_a, \sW)$ and we arrive at the desired presentation of 
$H^*_{N}(z,\sW)$ as $H^*(BN,\sW)$-algebra. The description of $H^*_{N}(z,\sW)[e^{-1}]$ follows directly from this and the presentation of $H^0(BN, \sW)[e^{-1}]$ given in Theorem~\ref{thm:Presentation}
\end{proof}

\begin{remark}\label{rem:Nontrivial0DimlCase} Take $b\in k^\times$ and let $a=b^2$,  $\sO=k[X]/X^2-a$, and let $z=\Spec\sO$. We let $N$ act on $z$ through $N/T_1=\<\bar\sigma\>$, with $\bar\sigma^*(X)=-X$. Then $\<a\>=1$, so the ideal $I_a:=\ker(\times(1-\<a\>){\colon}W(k)\to W(k)$ is the unit ideal. Moreover, we have
\begin{multline*}
H^*_N(z,\sW)=H^*_{T_1}(\Spec k, \sW)\cong W(k)[e]/(e)\\\cong
W(k)[e, x,y]/(x-\<a\>,y^2-2(\<1\>-\<a\>),I_a\cdot y, I_a\cdot e),
\end{multline*}
so the description of  $H^*_N(\Spec k[X]/(X^2-a),\sW)$ given by Theorem~\ref{thm:Nontrivial0DimlCase} in case of $a\in k^\times$, not a square, extends to the case of $a$ a square.
\end{remark}

We complete this section by computing the push-forward map $H^*_{N}(z,\sW)\to H^*(BN,\sW)$.

\begin{corollary} \label{cor:Nontrivial0DimlCasePushforward} Let $z\in \Sch^N/k$ be as in Theorem~\ref{thm:Nontrivial0DimlCase}. Let $\pi{\colon}z\to \Spec k$ be the structure map, where $N$ acts trivially on $\Spec k$. Then with respect to the presentation of $H^*_{N}(z,\sW)$ given in Theorem~\ref{thm:Nontrivial0DimlCase} and $H^*_N(\Spec k,\sW)=H^*(BN,\sW)$ given in Theorem~\ref{thm:Presentation}(2), the map
\begin{multline*}
W(k)[e, x][y]/(x-\<a\>,y^2-2(\<1\>-\<a\>), I_a\cdot y, I_a\cdot e)\\\xrightarrow{\pi_*}
W(k)[e, x]/(x^2-1, (1+x)e)
\end{multline*}
is the $W(k)[e]$-module map given by
\begin{align*}
\pi_*(1)&=\<2\>+\<2a\>x\\
\pi_*(e)&=(\<2\>-\<2a\>)\cdot e\\
\pi_*(y)&=0.
\end{align*}
After inverting $e$, this is the $W(k)[e]$-module map $\bar{W}(k)[e, y, e^{-1}]\to W(k)[e,e^{-1}]$ given by
\begin{align*}
\pi_*(1)&=<2>-<2a>\\
\pi_*(e)&=(<2>-<2a>)\cdot e\\
\pi_*(y)&=0
\end{align*}
\end{corollary}

\begin{proof} The fact that $\pi_*$ is a $W(k)[e]=H^*(\BSL_2, \sW)$-module map is the projection formula; this implies that $\pi_*$ is determined by its values on $1, e$ and $y$.

Inverting $e$ gives the relation $x=-1$, so the second set of formulas follows from the first set.

To avoid notational misunderstandings, we temporarily write $1', x', e', y'$ for the images in $H^*_{N}(z,\sW)$ of $1, x,e\in H^*(BN, \sW)$ and $y\in H^*(BN, \sW(\gamma))$ under the pullback maps $H^*(BN, \sW)\to H^*_{N}(z,\sW)$, $H^*(BN, \sW(\gamma))\to H^*_{N}(z,\sW)$

The map $\pi{\colon}z\to \Spec k$ corresponds to the $\SL_2$-equivariant projection $\pi_a{\colon}X_a\to \Sym^2(F)\setminus\Sym(F)$, which realizes $X_a$ as the \'etale double cover of  $\P^2\setminus D$ defined by $T_3^2=aQ$. This implies that $\pi_*\circ\pi^*$ is multiplication by the corresponding trace form, namely $\<2\>+\<2aQ\>$.  Since 
\[
1'=\pi^*(1),\ x'=\pi^*([\<Q\>]),\ e'=\pi^*(e), \ y'=\pi^*([\tilde{Q}]), 
\]
and $xe=-e$, $x=[\<Q\>]$ in $H^*(BN,\sW)$,  this gives the formulas for $\pi_*(1')$ and $\pi_*(e')$. 

To compute $\pi_*(y')$, we diagonalize the restriction of $\pi_a^*\tilde{Q}$ to $T_0T_3\neq0$,  as in proof of  Proposition~\ref{prop:EquivWittCohXaComp}, 
\[
\tilde{q}_a= (\<-1\>+\<a\>)\<t_3\>.
\]
Since $X_a\cap (T_0\neq0)$ is defined by $t_3^2=a(Q/T_0^2)$, we have $\pi_*(\<t_3\>)=\<1\>+\<-1\>=0\in W(k)$, so 
\[
\pi_*(\tilde{q})=  (\<-1\>+\<a\>)\cdot \pi_*(\<t_3\>)=0.
\]
\end{proof}

\section{Spectral sequences and Euler classes}\label{sec:SpectSeqEulerClass}

In this section we compute the $\SL_2$-equivariant cohomology of $\P(\Sym^2F)\setminus \P(F)$; since $N\backslash \SL_2$ is $\SL_2$-equivariantly isomorphic to $\P(\Sym^2F)\setminus \P(F)$, this also gives a computation of $H^*(BN, \sW)$ and  $H^*(BN, \sW(\gamma))$ ($\gamma$ the generator of $\Pic(BN)\cong \Z/2$).  We had previously computed $H^*(BN, \sW)$ and  $H^*(BN, \sW(\gamma))$ in \cite[Proposition 5.3]{LevineBG}, but our description there contained an error, in that we ignored a possible contribution from  $H^0(BN, \sW(\gamma))$ (which we will see here is non-zero). We computed in {\it loc. cit.}   the ring structure on $H^*(BN, \sW)$ and the $H^*(\BSL_2, \sW)$-module structure on $H^{*\ge2}(BN, \sW(\gamma))$, which we extend here to a computation of the  ring structure of $H^*(BN, \sW)\oplus H^*(BN, \sW(\gamma))$ (see Theorem~\ref{thm:CohBN} for the precise statement); we have used this in our computation of $H^*_{\SL_2}(X_a, \sW)$ given in the previous section.

As a first step,  we identify the Euler class of the universal rank two bundle with the value of a differential in a Leray spectral sequence. Before we do this, we pause to fix our conventions for the total complex of a double complex and the spectral sequence for a double complex, so that we can keep track of signs properly.

Recall that a {\em double complex} $(K^{**}, \del_1, \del_2)$ is a bigraded object $K^{**}=\oplus_{a,b}K^{a,b}$ in some additive category, with differentials $\del_1$ of bi-degree $(1,0)$ and $\del_2$ of bi-degree $(0,1)$ satisfying
\[
\del_1^2=0,\ \del_2^2=0,\ \del_1\circ\del_2=\del_2\circ\del_1.
\]
For simplicity, we will assume that for each $n\in \Z$, there are only finitely many indices $(a,b)$ with $a+b=n$ and $K^{a,b}\neq0$. The total complex $(\Tot^* K, d)$ is the complex with
\[
\Tot^nK:=\oplus_{a+b=n}K^{a,b}
\]
and differential $d^n{\colon}\Tot^nK\to \Tot^{n+1}K$ defined as 
\[
d^n=\sum_{a+b=n}(\del_1^a{\colon}K^{a,b}\to K^{a+1,b})+((-1)^a\del_2^b{\colon}K^{a,b}\to K^{a,b+1})
\]

$(\Tot^* K, d)$ comes with two descending filtrations by subcomplexes $F_I^p\Tot^* K$, 
$F_{II}^q\Tot^* K$ with
\[
F_I^p\Tot^nK=\oplus_{a+b=n, a\ge p}K^{a,b},\ F_{II}^q\Tot^n K=\oplus_{a+b=n, b\ge q}K^{a,b}
\]
giving rise to spectral sequences
\[
E^{a,b}_1(I):=H^{a+b}(\gr_{F_I}^a\Tot^* K)\Rightarrow H^{a+b}\Tot^*K
\]
\[
E^{a,b}_1(II):=H^{a+b}(\gr_{F_{II}}^a\Tot^* K)\Rightarrow H^{a+b}\Tot^*K
\]
These are strongly convergent if for example there are integers $a_0, b_0$ with $K^{a,b}=0$ if $a<a_0$ or if $b<b_0$. 

We will mainly be using the filtration $F_I$, and we give here an explicit description of the $d_1$ and $d_2$ differentials. The differential in the complex $\gr_{F_I}^a\Tot^* K$ is $(-1)^a\del_2$ and the map $d_1^{a,b}{\colon}H^{a+b}(\gr_{F_I}^a\Tot^* K)\to H^{a+b+1}(\gr_{F_I}^{a+1}\Tot^* K)$ is induced by $\del_1$. If we have an element $x\in \ker(d_1^{a,b})$, then lift $x$ to $\tilde{x}\in K^{a,b}$ with $\del_2(\tilde{x})=0$,  and there is a $y\in K^{a+1, b-1}$ with
\[
\del_1^{a,b}(\tilde{x})=(-1)^a\del_2^{a+1, b-1}(y).
\]
Then the image  of $\del_1^{a+1,b-1}(y)$ in $H^{a+b+1}(\gr_{F_I}^{a+2}\Tot^* K)$ is in $\ker(d_1^{a+2,b-1})$, giving
\[
d_2^{a,b}(x):=\overline{\del_1^{a+1,b-1}(y)}\in E_2^{a+2,b-1}(I).
\]

\begin{ex}\label{ex:CechDoubleComplex} Let $T$ be a topological space, $\sU=\{U_i, i\in I\}$ an open cover and $(\sC^*, d_{\sC})$ a complexes of presheaves (with values in some abelian category) on $T$. We have the associated double complex $(K^{**}, \del_1, \del_2)$, with $(K^{*,b}, \del_1)$ the \v{C}ech complex $\check{C}^*(\sU, \sC^b)$ and $\del_2$ induced by the differential $d_{\sC}$. Letting $\sH^b(\sC^*)$ be the cohomology presheaf, $E_1^{**}$ is the \v{C}ech complex
\[
E_1^{a,b}(I)=\check{C}^a(\sU, \sH^b(\sC^*))
\]
and $E_2^{**}$ is the \v{C}ech cohomology
\[
E_2^{a,b}(I)=\check{H}^a(\sU, \sH^b(\sC^*)).
\]
\end{ex}

\begin{construction}\label{const:Cocycle}
For $Y$ a smooth $k$-scheme, and $L$ a line bundle on $Y$,  we have the  Rost-Schmid complex with coefficients in the twisted Witt sheaf $\sW(L)$, $C_{RS}^*(Y, \sW(L))$, and the associated presheaf of complexes on $Y_\Nis$,  $\sC_{RS}^*(\sW(L))$, as described in \S\ref{subsec:RS}.  We use Notation~\ref{Not:RSNotation}.

Let $p{\colon}E\to X$ be a rank two vector bundle on some smooth $k$-scheme $X$. Suppose we have a trivializing Zariski open cover $\sU=\{U_i\mid i\in I\}$, with vector bundle isomorphisms $\psi_i{\colon}p^{-1}(U_i)\xrightarrow{\sim} U_i\times \A^2$, and let $\{\psi_{ij}:=\psi_{i|U_{ij}}\circ\psi^{-1}_{j|U_{ij}}\in \GL_2(U_{ij})\}_{ij}$ be the corresponding cocycle.

We will assume that $H^1(U_{ij}, \sW)=0$ for all $i\neq j$. If in addition we  assume there is a given isomorphism $\det E\cong \sO_X$, and that the $\psi_i$ respect this isomorphism, then  the $\psi_{ij}$ are sections of $\SL_2$ over $U_{ij}$.  

In general, we have an explicit representation of the Thom class of $E$ in the Rost-Schmid complex $C_{RS}^*(E,\sW(\det^{-1}E))$, as follows. Let $s{\colon}X\to E$ be the zero-section. Then the Thom class $\th(E)\in H^2_{s(X)}(E, \sW(\det^{-1}E))$ is represented by $\th_{RS}(E):=(\<1\>\text{ on }s(X))$ in  $C^2_{RS, s(X)}(E,\sW(\det^{-1}E))$.  Indeed,
\[
C^n_{RS, s(X)}(E, \sW(\det^{-1}E))=\oplus_{x\in E^{(n)}\cap s(X)}W(k(x), \Lambda^E_x\otimes \det^{-1}E))
\]
and for $x\in E^{(n)}\cap s(X)$, the fact that the normal bundle of $s(X)$ in $E$ is $E$ gives us  a canonical isomorphism $\Lambda^E_x\otimes \det^{-1}E\cong \Lambda^X_x$. This in turn gives the isomorphism of Rost-Schmid complexes
\[
s^{RS}_*{\colon}C^*_{RS}(X, \sW(\det^{-1}E))\to C^*_{RS, s(X)}(E, \sW(\det^{-1}E))[2]
\]
representing the Thom isomorphism $s_*{\colon}H^*(X, \sW)\to H^{*+2}_{s(X)}(E, \sW(\det^{-1}E))$. As the Thom class $\th(E)\in H^2_{s(X)}(E, \sW(\det^{-1}E))$ is by definition $s_*(\<1\>)$, we see that  
$\th(E)$ is represented by $s^{RS}_*(\<1\>\text{ on }X)=\<1\>\text{ on }s(X)$ in $C^2_{RS, s(X)}(E, \sW(\det^{-1}E))$, as claimed.  

It follows that $p^*(e(E))$ is represented by the image of $\th_{RS}(E)$ in the Rost-Schmid complex without support, that is,
\[
p^*(e(E))=[\<1\>\text{ on }s(X)]\in H^2(\sC^*_{RS}(E, \sW(\det^{-1}E)).
\]
This follows from the identity $e(E)=s^*\th(E)$, and the fact that $p^*s^*=\id$ on $H^*(E,\sW(\det^{-1}E))$. Passing from $e(E)$ to $p^*e(E)$ does not lose any information, since by the $\A^1$-homotopy invariance of $H^*(-, \sW(-))$, the map 
\[
p^*{\colon}H^*(X,\sW((\det^{-1}E))\to H^*(E, \sW(\det^{-1}E))
\]
is an isomorphism.

Next, given our affine cover $\sU=\{U_i\}$ of $X$, trivializing $E$ as above, we give a translation into the \v{C}ech double complex for the cover $p^{-1}\sU:=\{p^{-1}(U_i)\}$ of $E$ with coefficients in the Rost-Schmid complex for $\sW(\det^{-1}E)$. The isomorphism $\psi_i{\colon}p^{-1}(U_i)\xrightarrow{\sim}U_i\times\A^2$ gives us the trivialized sub-bundle $\A^1_i:=\psi^{-1}_i(\A^1_{U_i}\times 0)\subset p^{-1}(U_i)$. Let $p_1, p_2{\colon}\A^2\to \A^1=\Spec k[t]$ be the two projections, giving us the two coordinates $t_j=p_j^*t$, $j=1,2$ on $U_i\times\A^2$ and the corresponding coordinates $t_j^{(i)}:=\psi_i^*(t_j)$ on $p^{-1}(U_i)$. Let $\bar{t}^{(i)}_1$ be the restriction of $t_1^{(i)}$ to $\A^1_i$. The ideal defining $\A^1_i$ in $p^{-1}(U_i)$ is generated by $t^{(i)}_2$, giving us the generator $\del/\del t^{(i)}_2$ of  $\Lambda^{p^{-1}(U_i)}_{\eta}$, where $\eta$ is the generic point of $\A^1_i$. Similarly, we may consider $dt_1^{(i)}\wedge dt_2^{(i)}$ as a basis for $\det^{-1}E$ on $p^{-1}(U_i)$, dual to $\del/\del t_1^{(i)}\wedge \del/\del t_2^{(i)}$.

Putting this all together gives us the element
\[
\phi_i:= \<\bar{t}^{(i)}_1\>\otimes (\del/\del t^{(i)}_2\otimes dt_1^{(i)}\wedge dt_2^{(i)})\text{ on }\A^1_i \in
 C_{RS}^1(p^{-1}(U_i), \sW(\det^{-1}E)).
 \]
 We compute the boundary $\del_{RS}\phi_i\in  C_{RS}^2(p^{-1}(U_i), \sW(\det^{-1}E))$ as
 \begin{align*}
 \del_{RS}\phi_i&=(\del_{\bar{t}^{(i)}_1}\otimes \del/\del \bar{t}^{(i)}_1)( \<\bar{t}^{(i)}_1\>\otimes (\del/\del t^{(i)}_2\otimes dt_1^{(i)}\wedge dt_2^{(i)})\text{ on }\A^1_i)\\
 &=\<1\>\otimes  (\del/\del \bar{t}^{(i)}_1\wedge \del/\del t^{(i)}_2\otimes dt_1^{(i)}\wedge dt_2^{(i)})\text{ on }s(U_i)\\
 &=\<1\>\text{ on }s(U_i),
 \end{align*}
following our conventions outlined in Remark~\ref{rem:RSMultilinearAlgebra}.

If we apply the \v{C}ech coboundary $\check{\delta}^0$ to $\{\phi_i\}_i  \in \check{C}^0(p^{-1}(\sU), \sC_{RS}^1(\sW(\det^{-1}E))$, we get the
element
$\check{\delta}^0(\{\phi_i\}_i)\in \check{C}^1(p^{-1}(\sU), \sC_{RS}^1(\sW(\det^{-1}E))$ with 
\[
\del_{RS}(\check{\delta}^0(\{\phi_i\}_i))_{ij}=0
\] 
 in $C_{RS}^2(p^{-1}(U_{ij}), \sW(\det^{-1}E))$, for all $i,j$. The cohomology of the Rost-Schmid complex $C_{RS}^*(-, \sW)$ computes the $\sW$-cohomology on smooth $k$-schemes, so by our assumption $H^1(U_{ij}, \sW)=0$ and homotopy invariance of $H^*(-,\sW)$, each  $\check{\delta}^0(\{\phi_i\}_i)_{ij}$ is a Rost-Schmid coboundary. Thus, for each $i,j$, there is a   $\lambda_{ij}\in W(k(p^{-1}(U_{ij})), \det^{-1}E)$ with
\[
\del_{RS}(\{\lambda_{ij}\}_{ij})=\check{\delta}^0(\{\phi_i\}_i)_{ij}\in C_{RS}^1(p^{-1}(U_{ij}), \sW(\det^{-1}E)).
\]

Letting $\{\lambda_{ijk}\}_{ijk}:=\check{\delta}^1(\{\lambda_{ij})\}_{ij})$, we see that  each $\lambda_{ijk}$ has zero Rost-Schmid coboundary, and thus determines an element $\lambda_{ijk}\in  H^0(p^{-1}(U_{ijk}), \sW(\det^{-1}E))$, defining a \v{C}ech cocycle $\{\lambda_{ijk}\}_{ijk}\in \check{C}^2(p^{-1}(\sU), \sW(\det^{-1}E))$.  Finally, assuming we have chosen an isomorphism $\det E\cong \sO_X$, we consider $\{\lambda_{ijk}\}_{ijk}$ as a  \v{C}ech cocycle in  $\check{C}^2(p^{-1}(\sU), \sW)$.

To summarize, we have $\<1\>\text{ on }s(X)$ in $C^2_{RS}(E, \sW(\det^{-1}E)))$ mapping to  the cocycle $\{\<1\>\text{ on }s(U_i)\}_i$ in $\check{C}^0(p^{-1}\sU, \sC^2_{RS}(\sW(\det^{-1}E)))$ under the augmentation 
\[
\epsilon_1{\colon}C^*_{RS}(E, \sW(\det^{-1}E)))\to \check{C}^0(p^{-1}\sU, \sC^*_{RS}(\sW(\det^{-1}E))), 
\]
the cochain $\{\phi_i\}_i\in  \check{C}^0(p^{-1}\sU, \sC^1_{RS}(\sW(\det^{-1}E)))$  with $\del_{RS}(\{\phi_i\}_i)=\{\<1\>\text{ on }s(U_i)\}_i$, the cochain $\{\lambda_{ij}\}_{ij}\in \check{C}^1(p^{-1}\sU, \sC_{RS}^0(\sW(\det^{-1}E)))$ with $\del_{RS}(\{\lambda_{ij}\}_{ij})=\check{\delta}^0(\{\phi_i\}_i)$ in $\check{C}^1(p^{-1}\sU, \sC_{RS}^2(\sW(\det^{-1}E)))$ and the cocycle $\{\lambda_{ijk}\}_{ijk}\in \check{C}^2(p^{-1}\sU, \sW(\det^{-1}E))$ mapping to $\check{\delta}^1(\{\lambda_{ij}\}_{ij})$ in $\check{C}^2(p^{-1}\sU, \sC^0_{RS}(\det^{-1}E))$ under the augmentation 
\[
\epsilon_2{\colon}
\check{C}^*(p^{-1}\sU, \sW(\det^{-1}E))\to \check{C}^*(p^{-1}\sU, \sC^0_{RS}(\det^{-1}E)).
\]
\end{construction}

Let $U_{2,n}=E_{n-2}\SL_2$ be the scheme of $2\times n$ matrices of rank $2$ and $U_{2,\infty}=\colim\, U_{2,n}$, with the evident left $\SL_2$-action. We use $U_{2,\infty}$ as our model for $E\SL_2$ and $\SL_2\backslash U_{2,\infty}$ for $\BSL_2$. The tautological vector bundle  $\pi{\colon}E_2\to \BSL_2$  is the vector bundle associated to the fundamental representation $\SL_2\subset \GL_2$, in particular, we have a canonical  trivialization of $\det E_2$.

We have the approximation $B_n\SL_2:=\SL_2\backslash U_{2,n+2}$ to $\BSL_2$, and the corresponding restriction map
\[
\alpha_n^i{\colon}H^i(\BSL_2, \sW)\to H^i(B_nSL_2, \sW)
\]
We concentrate on $X:=B_1\SL_2=\SL_2\backslash U_{2,3}$. 

For a matrix  $A\in U_{2,3}$, let $|A|_{ij}$ denote the determinant of the 2 by 2 submatrix of $A$ with columns $i, j$. Sending a matrix $A\in U_{2,3}$ to the triple $(|A|_{23}, -|A|_{13}, |A|_{12})$ gives an isomorphism of 
$X$ with $\A^3\setminus\{0\}$. We use coordinates $y_0, y_1, y_2$ on $\A^3$, and let $U_i\subset \A^3\setminus\{0\}$ be the open subscheme $y_i\neq0$. This gives the affine open cover $\sU=\{U_i\}_{i=0,1,2}$ of $\A^3\setminus\{0\}=\SL_2\backslash U_{2,3}$. We note that $U_{ij}\cong \G_m\times\G_m\times\A^1$ for $i\neq j$, and thus we have $H^1(U_{ij}, \sW)=0$. 

The quotient $\GL_2\backslash U_{2,3}$ is isomorphic to $\P^2$, by sending $A\in U_{2,3}$ to $[|A|_{23}: -|A|_{13}: |A|_{12}]$; this identifies $\P^2$ with $\Gr(2,3)$, that is, our $\P^2$ is the dual $\P^2$ of linear planes in $\A^3$. The induced quotient map $X\to \P^2$ is the standard quotient map $\A^3\setminus\{0\}\to \P^2$ writing $\P^2=\G_m\backslash (\A^3\setminus\{0\})$.

Let $\alpha{\colon}X\to \BSL_2=\colim_nB_n\SL_2$ denote the canonical map. We have the universal bundle $E_2\to \BSL_2$, and we let $E_{2,3}\to X$ denote its restriction to $X$.

\begin{lemma}\label{lem:RestrictionToX} The restriction map
\[
\alpha_1^i{\colon}H^i(\BSL_2, \sW)\to H^i(X, \sW)
\]
is an isomorphism for $0\le i\le 3$.
\end{lemma}

\begin{proof} From the localization sequence for $i{\colon}\{0\}\hookrightarrow\A^3$, $j{\colon}\A^3\setminus\{0\}\hookrightarrow \A^3$, 
\begin{multline*}
\ldots\to H^n(\A^3_k,\sW)\xrightarrow{j^*} H^n(\A^3_k\setminus\{0\}, \sW)
\\\xrightarrow{\delta} H^{n-2}(\Spec k, \sW)\xrightarrow{i_*} H^{n+1}(\A^3_k, \sW)\to\ldots,
\end{multline*}
one sees that $H^i(X, \sW)=0$ for $i\neq 0, 2$,  $H^i(X, \sW)\cong W(k)$ for $i=0,2$, and that the pullback from $\Spec k$ induces the isomorphism $W(k)\cong H^0(X, \sW)$. Since $H^*(\BSL_2, \sW)=W(k)[e]$, with $e=e(E_2)\in H^2(\BSL_2, \sW)$, we see that $\alpha_1^i$ is an isomorphism for $i=0,1, 3$, and  to complete the proof, it suffices to show that $e(E_{2,3})$ is a $W(k)$-basis for $H^2(X, \sW)$.  

The bundle $E_2\to \BSL_2$ is induced by the right $\SL_2$-representation $F$, so the same holds for the bundle $E_{2,3}\to X$. We claim that via the map $q{\colon}X\to \P^2$, $E_{2,3}$ is isomorphic to $q^*\Omega_{\P^2}(1)$. To see this, we can identify $\P^2$ with $P\backslash\GL_3$, where $P$ is the parabolic subgroup-scheme
\[
P=\left\{\begin{pmatrix}*&*&*\\0&*&*\\0&*&*\end{pmatrix}\in \GL_3\right\}.
\] 
We have the exact sequence of right $P$-representations $0\to k^2\to k^3\to k\to0$. where $k^3$ is the restriction to $P$ of the standard  right $\GL_3$-representation on row vectors, $k^2$ is the subspace of vectors $(0,*,*)$ and $k$ is the quotient. Since the $P$-representation on $k^3$ extends to a $\GL_3$-representation, the induced locally free sheaf on $\P^2$ is just $\sO_{\P^2}^3$, and it is standard that the invertible sheaf induced by the $P$-representation on $k$ is $\sO_{\P^2}(1)$, with induced surjection $\sO_{\P^2}^3\to \sO_{\P^2}(1)$ the usual one. Comparing with the Euler sequence, we identify the locally free sheaf associated to the $P$-representation on $k^2$ with $\Omega_{\P^2}(1)$. On the other hand, projecting $\GL_3$ to $U_{2,3}$ by sending a matrix to its last two rows induces an isomorphism of $P\backslash \GL_3$ with $\GL_2\backslash U_{2,3}$ and identifies the bundle $k^2\times^P\GL_3$ with $F\times^{\GL_2}U_{2,3}$. Pulling back to $X$ we see that the bundle $E_{2,3}=F\times^{\SL_2}U_{2,3}$ on $X$ is $q^*\Omega_{\P^2}(1)$. 

The pushforward map  $\pi_{\P^2*}{\colon}H^2(\P^2, \sW(\omega_{\P^2}))\to H^0(\Spec k, \sW)=W(k)$ is an isomorphism. By the motivic Gau{\ss}-Bonnet theorem \cite[Theorem 4.6.1]{DJK} (see also \cite[Theorem 4.1]{LevineEuler}), and Hoyois' computation \cite[Example 1.7]{HoyoisTrace} of $\chi(\P^2/k)\in \GW(k)$ as $\<1\>+H$, the Euler class $e(T_{\P^2})\in H^2(\P^2, \sW(\omega_{\P^2}))$ has pushforward $1\in W(k)$. Thus $e(T_{\P^2})$ is a $W(k)$-basis for $H^2(\P^2, \sW(\omega_{\P^2}))\cong W(k)$. By \cite[Theorem 8.1]{LevineEuler}, we have $e(T_{\P^2})=e(\Omega_{\P^2})$, after applying the canonical isomorphism 
\[
H^2(\P^2, \sW(\omega_{\P^2}))\cong H^2(\P^2, \sW(\omega^{-1}_{\P^2})). 
\]
Finally, by \cite[Theorem 10.1]{LevineEuler}, we have 
$e(\Omega_{\P^2})\cong e(\Omega_{\P^2}(1))$ after applying the canonical isomorphism $H^2(\P^2, \sW(\omega^{-1}_{\P^2}))\cong H^2(\P^2, \sW(\omega^{-1}_{\P^2}(-2)))$. Since $\omega^{-1}_{\P^2}(-2)\cong \sO_{\P^2}(1)$, we may consider $e(\Omega_{\P^2}(1))$ as a $W(k)$-basis for $H^2(\P^2, \sW(\sO_{\P^2}(1)))$. 

Since $q{\colon}\A^3\setminus\{0\}\to \P^2$ is the principal $\G_m$-bundle associated to the line bundle $\sO_{\P^2}(1)$, the canonical section $t$ of $q^*\sO_{\P^2}(1)$ gives an isomorphism 
\[
H^2(X, \sW(q^*\sO_{\P^2}(1)))\cong H^2(X, \sW). 
\]
On the other hand, the fact that, for a field $L$, we have
\[
H^i(\A^1_L\setminus\{0\}, \sW)=\begin{cases} W(L)&\text{ for }i=0\\0&\text{ for }i>0\end{cases}
\]
shows that 
\[
R^iq_*\sW(q^*\sO_{\P^2}(1))=\begin{cases}  \sW(\sO_{\P^2}(1))&\text{ for }i=0\\0&\text{ for }i>0.\end{cases}
\]
Using this, the Leray spectral sequence for $q$ shows that
\[
q^*{\colon}H^i(\P^2, \sW(\sO_{\P^2}(1)))\to H^i(X, \sW(q^*\sO_{\P^2}(1)))
\]
is an isomorphism. Thus $e(q^*\Omega_{\P^2}(1))\in H^2(X, \sW(q^*\sO_{\P^2}(1)))$ is a $W(k)$-basis for $H^2(X, \sW(q^*\sO_{\P^2}(1)))\cong W(k)$, so we conclude that $e(E_{2,3})$ is a $W(k)$-basis for $H^2(X, \sW)$, as claimed.

\end{proof}

We apply Construction~\ref{const:Cocycle} to  $X:=\SL_2\backslash U_{2,3}$ with its open  open cover $\sU=\{U_0, U_1, U_2\}$,   and the tautological bundle $p{\colon}E_{2,3}\to X$ with the canonical trivialization of $\det E_{2,3}$.  This gives us the cocycle $\{\lambda_{ijk}\}_{ijk} \in \check{C}^2(p^{-1}(\sU), \sW)$ and the Euler class $e(E_{2,3})\in H^2(X, \sW)$.
 
 \begin{lemma}\label{lem:CechEulerClass}  Under the canonical map
\[
\check{H}^2(p^{-1}\sU, \sW)\to H^2(E_{2,3}, \sW)\cong H^2(X, \sW)
\]
the cocycle $\{\lambda_{ijk}\}_{ijk} \in \check{C}^2(p^{-1}(\sU), \sW)$ maps to $-e(E_{2,3})$.
\end{lemma}
 
\begin{proof}  Consider the double complex $\check{C}^*(p^{-1}(\sU), \sC^*_{RS}(\sW))$ of Example~\ref{ex:CechDoubleComplex}. Since  the Rost-Schmid complex is a flasque resolution of $\sW$, one  shows using the spectral sequence for the  filtration $F_{II}$ that the augmentation $\epsilon_1{\colon}C^*_{RS}(E_{2,3}, \sW)\to \Tot\check{C}^*(p^{-1}(\sU), \sC^*_{RS}(\sW))$ induces an isomorphism on cohomology. Similarly, the augmentation $\epsilon_2{\colon}\check{C}^*(p^{-1}(\sU), \sW)\to 
\Tot\check{C}^*(p^{-1}(\sU), \sC^*_{RS}(\sW))$ induces the natural map $\check{H}^n(p^{-1}(\sU), \sW)\to H^n(E_{2,3}, \sW)$.

The last  paragraph in Construction~\ref{const:Cocycle}  shows that    $\epsilon_2(\{\lambda_{ijk}\}_{ijk})+\epsilon_1(\<1\>\text{ on }s(X)))$ is a coboundary in  $\Tot\check{C}^*(p^{-1}(\sU), \sC^*_{RS}(\sW))$. Indeed,
\begin{align*}
d_\Tot(\{\phi_i\}_i+\{\lambda_{ij}\}_{ij})&= 
\del^{0,1}_2(\{\phi_i\}_i)+\del^{0,1}_1(\{\phi_i\}_i)
-\del^{1,0}_2(\{\lambda_{ij}\}_{ij})+\del^{1,0}_1(\{\lambda_{ij}\}_{ij})\\
&= 
\del_{RS}(\{\phi_i\}_i)+\check{\delta}^0(\{\phi_i\}_i)
-\del_{RS}(\{\lambda_{ij}\}_{ij})+\check{\delta}^1(\{\lambda_{ij}\}_{ij})\\
&= \{\<1\>\text{ on }s(U_i)\}_i+\epsilon_2(\{\lambda_{ijk}\}_{ijk}\})\\
&=\epsilon_1(\<1\>\text{ on }s(X))+\epsilon_2(\{\lambda_{ijk}\}_{ijk}\}).
\end{align*}
Since   $(\<1\>\text{ on }s(X))$ in $C^2_{RS}(E_{2,3}, \sW)$ represents $p^*e(E_{2,3})$ in $H^2(E_{2,3}, \sW)$, this proves the Lemma.
\end{proof}

Let $Y$ be a smooth $k$-scheme with right $\SL_2$-action and let  $p_Y{\colon}Y\times^{\SL_2} U_{2,3}\to \A^3\setminus\{0\}$ be the projection induced by the structure map $\pi_Y{\colon}Y\to\Spec k$. Let $L$ be an $\SL_2$-linearized invertible sheaf on $Y$, giving the invertible sheaf $\sL$ on $Y\times^{\SL_2}\sU_{2,3}$ by descent. Let $\sU=\{U_i\mid i=0,1,2\}$ be the affine open cover of $\A^3\setminus\{0\}$ defined above, giving us the cover $p_Y^{-1}\sU$ of $Y\times^{\SL_2}\sU_{2,3}$ by pull-back. As before, we have the double complex $\check{C}^*(p_Y^{-1}\sU, \sC^*_{RS}(\sW(\sL)))$, that is,   the \v{C}ech complex with values in the presheaf of Rost-Schmid complexes for $\sW(\sL)$ , with $p,q$ term $\check{C}^p(p_Y^{-1}\sU, \sC^q_{RS}(\sW(\sL))$.  

If we take an injective resolution on $(Y\times^{\SL_2}\sU_{2,3})_\Nis$, $\sC^*_{RS}(\sW(\sL))\to \sI^*$, the induced map of double complexes 
\[
\check{C}^*(p_Y^{-1}\sU, \sC^*_{RS}(\sW(\sL)))\to \check{C}^*(p_Y^{-1}\sU, \sI^*)
\]
gives a map in the derived category $\check{C}^*(p_Y^{-1}\sU, \sC^*_{RS}(\sW(\sL)))\to R\Gamma(Y, Rp_*\sW(\sL))$, and gives a map $\epsilon_\sU{\colon}E^{**}_*(\text{\v{C}ech}) \to E_*^{**}(\text{Leray})$
of the spectral sequence for the filtration $F_I$,
\[
E_2^{p,q}(\text{\v{C}ech}):=\check{H}^p(p_Y^{-1}\sU, \sH^q\sC^*_{RS}(\sW(\sL))\Rightarrow H^{p+q}(\Tot \check{C}^*(p_Y^{-1}\sU, \sC^*_{RS}(\sW(\sL))),
\]
to  the Leray spectral sequence,
\[
E_2^{p,q}(\text{Leray})=H^p(\A^3\setminus\{0\}, R^qp_{Y*}\sW(\sL))\Rightarrow H^{p+q}(Y\times^{\SL_2}\sU_{2,3}, \sW(\sL)).
\]

Let $\sH^n(Y, \sW(\sL))$ be the Nisnevich sheaf on $\Sm/k$ associated to the presheaf $U\mapsto H^n(Y\times_kU, \sW(p_1^*L))$. We note that  $\SL_2$-action on $\sH^n(Y, \sW(\sL))$ induced by the $\SL_2$-action on $Y$ is trivial, by Remark~\ref{rem:SL2Inv}. 

Sending $(x,y,z)\in U_0$ to 
\[
s_0(x,y,z):=\begin{pmatrix}-y&x&0\\-z/x&0&1\end{pmatrix} 
\]
gives a section $s_0{\colon}U_0\to U_{2,3}$ of $U_{2,3}\to \A^3\setminus\{0\}$; we have similarly defined sections over $U_1, U_2$,
\[
s_1(x,y,z):= \begin{pmatrix}-y&x&0\\0&-z/y&1\end{pmatrix},\ 
 s_2(x,y,z):= \begin{pmatrix}z&0&-x\\0&1&-y/z\end{pmatrix}.
\]
These define trivializations $\Psi_i$ of the $\SL_2$-principal bundle $U_{2,3} \to \A^3\setminus\{0\}$ over $U_i$. The $\Psi_i$ induce trivializations $\psi_i$ of $E_{2,3}\to \A^3\setminus\{0\}$ over $U_i$, with cocycle $\{\psi_{ij}\}_{ij}$.

\begin{lemma}\label{lem:CechLerayIso}  Suppose  that   $\sH^q(Y, \sW(\sL))$ is a free, finitely generated  $\sW$-module for each $q\ge0$.  Then the natural maps
\begin{equation}\label{eqn:CechLerayIso1}
\check{H}^p(p_Y^{-1}\sU, \sH^q\sC^*_{RS}(\sW(\sL))\to H^p(\A^3\setminus\{0\}, R^qp_{Y*}\sW(\sL)),
\end{equation}
and
\begin{equation}\label{eqn:CechLerayIso2}
H^{p+q}(\Tot \check{C}^*(p_Y^{-1}\sU, \sC^*_{RS}(\sW(\sL)))\to H^{p+q}(Y\times^{\SL_2}U_{2,3}, \sW(\sL))
\end{equation}
are isomorphisms and $\epsilon_\sU{\colon}E^{**}_*(\text{\v{C}ech}) \to E_*^{**}(\text{Leray})$ is an isomorphism of  spectral sequences.
\end{lemma}

\begin{proof}  Since the $\SL_2$ action on $\sH^*(Y, \sW(\sL))$ is trivial, we have a canonical isomorphism of sheaves on $(\A^3\setminus\{0\})_\Nis$, $\sH^q(Y, \sW(\sL))\cong R^qp_{Y*}(\sW(\sL))$, for each $q\ge0$, giving us the map \eqref{eqn:CechLerayIso1}.

We now show that the map \eqref{eqn:CechLerayIso1} is an isomorphism for all $p,q$. Indeed, it follows from our assumption on $Y$ that for each $q\ge0$
 there is an integer $r_q\ge0$ and an isomorphism of Nisnevich sheaves on $\A^3\setminus\{0\}$, $R^qp_{Y*}\sW(\sL)\cong  \sW^{r_q}$.

 For each subset $I\subset \{0,1,2\}$, the corresponding open subset $U_I\subset \A^3\setminus\{0\}$ is isomorphic to $\A^{3-|I|}\times(\A^1\setminus\{0\})^{|I|}$. Moreover, using the sections $s_i{\colon}U_i\to U_{2,3}$, we have isomorphisms $p_Y^{-1}(U_i)\cong Y\times U_i$ for $i=0,1,2$.  Multiplication by the algebraic Hopf map $\eta$ is invertible on $\sW$, and after inverting $\eta$ and $-\wedge S^1$, $\A^1\setminus\{0\}$ is isomorphic to $\Spec k\amalg \Spec k$. Thus, 
\[
H^i(U_I, R^qp_{Y*}\sW(\sL))=H^i(U_I, \sW^{r_q})=0
\]
for $i>0$,   that is, $\sU$ is a Leray cover for the sheaves $R^qp_{Y*}\sW(\sL)$. This shows that the map \eqref{eqn:CechLerayIso1} is an isomorphism. 

As the map $\epsilon_\sU$ on the $E_2$-terms is the isomorphism \eqref{eqn:CechLerayIso1},  $\epsilon_\sU$ is an isomorphism of spectral sequences.  Since both spectral sequences are strongly convergent, it follows that \eqref{eqn:CechLerayIso2} is an isomorphism.  
\end{proof}

\begin{exs}\label{exs:CechExamples} Here we give some examples of $Y\in\Sm^{\SL_2}/k$ satisfying the hypotheses in Lemma~\ref{lem:CechLerayIso}.

We may take $Y=F\setminus\{0\}$ with the right $\SL_2$-action induced from that on $F$. Forgetting the $\SL_2$-action, $Y=\A^2_k\setminus\{0\}$, and the localization sequence for the open immersion $Y\times U\hookrightarrow \A^2\times U$ with closed complement $0\times U\hookrightarrow \A^2\times U$,
 \begin{multline*}
\ldots\to H^n(\A^2\times U, \sW)\xrightarrow{j^*} H^n(Y\times U, \sW)\\
\xrightarrow{\delta}H^{n-1}(U, \sW)\to  H^{n+1}(\A^2\times U, \sW)\to\ldots,
\end{multline*}
 gives isomorphisms $\sH^n(Y, \sW)\cong  \sW$ for $n=0,1$ and shows that $\sH^n(Y, \sW)=0$ for $n>1$.   

We may also take $Y=\P^1(F)$, again with the right $\SL_2$ action induced by the $\SL_2$-action on $F$. Since the (unstable)  Hopf map $\eta{\colon}\A^2\setminus\{0\}\to \P^1$ acts invertibly on $H^*(-, \sW)$, $\eta$ induces an $\SL_2$-equivariant isomorphism
\[
\eta^*{\colon}\sH^n(\P^1(F), \sW)\to \sH^n(F\setminus\{0\}, \sW)
\]

As a third example, the sheaf $\sH^n(\P^1(F), \sW(\sO_{\P^1(F)}(1)))$ is the zero sheaf for all $n\ge0$.  

In all these cases,   the fact that $\SL_2$ acts trivially on $\sH^n(Y,\sW)$ (Remark~\ref{rem:SL2Inv}) gives us a canonical isomorphism $R^np_{Y*}\sW\cong \sH^n(Y,\sW)$ for each $n\ge0$.
\end{exs}

\begin{lemma}\label{lem:DifferentialComp1} Let $Y\in \Sm/k$ be either\\[5pt]
a) $Y=F\setminus\{0\}$\\[2pt]
b) $Y=\P(F)$,\\[5pt]
endowed with the corresponding right $\SL_2$-action, and let $p_Y{\colon}Y\times^{\SL_2}U_{2,3}\to \A^3\setminus\{0\}$ be the map induced by $Y\to \Spec k$. In both cases, we have defined isomorphisms of sheaves on $\A^3\setminus\{0\}$,
\[
\sW\xrightarrow{\tau}  R^1p_{Y*}\sW
\]
and
\[
\sW\xrightarrow{p_Y^*}p_{Y*}\sW.
\]
Let $\phi^0\in H^0(\A^3\setminus\{0\}, Rp_{Y*}^1\sW)$ be the global section $\tau(1)$. Then in the Leray spectral sequence
\[
E_2^{p,q}=H^p(\A^3\setminus\{0\}, Rp_{Y*}^q(\sW))\Rightarrow H^{p+q}(Y\times^{\SL_2}U_{2,3}, \sW), 
\]
we have 
\[
d_2^{0,1}(\phi^0)=p_Y^*e(E_{2,3})\in H^2(\A^3\setminus\{0\}, p_{Y*}\sW)\cong H^2(\A^3\setminus\{0\}, \sW). 
\]
\end{lemma}
 
\begin{proof}  Using the various isomorphisms induced by the unstable Hopf map $\eta{\colon}F\setminus\{0\}\to \P(F)$, we reduce to the case $Y=F\setminus\{0\}$. In this case,  $Y\times^{\SL_2}U_{2,3}$ is the punctured tautological bundle $E_{2,3}\setminus \{0_E\}$, and $p_Y$ is   the  map $p_0{\colon}E_{2,3}\setminus \{0_E\}\to \A^3\setminus\{0\}$ induced from the projection $p{\colon}E_{2,3}\to \A^3\setminus\{0\}$ by restriction. We also have a canonical trivialization of $\det E_{2,3}$. 

By Lemma~\ref{lem:CechLerayIso}, we may replace the Leray spectral sequence with the \v{C}ech spectral sequence $E_*^{**}(\text{\v{C}ech})$ for the filtration $F_I$.  

Using the Rost-Schmid resolution of $\sW$, we use Construction~\ref{const:Cocycle} to represent  the differential $d_2^{0,1}$  in \v{C}ech cohomology.  We retain the notation from Construction~\ref{const:Cocycle}  and the proof of Lemma~\ref{lem:CechEulerClass}. 

Following our conventions for the spectral sequence of a double complex, to compute $d_2^{0,1}$, one starts with a class $\{\tau_i\}_i\in  \check{C}^0(p^{-1}_0\sU, \sC_{RS}^1(\sW))$ with $\del_2^{0,1}(\{\tau_i\}_i)=\{\del_{RS}(\tau_i)\}_i=0$ and with $\check{\delta}^0(\{\tau_i\}_i)=\del^{0,1}_1(\{\tau_i\}_i)$ a Rost-Schmid coboundary:
\[
\del^{0,1}_1(\{\tau_i\}_i)=\del_{RS}(\{\mu_{ij}\}_{ij}):=-\del^{1,0}_2(\{\mu_{ij}\}_{ij}).
\]
In particular,  the cochain $\{\tau_i\}_i$ defines a \v{C}ech cocycle 
\[
\{\bar\tau_i\}_i\in \check{C}^0(p^{-1}_0\sU, \sH^1\sC^*_{RS}(E_{2,3}\setminus\{0\}, \sW)). 
\]
Let $\bar\tau\in H^0(\A^3\setminus \{0\}, \sH^1(E_{2,3}\setminus\{0\}, \sW))$ be the corresponding class. 
Then $d_2^{0,1}(\bar\tau)$ is represented by the \v{C}ech cocycle $\del^{1,0}_1(-\{\mu_{ij}\}_{ij})=-\check{\delta}^1(\{\mu_{ij}\}_{ij})$. 

We recall the class $\phi_i\in \sC_{RS}^1(p^{-1}(U_i), \sW)$ from Construction~\ref{const:Cocycle}.  We take $\tau_i$ to be the restriction $\phi^0_i$ of $\phi_i$ to $p_0^{-1}(U_i)= E_{2,3}\setminus\{0\}\cap p^{-1}(U_i)$. Since $\del_{RS}(\phi_i)=\<1\>_{s(U_i)}$, and we are removing the 0-section from $E_{2,3}$,  we have $\del_{RS}(\phi^0_i)=0$. Moreover, since $\del_{RS}(\phi_i)=\<1\>_{s(U_i)}$, one sees from the localization sequence for the open immersion $E_{2,3}\setminus\{0_E\}\hookrightarrow E_{2,3}$ with closed complement $s(\A^3\setminus\{0\})$, that  the cocycle
 $\{\phi^0_i\}_i\in \check{C}^0(p_0^{-1}(\sU), \sC^1_{RS}(\sW))$ is a representative of  $\phi^0\in H^0(\A^3\setminus\{0\}, R^1p_{Y*}\sW)$. 
 
Replace the cochains $\{\lambda_{ij}\}_{ij}$ and $\{\lambda_{ijk}\}_{ijk}$ from the proof of Lemma~\ref{lem:CechEulerClass} with their respective restrictions $\{\lambda^0_{ij}\}_{ij}$ and $\{\lambda^0_{ijk}\}_{ijk}$  to $\check{C}^*(p_0^{-1}\sU,-)$. We let $\bar\lambda\in \check{H}^2(p_0^{-1}\sU, \sH^0\sC^*_{RS}(\sW))$ be the class represented by $\{\lambda^0_{ijk}\}_{ijk}$.

 Following the last paragraph of Construction~\ref{const:Cocycle}, we have identities
\[
-\del_2^{1,0}(\{\lambda_{ij}\}_{ij})=\del_{RS}(\{\lambda_{ij}\}_{ij})=\check{\delta}^0(\{\phi_i\}_i),
\]
and
\[
\del_1^{1,0}((\{\lambda_{ij}\}_{ij}))=\check{\delta}^1((\{\lambda_{ij}\}_{ij}))=\{\lambda_{ijk}\}_{ijk}.
\]
These with fact that $\{\phi^0_i\}_i$ represents $\phi^0$ 
show that $d^{0,1}_2(\phi^0)=-\bar\lambda$. On the other hand, the restriction map 
$p_*(\sH^0\sC^*_{RS}(\sW)_{|E_{2,3}})\to p_{0*}(\sH^0\sC^*_{RS}(\sW)_{|E_{2,3}\setminus\{0\}})$ is an isomorphism, as is  the pullback $p_0^*{\colon}\sW\to p_{0*}\sW$ (all considered as presheaves on $(\A^3\setminus\{0\})_\Nis$); in addition the augmentation
\[
\sW\to \sH^0\sC^*_{RS}(\sW)
\]
is also an isomorphism (of presheaves on $\Sm/k_\Nis$).

 Thus,  considering $d^{0,1}_2(\phi^0)$ as an element in 
$\check{H}^2(\sU, p_{0*}\sW)$ via these isomorphisms, we see that $d^{0,1}_2(\phi^0)\in 
\check{H}^2(\sU, p_{0*}\sW)$ is represented by $\{-\lambda_{ijk}\}_{ijk}$. By Lemma~\ref{lem:CechEulerClass}, this class is exactly $p_0^*e(E_{2,3})$, via the isomorphism $\check{H}^2(\sU, p_{0*}\sW)\cong H^2(\A^3\setminus\{0\}, p_{0*}\sW)$. Since $p_0^*{\colon}H^2(\A^3\setminus\{0\},\sW)\to H^2(\A^3\setminus\{0\}, p_{0*}\sW)$ is an isomorphism, the lemma is proved.
\end{proof} 
 
 \begin{corollary}\label{lem:LerayDifferentialComputation} Let $Y=F\setminus\{0\}$ or $\P(F)$, giving $p_Y{\colon}Y\times^{\SL_2}E\SL_2\to \BSL_2$. In the Leray spectral sequence 
  \[
E_2^{p,q}=H^p(\BSL_2, R^qp_{Y*}(\sW))\Rightarrow H^{p+q}(Y\times^{\SL_2}E\SL_2, \sW).
\]
and with respect to our choice of isomorphisms
\[
p_Y^*{\colon}\sW\xrightarrow{\sim} p_{Y*}(\sW),\ \tau{\colon}\sW\xrightarrow{\sim} R^1p_{Y*}(\sW),
\]
described above, the differential
\[
d_2^{0,1}{\colon}H^0(\BSL_2,  \sW)\to
H^2(\BSL_2, \sW)
\]
 satisfies
\[
d^{0,1}_2(\<1\>_{\BSL_2})=e(E_2)\in H^2(\BSL_2, \sW).
\]
\end{corollary}

\begin{proof} The class $\phi^0$ of Lemma~\ref{lem:DifferentialComp1} is the restriction of $\tau(\<1\>_{\BSL_2})$ to $H^0(\A^3\setminus\{0\}, \sW)$ and $d^{0,1}_2(\phi^0)$ is the restriction of $d^{0,1}_2(\<1\>_{\BSL_2})$ to $H^2(\A^3\setminus\{0\}, \sW)$ via the restriction map $\alpha_1^*{\colon}H^*(\BSL_2, \sW)\to H^*(\A^3\setminus\{0\})$. By Lemma~\ref{lem:RestrictionToX}, Lemma~\ref{lem:DifferentialComp1} yields the result.
\end{proof}

We recall  the notations introduced in Notation~\ref{not:QuadraticForms}. In particular,  giving $F$ (with the standard right $\SL_2$-action)  the standard basis $e_0, e_1$, we identify $\Sym^2F^\vee=H^0(\P(\Sym^2F), \sO(1))$ with the space of quadratic forms on $F$, with the linear forms $T_0, T_1, T_2$ on $\P^2=\P(\Sym^2F)$ being the basis dual to $e_0^2, e_0e_1, e_1^2$.  Similarly, we have the basis $X_0, X_1$ for $H^0(\P(F), \sO(1))$ dual to $e_0, e_1$. 

The squaring map $\sq{\colon} F\to  \Sym^2F$ is given by
 \[
 sq(x_0e_0+x_1e_1)=x_0^2e_0^2+2x_0x_1e_0e_1+x_1^1e_1^2
 \]
 so $\sq{\colon}\P(F)\to \P(\Sym^2F)$ satisfies
 \[
 sq^*(T_0)=X_0^2,\ sq^*(T_1)=2X_0X_1,\ sq^*(T_2)=X_1^2.
 \]
The map $\sq{\colon}\P(F)\to \P(\Sym^2F)$ is a closed immersion, with image curve $D\subset \P^2$ with defining polynomial $Q:=T_1^2-4T_0T_2$.

The map $\pi{\colon}\SL_2\to \P(\Sym^2F)=\P^2$ defined by
 \[
 \pi\left(\begin{pmatrix}x_0&x_1\\y_0&y_1\end{pmatrix}\right):=[x_0y_0: x_0y_1+x_1y_0:x_1y_1]
 \]
 gives isomorphisms of right $\SL_2$-homogeneous spaces 
 \[
 N\backslash\SL_2\cong \P(\Sym^2F)\setminus \P(F)\cong\P^2\setminus D.
 \]
 
For $W\in \Sm/k$, we write $H^*(W, \sW(\omega))$ for $H^*(W, \sW(\omega_W))$, to shorten the notation. For the rest of this section, we take $Y:=\P(\Sym^2F)\setminus \P(F)$. By the above identifications, we have a natural isomorphism $H^*_{\SL_2}(Y,  \sW)\cong H^*(BN, \sW)$. 

 The following lemma shows that $H^*_{\SL_2}(Y,  \sW(\omega))$ computes $H^*(BN, \sW(\gamma))$. 

\begin{lemma}\label{lem:SignOmega} Let $\gamma$ be the invertible sheaf on $BN$ corresponding to the sign representation $N\to N/T_1=\Z/2= \{\pm1\}$. Then there is a canonical isomorphism $H^*_{\SL_2}(Y,  \sW(\omega))\cong 
H^*(BN, \sW(\gamma))$.
\end{lemma}

\begin{proof} We have the $\SL_2$-equivariant identifications
\[
N\backslash \SL_2\cong \<\bar{\sigma}\>\backslash(T_1\backslash\SL_2)\cong 
\<\bar{\sigma}\>\backslash(\P^1\times\P^1\setminus \Delta)\cong Y,
\]
where $\bar{\sigma}$ is the image of $\sigma$ in $N/T_1$, acting on $\P^1\times\P^1$ by the exchange of factors. Via these identifications, $\gamma$ corresponds to the invertible sheaf on $Y$ induced by twisting the canonical $\bar{\sigma}$-action on $\sO_{\P^1\times\P^1\setminus \Delta}$ by the sign representation, and $\omega_Y$ is similarly induced by the canonical action of $\bar{\sigma}$ on $\omega_{\P^1\times\P^1\setminus \Delta}$. 

Let $p{\colon}{\P^1\times\P^1\setminus \Delta}\to Y\subset \P(\Sym^2F)$ be the projection; as we have mentioned above, this identifies ${\P^1\times\P^1\setminus \Delta}$ with the cover of $Y=\P^2\setminus V(T_1^2-4T_0T_2)$ defined by $T_3^2=T_1^2-4T_0T_2$. We have the global $\SL_2$-invariant generating section $p^*\Omega_2/T_3^2$ of $\omega_{\P^1\times\P^1\setminus \Delta}$, defining an $\SL_2$-equivariant isomorphism
\[
\times p^*\Omega_2/T_3^2{\colon}\sO_{\P^1\times\P^1\setminus \Delta}\xrightarrow{\sim}\omega_{\P^1\times\P^1\setminus \Delta}.
\]
Since $p^*\Omega_2$ is invariant under $\bar{\sigma}$ and $\bar{\sigma}^*(T_3)=-T_3$,   $\times p^*\Omega_2/T_3^2$ intertwines the sign-twisted action of $\bar{\sigma}$ on $\sO_{\P^1\times\P^1\setminus \Delta}$ with the canonical action on $\omega_{\P^1\times\P^1\setminus \Delta}$.
\end{proof}

We now apply the constructions of the first part of this section to compute $H^*_{\SL_2}(Y,  \sW)$ and  $H^*_{\SL_2}(Y,  \sW(\omega))$.

Following \eqref{eqn:OrOmega} and Remark~\ref{rem:SLInvOrientation}, we have the $\SL_2$-invariant oriented generating section $\Omega_2$ of $\omega_{\P(\Sym^2F)}(3)$,
and the $\SL_2$-invariant  oriented generating section $\Omega_1$ of $\omega_{\P(F)}(2)$. Following Lemma~\ref{lem:OriGen}, we use  $\Omega_1$ to define our chosen orientation $\omega_{\P(F)}\xrightarrow{\sim}\sO_{\P(F)}(-1)^{\otimes 2}$

Let $\sV$ be the $\SL_2$-linearized sheaf on $Y$, $\sV:=\sO_Y(-2)\otimes_kF$. We have the $\SL_2$-invariant, non-degenerate quadratic form on $\sV$ \eqref{eqn:TildeQ}
\[
\tilde{Q}{\colon}\sV\to \omega_{Y/k}
\]
defined for local sections $x_0, x_1$ of $\sO_Y(-2)$ by
\[
\tilde{Q}(x_0e_0+x_1e_1)=2T_0x_0^22+2T_1x_0x_1+2T_2x_1^2\otimes\Omega_2.
\]
We have the corresponding element $\tilde{Q}\in H^0(Y, \sW(\omega))$ and since $\tilde{Q}$ is an $\SL_2$-invariant, this lifts to  a class $[\tilde{Q}]\in H^0_{\SL_2}(Y, \sW(\omega_Y))$. Similarly, $Q$ defines an $\SL_2$-invariant quadratic form $\<Q\>$ on $\sO_Y(-1)$, giving the element $[\<Q\>]\in  H^0_{\SL_2}(Y, \sW)$.

The elements $[\<Q\>]\in H^0_{\SL_2}(Y, \sW)$ and $[\tilde{Q}]\in H^0_{\SL_2}(Y, \sW(\omega_Y))$ define the elements $x\in H^0(BN, \sW)$ and $y\in H^0(BN, \sW(\gamma))$, via canonical isomorphisms $\phi_0{\colon}H^0_{\SL_2}(Y,\sW)\xrightarrow{\sim} H^0(BN, \sW)$ and $\phi_0'{\colon}H^0_{\SL_2}(Y,\sW(\omega_Y))\xrightarrow{\sim} H^0(BN, \sW(\gamma))$ mentioned in the statement of Theorem~\ref{thm:Presentation2}. 

Next, we recall some facts from \cite{LevineBG}. 
In \cite[Proposition 5.5]{LevineBG}, we have computed $H^*(BN\, \sW)=H^*_{\SL_2}(Y, \sW)$ as the $H^*(\BSL_2,\sW)=W(k)[e]$-algebra $W(k)[e, x]/((1+x)e)$, with $x$ corresponding to  $[\<Q\>]\in H^0_{\SL_2}(Y, \sW)$.

To compute $H^*_{\SL_2}(Y, \sW(\omega))$, we use a localization sequence, first computing 
$H^*_{\SL_2}(\P(\Sym^2F), \sW(\omega))$ and $H^*_{\SL_2}(\P(F), \sW(\omega))$. 

\begin{lemma}\label{lem:SL2TwistCohP2} $H^*_{\SL_2}(\P(\Sym^2F), \sW(\omega))$ is a free  $H^*(\BSL_2, \sW)=W(k)[e]$ module with basis $e_{\SL_2}(T_{\P(\Sym^2F)})$ (in degree 2).
\end{lemma}

\begin{proof}
By Lemma~\ref{lem:PnCoh}, we have $H^p(\P(\Sym^2F), \sW(\omega)))=0$ for $p\neq 2$, and  the projection 
\[
\pi_{\P(\Sym^2F)*}{\colon}H^2(\P(\Sym^2F), \sW(\omega)\to H^0(\Spec k, \sW)=W(k)
\]
is an isomorphism. 
On the other hand, by \cite[Example 2.6(1)]{LevineEuler},   the motivic Euler characteristic $\chi(\P(\Sym^2F)/k)\in \GW(k)$ is 
\[
\chi(\P(\Sym^2F)/k)=\<1\>+H\in \GW(k).
\]
We have the Euler class $\tilde{e}:=e(T_{\P(\Sym^2F)})\in H^2(\P(\Sym^2F), \sW(\omega))$, and by the  motivic Gau{\ss}-Bonnet theorem \cite[Theorem 2.6.1]{DJK}, \cite[Theorem 5.3]{LR}, the  image of $\chi(\P(\Sym^2F)/k)$ in $W(k)$ is $\pi_{\P(\Sym^2F)*}(\tilde{e})$. Thus  $H^2(\P(\Sym^2F), \sW(\omega))$ is the free $W(k)$-module with basis $\tilde{e}$.
 
Since $\SL_2$ acts trivially on $H^*(\P(\Sym^2F), \sW(\omega))$ (Remark~\ref{rem:SL2Inv}), the  Leray spectral sequence
 for $H^*_{\SL_2}(\P(\Sym^2F), \sW(\omega))$ has 
 \[
 E_2^{p,q}=\begin{cases}0&\text{ for }q\neq2\\ H^p(\BSL_2, \sW)\cdot \tilde{e}
 &\text{ for }q=2,
 \end{cases}
 \]
 and  thus degenerates at $E_2$, and yielding the description of $H^*_{\SL_2}(\P(\Sym^2F), \sW(\omega))$ as given in the statement of the Lemma.
 \end{proof}

Let $t:=X_1/X_0$ be our standard coordinate on $\P(F)$.  By Lemma~\ref{lem:OriGen},
$dt$ is an oriented generating section of  $\omega_{\P(F)}$ on $\P(F)\setminus \{[0,1]\}$, and $-d(1/t)=t^{-2}dt$ is an oriented generating  section on $\P(F)\setminus \{[0:1]\}$.  The corresponding orientation $\rho{\colon}\omega_{\P(F)}\to \sO_{\P(F)}(-1)^{\otimes 2}$ gives us the well-defined $\SL_2$-equivariant isomorphism
\[
-\otimes dt{\colon}H^0(\P(F), \sW)\to H^0(\P(F), \sW(\omega)).
\]
Composing with the isomorphism $\pi_{\P(F)}^*{\colon}W(k)\to H^0(\P(F), \sW)$ gives us our canonical choice of isomorphism $\pi_{\P(F)}^*\otimes dt{\colon}W(k)\to H^0(\P(F), \sW(\omega_{\P(F)}))$.

We also have the isomorphism $i_{[1:0]*}{\colon}W(k)\to H^1(\P(F), \sW(\omega))$. 

We have the localization sequence for $\P(F)\subset \P(\Sym^2F)$
\begin{multline}\label{multline:LocSeq}
\ldots\to H^{p-1}(\P(F), \sW(\omega))\xrightarrow{i_{\P(F)*}} H^p(\P(\Sym^2F), \sW(\omega))\\\xrightarrow{j^*}
 H^p(Y, \sW(\omega))\xrightarrow{\del}H^{p}(\P(F), \sW(\omega)\to\ldots . 
 \end{multline}
 and the corresponding version in $\SL_2$-equivariant cohomology (see Lemma~\ref{lem:HtpyModBound}  and Lemma~\ref{prop:BMLocalization}).

\begin{lemma} The restriction map 
 \[
j^*{\colon}H^p_{\SL_2}(\P(\Sym^2F), \sW(\omega))\xrightarrow{\sim}H^p_{\SL_2}(Y, \sW(\omega))
\]
is an isomorphism for $p\ge 1$ and the boundary map 
\[
H^0_{\SL_2}(Y, \sW(\omega))\xrightarrow{\del} H^0_{\SL_2}(\P(F), \sW(\omega_{\P(F)}))\cong W(k)
\]
is an isomorphism.
\end{lemma}
\begin{proof} By \cite[Lemma 5.1]{LevineBG}, we have
\[
H^*_{\SL_2}(\P(F), \sW(\omega_{\P(F)}))=W(k),
\]
concentrated in degree 0 and $H^p_{\SL_2}(\P(\Sym^2F), \sW(\omega))=0$ for $p<2$ (Lemma~\ref{lem:SL2TwistCohP2}). We then use the $\SL_2$-equivariant version of the  localization sequence \eqref{multline:LocSeq}  to  finish the proof.
\end{proof}

The next computation will help us understand the Leray spectral sequence for the projection  $\pi_Y{\colon}Y\times^{\SL_2}E\SL_2\to \BSL_2$.

\begin{lemma}\label{lem:WittCohP2-D} Take $U\in \Sm/k$.\\[2pt]
1. $H^q(Y\times U, \sW(\omega_Y))=0$ for $q\neq 0$ and $H^0(Y\times U, \sW(\omega_Y))$ is a free $H^0(U, \sW)$-module with generator $[\tilde{Q}]$.
\\[2pt]
2. Let $\delta{\colon}H^0(Y, \sW(\omega_Y))\to H^0(\P(F), \sW(\omega_{\P(F)}))$ be the boundary in the localization sequence for $\sq{\colon}\P(F)\to \P(\Sym^2F)$. Then
\[
\delta(\tilde{Q})=-1\otimes dt\in H^0(\P(F), \sW(\omega_{\P(F)})).
\]
\end{lemma}

\begin{proof}  
 For (1), we have the localization sequence for the closed immersion $\sq{\colon}\P(F)\times U\to \P(\Sym^2F)\times U$ with open complement $j{\colon}Y\times U\to \P(\Sym^2F)\times U$
\begin{multline*}
\ldots\to H^{q-1}(\P(F)\times U,\sW(\omega_{\P(F)}))\xrightarrow{sq_*}  H^q(\P(\Sym^2F)\times U, \sW(\omega_{\P(\Sym^2F)}))\\\xrightarrow{j^*}
H^q(Y\times U, \sW(\omega_Y))\xrightarrow{\del} H^{q}(\P(F)\times U,\sW(\omega_{\P(F)}))\to\ldots
\end{multline*}
We have  $H^0(\P(F)\times U,\sW(\omega_{\P(F)}))=H^1(\P(F)\times ,\sW(\omega_{\P(F)}))=H^0(U, \sW)$. Also,  $H^*(\P(\Sym^2F)\times U, \sW(\omega_{\P(\Sym^2F)}))$ is $H^0(U, \sW)$, supported in degree 2, and  $sq_*{\colon}H^1(\P(F)\times U,\sW(\omega_{\P(F)}))\to H^2(\P(\Sym^2F)\times U, \sW(\omega_{\P(\Sym^2F)}))$ is an isomorphism (Lemma~\ref{lem:PnCoh}). Thus (1) follows, once we know that the boundary  \[
\del{\colon}H^0(Y\times U, \sW(\omega_Y))\to H^0(\P(F)\times U,\sW(\omega_{\P(F)}))
\]
is an isomorphism. But this follows from (2), using Lemma~\ref{lem:PnCoh} again.

For (2), we identify $\P(\Sym^2F)$ with $\P^2$ and $\P(F)$ with $D$, and use our isomorphism $\phi{\colon}\P^1\to D$, $\phi([x_0,x_1])=[x_0^2, 2x_0x_1, x_1^2]$ as before. We have affine coordinates  $t_i=T_i/T_0$ on $\P^2\setminus\{T_0=0\}$. We use our orientation on $D$ to define an isomorphism $H^0(D, \sW(\omega))\cong H^0(D, \sW)$; recall that $(1/2)dt_1$ is a  oriented section of $\omega_D$, generating over $D\setminus \{[0,0,1]\}$, since $\phi^*(1/2)dt_1)=dt$.  

Since the restriction map 
$H^0(D, \sW) \to W(k(D))$ is injective,  we may restrict $\tilde{Q}$ to $\Spec \sO_{Y, D}$ and use the boundary map $\delta{\colon}H^0(\Spec \sO_{Y,D}, \sW(\omega_Y))\to W(k(D))$ to compute  $\delta(\tilde{Q})$. In particular, we may use the generating section $T_0$ of $\sO_Y(1)$ in a neighborhood of the generic point of $D$ to give oriented generator $T_0^{-2}$ for $\sO_Y(-2)$. Letting $t_i=T_i/T_0$, and noting that $\Omega_2/T_0^3=dt_1dt_2$,  this rewrites $\tilde{Q}$ as the quadratic form $\tilde{q}{\colon}\sO_{Y\setminus V(T_0)}^2\to \omega_{Y\setminus V(T_0)}$, 
\[
\tilde{q}(x_0, x_1)=x_0^2+t_1x_0x_1+t_2x_1^2\otimes 2dt_1dt_2.
\]
We may then diagonalize $\tilde{q}$, giving the equivalent form 
\[
\tilde{q}\sim (\<1\>+\<4t_2-t_1^2\>)\otimes 2dt_1dt_2.
\]

Let $g=4t_2-t_1^2$. To compute  the boundary, we may use the ``coordinatized'' version
\[
\del=\del_g\otimes \del/\del g.
\]
appearing the the Rost-Schmid complex. We then have
\begin{align*}
\del(\tilde{q})&=(\del_g\otimes \del/\del g)((\<1\>+\<g\>)\otimes 2dt_1dt_2)\\
&=\<1\>\otimes (\del/\del g, 2dt_1dt_2)\\
&=\<1\>\otimes (\del/\del g, (-1/2)dgdt_1)\\
&=\<1\>\otimes   (-1/2)dt_1)\\
&=-1\otimes (2dt_1)\\
&=-1\otimes dt
\end{align*}
\end{proof}

Recall that $\tilde{Q}$ is an $\SL_2$-invariant quadratic form, defining the element $[\tilde{Q}]\in H^0_{\SL_2}(Y, \sW(\omega_Y))$ mentioned above. 

Recall from Remark~\ref{rem:SL2Inv} that $\SL_2$ acts trivially on $H^*(Y\times_kU, \sW(\omega))$ for any $U\in \Sm/k$.  Thus, the sheaf $R^q\pi_{Y*}\sW(\omega)$ on $\BSL_2$ is the Nisnevich sheaf  $\sH^q(Y, \sW(\omega_Y))$ associated to $U\mapsto H^q(Y\times U,\sW(\omega_Y)))$.

\begin{lemma}\label{lem:EdgeHom} 1. The edge homomorphism
\[
H^2(\BSL_2, \sH^0(Y, \sW(\omega_Y)))\to H^2_{\SL_2}(Y, \sW(\omega_Y))
\]
with respect to  the Leray spectral sequence
\[
E_2^{p,q}=H^p(\BSL_2, \sH^q(Y, \sW(\omega_Y)))\Rightarrow H^{p+q}_{\SL_2}(Y, \sW(\omega_Y))
\]
is an isomorphism.
\\[2pt]
2. Using $[\tilde{Q}]\in H^0(Y, \sW(\omega_Y))$ as $\sW$-generator for the sheaf 
$\sH^0(Y, \sW(\omega_Y)))\cong \sW$ on $\BSL_2$, we have the isomorphism
\[
H^2(\BSL_2, \sH^0(Y, \sW(\omega_Y)))=[\tilde{Q}]\cdot H^2(\BSL_2, \sW)=W(k)\cdot ([\tilde{Q}]\cdot e)
\]
3. Under the isomorphisms of (1) and (2), the Euler class $e_{\SL_2}(T_Y)\in  H^2_{\SL_2}(Y, \sW(\omega_Y))$ gets sent to $[\tilde{Q}]\cdot e$. 
\end{lemma}

\begin{proof} By Lemma~\ref{lem:WittCohP2-D},  $\sH^q(Y, \sW(\omega_Y))$ is the sheaf $\sW$, generated by the class $[\tilde{Q}]$,  for $q=0$, and is zero for $q\neq0$. This proves  (1) and (2).

For (3), let $j{\colon}Y\to \P(\Sym^2 F)$ be the inclusion, so we have
\[
e_{\SL_2}(T_Y)=j^*e_{\SL_2}(T_{\P(\Sym^2F)}).
\]
Since the inclusion $\alpha{\colon} \A^3\setminus\{0\}\to \BSL_2$ induces an isomorphism on $H^p(-, \sW)$ for $p\le 3$ (Lemma~\ref{lem:RestrictionToX}), we may restrict to the corresponding bundles over $\A^3\setminus\{0\}$ via  $\alpha$.

Let $p_{\P(\Sym^2F)}{\colon}\P(\Sym^2F)\times^{\SL_2}U_{2,3}\to \A^3\setminus\{0\}$ be the projection induced by the structure map $\P(\Sym^2F)\to \Spec k$.  Under the isomorphism 
\[
H^2(\P(\Sym^2F)\times^{\SL_2}U_{2,3}, \sW(\omega))\cong
H^2(\Tot(\check{C}^*(p_{\P(\Sym^2F)}^{-1}\sU, \sC^*_{RS}(\sW(\omega)))))
\]
of Lemma~\ref{lem:CechLerayIso}, the Euler class $\alpha^*e_{\SL_2}(T_{\P(\Sym^2F)})$ is represented by a collection of elements $\{e_{p,q}\}_{p+q=2}$,
\[
e_{p,q}\in \check{C}^p(\sU, \sC^q_{RS,\ \P(\Sym^2F)}(\sW(\omega))),
\]
where $\sC^q_{RS,\ \P(\Sym^2F)}(\sW(\omega)))$ is the presheaf on $\A^3\setminus\{0\}$
\[
\sC^q_{RS,\ \P(\Sym^2F)}(\sW(\omega))(U):= C^q_{RS}(p_{\P(\Sym^2F)}^{-1}(U), \sW(\omega))).
\]

Similar to (1) and (2), it follows from Lemma~\ref{lem:PnCoh} that the Leray spectral sequence
\[
E_2^{p,q}=H^p(\A^3\setminus\{0\}, \sH^q(\P(\Sym^2F), \sW(\omega_Y)))\Rightarrow H^{p+q}(\P(\Sym^2F)\times^{\SL_2}U_{2,3}, \sW(\omega_Y))
\]
degenerates at $E_2$ and the edge homomorphism 
\[
 H^2(\P(\Sym^2F)\times^{\SL_2}U_{2,3}, \sW(\omega_Y))\to H^0(\A^3\setminus\{0\}, \sH^q(\P(\Sym^2F), \sW(\omega_Y)))
 \]
  is an isomorphism. 

Writing $e_{0,2}$ as the family $\{e_{0,2}^i\}_i$ with  $e_{0,2}^i\in C_{RS}^2(U_i\times \P(\Sym^2F),  \sW(\omega))$ and with $\del_{RS}(e_{0,2}^i)=0$. Via the isomorphisms (see Lemma~\ref{lem:PnCoh})
\begin{multline*}
H^0(\A^3\setminus\{0\}, R^2p_{\P(\Sym^2F)*}\sW(\omega))\cong H^0(\A^3\setminus\{0\}, \sH^2(\P(\Sym^2F), \sW(\omega)))\\\cong H^2(\P(\Sym^2F), \sW(\omega))\cong W(k)
\end{multline*}
we see that the classes $[e_{0,2}^i]\in H^2(p_{\P(\Sym^2F)}^{-1}(U_i), \sW(\omega))$ agree on overlaps, defining an element $[e_{0,2}]\in H^0(\A^3\setminus\{0\}, \sH^2(\P(\Sym^2F), \sW(\omega)))\cong H^2(\P(\Sym^2F), \sW(\omega))$. In addition, $[e_{0,2}]$ determines 
$e_{\SL_2}(T_{\P(\Sym^2F)})$ and is  equal to the (non-equivariant) Euler class $e(T_{\P(\Sym^2F)})$. 

Conversely, we can construct a collection $\{e_{p,q}\}$ representing $e(T_{\P(\Sym^2F)\times^{SL_2}U_{2,3}})$ by starting with any chosen representative of $e(T_{\P(\Sym^2F)})$, and pulling back to $p_{\P(\Sym^2F)}^{-1}(U_i)\cong \P(\Sym^2F)\times U_i$. The vanishing of $\sH^q(\P(\Sym^2F), \sW(\omega))$ for $q\neq2$ says that we can then find the necessary elements $e_{1,1}, e_{2,0}$ to build a cocycle in $\Tot^2(\check{C}^*(p_{\P(\Sym^2F)}^{-1}\sU, \sC^*_{RS}(\sW(\omega))))$, whose image in $H^2(\Tot^*(***))$ will automatically represent $\alpha^*e_{\SL_2}(T_{\P(\Sym^2F)})$. 

Having done this, we can then apply the restriction map for $j{\colon}Y\to \P(\Sym^2F)$, giving the collection
\[
\{j^*e_{p,q}\}_{p+q=2}\in \Tot^2(\check{C}^*(p_Y^{-1}\sU, \sC^*_{RS}(\sW(\omega))))
\]
representing $j^*\alpha^*e_{\SL_2}(T_{\P(\Sym^2F)})=\alpha^*j^*e_{\SL_2}(T_{\P(\Sym^2F)})=\alpha^*e_{\SL_2}(T_Y)$. Now using (1) and (2), in particular, the vanishing of $\sH^q(Y, \sW(\omega))$ for $q\neq 0$, we see that the edge homomorphism 
\[
H^2(\A^3\setminus\{0\}, \sW)=H^2(\A^3\setminus\{0\}, \sH^0(Y, \sW(\omega)))\to
H^2(Y\times^{\SL_2}U_{2,3}, \sW(\omega))
\]
is an isomorphism, and that, after modifying $\{j^*e_{p,q}\}$ by a boundary to give a new collection $\{\tilde{e}_{p,q}\}$ with $\tilde{e}_{0,2}=\tilde{e}_{1,1}=0$, if necessary, $\tilde{e}_{2,0}$  maps to $\alpha^*e_{\SL_2}(T_Y)$. We then prove (3) by examining $\tilde{e}_{2,0}$ for a suitable choice of the family $\{e_{p,q}\}$.  We now carry out this program. 

By the Gau{\ss}-Bonnet theorem for the $W(k)$-valued Euler characteristic $\chi^\sW(-/k)$ \cite[Theorem 1.5]{LR}, and the computation of $\chi(\P^2/k)=\<1\>+H\in \GW(k)$ by Hoyois \cite[Example 1.7]{HoyoisTrace}, we have
\[
\pi_{\P(\Sym^2F)*}(e(T_{\P(\Sym^2F)}))=\chi^\sW(\P^2/k)=\<1\>
\]
in $W(k)$. In addition, 
\[
\pi_{\P(\Sym^2F)*}{\colon}H^2(\P(\Sym^2F),\sW(\omega)))\to W(k)
\]
is the isomorphism we use to identify  $H^2(\P(\Sym^2F),\sW(\omega)))$ with $W(k)$ from Lemma~\ref{lem:PnCoh}. Moreover,   $\pi_{\P(\Sym^2F)*}$ is inverse to the isomorphism 
\[
i_{p*}{\colon}W(k)\to H^2(\P(\Sym^2F),\sW(\omega)))
\]
for the inclusion $i_p{\colon}\Spec k\to \P(\Sym^2F)$ for an arbitrary $k$-point $p$ of $\P(\Sym^2F)$; these are canonical push-forwards, so do not require any orientation choices.

 Thus, we may take $p=[1,0,0]\in \P(\Sym^2F)$, giving the representative
\[
\<1\>\otimes \del/\del t_1\wedge \del/\del t_2\otimes dt_1dt_2\text{ on }p
\]
for $e(T_{\P(\Sym^2F)})$ in the Rost-Schmid complex for $\sW(\omega)$ on $\P(\Sym^2F)$, where $t_i:=T_i/T_0$. We then take $e_{0,2}^i\in C^2_{RS}(U_i\times \P(\Sym^2F),\sW(\omega))$ to be $<1\>\otimes \del/\del t_1\wedge \del/\del t_2\otimes dt_1dt_2\text{ on }U_i\times p$.

Write $i_{\P(F)}$ for the squaring map $\sq{\colon}\P(F)\to \P(\Sym^2F)$. Note that $p=i_{\P(F)}([1,0])\in i_{\P(F)}(\P(F))\subset \P(\Sym^2F)$, and  under our isomorphism of $\P(F)$ with $\P^1$, we have
\[
\<1\>\otimes \del/\del t_1\wedge \del/\del t_2\otimes  dt_1dt_2\text{ on }p=i_{\P(F)*}(\<1\>\otimes\del/\del (x_1/x_0)\otimes d(x_1/x_0)\text{ on }[1,0]).
\]
We take  
\[
e^i_{\P(F), 0,1}=\<1\>\otimes d(x_1/x_0)\text{ on }U_i\times (1,0),
\]
so that $i_{\P(F)*}(e^i_{\P(F), 0,1})=e_{0,2}^i$.

The term $e_{1,1}$ is a collection of elements $\{e_{1,1}^{ij}\in C^1_{RS}(U_{ij}\times \P(\Sym^2F), \sW(\omega))\}_{ij}$ with 
\[
\del^1_{RS}(\{e_{1,1}^{ij}\}_{ij})=\check{\delta}^0(\{e_{0,2}^i\}_i),
\]
where $\check{\delta}^0, \del^1_{RS}$ are the \v{C}ech, respectively, Rost-Schmid, coboundaries.  Since the canonical push-forward $i_{\P(F)*}$ gives a map of the \v{C}ech-Rost-Schmid double complex, we can solve this equation by taking  $e_{1,1}^{ij}=i_{\P(F)*}(e_{\P(F), 1,0}^{i,j})$, where the $e_{\P(F), 1,0}^{i,j}\in C^0_{RS}(U_{ij}\times \P(F), \sW(\omega_{\P(F)}))$ solve the corresponding equation for $\{e_{\P(F), 0,1}^i\}_i$, that is
\[
\del^0_{RS}(\{e_{\P(F), 1,0}^{ij}\}_{ij})=\check{\delta}^0(\{e_{\P(F), 0,1}^i\}_i).
\]
The $e_{\P(F), 1,1}^{ij}$ exist since $\SL_2$ acts trivially on $H^1(\P(F), \sW(\omega_\P(F)))$ (Remark~\ref{rem:SL2Inv}), so for each $ij$ , 
$\check{\delta}^0(\{e_{\P(F), 0,1}^i\}_i)_{ij}$ goes to zero in $H^1(U_{ij}\times\P(F),\sW(\omega_{\P(F)}))=H^1(C^*_{RS}(U_{ij}\times \P(F), \sW(\omega_{\P(F)})))$.

Finally,  $e_{2,0}$ is an element $e_{2,0}^{012}\in C^0_{RS}(U_{012}\times \P(\Sym^2F), \sW(\omega))$ with
\[
-\del^0_{RS}(\{e_{2,0}^{012}\}_{012})=\check{\delta}^1(\{e_{1,1}^{ij}\}_{ij}).
\]
Letting $e_{\P(F), 2,0}^{012}\in H^0(U_{012}\times \P(F), \sW(\omega_{\P(F)}))$ be the \v{C}ech coboundary 
$\check{\delta}^1(\{e_{\P(F), 1,0}^{ij}\}_{ij})$,   we have
\[
\check{\delta}^1(\{e_{1,1}^{i,j}\}_{ij})=i_{\P(F)*}(e^{012}_{\P(F), 2,0}).
\]
so 
\begin{equation}\label{eqn:DeltaId}
-\del^0_{RS}(e_{2,0}^{012})=i_{\P(F)*}(e^{012}_{\P(F), 2,0}).
\end{equation}

Letting $K$ be the function field of $U_{012}\times \P(\Sym^2F)$, $C^0_{RS}(U_{012}\times \P(\Sym^2F), \sW(\omega))$ is by definition $W(K;\omega)$, and the Rost-Schmid boundary map is
\[
W(K;\omega)\xrightarrow{\del^0_{RS}=\oplus_x\del^0_{RS,x}}\oplus_{x\in (U_{012}\times \P(\Sym^2F))^{(1)}} W(k(x);\omega\otimes \mathfrak{m}_x/\mathfrak{m}_x^2).
\]
The identity \eqref{eqn:DeltaId} says that $\del^0_{RS,x}=0$ except for $x$ the generic point $\eta$ of $U_{012}\times\P(F)$ and
\[
\del^0_{RS,\eta}(e_{2,0}^{012})=-e^{012}_{\P(F), 2,0}
\in W(k(U_{012}\times\P(F));\omega_{\P(F)}).
\]
This implies that $e_{2,0}^{012}$ extends uniquely to an element $e_{012}\in H^0(U_{012}\times(\P(\Sym^2F)\setminus\P(F), \sW(\omega))$ and
\begin{equation}\label{eqn:BoundaryEqn}
\del_{\P(F)}(e_{012})=-e^{012}_{\P(F), 2,0}\in H^0(U_{012}\times\P(F), \sW(\omega_{\P(F)})),
\end{equation}
where
\[
\del_{\P(F)}{\colon}H^0(U_{012}\times(\P(\Sym^2F)\setminus\P(F), \sW(\omega))\to 
H^0(U_{012}\times\P(F), \sW(\omega_{\P(F)}))
\]
is the coboundary in the localization sequence for $U_{012}\times\P(F)\subset U_{012}\times\P(\Sym^2F)$.

From our construction, the steps we used to construct $e_{\P(F), 2,0}$ show that the class $[e_{\P(F), 2,0}]\in \check{H}^2(\A^3\setminus\{0\}, \sH^0(\P(F), \sW(\omega)))$ represented by $e^{012}_{\P(F), 2,0}$ is exactly $d_2^{0,1}(i_{[1,0]*}(\<1\>))$. By Lemma~\ref{lem:LerayDifferentialComputation}, we have
\[
[d_2^{0,1}(i_{(1,0)*}(\<1\>))]=\alpha^*e(E_2)\otimes d(x_1/x_0),
\]
and by Lemma~\ref{lem:WittCohP2-D} and \eqref{eqn:BoundaryEqn}, we may take
\[
e_{2,0}^{012}:=[\tilde{Q}]\cdot \tilde{e},
\]
where $\tilde{e}\in  \check{C}^2(\A^3\setminus\{0\},\sW)$ represents $\alpha^*(e(E_2))$. 

Now we apply $j^*$. By construction, we have $j^*(e_{0,2})=0$ and $j^*(e_{1,1})=0$, so by our discussion above, we have $[j^*(e_{2,0})]$ defines an element of $H^2(\A^3\setminus\{0\}, \sW)$ representing $\alpha^*e_{\SL_2}(T_Y)$. Since $e_{2,0}^{012}:=[\tilde{Q}]\cdot \tilde{e}$, we find that $\alpha^*e_{\SL_2}(T_Y)=[\tilde{Q}]\cdot \alpha^*e(E_2)$,  so by Lemma~\ref{lem:RestrictionToX}, 
$e_{\SL_2}(T_Y)=[\tilde{Q}]\cdot(E_2)$, completing the proof of (3).
\end{proof} 

We now have the information we need to  compute products for our generators for $H^*_{\SL_2}(Y, \sW)\oplus H^*_{\SL_2}(Y, \sW(\omega))$.  

\begin{proposition}\label{prop:MultiplicativeStructure}  With respect to the canonical isomorphism  $\sW(\omega^{\otimes 2})\cong \sW$, we have the following identities.\\[5pt]
 1.  $[\tilde{Q}]^2=2\cdot(1-[\<Q\>])$ in $H^0_{\SL_2}(Y, \sW)$. \\[2pt]
 2. $[\<Q\>]\cdot[\tilde{Q}]=-[\tilde{Q}]$ in  $H^0_{\SL_2}(Y, \sW(\omega))$. \\[2pt]
3.  $[\tilde{Q}]\cdot e_{\SL_2}(T_{X_a})=4e$ in $H^2_{\SL_2}(Y, \sW)$. \\[2pt]
4. $(e_{\SL_2}(T_{X_a}))^2=4e^2$ in $H^2_{\SL_2}(Y, \sW)$.
\end{proposition}

\begin{proof}  (1) Since $H^0_{\SL_2}(Y, \sW)=H^0(Y, \sW)$, we can forget the $\SL_2$-action and then pass to any open subset of $Y$  to compute $[\tilde{Q}]^2$. After inverting $T_0$, we diagonalize $\tilde{Q}$ as 
\[
\tilde{q} = (\<1\>-\<Q/T_0^2\>)\otimes\frac{2\Omega_2}{T_0^3}=(\<1\>-\<Q\>)\otimes\frac{2\Omega_2}{T_0^3},
\]
so after applying the isomorphism $\sW(\omega^{\otimes 2})\cong \sW$, and recalling that $\<Q\>^2=1$, we have
\[
[\tilde{Q}]^2= (1-[\<Q\>])^2=2\cdot(1-[\<Q\>])
\]
The proof of (2) is similar: $\<Q/T_0^2\>\cdot \tilde{q}=-\tilde{q}$,
whence $[\<Q\>]\cdot[\tilde{Q}]=-[\tilde{Q}]$.

For (3), we have already shown that $e_{\SL_2}(T_Y)=[\tilde{Q}]\cdot e$, so 
\[
[\tilde{Q}]\cdot e_{\SL_2}(T_Y)=[\tilde{Q}]^2\cdot e =2\cdot(1-[\<Q\>])e
\]
The element $[\<Q\>]\in H^0_{\SL_2}(Y, \sW)$ is by \cite[Proposition 5.5]{LevineBG} exactly the element $x$ of Theorem~\ref{thm:Presentation}, hence  $[\<Q\>]\cdot e=- e$. Thus we have $[\<Q\>]\cdot e_{\SL_2}(T_{X_a})=4e$. 

(4) is similar: $e_{\SL_2}(T_Y)=[\tilde{Q}]\cdot e\Rightarrow$
\[
 (e_{\SL_2}(T_Y))^2=[\tilde{Q}]^2\cdot e\cdot e= 4e^2
 \]
by (3). 
\end{proof}

Consider $H^*(BN, \sW)\oplus H^*(BN, \sW(\gamma))$ as a $\Z\times \Z/2$-graded $W(k)$-algebra via the isomorphism $\sW(\gamma^{\otimes 2})\cong \sW$. Recall that we identify $BN$ with $Y\times^{\SL_2}E\SL_2$; let $p{\colon}BN\to \BSL_2$ correspond to the projection $Y\to \Spec k$. The quadratic forms $\<Q\>$ and $\tilde{Q}$ define the elements $[\<Q\>]\in H^0(BN, \sW)$, $[\tilde{Q}]\in H^0(BN, \sW(\gamma))$. We have the Euler class  $e(E_2)\in H^*(\BSL_2, \sW)$ and  the Euler class $e_{\SL_2}(T_Y)\in H^2(BN, \sW(\gamma))$. 

\begin{theorem}\label{thm:CohBN} Give variables $x, e, y, \tilde{e}$ bi-degrees $(0,0), (1,0), (0,1)$ and $(1,1)$ in $\Z\times\Z/2$, respectively. Let $I\subset W(k)[x,e,y,\tilde{e}]$ be the ideal
\[
I:=(x^2-1, (1+x)e, ye-\tilde{e}, (1+x)y, y^2+2x-2).
\]
Then the assignments 
\[
x\mapsto [\<Q\>],\  e\mapsto p^*e(E_2),\ y\mapsto [\tilde{Q}],\ \tilde{e}\mapsto e_{\SL_2}(T_Y)
\]
define an isomorphism of $\Z\times \Z/2$-graded $W(k)$-algebras
\[
W(k)[x,e,y,\tilde{e}]/I\xymatrix{\ar[r]^\alpha_\sim&}
H^*(BN, \sW)\oplus H^*(BN, \sW(\gamma)).
\]
\end{theorem}

\begin{proof} The $W(k)$-module $W(k)[x,e,y,\tilde{e}]/I$
 is the free $W(k)$-module with basis 
$\{e^n\}_{n\ge0}\cup \{e^n\cdot \tilde{e}\}_{n\ge0}\cup \{x,y\}$.
By Theorem~\ref{thm:Presentation}, $ \{[\<Q\>]\}\cup\{p^*e_{\SL_2}(F)^n\}_{n\ge0}$ is a basis of the   free $W(k)$-module $H^*(BN, \sW)$,  and the free $W(k)$-module 
  $H^*(BN, \sW(\gamma))$ has basis $\{[\tilde{Q}]\}\cup  \{p^*e_{\SL_2}(F)^n\cdot e_{\SL_2}(T_Y)\}_{n\ge0}$. Thus the map $\alpha$ is an isomorphism of  $\Z\times \Z/2$-graded $W(k)$-modules. By Theorem~\ref{thm:Presentation} and Proposition~\ref{prop:MultiplicativeStructure}, $\alpha$ is a $W(k)$-algebra isomorphism. 
\end{proof}

 \end{document}